%% file: thesis.tex
\documentclass{ut-thesis}

\usepackage{amssymb,amsfonts, amsmath, amsthm, graphicx, enumerate, upgreek}
\usepackage{mathtools, etoolbox, bm}
\usepackage[table]{xcolor}
\usepackage{tikz}
\usepackage[all,arc]{xy}
\usepackage[utf8]{inputenc}
\usepackage[english]{babel}

\usepackage[pagebackref, hidelinks]{hyperref} 

\newtheorem{thm}{Theorem}[section]
\newtheorem*{thm*}{Theorem}
\newtheorem{thmx}{Theorem}
\newtheorem{thmy}{Theorem}
\newtheorem{cor}[thm]{Corollary}
\newtheorem*{cor*}{Corollary}
\newtheorem{prop}[thm]{Proposition}
\newtheorem*{prop*}{Proposition}
\newtheorem{lem}[thm]{Lemma}
\newtheorem*{lem*}{Lemma}
\newtheorem{conj}[thm]{Conjecture}
\newtheorem*{conj*}{Conjecture}
\newtheorem{quest}[thm]{Question}
\newtheorem*{quest*}{Question}

\theoremstyle{definition}
\newtheorem{defn}[thm]{Definition}

\newtheorem{con}[thm]{Construction}
\newtheorem{exmp}[thm]{Example}

\newtheorem{assump}[thm]{Assumption} 

\theoremstyle{remark}
\newtheorem{rem}[thm]{Remark}

\newtheorem*{ack*}{Acknowledgements}

\newcommand{\modsf}{\mathsf{mod}\ \!} 
\newcommand{\Modsf}{\mathsf{Mod}\ \!} 
\newcommand{\grmod}{\mathsf{grmod}\ \!} 
\newcommand{\Grmod}{\mathsf{Grmod}\ \!} 
\newcommand{\qgr}{\mathsf{qgr}\ \!} 
\newcommand{\QGr}{\mathsf{QGr}\ \!} 
\newcommand{\coh}{\mathsf{coh}\ \!} 
\newcommand{\qcoh}{\mathsf{QCoh}\ \!} 

\newcommand{\N}{\mathbb{N}} 
\newcommand{\Z}{\mathbb{Z}} 
\newcommand{\Q}{\mathbb{Q}} 
\newcommand{\C}{\mathbb{C}} 

\newcommand{\D}{\mathcal{D}} 
\newcommand{\T}{\mathcal{T}} 
\newcommand{\CC}{\mathcal{C}} 
\newcommand{\A}{\mathcal{A}} 
\newcommand{\B}{\mathcal{B}} 
\newcommand{\Se}{\mathbb{S}} 
\newcommand{\perf}{\mathsf{perf}} 
\newcommand{\thick}{\mathsf{thick}} 

\newcommand{\m}{\mathfrak{m}} 
\newcommand{\depth}{{\rm depth}\ \!} 
\newcommand{\edim}{{\rm edim}\ \!} 
\newcommand{\pdim}{\mathsf{pdim}\ \!} 
\newcommand{\idim}{\mathsf{idim}\ \!} 
\newcommand{\gldim}{\mathsf{gldim}\ \!} 

\newcommand{\rad}{\mathsf{rad}} 

\newcommand{\spec}{\mathsf{spec}\ \!} 
\newcommand{\proj}{\mathsf{proj}\ \!} 
\newcommand{\OO}{\mathcal{O}} 
\newcommand{\E}{\mathcal{E}} 
\newcommand{\F}{\mathcal{F}} 
\newcommand{\G}{\mathcal{G}} 
\newcommand{\ZZ}{\mathcal{Z}} 

\newcommand{\X}{\mathbb{X}} 

\newcommand{\U}{\mathcal{U}} 

\degree{Doctor of Philosophy}
\department{Mathematics}
\gradyear{2018}
\author{Vincent G\'elinas}
\title{Contributions to the Stable Derived Categories of Gorenstein Rings} 

\begin{document}

\begin{preliminary}

\maketitle


\begin{abstract}
	The stable derived category ${\rm D}_{sg}(R)$ of a Gorenstein ring $R$ is defined as the Verdier quotient of the bounded derived category $ {\rm D}^b(\mathsf{mod}\ \! R) $ by the thick subcategory of perfect complexes, and was introduced by Ragnar-Olaf Buchweitz as a homological measure of the singularities of $R$. This thesis contributes to its study, centered around representation-theoretic and Koszul duality aspects.	

	In Part I, we first complete (over $\mathbb{C}$) the classification of homogeneous complete intersection isolated singularities $R$ for which the graded stable derived category ${\rm D}^{\mathbb{Z}}_{sg}(R)$ (respectively, $ {\rm D}^b(\mathsf{coh}\ \! X) $ for $X = \mathsf{proj}\ \! R$) contains a tilting object. For this, we prove the existence of a full strong exceptional collection of vector bundles on a $2n$-dimensional smooth complete intersection of two quadrics $X = V(Q_1, Q_2) \subseteq \mathbb{P}^{2n+2}$, building on work of Kuznetsov. We then use results of Buchweitz-Iyama-Yamaura to describe the indecomposable objects in ${\rm D}_{sg}^{\mathbb{Z}}(R_Y)$ for $R_Y$ the homogeneous coordinate rings of $4$ points on $\mathbb{P}^1$ and $4$ points on $\mathbb{P}^2$ in general position, and classify their Betti tables. 

	In Part II, for $R$ a Koszul Gorenstein algebra, we study a natural pair of full subcategories whose intersection $\mathcal{H}^{\mathsf{lin}}(R) \subseteq {\rm D}_{sg}^{\mathbb{Z}}(R)$ consists of modules with eventually linear projective resolutions. We prove that such a pair forms a bounded t-structure if and only if $R$ is absolutely Koszul in the sense of Herzog-Iyengar, in which case there is an equivalence of triangulated categories ${\rm D}^b(\mathcal{H}^{\mathsf{lin}}(R)) \cong {\rm D}_{sg}^{\mathbb{Z}}(R)$. We then relate the heart to modules over the Koszul dual $R^!$. This extends the Bernstein-Gel'fand-Gel'fand correspondence to any pair of Koszul dual algebras $(R, R^!)$ with $R$ absolutely Koszul Gorenstein, applying in particular to the coordinate ring of elliptic normal curves of degree $\geq 4$ and to the anticanonical model of del Pezzo surfaces of degree $\geq 4$. We then relate our results to conjectures of Bondal and Minamoto on the graded coherence of Artin-Schelter regular algebras and higher preprojective algebras; we characterise when these conjectures hold in a restricted setting, and give counterexamples to both in all dimension $\geq 4$.

\end{abstract}


\begin{dedication}
       To the memory of Ragnar-Olaf Buchweitz.\\
	To my Mary.
\end{dedication}

\newpage  


\begin{acknowledgements}

This entire thesis is dedicated to Ragnar-Olaf Buchweitz, who, as a broadly knowledgeable and deeply creative mathematician, as a contagiously enthusiastic lecturer, and as the kindest man I have ever met, has always set the standard to which one should aspire. \bigskip 

	That I have made it all the way here is largely the result of having been surrounded by absolutely outstanding people, to whom thanks are due but for which words hardly seem enough.\bigskip

	To Ben, who along with Ragnar has taught me all of the mathematics that I know. \bigskip

	To \"Ozg\"ur, who is always high in spirit, for all the good times, for being the best conference roommate, and for being a solid friend all this time. \bigskip 

	To Louis-Philippe, for being amazingly kind, thoughtful, and also for smooth driving. \bigskip

	To Jemima, for all the constant help and for taking care of all of us. \bigskip

	To the inhabitants of the math house, for some of the best food I've had in grad school. \bigskip

	To my parents and my sister, for the constant support, and without whom I would never have made it here.\bigskip

	I owe an immense debt to my supervisors Colin and Joel, who have been extremely generous with their time, advice and mentorship. The same goes for John Scherk, Paul Selick, Marco Gualtieri and Matilde Marcolli, for agreeing to serve on my committee, for the help on many occasions and for the solid advice that I have received from them. \bigskip 
	
	I have benefitted greatly from discussions with Srikanth Iyengar during a very pleasant visit to the University of Utah. The second part of this thesis is largely a result of my time spent there.\bigskip 

	Finally, I owe everything to Mary McDonald, who pushes me and challenges me, who makes our lives comfortable and interesting, who constantly gets me to try new things, who shares her interests with me, and who has always been the best partner for exploring the road ahead. 

\end{acknowledgements}

\tableofcontents




\end{preliminary}


\include{chapters/introduction/introduction}

\part{MCM Modules and Representation Theory of Algebras}
\include{chapters/tilting/tilting}

\include{chapters/classifications/classifications}

\part{Absolutely Koszul Gorenstein algebras}
\include{chapters/fano/fano}

\include{chapters/tstructure/tstructure}

\include{chapters/applications/applications}

\appendix
\include{chapters/appendix/appendix}

\bibliographystyle{plain}
\bibliography{thesis}


\end{document}

%% file: chapters/introduction/introduction.tex
\chapter{Introduction}
This thesis is concerned with representation-theoretic and homological aspects of the stable derived categories of Gorenstein rings, mostly centered around questions of derived Morita theory, Koszul duality and tame classification problems. Stable derived categories were originally introduced by Ragnar-Olaf Buchweitz and later independently by D. Orlov under the name of \emph{singularity categories}. The stable derived category of a Noetherian ring $R$ is defined as the Verdier quotient
\[
	{\rm D}_{sg}(R) = {\rm D}^b(\modsf R)/{\rm D}^{\perf}(\modsf R)
\]
of the bounded derived category of $R$ by the thick subcategory of perfect complexes, meaning complexes with finite projective resolutions. When $R$ is commutative, the Auslander-Buchsbaum-Serre characterisation of regular rings as the commutative rings of finite global dimension shows that the singularity category ${\rm D}_{sg}(R)$ vanishes precisely for $R$ regular, and thus in general provides a measure of the singularities of $\spec R$. 

A two-sided Noetherian, not necessarily commutative ring $R$ is Gorenstein if both injective dimensions $\idim( _R R) < \infty$ and $\idim(R_R) < \infty$ are finite, in which case those dimensions are equal and define the Gorenstein dimension of $R$. When $R$ is Gorenstein, a foundational theorem of Buchweitz identifies ${\rm D}_{sg}(R) \cong \underline{{\rm MCM}}(R)$ with the stable module category of maximal Cohen-Macaulay (MCM) $R$-modules, whose study interpolates singularity theory, commutative algebra and representation theory. When $R$ is graded, we write ${\rm D}_{sg}^{\Z}(R)$ for the version involving graded modules. These triangulated categories have turned out to be ubiquitous in the last 30 years, and routinely appear in algebra and geometry. We list some notable examples: 
\begin{enumerate}[i.]
	\item The homotopy category of matrix factorizations of $f \in S$ for $S$ regular is equivalent to the stable derived category of the hypersurface singularity $S/f$, that is $\underline{{\rm MF}}(S, f) \cong \underline{{\rm MCM}}(S/f)$ by a theorem of Eisenbud. 
	\item The stable derived category of a self-injective algebra $\Lambda$ is equivalent to its stable module category $\underline{\mathsf{mod}}\ \! \Lambda$. This applies in particular to the group algebra $\Lambda = kG$ of a finite group $G$ over a field of characteristic $p > 0$ dividing the order of $G$.
	\item The stable derived category of the group ring $SG$ for $S$ regular and $G$ finite is equivalent to the stable category of lattices over $SG$. 
	\item Stable derived categories arise on the opposite side of Koszul duality equivalences with the bounded derived categories $ {\rm D}^b(\coh X) $ of coherent sheaves on projective varieties $X$, as in the Bernstein-Gel'fand-Gel'fand (BGG) correspondence 
		\[
			{\rm D}^{\Z}_{sg}(\wedge V^*) \cong {\rm D}^b(\coh \mathbb{P}(V)) 
		\]
		and its generalisations to complete intersections of quadrics by Buchweitz and Kapranov. 
	\item Let $X \subseteq \mathbb{P}^n$ be an arithmetically Gorenstein projective variety with homogeneous coordinate ring $R_X$, meaning that the affine cone $\spec R_X$ over $X$ has Gorenstein singularities. Then there is an adjoint pair ${\rm D}_{sg}^{\Z}(R_X) \leftrightarrows {\rm D}^b(\coh X)  $ between the graded stable derived category of $R_X$ and the derived category of coherent sheaves on $X$, which, by a theorem of Orlov, is always an equivalence when $X$ is Calabi-Yau, and in general embeds one category inside the other according to the sign of the twist on the canonical module $\omega_{R_X} \cong R_X(a)$ of $R_X$. 
	\item Let $\Lambda$ be a finite dimensional algebra over a field, with finite global dimension $\gldim \Lambda < \infty$. By Happel's Theorem, one can associate to $\Lambda$ a graded Gorenstein finite dimensional algebra $T(\Lambda)$, its \emph{trivial extension algebra}, for which we always have a triangulated equivalence $ {\rm D}^b(\modsf \Lambda) \cong {\rm D}^{\Z}_{sg}(T(\Lambda))$. 
\end{enumerate} 

The present thesis makes contributions to various aspects of the study of ${\rm D}_{sg}(R)$ for $R$ Gorenstein. The work is separated into two independent sections, motivated by different ways of presenting ${\rm D}_{sg}(R)$. Each comes with applications to problems in related fields. In particular, this thesis pulls ideas from commutative algebra, algebraic geometry, representation theory of associative algebras and noncommutative algebraic geometry, which are tied together by the use of triangulated categories and can be fruitfully studied by means of the stable derived category.

\begin{enumerate}[I.]
	\item In the first part, we study the \emph{tilting problem}, which asks for abstract realisations
		\[
		{\rm D}^{\Z}_{sg}(R) \cong {\rm D}^{\perf}(\Lambda) 
		\]
	of the graded stable derived categories of graded Gorenstein $k$-algebras as the perfect derived category ${\rm D}^{\perf}(\Lambda) \subseteq {\rm D}^b(\modsf \Lambda)$ of a finite dimensional $k$-algebra $\Lambda$, with ${\rm D}^{\perf}(\Lambda) = {\rm D}^b(\modsf \Lambda)$ when $\gldim \Lambda < \infty$. This gives a strong handle on the structure of the triangulated category ${\rm D}^\Z_{sg}(R)$, giving an understand of its numerical invariants such as its Grothendieck group $K_0$, or in good cases producing classifications of the indecomposable objects in $ {\rm D}^\Z_{sg}(R) $ when the representation theory of $\Lambda$ is well-understood. Producing such equivalences amounts to finding a \emph{tilting object} in ${\rm D}^{\Z}_{sg}(R)$, from which the equivalence arises through Morita theory. We study existence questions and obstructions to such equivalences, and apply tilting theory to classification problems over tame curve singularities. 
	\item In the second part, we study the stable derived categories of Koszul Gorenstein algebras. We relate questions of existence of t-structures of `Koszul type' to classical rationality problems in commutative algebra. We then show that the Koszul Gorenstein algebras satisfying the statement of the Bernstein-Gel'fand-Gel'fand correspondence
		\[
			D_{sg}^{\Z}(R) \cong {\rm D}^b(\qgr R^!) 
		\]
		are precisely the absolutely Koszul algebras of Herzog-Iyengar. This broadly extends the class of algebras for which the correspondence holds, and we study the equivalence of triangulated categories in many new cases, as well as produce explicit counterexamples where such an equivalence fails to exist. We further relate the matter to structural questions in noncommutative algebraic geometry: we characterise the graded coherence of a class of Artin-Schelter regular algebras and higher preprojective algebras by the absolute Koszulity of an associated Gorenstein Koszul algebras. This leads us to construct examples of non-coherent Artin-Schelter regular algebras and non-coherent higher preprojective algebras in all global dimension $\geq 4$, thereby giving counterexamples to conjectures of A. Bondal and H. Minamoto, respectively. 
\end{enumerate}

\input{chapters/introduction/overview}
\input{chapters/introduction/prelim}
\input{chapters/introduction/background}

%% file: chapters/introduction/overview.tex
\section{Overview of the thesis}
\subsection*{Chapter 1}
We begin with various preliminaries on standard notation, well-known constructions and basic notions to be taken for granted throughout this thesis. We then give a short background exposition of the theory of maximal Cohen-Macaulay (MCM) modules over a Gorenstein ring, beginning with classical local commutative algebra, continuing on to Buchweitz's manuscript \cite{BuchweitzManuscript}, and ending with Orlov's semiorthogonal decomposition theorem. \bigskip

For the purpose of this thesis, the importance of the class of MCM modules over a Gorenstein ring $R$ is encapsuled in the following two results of Buchweitz and Orlov.\bigskip 

For any complex of $R$-module with bounded cohomology $N \in {\rm D}^b(\modsf R) $, there is an MCM module $M = N^{st}$, unique up to stable isomorphism, which becomes isomorphic to $N$ in the Verdier quotient ${\rm D}_{sg}(R) := {\rm D}^b(\modsf R) / {\rm D}^{\perf}(R)$ by complexes of finite projective dimension. Moreover, MCM modules can be characterised as those modules $M$ admitting a two-sided projective resolution $C_*$ 
\[
	\xymatrix@R=5pt{ \cdots \ar[r] & C_1 \ar[r] & C_0 \ar[dr] \ar[rr] && C_{-1} \ar[r] & C_{-2} \ar[r] & \cdots \\
	&&& M \ar[ur] } 
\]
called its \emph{complete} resolution. The homotopy theory of the resulting acyclic complexes $C_*$ then describes the triangulated category ${\rm D}_{sg}(R)$, which is the content of Buchweitz's Theorem. \bigskip

Given a projective variety $X$ and an ample line bundle $\mathcal{L}$ with section ring 
\[
	R_X = R_{X, \mathcal{L}} = \bigoplus_{n \geq 0} {\rm H}^0(X, \mathcal{L}^{\otimes n})
\]
a graded Gorenstein ring, the structure of graded MCM modules over $R_X$ reveals a surprising amount of the geometry of $X$. Letting $M$ be a graded MCM $R_X$-module with complete resolution $C_*$, and killing free summands in $C_*$ generated in degree below some fixed cut-off $i \in \Z$, one obtains a right bounded complex $F = C_{[\geq i]}$ of free $R_X$-modules with bounded cohomology 
\[
	\xymatrix@R=5pt{ \cdots \ar[r] & F_{n+3} \ar[r] & F_{n+2} \ar[r] & F_{n+1}  \ar[r] & F_{n} \ar[r] & 0 }
\]
which sheafifies to a complex of coherent sheaves $\F = \widetilde{F} \in {\rm D}^b(\coh X) $ with bounded cohomology. By Orlov's Theorem, the resulting functor 
\[
	{\rm D}_{sg}^\Z(R_X) \to {\rm D}^b(\coh X) 
\]
is an equivalence of triangulated categories whenever $X$ is Calabi-Yau, and otherwise induces an equivalences upon taking Verdier quotient by a certain exceptional collection. 

\subsection*{Chapter 2}
In Chapter 2, we review historical background and definitions concerning tilting theory for derived categories of coherent sheaves $ {\rm D}^b(\coh X) $ and graded singularity categories $ {\rm D}_{sg}^\Z(R)$. Roughly speaking, a tilting object $T \in \T$ in a triangulated category $\T$, with $\T$ idempotent complete and of algebraic origin, is the triangulated category analog of a small projective generator in Morita Theory, in that it induces an equivalence of triangulated categories
\[
	{\rm RHom}(T, -): \T \xrightarrow{\cong} {\rm D}^{\perf}(\Modsf {\rm End}_\T(T)).
\]
We review briefly the important classes of Gorenstein algebras for which ${\rm D}_{sg}^\Z(R)$ is known to contain a tilting object, with special attention to results of Yamaura and Buchweitz-Iyama-Yamaura which guarantee this in dimension $\leq 1$, under mild hypotheses. In particular, this applies to connected graded Artinian self-injective algebras (dimension zero) and reduced Gorenstein curve singularities (dimension one). \bigskip

We then turn our attention to dimension $\geq 2$. Restricting to the case of commutative Gorenstein algebras, always assumed Noetherian, we prove the following results. Recall that $R$ is called standard graded if $R = k \oplus R_1 \oplus R_2 \oplus \cdots$ is generated by $R_1$ over $R_0 = k$. 
\begin{prop*}[Prop. \ref{hodgeobstructionmcm}] Let $k$ be an algebraically closed field and let $R$ be a standard graded Gorenstein $k$-algebra of $\dim R \geq 2$, having at most isolated singularities. Let $X = \proj R$. Assume that ${\rm D}^{\Z}_{sg}(R)$ admits a tilting object $T$. Then:
	\begin{enumerate}[i)]
	\item ${\rm H}^q(X, \OO_X) = 0$ for $q > 0$. 
	\item The $a$-invariant of $R$ satisfies $a < 0$. In particular $X$ is a Fano variety. 
\end{enumerate}
Moreover, if $k$ has characteristic zero, then $i)$ may be strenghthened to the vanishing of Hodge numbers:
\begin{enumerate}[i')]
\item ${\rm H}^{q}(X, \Omega^p_X) = 0$ for $p \neq q$.
\end{enumerate}
\end{prop*}

\begin{thm*}[Thm. \ref{hodgeconverse}] Let $k$ be an algebraically closed field of characteristic zero. Let $R$ be a standard graded complete intersection with at most isolated singularities, and let $X = \proj R$. The following are equivalent:
	\begin{enumerate}[1)]
		\item ${\rm D}^{\Z}_{sg}(R)$ admits a tilting object $T$.
		\item $ {\rm D}^b(\coh X) $ admits a tilting object $\E$.
		\item We have $H^q(X, \Omega^p_X) = 0$ for $p \neq q$. 
	\end{enumerate}
	Moreover, when $\dim R \geq 2$ and $R$ is non-regular, this is equivalent to 
\begin{enumerate}[1)] 
\setcounter{enumi}{3}
	\item $X$ is one of the following: 
			\begin{enumerate}[a)]
			\item A smooth quadric hypersurface.
			\item A smooth $2n$-dimensional complete intersection of two quadrics.
			\item A smooth cubic surface.
		\end{enumerate}
	\end{enumerate}
Lastly, when these equivalent conditions hold, the tilting object $\E \in {\rm D}^b(\coh X)$ can always be chosen to come from a full strong exceptional collection of vector bundles. 
\end{thm*} 
The first proposition should be expected if not folklore amongst experts. Likewise many of the implications in the above theorem are known, and follow from work of Beilinson, Kapranov, Kuleshov-Orlov, Kuznetsov, Rapoport and Buchweitz-Iyama-Yamaura. Our contribution is to treat the remaining case of the $2n$-dimensional intersection of two quadrics $X = V(Q_1, Q_2) \subseteq \mathbb{P}^{2n+2}$, in particular proving that $ {\rm D}^b(\coh X) $ admits a full strong exceptional collection of vector bundles, which is compatible with an analogous collection in the graded singularity category $ {\rm D}_{sg}^\Z(R)$. 
\subsection*{Chapter 3}
In Chapter 3, we explore applications of tilting theory to reduced curve singularities, where we apply the recent tilting result of Buchweitz-Iyama-Yamaura to the classification problem for indecomposable MCM modules over certain tame curve singularities. Namely, we study the homogeneous coordinate rings $R_Y$ of $4$ points on $\mathbb{P}^2$ in general position and $4$ points on $\mathbb{P}^1$, respectively. To each algebra $R_Y$, the theorem of Buchweitz-Iyama-Yamaura associates a quiver path algebra with relations $kQ/I$ along with an equivalence of triangulated categories
\[
	{\rm D}_{sg}^{\Z}(R_Y) \cong {\rm D}^b(\modsf kQ/I).
\]
In the case of $4$ points $Y \subseteq \mathbb{P}^2$ in general position, the algebra is given by the path algebra $kQ$ of the $4$-subspace quiver
	{\tiny{
			\[
	\xymatrix@R=5pt@C=20pt{&\bullet\\
		&\bullet\\
		\bullet \ar[uur]^{p_1} \ar[ur]_{\ p_2} \ar[dr]^{\ p_3} \ar[ddr]_{p_4}& \\
		&\bullet\\
		&\bullet
	}
	\]
}}

\noindent To $4$ points $Y \subseteq \mathbb{P}^1$ on the projective line, one attaches the `Squid' path algebra $Sq(2,2,2,2;\lambda) = kQ/I$ with quiver 
{\tiny{
		\[
	\xymatrix@R=5pt@C=20pt{&&\bullet\\
		&&\bullet\\
		\bullet \ar@<.5ex>[r]^y \ar@<-.5ex>[r]_x &\bullet \ar[uur]^{p_1} \ar[ur]_{\ p_2} \ar[dr]^{\ p_3} \ar[ddr]_{p_4}& \\
		&&\bullet\\
		&&\bullet
	}
\]
}}
and relations $p_i l_i(x,y) = 0$ for $i = 1,2,3,4$, where $l_i$ is the linear form cuting out the $i$-th point in $Y$. \bigskip

Both algebras have derived tame representation type, meaning that the classification problem for indecomposables in $ {\rm D}^b(\modsf kQ/I) $ is of mild complexity. In fact the representation theory of both algebras is well-understood. Our contribution is to studying the analogous classification in the graded singularity category ${\rm D}_{sg}^{\Z}(R_Y)$: in the first case we show that one recovers the regular components of $\modsf kQ$ over the $4$-subspace quiver from the geometry of the pencil of conics through $4$ points $Y \subseteq \mathbb{P}^2$. We then deduce the classification of indecomposable MCM modules over $R_Y$ and write down the Betti tables of indecomposable complete resolutions. \bigskip

In the second case of $R_Y$ the coordinate ring of $4$ points $Y \subseteq \mathbb{P}^1$, we study the classification problem for MCM modules over $R_Y$ by reducing the problem to the weighted projective line $\X = \mathbb{P}^1(2,2,2,2;\lambda)$ of genus one, which is derived equivalent to the tubular algebra $Sq(2,2,2,2;\lambda)$. We discuss the role of the induced braid group action on the singularity category, write down the indecomposable MCM modules corresponding to the simple torsion sheaves on $\X$, and give a complete classification of the Betti tables of indecomposable graded MCM modules over $R_Y$. 

\subsection*{Chapter 4}
Chapter 4 is largely expository but contains many of the key ideas and setup required for later chapters. We begin with recalling the basics of a class of finite dimensional algebras $\Lambda$ called \emph{Fano algebras}. The prototypical example of a Fano algebra arises as follows: let $X$ be a smooth projective Fano variety with a tilting bundle $\E$, and set $\Lambda = {\rm End}_X(\E)$. We obtain an equivalence 
\[
	{\rm RHom}(\E, -): {\rm D}^b(\coh X) \xrightarrow{\cong} {\rm D}^b(\modsf \Lambda).
\]
One may take here $X = \mathbb{P}^n$ and $\Lambda = {\rm End}(\bigoplus_{i=0}^n \OO(i))$ the $n$-th Beilinson algebra. Further assume that $\dim X = \gldim \Lambda =: n$, which holds for the example above. Let
\[
	\Se_\Lambda = - \otimes_\Lambda^{\mathbb{L}} D\Lambda 	
\]
be the Serre functor for $ {\rm D}^b(\modsf \Lambda) $ where $D\Lambda = {\rm Hom}_k(\Lambda, k)$, and set $\Se_n = \Se_\Lambda \circ [-n]$ for the $n$-shifted Serre functor on $ {\rm D}^b(\modsf \Lambda) $, which corresponds under the above equivalence to $- \otimes \omega_X$. Since $X$ is Fano, using ampleness of $\omega_X^{-1}$ it is easy to see that the equivalence sends the subcategory $\coh X \subseteq {\rm D}^b(\coh X) $ onto the subcategory $\mathcal{H}(\Lambda) \subseteq {\rm D}^b(\modsf \Lambda) $ of \emph{asymptotic modules} 
\[
	\mathcal{H}(\Lambda) := \{ M \in {\rm D}^b(\modsf \Lambda) \ | \ \Se_n^{-m}(M) \in \modsf \Lambda \subseteq {\rm D}^b(\modsf \Lambda) \textup{ for all } m \gg 0 \}
\]
and so we recover the coherent sheaves on $X$ via a canonically defined subcategory of $ {\rm D}^b(\modsf \Lambda) $. \bigskip

The class of Fano algebras, as introduced by H. Minamoto, precisely captures those finite dimensional algebras which come from Fano varieties, at least if one allows the varieties to be noncommutative. More precisely, by a beautiful theorem of Minamoto, any Fano algebra $\Lambda$ satisfies 
\[
	{\rm D}^b(\modsf \Lambda) \cong {\rm D}^b(\qgr \Pi) 
\]
for some graded algebra $\Pi$, where $\qgr \Pi$ is the category of finitely presented modules modulo torsion modules, thought of as a category of coherent sheaves of a noncommutative variety, and furthermore $\qgr \Pi$ corresponds to $\mathcal{H}(\Lambda)$ under this equivalence. The graded algebra $\Pi$ is moreover \emph{coherent}, a weakening of the Noetherian condition which is required for the whole machinery to work. The algebra $\Pi$ is constructed out of $\Lambda$ by a process reminescent of the anticanonical ring $\bigoplus_{n \geq 0} {\rm H}^0(X, \omega_X^{-n})$ of a smooth variety $X$.\bigskip

Fano algebras belong to the larger class of \emph{almost} Fano algebras, which are finite dimensional algebras containing enough data to construct $\Pi$ canonically, but without the a priori knowledge that $\Pi$ is a graded coherent algebra. Almost fano algebras arise much more often in practice, and it is important to determine whether they are actually Fano. \bigskip

The class of Fano algebras contains all representation infinite quiver path algebras $kQ$, meaning that $Q$ is a finite acyclic non-Dynkin quiver $Q$, and $\Pi = \Pi(Q)$ is the classical preprojective algebra of $Q$. Similarly, the class of almost Fano algebras contains the subclass of \emph{higher representation infinite algebras} introduced by Herschend-Iyama-Oppermann, and $\Pi = \Pi(\Lambda)$ is the associated higher preprojective algebra. The missing ingredient for such algebras to be Fano is to establish the graded coherence of $\Pi(\Lambda)$. As such, central to this storyline is the conjecture:
\begin{conj*}[Minamoto, Conjecture \ref{Minamotoconjecture}] Higher preprojective algebras $\Pi(\Lambda)$ of higher representation infinite algebras are always graded coherent. 
\end{conj*}
This conjecture appears to be partially modelled on an older conjecture of A. Bondal, concerning the class of Artin-Schelter regular algebras, which share many similarities with $\Pi(\Lambda)$. 
\begin{conj*}[Bondal, Conjecture \ref{Bondalconjecture}] Artin-Schelter regular algebras are always graded coherent. 
\end{conj*}

In this chapter, we will see how to associate, to any Koszul Frobenius algebras $A$, a pair $(E, \Pi)$ consisting of an Artin-Schelter regular algebra $E$ and the higher preprojective algebra $\Pi = \Pi(\Lambda)$ of a higher representation infinite algebra $\Lambda$. That one can do this is not new, but this doesn't appear to have been exploited as much as one would expect. More interestingly, we will see how to characterise the coherence of $E$ and $\Pi$ in terms of the graded singularity category ${\rm D}_{sg}^{\Z}(A)$, and therefore characterise the algebras $A$ giving rise to a pair of coherent graded algebras $(E, \Pi)$. Our main results will then be:
\begin{thm*}[Thm. \ref{summary6.2}] The following are equivalent:
	\begin{enumerate}[i)]
	\item $\Pi$ is coherent.
	\item $E$ is coherent.
	\item $A$ is absolutely Koszul in the sense of Herzog-Iyengar. 
\end{enumerate}
\end{thm*}
This last condition is well-studied and holds for many classes of Koszul algebra. However we will also construct Koszul Frobenius algebras which are not absolutely Koszul, and so obtain: 

\begin{thm*}[Thm. \ref{abscounterexamples}]
	There are Artin-Schelter regular algebras $\{ E_d\}_{d \geq 4}$ and higher preprojective algebras $\{ \Pi_d \}_{d \geq 4}$ in each global dimension $d \geq 4$ which fail to be graded coherent. 
\end{thm*}
This produces a counterexample to the conjectures of Minamoto and Bondal in all global dimensions $d \geq 4$, although the counterexamples are somehow isolated and don't seem representative of typical behavior. Since both conjectures hold in global dimension $d \leq 2$, this leaves open the case of $d = 3$.\bigskip 

The proofs of the above theorems appear only in Chapter 6 after the relevant machinery has been setup in Chapter 5. 
\subsection*{Chapter 5}
Chapter 5 forms the technical core of this thesis. Initially motivated by ideas of Chapter 4, this chapter can be seen as a long digression and can be read mostly independently from the previous ones. \bigskip

Given a Koszul Gorenstein algebra $R$, we are lead to consider a natural candidate for a t-structure on $ {\rm D}_{sg}^{\Z}(R) $, namely a natural pair of full subcategories $({\rm D}^{\leq 0}_{sg}(R), {\rm D}^{\geq 0}_{sg}(R))$ with intersection $\mathcal{H}^{\mathsf{lin}}(R) := {\rm D}^{\leq 0}_{sg}(R) \cap {\rm D}^{\geq 0}_{sg}(R) \subseteq {\rm D}_{sg}^{\Z}(R)$ the subcategory of modules with eventually linear minimal projective resolution 
\[
	\mathcal{H}^{\mathsf{lin}}(R) := \{ M \in {\rm D}_{sg}^{\Z}(R) \ | \ \beta_{i,j}(M) = 0 \textup{ for all } i \neq j \textup{ whenever } i \gg 0 \}.
\]
Recall that $\beta_{i,j}(M)$ is the number of generators of degree $j$ in the $i$-th term of the minimal projective resolution $P_* \xrightarrow{\sim} M$, and will be independent of representative of the isomorphism class of $M \in {\rm D}_{sg}^\Z(R)$ so long as $i \gg 0$. Pairs of the form $( {\rm D}^{\leq 0}_{sg}(R), {\rm D}^{\geq 0}_{sg}(R))$ often arise through Koszul duality equivalence as the t-structures pulled back from the standard t-structure attached to the Koszul dual $R^! = {\rm Ext}^*_R(k,k)$. For instance, under the (contravariant) BGG correspondence
\[
	{\rm D}^\Z_{sg}(\wedge V^*)^{op} \cong {\rm D}^b(\coh \mathbb{P}(V)) 
\]
one sees that the standard t-structure on $ {\rm D}^b(\coh \mathbb{P}(V)) $ gives rise to $( {\rm D}^{\leq 0}_{sg}(\wedge V^*), {\rm D}^{\geq 0}_{sg}(\wedge V^*))$ and so induces a (contravariant) equivalence of abelian categories
\[
	\mathcal{H}^{\mathsf{lin}}(\wedge V^*)^{op} \cong \coh \mathbb{P}(V).
\]
However, for a general Koszul Gorenstein algebra $R$ with Koszul dual $R^! = {\rm Ext}_R^*(k, k)$, it isn't a priori clear that ${\rm D}_{sg}^{\Z}(R)$ takes part in an equivalence of the type above, and so whether $({\rm D}^{\leq 0}_{sg}(R), {\rm D}^{\geq 0}_{sg}(R))$ should form a t-structure at all. \bigskip

It is well-known that an appropriate equivalence exists whenever $R$ is Artinian and $R^!$ Noetherian, as in the classical situation of $(R, R^!) = (\wedge V^*, {\rm Sym} V)$. It isn't hard to see that the same hold for $R$ Artinian if we only require that $R^!$ be coherent as a graded algebra. Much more difficult is the extension to pairs of Noetherian Koszul Gorenstein algebras $(R, R^!)$, as is done in by Buchweitz \cite[Appendix]{BEH} in a beautiful tour de force. Buchweitz studies in particular the case where one of $R, R^!$ is commutative, say $R$. Requiring that $R^! = {\rm Ext}^*_R(k, k)$ be Noetherian is then such a strong constraint that it forces $R$ to be given by a complete intersection of quadrics $R = k[ \underline{x}]/( \underline{q})$, by an important theorem of B{\o}gvad-Halperin. No improvement beyond the case of complete intersections of quadrics has been obtained since the work of Buchweitz \cite{BEH} in 1987. 

In this chapter, we will obtain a complete characterisation of the Koszul Gorenstein algebras $R$ for which the natural candidate $( {\rm D}^{\leq 0}_{sg}(R), {\rm D}^{\geq 0}_{sg}(R))$ forms a bounded t-structure on ${\rm D}_{sg}^{\Z}(R)$. Our main results are as follows: 

\begin{thm*}[Theorem \ref{thmA}] Let $R$ be a Koszul Gorenstein algebra. The following are equivalent: 
	\begin{enumerate}[i)]
	\item $({\rm D}^{\leq 0}_{sg}(R), {\rm D}^{\geq 0}_{sg}(R))$ forms a bounded t-structure for $ {\rm D}_{sg}^{\Z}(R) $.
	\item $R$ is absolutely Koszul in the sense of Herzog-Iyengar. 
\end{enumerate}
Furthermore, when these equivalent conditions hold, the natural realisation functor
\[
	\mathsf{real}:{\rm D}^b(\mathcal{H}^{\mathsf{lin}}(R)) \xrightarrow{\cong} {\rm D}_{sg}^{\Z}(R)
\]
is an equivalence of triangulated categories which extends the inclusion on $\mathcal{H}^{\mathsf{lin}}(R)$. 
\end{thm*}

Note that the theorem makes no mention of the Koszul dual $R^! = {\rm Ext}_R^*(k, k)$, and so in particular makes no demand on the structure of the Ext algebra. The role of $R^!$ is relegated to giving a convenient presentation of the abelian category $\mathcal{H}^{\mathsf{lin}}(R)$, as in the next theorem. 

\begin{thm*}[Theorem \ref{thmB}] Let $R$ be a Koszul Gorenstein algebra. If $R$ is absolutely Koszul, then the graded algebra $E = (R^!)^{op} = {\rm Ext}_R^*(k, k)^{op}$ is coherent, and there is a contravariant equivalence of abelian categories
	\[
		\mathcal{H}^{\mathsf{lin}}(R)^{op} \xrightarrow{\cong} \qgr E
	\]
	given by $M \mapsto {\rm Ext}^*_R(M, k)$, where $\qgr E$ is the Serre quotient of the category of finitely presented graded $E$-modules modulo the subcategory of finite length modules. \bigskip 

	\noindent Conversely, if $R$ is Artinian and $E = {\rm Ext}_R^*(k, k)^{op}$ is coherent, then $R$ is absolutely Koszul. 
\end{thm*}

Putting those two theorems together, we deduce the general Bernstein-Gel'fand-Gel'fand correspondence. 
\begin{thm*}[Theorem \ref{thmC}] Let $R$ be an absolutely Koszul Gorenstein algebra with $E = (R^!)^{op} = {\rm Ext}_R^*(k, k)^{op}$. Then there is a contravariant equivalence of triangulated categories
	\[
		{\rm D}^\Z_{sg}(R)^{op} \cong {\rm D}^b(\qgr E)
	\]
	such that the bounded t-structure $( {\rm D}^{\leq 0}_{sg}(R), {\rm D}^{\geq 0}_{sg}(R))$ arises as the pullback of the standard t-structure on the right-hand side. 
\end{thm*}

\subsection*{Chapter 6}
Finally, Chapter 6 consists of assorted corollaries, applications and worked out examples arising out of Chapter 4 and 5. \bigskip

We first spend some time collecting from the literature examples and classes of Koszul Gorenstein algebras which are known to be absolutely Koszul, and so to which our results apply. Of note, all Koszul Gorenstein algebras of codimension $\leq 4$ are automatically absolutely Koszul by results of Avramov-Kustin-Miller and Herzog-Iyengar. Moreover, many interesting projective varieties admit embeddings $X \subseteq \mathbb{P}^n$ with absolutely Koszul Gorenstein homogeneous coordinate ring $R_X$. In particular we take a close look at the cone over an elliptic normal curve $E \subseteq \mathbb{P}^{d-1}$ of degree $d \geq 4$, at the anticanonical model of a smooth del Pezzo $X_d \subseteq \mathbb{P}^d$ of degree $d \geq 4$, and at the canonical embedding $C \subseteq \mathbb{P}^{g-1}$ of a non-hyperelliptic smooth projective curve of genus $g \geq 3$, which is assumed to be cut-out by quadrics. \bigskip 

For the homogeneous coordinate ring $R_{E,d}$ of the elliptic normal curve $E \subseteq \mathbb{P}^{d-1}$, we obtain from the BGG correspondence an equivalence
\[
{\rm D}^b(\coh E) \cong {\rm D}^b(\mathcal{H}^{\mathsf{lin}}(R_{E,d})).
\]
The induced heart $\mathcal{H}^{\mathsf{lin}}(R_{E,d}) \subseteq {\rm D}^b(\coh E) $ is also hereditary and admits a simple description by results of Pavlov \cite{Pavlov}. Recall that the derived category of the elliptic curve $E$ admits a fully faithful action by the braid group ${\rm B}_3$ on three strands, generated by Thomas-Seidel twists along the spherical objects $\OO_E$ and the structure sheaf $k(x)$ of a closed point $x \in E$. For each $d \geq 4$, one can associate a positive braid element $\sigma = \sigma_d \in { \rm B}_3 \subseteq {\rm Aut}({\rm D}^b(\coh E) )$ acting on the derived category of $E$, which one normalises to an antoequivalence $\upgamma := \sigma[-1]$. The full subcategory 
	\[
		\mathcal{H}^{\upgamma}(E) := \{ \F \in {\rm D}^b(\coh E) \ | \ {\rm H}^n(E, \upgamma^j \F) = 0 \textup{ for } n \neq 0 \textup{ whenever } j \gg 0 \}
	\]
	then corresponds to $\mathcal{H}^{\mathsf{lin}}(R_{E,d})$ under the above equivalence. \bigskip 

For the anticanonical model $R_{X_d}$ of the del Pezzo surface $X_d \subseteq \mathbb{P}^{d}$, using a theorem of Happel we obtain that $\mathcal{H}^{\mathsf{lin}}(R_{X_d})$ must be derived equivalent to the representation category $\modsf kQ$ of a finite acyclic quiver or the category of coherent sheaves $\coh \X$ over a weighted projective line in the sense of Geigle-Lenzing. \bigskip 

Afterwards, we then construct a sequence of Artinian Koszul Gorenstein algebras $\{R_n\}_{n \geq 4}$ which fail to be absolutely Koszul, building on work of J.-E. Roos. This is directly applied to construct the counterexamples to the coherence conjectures discussed in Chapter 4. \bigskip 

Lastly, we end this thesis with a discussion of various philosophical points behind the results obtained and lay out some conjectures. 

%% file: chapters/introduction/prelim.tex
\section{Preliminaries}

Throughout this thesis, $k$ will stand for a fixed choice of field. A graded object $X$ in a category $\CC$ consists of a sequence of objects $\{ X_i \}_{i \in \Z}$ in $\CC$. A graded algebra $R$ over $k$ will always stand for a non-negatively graded $k$-module, meaning that $R_i = 0$ for $i < 0$, with products $R_i \otimes R_j \to R_{i+j}$, and abusing notation we will typically conflate $R$ with its direct sum totalisation and write $R = \bigoplus_{i \geq 0} R_i$. We say that $R$ is connected if $R_0 = k$, and standard graded if $R$ is generated in degree one over $R_0$, that is $R = R_0[R_1]$. Similarly we will often write $M = \bigoplus_{i \in \Z} M_i$ for a graded module over a graded algebra. When done this way, elements $r \in R$ or $x \in M$ are understood to be homogeneous, and we write $|r|$, $|x|$ for their degree. Outside of specific parts at the beginning of Chapter 2, graded algebras will be generated in degree one over $R_0$, which will often be $k$ but sometimes a finite dimensional semisimple $k$-algebra. \bigskip

Given a ring $S$, we write $\Modsf S$ for the category of right $S$-modules, and $\modsf S$ for the full subcategory of finitely presented $S$-modules, which is an abelian category whenever $S$ is right Noetherian (more generally right coherent). Similarly we write $\Grmod S$ and $\grmod S$ for the corresponding categories of graded modules. All module-theoretic notions will always implicitly refer to right modules, and we identify the category of left modules with right modules over the opposite ring $S^{op}$.\bigskip 

A complex $X = (X, d)$ in an abelian category $\A$ will mean a graded object $X = \{X_i \}_{i \in \Z}$ equipped with a differential $d_i: X_i \to X_{i-1}$, meaning that $d^2 = 0$. We always write complexes from left to right
\[
	\cdots \to X_{i+1} \xrightarrow{d_{i+1}} X_i \xrightarrow{d_i} X_{i-1} \to \cdots	 
\]and we define its homology by ${\rm H}_i(X) = {\rm ker}(d_i)/{\rm im}(d_{i-1})$. We will often write complexes cohomologically by setting $X^i = X_{-i}$, in which case the differential increases degree. The suspension $X[1]$ of a complex (sometimes written $\Sigma X$) is defined by $(X[1])^i = X^{i+1}$, or equivalently $(X[1])_i = X_{i-1}$, and we define $X[n]$ analogously for any $n \in \Z$. A chain-map $f: X \to Y$ of complexes consists maps $f_i: X_i \to Y_i$ commuting with the differential. A chain-map $f: X \to Y$ is called nullhomotopic if it is of the form $f = d_Y \circ h + h \circ d_X$ for a morphism of graded object $h: X \to Y[-1]$. We say that $f: X \to Y$ is a quasi-isomorphism if the induced map on homology $f_*: {\rm H}_*(X) \xrightarrow{\cong} {\rm H}_*(Y)$ is an isomorphism. The mapping cone of $f: X \to Y$ is defined as the complex
\[
	{\rm Cone}(f)_n := X_{n-1} \oplus Y_n
\]
with differential $\partial(x + y) = f(x) - d(x) + d(y)$. \bigskip

Given two graded objects $X, Y$ in $\A$, their graded Hom abelian group ${\rm Hom}(X, Y)$ is defined by 
\[
	{\rm Hom}^n(X, Y) = \prod_{i \in \Z} {\rm Hom}(X^i, Y^{i+n}). 
\]
When $X, Y$ are complexes, ${\rm Hom}(X, Y)$ inherits a differential by $\partial(f) = d_Y \circ f - (-1)^{|f|} f \circ d_X$. When $\A$ admits a tensor product, the tensor product $X \otimes Y$ is the graded object with components 
\[
	(X\otimes Y)^n = \bigoplus_{p+q=n} X^p \otimes Y^q.
\]
This inherits the structure of a complex from $X, Y$ by setting $d_{X \otimes Y} = d_X \otimes 1 + 1 \otimes d_Y$. \bigskip 

We define the following categories:
\begin{enumerate}[i.]
	\item $C(\A)$ is the category of complexes over $\A$, with morphisms given by chain-maps. 
	\item $\mathcal{K}(\A)$ is the homotopy category of complexes over $\A$, which is obtained from $C(\A)$ by quotienting out nullhomotopic morphisms.
	\item ${\rm D}(\A)$ is the derived category of $\A$, obtained from $\mathcal{K}(\A)$ by inverting quasi-isomorphisms (see \cite{Weibel} or \cite{Huy} for details). 
	\end{enumerate}
	We let ${\rm D}^*(\A)$ for $* = \{-, +, b\}$ denote the full subcategories of ${\rm D}(\A)$ consisting of complexes with right bounded cohomology, left bounded cohomology or bounded cohomology, respectively. When $R$ is a ring, we write ${\rm D}(R) := {\rm D}(\Modsf R)$. When $R$ is right Noetherian (or right coherent) we denote ${\rm D}^b(R) := {\rm D}^b(\modsf R)$, and we write $ {\rm D}(X) := {\rm D}(\qcoh X)$ and ${\rm D}^b(X) := {\rm D}^b(\coh X)$ for any scheme $X$. \bigskip 

	In any triangulated category $\T$, for any set of objects $S \subseteq \T$ we define the thick closure $\thick(S) \subseteq \T$ to be the smallest triangulated category of $\T$ containing $S$ which is closed under finite direct sums and summands. When $\T$ is closed under arbitrary direct sums, we define the localising closure $\mathsf{loc}(S) \subseteq \T$ to be the smallest triangulated category of $\T$ containing $S$ closed under arbitrary direct sums; note that $\mathsf{loc}(S)$ is then automatically closed under direct summands by \cite[Prop. 1.6.8]{NeemanBook}. We write $[1]$ for the suspension in a triangulated category, and ${\rm Ext}^n_\T(X, Y) := {\rm Hom}_\T(X, Y[n])$. We say that a $k$-linear triangulated category is Ext-finite if $\dim_k \bigoplus_{n \in \Z}{\rm Ext}^n_\T(X, Y) < \infty$ for any $X, Y \in \T$. All triangulated categories in this thesis will be $k$-linear over our base field $k$, and typically will be Ext-finite.\bigskip 

	A t-structure $t = (\T^{\leq 0}, \T^{\geq 0})$ on $\T$ is a pair of full subcategories satisfying three defining axioms. Define $\T^{\leq n} = \T^{\leq 0} \circ [-n]$ and $\T^{\geq n} = \T^{\geq 0} \circ [-n]$. Then $t$ forms a t-structure on $\T$ if
	\begin{enumerate}[T1.]
		\item We have containments $\T^{\leq 0} \subseteq \T^{\leq 1}$ and $\T^{\geq 1} \subseteq \T^{\geq 0}$.
		\item We have ${\rm Hom}_\T(\T^{\leq 0}, \T^{\geq 1}) = 0$.
		\item For every object $X \in \T$, there are objects $X^{\leq 0} \in \T^{\leq 0}$ and $X^{\geq 1} \in \T^{\geq 1}$ fitting in a distinguished triangle
			\[
				X^{\leq 0} \to X \to X^{\geq 1} \to X^{\leq 0}[1].
			\]
\end{enumerate}
When $t$ forms a t-structure, there are truncation functors $\tau^{\leq n}: \T \to \T^{\leq n}$ and $\tau^{\geq n}: \T \to \T^{\geq n}$ which are adjoint to the respective inclusions 
\begin{align*}
	&\tau^{\leq n}: \T \leftrightarrows \T^{\leq n}: i_n\\
	&j_n: \T^{\geq n} \leftrightarrows \T : \tau^{\geq n}
\end{align*}
such that $\tau^{\leq 0}X = X^{\leq 0}$ and $\tau^{\geq 1}X = X^{\geq 1}$ above, and with the maps in the distinguished triangle coming from the counit and unit maps of the respective adjunctions. \bigskip

The triangulated category $\T = {\rm D}(\A)$ has a standard t-structure $({\rm D}^{\leq 0}, {\rm D}^{\geq 0})$, given by complexes with cohomology supported in degree $\leq 0$ and $\geq 0$, respectively. The truncation functors are given by 
\[
	\xymatrix{
		\tau^{\leq n}X: & \cdots \ar[r] & X^{n-1} \ar[r]^-{d_{n-1}} & {\rm ker}(d_n) \ar[r] & 0 \ar[r] & \cdots \\
		X: & \cdots \ar[r] & X^{n-1} \ar[r]^-{d_{n-1}} & X^n \ar[r]^-{d_n} & X^{n+1} \ar[r] & \cdots \\
		\tau^{\geq n}X: & \cdots \ar[r] & 0 \ar[r] & {\rm coker}(d_{n-1}) \ar[r]^-{d_n} & X^{n+1} \ar[r] & \cdots \\
	}
\] 

A differential graded (dg) algebra over $k$ is a complex $A = (A, d)$ with an associative multiplication $m: A^{\otimes 2} \to A$ which is a chain-map, or equivalently $d$ satisfies the Leibniz rule $d(ab) = (da)b + (-1)^{|a|}a(db)$. Every graded algebra can be thought of as a dg algebra with trivial differential.  A morphism of dg algebras $f: A \to B$ is a homomorphism of algebras which respects the differential, and a quasi-isomorphism of dg algebras $f: A \xrightarrow{\sim} B$ is a morphism which is a quasi-isomorphism of underlying complexes. We will say that two dg algebras $A$, $B$ are quasi-isomorphic if they are connected to each other by a zig-zag of quasi-isomorphisms, and denote this by $A \simeq B$. Moreover, $A$ is \emph{formal} if it is quasi-isomorphic to its cohomology algebra ${\rm H}^*(A)$, that is $A \simeq {\rm H}^*(A)$. \bigskip 

Finally, we will make some usage of the theory of Koszul algebras. Let $A = A_0 \oplus A_1 \oplus \cdots $ be a locally finite graded $k$-algebra generated by $A_1$ over $A_0$, and we assume that $A_0$ is semisimple. Let $P_* \xrightarrow{\cong} A_0$ be the minimal graded projective resolution of $A_0 = A/A_{\geq 1}$ over $A$. We say that $A$ is \emph{Koszul} if $P_i$ is generated in degree $i$ for all $i \geq 0$. \bigskip

When $A$ is Koszul, we denote by $A^! = {\rm Ext}^*_A(A_0, A_0)$ its Koszul dual algebra. One can see that $A^!$ is also Koszul with degree zero part $A^!_0 = A_0$ and that ${\rm Ext}^*_{A^!}(A^!_0, A^!_0) \cong A$. While we assume some familiarity with the class of Koszul algebras, perhaps at the level of \cite{BGS}, the two above facts will suffice for the entirety of this thesis. 
	

%% file: chapters/introduction/background.tex
\section{Background: Maximal Cohen-Macaulay modules over Gorenstein rings}
In this section, we review classical definitions and characterizations of MCM modules, and give an exposition of Buchweitz's manuscript. We refer the reader to \cite{BHBook, YosBook, KunzBook, BuchweitzManuscript} for complete references.

\subsection*{MCM modules over CM local rings}
Let $R = (R, \m, k)$ here denote a commutative Noetherian local ring (respectively, graded local ring) of finite Krull dimension $d$. All definitions and results have natural extensions to the case of graded modules over graded rings which this thesis is concerned with, and we will spell out the graded version of a result when appropriate. For an $R$-module $M$, a sequence of elements $x_1, \dots, x_n$ in $\m$ is a regular $M$-sequence if $x_{i+1}$ is a non zero divisor on $M/(x_1, \dots, x_i)M$ for all $i = 0, 1, \dots, n-1$.
\begin{defn}
	The depth of an $R$-module $M$, denoted $\depth M$, is the maximal length of regular $M$-sequences. 
\end{defn}
The depth of a module is always upper bounded by the Krull dimension of its support, that is $\depth M \leq \dim M$.
\begin{defn}
	A module $M$ is Cohen-Macaulay if $\depth M = \dim M$. A module $M$ is maximal Cohen-Macaulay if furthermore $\dim M = \dim R$, that is $\depth M = d$. 
\end{defn}
\begin{defn}
	A ring $R$ is Cohen-Macaulay if $R$ is maximal Cohen-Macaulay as an $R$-module. That is, $\m$ contains a regular $R$-sequence of length $d$. The ring $R$ is regular if $\m$ is actually generated by a regular sequence (then of length $d$).
\end{defn}

\begin{prop}\label{MCMregular} Let $M$ be a finite MCM module over a regular local ring $R$. Then $M$ is free.
\end{prop}
\begin{proof}
	Since $R$ is regular $\pdim M < \infty$ and the Auslander-Buchsbaum formula $\pdim M + \depth M = \depth R$ gives $\pdim M = 0$, and $M$ is free since $R$ is local.
\end{proof}
It follows that we may think of non-free MCM modules as measuring the singularities of $R$ in some way. Depth and the Cohen-Macaulay property are best recognized by cohomological criteria. To this end, let $X = \spec R$, and for any $R$-module $M$, by abuse of notation we denote by $M$ the associated quasi-coherent sheaf on $X$. Let $\Gamma_\m = \Gamma_\m(X, -): \qcoh X \to \Modsf R$ be the functor of global sections with support in $\m$, meaning $$\Gamma_\m(M) = \Gamma_\m(X, M) = \{ x \in M \ | \ \m^r \cdot x = 0 \textup{ for some } r \gg 0 \}.$$ $\Gamma_\m$ is left exact with total derived functor ${\rm R}\Gamma_\m$, and we denote by ${\rm H}^i_\m(-) = {\rm R}^i\Gamma_\m$ the $i$-th local cohomology functor. Local cohomology detects depth.

\begin{prop}[\cite{BHBook}]\label{depthlocalcoh} Let $M$ be a finite $R$-module. We have:
	\begin{enumerate}[i)]
	\item $\depth_R M = \inf \{i \ | \ {\rm H}^{i}_\m(M) \neq 0\}$.
	\item $M$ is Cohen-Macaulay if and only if ${\rm H}^i_\m(M) = 0$ for $i \neq \dim M$.
\end{enumerate}
\end{prop}

In the graded local case, say when $R = \bigoplus_{n \geq 0} R_n$ is a standard graded connected $k$-algebra and $\m = R_{+}$, we can reinterpret local cohomology of graded $R$-modules in terms of the projective $k$-scheme $\proj R$. Here let $X = \proj R$ and $ \widehat{X} = \spec R$ is the affine cone over $X$. Let $U \subseteq \widehat{X} $ be the punctured spectrum $ U = \widehat{X} \setminus \m $. For a graded module $M$, by abuse of notation we denote by $M$ the quasi-coherent sheaf on $ \widehat{X} $ and $\widetilde{M}$ denote the quasi-coherent sheaf on $X$.
\[
	0 \to \Gamma_\m(\widehat{X}, M) \to \Gamma( \widehat{X} , M) \xrightarrow{\mathsf{res}} \Gamma(U, M) \to 0
\]
giving rise to distinguished triangles
\[
	{\rm R}\Gamma_\m( \widehat{X} , M) \to {\rm R}\Gamma( \widehat{X} , M) \xrightarrow{\mathsf{res}} {\rm R}\Gamma(U, M) \to {\rm R}\Gamma_\m( \widehat{X}, M)[1].
\]
Quasi-coherent sheaf cohomology vanishes over the affine scheme $ \widehat{X} $. By \cite{KunzBook}, we have ${\rm R}\Gamma(U, M) \cong {\rm R}\Gamma_*(X, \widetilde{M} ) := \bigoplus_{n \in \Z} {\rm R}\Gamma(X, \widetilde{M} (n))$. Putting these together, we have:
\begin{prop} Let $M$ be a graded $R$-module.
	 \begin{enumerate}[i)]
	 \item There is an exact sequence of graded modules
		 \[
			 0 \to \Gamma_\m(M) \to M \xrightarrow{\mathsf{res}} \Gamma_*(X, \widetilde{M}) \to {\rm H}^1_\m(M) \to 0.
		 \]
	 \item We have natural isomorphisms of graded $R$-modules ${\rm H}^i_*(X, \widetilde{M}) \cong {\rm H}^{i+1}_\m(M)$ for all $i \geq 1$.
 \end{enumerate}
 \end{prop}

 \begin{cor}\label{ACMsheaf} Let $M$ be a finite graded $R$-module. The following are equivalent:
	 \begin{enumerate}[i)]
	 \item $M$ is maximal Cohen-Macaulay.
	 \item The map $\mathsf{res}$ is an isomorphism, and we have ${\rm H}^i(X, \widetilde{M}(n)) = 0$ for all $n \in \Z$ for all $i \neq 0, \dim X$.
 \end{enumerate}
 \end{cor}
 A coherent sheaf $\F$ corresponding to a finite MCM module $M$ is called Arithmetically Cohen-Macaulay (ACM).

 \begin{defn} Let $R$ be a local CM ring of dimension $d$. A canonical module $\omega_R$ for $R$ is a finite $R$-module such that
	 \[
		 {\rm Ext}_R^i(k, \omega_R) = \begin{cases} 0 & i \neq d \\ k & i=d. \end{cases}
	 \]
	 Such a module $\omega_R$ is unique up to isomorphism whenever it exists.
 \end{defn}

 \begin{prop}[Local Duality {\cite{BHBook}}]\label{localduality} Let $R$ be a local CM ring with a canonical module $\omega_R$. Let $E_R(k)$ be the injective hull of $k$. For any finite $R$-module $M$, we have natural isomorphisms of $R$-modules
	\[
		{\rm H}_\m^i(M) = {\rm Hom}_R({\rm Ext}^{d-i}_R(M, \omega_R), E_R(k)).
\]
\end{prop}

\begin{prop}[Graded Local Duality {\cite{BHBook}}]\label{grlocalduality} Let $R$ be a graded connected CM $k$-algebra with graded canonical module $\omega_R$. For any finite graded $R$-module $M$, we have natural isomorphisms of graded $R$-modules
	\[
		{\rm H}_\m^i(M) = {\rm Hom}_k({\rm Ext}^{d-i}_R(M, \omega_R), k).
\]
\end{prop}

\begin{prop}\label{MCMlocalcoh} A finite $R$-module $M$ is Cohen-Macaulay if and only if ${\rm Ext}^{d-i}_R(M, \omega_R) = 0$ for $i \neq \dim M$. In particular $M$ is MCM if and only if ${\rm Ext}^n_R(M, \omega_R) = 0$ for all $n > 0$.
\end{prop}

\begin{defn}
	Let $R$ be a CM local ring (respectively graded local). We say that $R$ is Gorenstein if $\omega_R := R$ is a canonical module  (respectively $\omega_R := R(a)$ for some $a \in \Z$). In the graded case, the integer $a$ is called the $a$-invariant.\footnote{The integer $-a$ is also called the Gorenstein parameter, and sometimes denoted by the same letter $a$. We will not follow this convention. }
\end{defn}

Note that by proposition \ref{grlocalduality}, we have $a = \max \{i \ | \ {\rm H}^d_\m(R)_i \neq 0 \}$. Gorenstein rings are ubiquitous \cite{Bass}, and have many characterizations and standard examples.

\begin{prop}[\cite{BHBook}]\label{fininjdim} Let $R$ be a local Noetherian ring of dimension $d$. The following are equivalent:
	\begin{enumerate}[i)]
	\item $R$ is Gorenstein.
	\item The injective dimension of $R_R$ is finite (then equal to $d$).
	\item $R/(\bf{x})$ is Gorenstein for some (and then all) regular $R$-sequence ${\bf x} = (x_1, \dots, x_n)$.
\end{enumerate}
\end{prop}

We say that $R$ is a complete intersection if its completion $\widehat{R}$ at $\m$ is isomorphic to $\widehat{R} \cong Q/(\bf{x})$ for $Q$ a regular local ring and ${\bf x} = (x_1, \dots, x_c) \subseteq \m_Q$ a regular sequence. We have strict implications:
\[
	\textup{regular} \implies \textup{complete intersection} \implies \textup{Gorenstein} \implies \textup{Cohen-Macaulay}.
\]
\begin{exmp}
	Let $R$ be of Krull dimension zero. The following are equivalent:
	\begin{enumerate}[i)]
	\item $R$ is Gorenstein.
	\item $R$ is self-injective.
	\item $R_R$ has simple socle.
\end{enumerate}
Moreover when $R$ is a graded connected $k$-algebra, $\omega_R = DR \cong R(a)$ with $a$ the degree of the socle element.
\end{exmp}

\begin{exmp}[Watanabe \cite{Wata}] 
	Let $V$ be a finite-dimensional vector space over an algebraically closed field $k$ and let $G \leq {\rm GL}(V)$ be a finite group with $\mathsf{char}\ k \nmid |G|$. Assume that $G$ is \emph{small}, meaning containing no pseudo-reflection. Then the invariant ring $R = k[V]^G$ is Gorenstein if and only if $G \leq {\rm SL}(V)$.
\end{exmp}

\subsection*{The Stable Derived Category of a Gorenstein Ring}
From now on, we will take rings $R$ to be two-sided Noetherian and modules are taken to be finitely generated, unless specified otherwise. Taking a cue from proposition \ref{fininjdim}, we call a possibly noncommutative, two-sided Noetherian ring $R$ Gorenstein (or Iwanaga-Gorenstein) if \begin{enumerate}[i)] 
\item $\idim(R_R) < \infty$
\item $\idim(_R R) < \infty$
	\end{enumerate}
	in which case both injective dimensions are equal, say to some integer $d$. We call $d$ the Gorenstein dimension of $R$. Gorenstein rings admit a natural duality given by 
	\[
		(-)^{\vee} := {\rm RHom}_{R}(- , R): {\rm D}^b(R)^{op} \to {\rm D}^b(R^{op}).
	\]
	\begin{lem}[{\cite[Lemma 5.3]{BurkeStev}}] The functor $(-)^{\vee}$ is a duality. That is, the natural map $X \to X^{\vee \vee}$ is an isomorphism for all $X \in {\rm D}^b(R)$. 
	\end{lem}

	\begin{defn} Let $R$ be Gorenstein. An $R$-module $M$ is maximal Cohen-Macaulay if $ {\rm Ext}^n_R(M, R) = 0$ for $n > 0$. Equivalently, $M^{\vee}$ is a module in $ {\rm D}^b(R^{op}) $.
	\end{defn}
	Over $R$ Gorenstein, for any MCM module $M$ it follows that $M^\vee = M^* := {\rm Hom}_R(M, R)$, and since $(M^*)^{\vee} = M^{\vee \vee} \cong M$ is an $R$-module, $M^*$ is an MCM module over $R^{op}$ and $M^{**} \cong M$. We denote by ${\rm MCM}(R)$ the full subcategory of MCM modules over a Gorenstein ring $R$, (respectively ${\rm MCM}^\Z(R)$ the category of graded MCM modules over a $\Z$-graded Gorenstein ring $R$). Collecting some standard properties, we have: 
	\begin{prop}Let $R$ be Gorenstein. The category ${\rm MCM}(R)$ has the following properties:
		\begin{enumerate}[i)]
		\item The category ${\rm MCM}(R)$ is closed under sums, summands and extensions. 
		\item MCM modules are closed under taking duals. 
		\item MCM modules are reflexive, that is $M^{**} \cong M$.
		\item Projective $R$-modules are MCM.
		\item The projective $R$-modules are injective objects in the category ${\rm MCM}(R) \subseteq \modsf R$. 
	\end{enumerate}
	\end{prop}

	We next define the stable module category. 
	\begin{defn}
		Let $R$ be any ring. The \emph{projectively stable}, or \emph{stable}, module category $\underline{\Modsf} R$ (resp. $\underline{\modsf} R$) has for objects all $R$-modules (resp. finitely generated $R$-modules), with morphisms given by $$\underline{ {\rm Hom}}_R(M, N) = {\rm Hom}_R(M, N)/\mathcal{P}(M, N)$$ where $\mathcal{P}(M, N)$ is the ideal of morphisms factoring through a projective object. 
	\end{defn}	
	\begin{defn} 
		Let $R$ be two-sided Noetherian and Gorenstein. The stable category of MCM modules $ \underline{{\rm MCM}}(R)$ is the full subcategory of $\underline{\modsf} R$ consisting of MCM modules. 
	\end{defn}

	Working over any Noetherian ring $R$, given a finite projective presentation $P_1 \xrightarrow{\partial_1} P_0 \to M$ of $M \in \grmod R$, we define the first syzygy $\Omega(M) = {\rm im}(\partial_1)$. The first syzygy of $M$ depends on the choice of presentation, but it is well-known (Schanuel's Lemma) that any two choice of presentations give rise to stably isomorphic syzygies, and that one obtains a well-defined functor $\Omega: \underline{\modsf} R \to \underline{\modsf} R$. When $R$ is Gorenstein, the functor $\Omega$ preserves the subcategory of MCM modules, and since $R$ has finite injective dimension, any syzygy module $\Omega^n(N)$ is MCM for $n \gg 0$ (note that $n \geq d$ is enough). By a fundamental result of Buchweitz, all MCM modules are actually of this form.

	\begin{prop}[Buchweitz] The following are equivalent over $R$ Gorenstein:
		\begin{enumerate}[i)]
		\item $M$ is an MCM module.
		\item For every $n \geq 0$, there is an $R$-module $N$ such that $M \cong \Omega^n(N)$.
	\end{enumerate}
	\end{prop}
	\begin{proof} Let $P \xrightarrow{\sim} M$ be a projective resolution and $Q \xrightarrow{\sim} M^*$ be a projective resolution of the dual. Since $M$ is MCM, so is $M^*$ and by cohomology vanishing ${\rm Ext}_{R^{op}}^{> 0}(M^*, R) = 0$ we have a quasi-isomorphism $M = M^{**} \xrightarrow{\sim} Q^*$. Composing the quasi-isomorphisms $P \xrightarrow{\sim} M = M^{**} \xrightarrow{\sim} Q^*$, we obtain an acyclic complex of projectives by taking the (shifted) cone	$C_* = {\rm Cone}(P \to Q^*)[-1]$
		\[
			\xymatrix@R0.6pc{\cdots \ar[r] & P_1 \ar[r] & P_0 \ar[rr] \ar[dr] && Q_0^* \ar[r] & Q_1^* \ar[r] & \cdots \\
			&&& M \ar[ur] }
		\]
		Truncating the complex $C_*$ sufficiently far to the right reveals $M$ as the $n$-th syzygy module of $N = {\rm coker}(C_{-n+1} \to C_{-n})$.
	\end{proof}

	Extracting a definition from the above proof, we set:
	\begin{defn}
		Let $R$ be Gorenstein. A complete resolution $C_* \to M$ of an MCM module $M$ is an infinite acyclic complex of finite projectives whose non-negative truncation resolves $M$. 
	\end{defn}

	As for projective resolutions, complete resolutions are unique up to homotopy and, when $R$ is local (or graded local), admit a minimal model which is unique up to non-canonical isomorphism. Let us write ${\rm proj}\ \! R \subseteq \modsf R$ for the full subcategory of finitely generated projectives. We denote by $\mathcal{K}_{ac}({\rm proj}\ \! R)$ the homotopy category of complete resolutions, or equivalently the homotopy category of acyclic complexes of finitely generated projectives. 
	\begin{prop}[Buchweitz {\cite{BuchweitzManuscript}}]\label{completeres}
		Let $R$ be Gorenstein. The functor $C_* \mapsto {\rm coker}(C_1 \to C_0)$ gives rise to an equivalence of categories $\mathcal{K}_{ac}({\rm proj}\ \! R) \xrightarrow{\cong} \underline{{\rm MCM}}(R)$. The inverse sends $M$ to its complete resolution.
	\end{prop}
	\subsubsection*{The triangulated structure on $ \underline{{\rm MCM}}(R)$}
	For the remainder of these sections, we will work entirely over a Gorenstein ring $R$. \bigskip 

	The homotopy category $\mathcal{K}_{ac}({\rm proj}\ \! R)$ is naturally triangulated, and we can pull back the triangulated structure onto $ \underline{{\rm MCM}}(R)$ under the equivalence of Prop. \ref{completeres}. Let us describe some of its main features: 
	\begin{enumerate}[i)]
	\item Suspension: We have $M[1] = \Sigma M = {\rm cosyz}_R(M)$ with inverse $M[-1] = \Omega M = {\rm syz}_R(M)$:
		\[
			\xymatrix@C1pc@R.6pc{\dots \ar[r] & C_1 \ar[rr] \ar[dr] && C_0 \ar[rr] \ar[dr] && C_{-1} \ar[rr] \ar[dr] && C_{-2} \ar[r] & \dots \\
		&& \Omega M \ar[ur] && M \ar[ur] && \Sigma M \ar[ur]}
		\]
	\item Distinguished triangles: The distinguished triangles in $ \underline{{\rm MCM}}(R)$ are the images of short exact sequences of MCM modules in ${\rm MCM}(R)$.
	\item Mapping cones: Given a map $f: M \to M'$, first embed $\iota: M \hookrightarrow C_{-1}$ into a projective module and define ${\rm Cone}(f) := {\rm Coker}(\iota \oplus f: M \hookrightarrow C_{-1} \oplus M')$. We have a distinguished triangle
		\[
			\xymatrix@R.8pc{ && C_{-1} \oplus M' \ar[dr] \\
				& M \ar[r]^{f} \ar[ur]^{\iota \oplus f} & M' \ar@{=}_{\rotatebox[origin=c]{90}{$\sim$}}[u] \ar[r] & {\rm Cone}(f) \ar[dr] \ar[r] & \Sigma M \\
			0 \ar[ur] & & && 0
}
		\]
		where the map ${\rm Cone}(f) \to \Sigma M$ is induced from the quotient $C_{-1} \twoheadrightarrow \Sigma M$.
\end{enumerate}

\subsubsection*{Tate cohomology} Every triangulated category carries a natural cohomology theory by taking Hom into a fixed object. In the present context, this gives rise to Tate cohomology. Slightly more generally, we define:

\begin{defn} Let $R$ be Gorenstein, let $N$ be an $R$-module and $M$ an MCM module, with complete resolution $C_* \to M$. For $n \in \Z$, we define the $n$-th Tate cohomology group of $M$ with coefficients in $N$ by
	\[
		\underline{{\rm Ext}}^n_R(M, N) := {\rm H}^n {\rm Hom}_R(C_*, N).
	\]
\end{defn}

\begin{prop}[{\cite{BuchweitzManuscript}}]
	Tate cohomology has the following properties.
	\begin{enumerate}[i)]
	\item We can calculate $ \underline{{\rm Ext}}_R^n(M, N)$ via:
		\[
			\underline{{\rm Ext}}^n_R(M, N) = \begin{cases} {\rm Ext}^n_R(M, N) & n \geq 1 \\
				\underline{ {\rm Hom}}_R(M, N) & n = 0 \\
				{\rm Tor}^R_{-1-n}(N, M^*) & n \leq -2
			\end{cases}
		\]
		and we have a short exact sequence
		\[
			0 \to \underline{{\rm Ext}}^{-1}_R(M, N) \to N \otimes_R M^* \xrightarrow{ev} {\rm Hom}_R(M, N) \to \underline{{\rm Ext}}^0_R(M, N) \to 0.
		\]
	\item We have $\underline{{\rm Ext}}^{n+1}_R(\Sigma M, N)  = \underline{{\rm Ext}}^n_R(M, N) = \underline{{\rm Ext}}_R^{n-1}(\Omega M, N)$.
	\item Let $N$ be a perfect module. Then $ \underline{{\rm Ext}}^n_R(M, N) = 0$ for all $n \in \Z$.
	\item Let $N$ be an MCM module with complete resolution $D_*$. Then the natural map $D_* \to N$ induces a quasi-isomorphism of Hom complexes
		\[
			{\rm Hom}_R(C_*, D_*) \xrightarrow{\sim} {\rm Hom}_R(C_*, N).
		\]
		Hence we have natural isomorphisms ${\rm H}^n {\rm Hom}_R(C_*, D_*) \cong \underline{{\rm Ext}}^n_R(M, N)$. In particular $${\rm H}^0 {\rm Hom}_R(C_*,D_*)  \cong \underline{ {\rm Hom}}_R(M, N).$$
\end{enumerate}
\end{prop}

\subsubsection*{Stable derived categories and Buchweitz's equivalence}
Next, still over a fixed Gorenstein ring $R$, consider the bounded derived category $ {\rm D}^b(R)$ and its subcategory of perfect complexes ${\rm D}^{\perf}(R) = \mathsf{thick}(R)$, meaning complexes quasi-isomorphic to a bounded complex of projective modules.
	\begin{defn}
		The \emph{stable derived category} of $R$, also called \emph{singularity category}, is the Verdier quotient
		\[
			{\rm D}_{sg}(R) = {\rm D}^b(R)/\mathsf{thick}(R).
		\]
		of the bounded derived category of $R$ by the subcategory of perfect complexes.
	\end{defn}
	The composition of the natural embedding and projection factors as
\[
	\xymatrix@R1pc{ {\rm MCM}(R) \ar[r] \ar[d] & {\rm D}^b(R) \ar[r] & {\rm D}_{sg}(R) \\
	 \underline{{\rm MCM}}(R) \ar[urr] && }
\]
\begin{thm}[Buchweitz {\cite{BuchweitzManuscript}}]\label{buchweitzmainthm}
	Let $R$ be Gorenstein. The induced functor $\underline{{\rm MCM}}(R) \xrightarrow{\cong} {\rm D}_{sg}(R)$ is an equivalence of triangulated categories.
	\end{thm}

Implicit in this theorem is the existence of an exact functor
\[
	(-)^{st}: {\rm D}^b(R) \twoheadrightarrow {\rm D}_{sg}(R) \cong \underline{{\rm MCM}}(R)
\]
which we call MCM approximation, or \emph{stabilization}. We can describe it as follows: let $F \in {\rm D}^b(R) $ be a complex with bounded cohomology and take a projective resolution $P_* \xrightarrow{\sim} F$. Since $P_*$ has bounded cohomology, its tail $P_{\geq n}$ is exact and resolves an MCM module for $n \gg 0$, say $M = {\rm Coker}(P_{n+1} \to P_n)$. Now, note that we can realise $P_*$ as the shifted cone $P_* = {\rm Cone}(P_{\geq n} \xrightarrow{\partial_n} P_{\leq n-1}[1])[-1]$, which gives rise to a distinguished triangle in $ {\rm D}^b(R) $
\[
	\xymatrix@R.8pc{P_* \ar[r] \ar@{=}^{\rotatebox[origin=c]{90}{$\sim$}}[d] & P_{\geq n} \ar@{=}^{\rotatebox[origin=c]{90}{$\sim$}}[d]  \ar[r] & P_{\leq n-1}[1] \ar@{=}[d] \ar[r] & P_*[1] \ar@{=}^{\rotatebox[origin=c]{90}{$\sim$}}[d]   \\
	F \ar[r] & M[n] \ar[r] & P_{\leq n-1}[1] \ar[r] & F[1]}
\]
with $P_{\leq n-1}$ perfect, and so $F^{st} = M[n] = {\rm cosyz}_R^n(M)$ in $ \underline{{\rm MCM}}(R)$. When $F = N$ is a module in $ {\rm D}^b(R) $, we have a stronger statement:
\begin{prop}[{\cite{BuchweitzManuscript}}] Let $R$ be Gorenstein and let $N \in \grmod R$. There is a short exact sequence of $R$-modules
	\[
		0 \to P \to N^{st} \to N \to 0
	\]
	$P$ perfect and $N^{st}$ MCM, unique up to stable isomorphism.
\end{prop}

\begin{cor} Tate cohomology is representable. That is, for any $R$-module $N$, we have natural isomorphisms
	\[
		\underline{{\rm Ext}}_R^n(M, N) \cong \underline{{\rm Hom}}_R(M, N^{st}[n])
	\]
	for all $M$ MCM and $n \in \Z$.
\end{cor}
In particular Tate cohomology forms a cohomological functor on $ \underline{{\rm MCM}}(R)$, in that any distinguished triangle gives rise to a long exact sequence of Tate cohomology groups.

\subsection*{Auslander-Reiten-Serre duality and Almost-Split sequences}
A standard concern of representation theory is the classification of indecomposables in various settings. This mostly make sense only in the context of Krull-Schmidt categories. Recall that a noncommutative ring $B$ is local if $B/ \mathsf{rad} B$ is a division ring. 

\begin{defn}
	Let $\CC$ be an additive $k$-linear category. We say that $\CC$ is Krull-Schmidt if every object has a direct sum decomposition $X = X_1 \oplus \cdots \oplus X_n$ into finitely many indecomposable objects, and indecomposable objects have local endomorphism rings. 
\end{defn}

Finite decomposition into indecomposables are essentially unique in a Krull-Schmidt category, in that we have $X \cong \bigoplus_{i = 1}^r X_i^{\oplus e_i}$ for unique indecomposables $\{X_i\}$, pairwise non-isomorphic, and unique multiplicities $e_i \in \N$. \bigskip

As this thesis is principally concerned with graded MCM modules, we mention the following, first noted by Iyama and Takahashi: 

\begin{prop}[Iyama-Takahashi {\cite[Prop. 1.4]{IT}}]\label{KSMCM} Let $R = R_0 \oplus R_1 \oplus \cdots $ be a non-negatively graded Gorenstein $k$-algebra, locally finite over $k$. Then ${\rm MCM}^{\Z}(R)$ and $ \underline{{\rm MCM}}^{\mathbb{Z}}(R)$ are Krull-Schmidt categories. 
\end{prop} 
\begin{proof} 
	Any finitely generated graded module $M$ over a Noetherian ring decomposes into finitely many indecomposables, and have finite dimensional endomorphism algebras ${\rm End}_{gr R}(M)$ by local finiteness. When $M$ is indecomposable, the finite dimensional algebra ${\rm End}_{gr R}(M)$ has for only idempotents $e = 0,1$, and thus must be local. This shows that ${\rm MCM}^{\Z}(R)$ is Krull-Schmidt. The category $ \underline{{\rm MCM}}^{\mathbb{Z}}(R)$ then inherits the Krull-Schmidt property from ${\rm MCM}^{\Z}(R)$ since idempotents in ${\rm \underline{End}}_{gr R}(M)$ come from idempotents in ${\rm End}_{gr R}(M)$, by idempotent lifting for Artinian rings. 
\end{proof} 

When $R$ is as in Prop. \ref{KSMCM}, so that ${\rm MCM}^{\Z}(R)$ and $ \underline{{\rm MCM}}^{\mathbb{Z}}(R)$ are Krull-Schmidt, we note that every graded MCM module $M$ has a unique decomposition $M = [M] \oplus F$ with $F$ projective and $[M]$ without projective summands. Moreover, when $X \in {\rm MCM}^{\Z}(R)$ is indecomposable non-projective, the kernel of ${\rm End}_{gr R}(X) \twoheadrightarrow {\rm \underline{End}}_{gr R}(X)$ is contained in the radical of the local ring ${\rm End}_{gr R}(X)$. One deduces that any two indecomposable non-projective MCM modules $X, Y$ which are isomorphic in the stable category $ \underline{{\rm MCM}}^{\mathbb{Z}}(R)$ must also be isomorphic in ${\rm MCM}^{\Z}(R)$. Comparing indecomposable summands of a general MCM module, we see that two graded MCM modules $M, M'$ in $ \underline{{\rm MCM}}^{\mathbb{Z}}(R)$ are stably isomorphic if and only if $[M]$ and $[M']$ are isomorphic in ${\rm MCM}^{\Z}(R)$. In this sense, the classification of MCM modules reduces at once to the classification of indecomposables in $\underline{{\rm MCM}}^{\mathbb{Z}}(R)$. \bigskip

We now review standard background that is common to algebraic geometry, commutative algebra and representation theory of Artin algebras. 

\begin{defn}
Let $\mathcal{T}$ be a triangulated Hom-finite $k$-linear category. A Serre functor for $\mathcal{T}$ is an exact autoequivalence $\mathbb{S}: \mathcal{T} \to \mathcal{T}$ equipped with natural isomorphisms
\[
	{\rm Hom}_{\mathcal{T}}(X, \mathbb{S}(Y)) \cong {\rm D}{\rm Hom}_{\mathcal{T}}(Y, X).
\]
When they exist, Serre functors are unique up to isomorphism.
\end{defn}

When $X$ is a smooth projective variety over $k$, the functor $\mathbb{S}_X(-) = - \otimes_{\mathcal{O}_X} \omega_X  [\textup{dim } X]$ is a Serre functor for ${\rm D}^b(X)$. It is quite remarkable that this extends to the stable category of MCM modules.
\begin{defn} A local (respectively, graded local) commutative ring $R = (R, \m, k)$ has isolated singularities if $R_\mathfrak{p}$ is regular for each $\mathfrak{p} \in \spec R \setminus \{\m\}$ (respectively, homogeneous primes $\mathfrak{p} \in \spec^* R \setminus \{\m \}$). 
\end{defn}

\begin{prop}[Auslander \cite{AuslanderIsolSing}, \cite{IT}] Let $R$ be a commutative Gorenstein, or graded Gorenstein, $k$-algebra. 
	\begin{enumerate}[i)]
	\item When $R$ is complete local with isolated singularities, then $ \underline{{\rm MCM}}(R)$ is Hom-finite over $k$ and has a Serre functor given by $\Se_R(-) = (-)[\dim R - 1]$. 
	\item When $R$ is graded connected with isolated singularities, then $ \underline{{\rm MCM}}^{\mathbb{Z}}(R)$ is Hom-finite over $k$ and has a Serre functor given by $\Se_R(-) = - \otimes_R \omega_R [\dim R - 1]$.
\end{enumerate}
\end{prop}
\begin{rem} Note that when $R$ is complete local above, we have $\omega_R \cong R$ and so the tensor $- \otimes_R \omega_R$ is implicit. In the graded case, $\omega_R \cong R(a)$ and so the $M \otimes_R \omega_R \simeq M(a)$ for $M \in \underline{{\rm MCM}}^{\mathbb{Z}}(R)$. \end{rem} 

\begin{rem}
	Note that when $R$ is a standard graded $k$-algebra, the punctured spectrum $\spec R \setminus \{ \m \}$ forms a $\mathbb{G}_m$-bundle over $X = \proj R$, and so $R$ has isolated singularities if and only if $X$ is smooth. We will use this fact implicitly throughout the thesis. 
\end{rem}

Serre functors were independently discovered by Auslander and Reiten in the guise of the  translate $\tau$, introduced in the context of stable module categories, see \cite{ARS}. In \cite{ReitVdBalmostsplit}, Reiten-Van den Bergh studied $\tau$ in the context of a Krull-Schmidt Hom-finite $k$-linear triangulated category $\mathcal{T}$. Let
\[
	\xi: X \to Y \to Z \xrightarrow{h} X[1]
\]
be a distinguished triangle in $\mathcal{T}$. We call $\xi$ an almost-split triangle if
\begin{enumerate}[i.]
\item X and Z are indecomposable,
\item $h \neq 0$,
\item if $W$ is indecomposable, then for every non-isomorphism $t: W \to Z$ we have $ht = 0$. Equivalently, we have a lift as in the diagram below
	\[
		\xymatrix@R1pc{ & & W \ar[d]^{t} \ar@{..>}[dl] \\
		X \ar[r] & Y \ar[r]& Z \ar[r]^{h}& X[1]. }
	\]
\end{enumerate}

\begin{prop}[Reiten-Van den Bergh, \cite{ReitVdBalmostsplit}] Assume that $k$ is algebraically closed. The following are equivalent over a Krull-Schmidt Hom-finite $k$-linear triangulated category $\mathcal{T}$: 
\begin{enumerate}[1.]
\item Each indecomposable $Z$ of $\mathcal{T}$ sits inside an almost-split triangle, say
	\[
		\tau Z \to Y \to Z \xrightarrow{h} \tau Z [1].
	\]
\item $\mathcal{T}$ admits a Serre functor $\mathbb{S}$.
\end{enumerate}
In this case $\tau = \mathbb{S}\circ [-1]$ and the map $Z \xrightarrow{h} \tau Z [1] = \mathbb{S}(Z)$ classifying the extension is Serre dual to the trace map ${\rm End}(Z) \to {\rm End}(Z)/\rad\ \!{\rm End}(Z) = k$.
\end{prop}

Now for $\mathcal{T}$ as above, the Auslander-Reiten quiver $\Gamma(\mathcal{T})$ is the quiver whose vertices are the isomorphism classes of indecomposables of $\mathcal{T}$ and arrows taken from a basis for ${\rm Irr}(X, Y)$, the space of equivalence classes of irreducible maps between indecomposables (see e.g. \cite{HappelBook} for details). The Auslander-Reiten quiver $\Gamma(\mathcal{T})$ is related to almost-split triangles as follows.
\begin{prop}[{\cite[4.8]{HappelBook}}]\label{Ringel} Let $X, M, Z$ be indecomposable objects in $\mathcal{T}$ and $X \to Y \to Z \to X[1]$ almost-split, with $Y = \bigoplus_{i=1}^r Y_{i}^{\oplus e_i}$ decomposed into pairwise non-isomorphic indecomposables. Then ${\rm Irr}(X, M) \neq 0$ if and only if $M \cong Y_i$ for some i, in which case $e_i = {\rm dim}\ {\rm Irr}(X, Y_i)$.
\end{prop}

\begin{rem} Given $R$ Gorenstein, the Auslander-Reiten quiver of ${\rm MCM}(R)$ was defined in \cite{YosBook} in terms of almost-split short exact sequences. It is immediate that almost-split sequences descend to almost-split triangles in $ \underline{{\rm MCM}}(R)$, and one obtains the Auslander-Reiten quiver of $ \underline{{\rm MCM}}(R)$ from that of ${\rm MCM}(R)$ by removing the vertex corresponding to $R$.
\end{rem}

\subsection*{Complete resolutions over complete intersections rings}

Let $Q$ be a regular local ring with element $f \in Q$, and set $R = Q/f$ the hypersurface ring. Complete resolutions of MCM modules over $R$ take a rather simple form using a construction of Eisenbud \cite{Eis}, \cite[Chp. 7]{YosBook}. \bigskip 

Let $M \in {\rm MCM}(R)$ be an MCM $R$-module, so that ${\rm depth}_R\ \! M = \dim R = \dim Q - 1$. The Auslander-Buchsbaum formula for the projective dimension of $M$ over $Q$ gives ${\rm pdim}_Q\ \! M + {\rm depth}_Q\ \! M = \dim Q$, since ${\rm depth}_Q\ \! M = {\rm depth}_R\ \! M = \dim Q - 1$ gives ${\rm pdim}_Q M = 1$. Hence we have a length two resolution 
\[
	\xymatrix{0 \ar[r] & F \ar[r]^-{A} & G \ar[r] & M \ar[r] & 0}
\]
with $F, G$ finite free modules over $Q$. Moreover, $M$ is annihilated by $f$ and so multiplication by $f$ on the deleted resolution must be nullhomotopic 
\[
	\xymatrix{ 
		0 \ar[r] & F \ar[d]_-{f} \ar[r]^-{A} & G \ar@{.>}[dl]_-{B} \ar[d]_-{f} \ar[r] & 0\\
		0 \ar[r] & F \ar[r]^-{A} & G \ar[r] & 0.
	}
\]
That is, we have factorisations $f \cdot {\rm id_G} = AB$ and $f \cdot {\rm id_F} = BA$. Since $f$ is a regular element, $A$ and $B$ become isomorphisms over $Q[f^{-1}]$ and so the free modules $F$ and $G$ have the same rank, say $r$. The tuple $(A, B)$ then forms a pair of $r \times r$ matrices in $Q$ whose product is the diagonal matrix $f \cdot {\rm I}_r$. Such a pair is called a matrix factorisation of $f$.\bigskip 

Since $f$ is regular, it isn't hard to see that the complex obtained from the above by applying $- \otimes_Q R$ and extending by $2$-periodicity is acyclic 
\[
	\cdots \to \overline{G} \xrightarrow{\overline{B}} \overline{F} \xrightarrow{\overline{A}} \overline{G} \xrightarrow{\overline{B}} \overline{F} \to \cdots 
\]
and yields a complete resolution of $M$ over $R$. It follows that all complete resolutions and MCM $R$-modules are described by matrix factorisations. We note that this extends to an equivalence of triangulated categories 
\[
	{\rm \underline{MF}}(Q, f) \cong {\rm \underline{MCM}}(R) 
\]
between the homotopy category (or stable category) of matrix factorisations, and the stable category of MCM $R$-modules (see \cite[Chp. 7]{YosBook} for details on this category). Lastly, all results above extend naturally to graded modules over graded hypersurface rings. \bigskip 

Finally, we note that the description of complete resolutions of MCM modules above was extended to arbitrary complete intersections $R = Q/(f_1, \cdots, f_c)$ by Buchweitz-Pham-Roberts in \cite{BuchweitzCI}.

\subsection*{Orlov's semiorthogonal decomposition Theorem} 
Finally we arrive at Orlov's Theorem, following \cite{Or1}. We will sketch the basic statement and refer to Appendix \ref{orlovappendix} and \cite{Or1} for the more general picture. Let $k$ be a field and $R$ a commutative, standard-graded, connected Gorenstein $k$-algebra (or more generally a noncommutative Artin-Schelter Gorenstein algebra, see Appendix \ref{orlovappendix}). \bigskip

In such a scenario, the projective scheme $X = \proj R$ and its derived category $ {\rm D}^b(\coh X) $ are closely related to the singularity category ${\rm D}_{sg}^\Z(R)$. Recall that by Serre's Theorem, we can reconstruct $\coh X \cong \qgr R$ as the quotient category of finitely generated graded $R$-modules $\grmod R$, modulo the Serre subcategory of torsion modules. For any choice of cut-off $i \in \Z$, denote by $\mathsf{grmod}_{\geq i}\ \! R$ the full abelian subcategory of graded $R$-modules with $M_{j} = 0$ for $j < i$. The quotient functor $\pi: \grmod R \twoheadrightarrow \coh X$ restricts to an essentially surjective exact functor $\pi_i: \mathsf{grmod}_{\geq i}\ \! R \twoheadrightarrow \coh X$. \bigskip

The functor $\pi_i$ admits a right adjoint $\Gamma_{\geq i}: \coh X \to \mathsf{grmod}_{\geq i}\ \! R$ given by
\[
	\Gamma_{\geq i}(\F) = \bigoplus_{n \geq i}{\rm \Gamma}(X, \F(n))
\]
with its natural graded $R$-module structure. Moreover, we have a natural isomorphism $\pi_i \circ \Gamma_{\geq i} \simeq {\rm id}$ and so $\Gamma_{\geq i}$ is fully faithful. This extends to an adjoint pair on derived categories 
\[
	\pi_i: 	{\rm D}^b(\mathsf{grmod}_{\geq i}\ \! R) \leftrightarrows {\rm D}^b(\qgr R): {\rm R}\Gamma_{\geq i}
\]
with ${\rm R}\Gamma_{\geq i}$ fully faithful and $\pi_i$ essentially surjective. \bigskip

We have another Verdier quotient ${\rm st}: {\rm D}^b(\grmod R) \twoheadrightarrow {\rm D}_{sg}^{\Z}(R) \cong \underline{{\rm MCM}}^{\mathbb{Z}}(R)$ given by stabilisation, or MCM approximation. Following an observation of Buchweitz, the stabilisation functor also admits a left adjoint. Given on an MCM module $M$, first pick a complete resolution $C$ 
\[
	\cdots \to C_{n+1} \to C_n \to C_{n-1} \to \cdots.
\]
Quotienting out the subcomplex $C_{[<i]} \subseteq C$ whose terms are given by summands generated in degree less than $i$, one obtains this way a complex $C_{[\geq i]} \in {\rm D}^b(\mathsf{grmod}_{\geq i}\ \! R) $ with at most bounded cohomology, which is functorial in $M$ and provides a left adjoint to stabilisation (see Appendix \ref{orlovappendix} and Prop. \ref{BuchweitzAdjoint} for details). With slight abuse of notation, we denote this by: 
\[
	(-)_{[\geq i]}: {\rm D}^{\Z}_{sg}(R) \rightleftarrows {\rm D}^b(\mathsf{grmod}_{\geq i}\ \! R): {\rm st} 
\]
Moreover, we have ${\rm st}\circ  (-)_{[\geq i]} \simeq {\rm id}$, and so $(-)_{[\geq i]}$ is also fully faithful. We may then compose both adjoints as shown below: 
\[
	\xymatrix{ & {\rm D}^b(\mathsf{grmod}_{\geq i}\ \!R) \ar@<.5ex>[dl]^-{st} \ar@<-.5ex>[dr]_-{\pi_i} & \\ {\rm D}_{sg}^{\Z}(R) \ar@<.5ex>[ur]^-{(-)_{[\geq i]}} \ar@<-.5ex>[rr]_{\Phi_i} && {\rm D}^b(\coh X) \ar@<-.5ex>[ll]_{\Psi_i} \ar@<-.5ex>[ul]_-{ {\rm R}\Gamma_{\geq i}} 
	}
\] 
This yields an adjoint pair $(\Phi_i, \Psi_i)$ with $\Phi_i = \pi_i \circ (-)_{[\geq i]}$ and $\Psi_i = {\rm st} \circ {\rm R}\Gamma_{\geq i}$. Note that while the categories $ {\rm D}^b(\coh X)$ and ${\rm D}^{\Z}_{sg}(R)$ do not depend on the resulting cutoff $i$, both functors $(\Phi_i, \Psi_i)$ do and will generally differ as $i$ varies. \bigskip 

The following is Orlov's semiorthogonal decomposition theorem in the current setting. We refer to \cite{Or1}, \cite{BurkeStev} or Appendix \ref{orlovappendix} for the more general statement and for the definition of semiorthogonal decomposition.  
\begin{thm}[{\cite[Thm. 2.5]{Or1}}] Let $R$ be a commutative, standard-graded, connected Gorenstein $k$-algebra with projective scheme $X = \proj R$, and $a$-invariant $a \in \Z$. The above functors and triangulated categories are related as follows:
	\begin{enumerate}[i)]
	\item (Fano case) if $a < 0$, there is a semiorthogonal decomposition
		\[
			{\rm R}\Gamma_{\geq i}\left({\rm D}^b(\coh X)\right) = \big\langle R(-i + a +1), R(-i+a+2), \dots, R(-i), {\rm D}^{\Z}_{sg}(R)_{[\geq i]} \big\rangle.
		\]
		Applying $\pi_i$, this descends to a semiorthogonal decomposition 
		\[
		{\rm D}^b(\coh X) = \big\langle  \OO_X(-i + a + 1), \OO_X(-i+a+2), \dots, \OO_X(-i), \Phi_i{\rm D}^{\Z}_{sg}(R)) \big\rangle.
		\]
	\item (Calabi-Yau case) if $a = 0$, the essential images of both embeddings in $ {\rm D}^b(\mathsf{grmod}_{\geq i}\ \! R)$ are equal 
		\[
			{\rm D}^{\Z}_{sg}(R)_{[\geq i]} = {\rm R}\Gamma_{\geq i}({\rm D}^b(\coh X))  
		\]
		hence $(\Phi_i, \Psi_i)$ give inverse equivalences
		\[
			\Phi_i: {\rm D}_{sg}^{\Z}(R) \cong {\rm D}^b(\coh X) : \Psi_i.
		\]
	\item (General type case) if $a > 0$, there is a semiorthogonal decomposition
		\[
			{\rm D}^{\Z}_{sg}(R)_{[\geq i]} = \big\langle k(-i), k(-i-1), \dots, k(-i-a+1), {\rm R}\Gamma_{\geq i+a} {\rm D}^b(\qgr R) \big\rangle.
		\]
		Applying the stabilisation ${\rm st}$, this descends to a semiorthogonal decomposition
		\[
			{\rm D}^{\Z}_{sg}(R) = \big\langle k^{st}(-i), k^{st}(-i-1), \dots, k^{st}(-i-a+1), \Psi_{i+a} {\rm D}^b(\coh X) \big\rangle.
		\]
\end{enumerate}
\end{thm}

%% file: chapters/tilting/tilting.tex
\chapter{Graded Gorenstein rings with tilting MCM modules}

\input{chapters/tilting/known}
\input{chapters/tilting/agv}
\input{chapters/tilting/ci}

%% file: chapters/tilting/known.tex
\section{Some history: Tilting theory, exceptional singularities and low dimension}


A classical problem of algebraic geometry asks for the realisation of a vector bundle $\E$ on projective space $\mathbb{P} = \mathbb{P}(V)$ over a field $k$ as the zeroth cohomology sheaf of a complex of `known' vector bundles
\[
	0 \to \E_m \to \dots \to \E_0 \to \dots \to \E_{-m} \to 0
\]
with cohomology concentrated in degree $0$ and where $\E_i$ a sum of bundles of the form $\OO(j)$, $\wedge^i\T_{\mathbb{P}}(j)$ or $\Omega^i_{\mathbb{P}^2}(j)$ for $j \in \Z$. This problem was solved in its general form by Beilinson \cite{Bei}, who proved that every coherent sheaf $\F$ over $\mathbb{P}^n$ was quasi-isomorphic to a complex of the form 
\[
	0 \to \E_n \to \dots \to \E_0 \to \dots \to \E_{-n} \to 0
\]
with general term
\[
	\E_i = \bigoplus\limits_{j} {\rm H}^{i+j}(\mathbb{P}^n, \F(-j))\otimes_k \Omega_{\mathbb{P}^n}^{j}(j)
\]
as well as to one with terms 
\[
	\E_i = \bigoplus\limits_{j} {\rm H}^{i-j}(\mathbb{P}^n, \F \otimes \Omega^j(j)) \otimes_k \OO_{\mathbb{P}^n}(-j).
\] 
These complexes are known as the Beilinson monads (see \cite{EFS} for a constructive approach). Beilinson's approach took a derived category perspective: he showed that the sequences of sheaves 
\[
	\left(\Omega^n(n), \Omega^{n-1}(n-1), \dots, \Omega^1(1), \OO \right)
\]
\[
	\left(\OO(-n), \OO(-n+1), \cdots, \OO(-1), \OO \right)
\]
form full, strong exceptional collections in ${\rm D}^b(\mathbb{P}^n)$, meaning that upon setting $E_i = \Omega^{n-i}(n-i)$ (respectively $E_i = \OO(-n+i)$), for $0 \leq i,j \leq n$ and $l \geq 1$ we have
\begin{align*}
	{\rm Hom}(E_i, E_j) = \begin{cases} 0, & i > j \\ k, & i=j\end{cases} && {\rm Ext}^l(E_i, E_j) = 0 \textup{ for all } i, j
\end{align*}
and that ${\rm D}^b(\mathbb{P}^n)$ is the smallest triangulated subcategory closed under direct summands containing the $\{E_i\}$. Letting $T = \bigoplus_{i = 0}^n E_i$ and $\Lambda = {\rm End}(T)$, it follows \cite[Thm. 1.8]{BuHi} that we have inverse equivalences of triangulated categories
\begin{eqnarray}\label{Beilinsontilting} 
	{\rm RHom}(T, -): {\rm D}^b(\mathbb{P}^n) \leftrightarrows {\rm D}^b(\Lambda): - \otimes^{\mathbb{L}}_{\Lambda} T
\end{eqnarray}
onto the derived category of a finite dimensional noncommutative $k$-algebra $\Lambda$. This sends the $E_i$ onto the indecomposable projective $\Lambda$-modules, and pulls back a general complex of $\Lambda$-modules to a complex whose terms are sums of $\Omega^{j}(j)$ (respectively $\OO(-j)$). \bigskip

Equivalences of the form (\ref{Beilinsontilting}) for smooth projective varieties $X$ have been heavily studied since, and they are always induced from a special object $T \in {\rm D}^b(X) $. It isn't necessary that the indecomposable summands of $T = \bigoplus_{i}T_i$ form a full strong exceptional collection. Rather, the slightly weaker condition is that $T$ be a \emph{tilting} object.\bigskip 

Let $\T$ be a triangulated $k$-linear category. For $\F \in \T$, we denote by $\thick(\F) \subseteq \T$ the smallest triangulated subcategory containing $\F$ closed under direct summands. We say that $T \in \T$ is \emph{tilting} if it is a classical generator for $\T$ with no non-trivial self-extensions, that is: 
\begin{enumerate}[1)]
\item (Generating) We have $\thick(T) = \T$, that is $\T$ is the smallest triangulated subcategory containing $T$ and closed under summands. 
\item (No self-extensions) We have ${\rm Hom}_{\T}(T, T[n]) = 0$ for $n \neq 0$. 
\end{enumerate} 
Other variants of the definition are in use but the above suffices for our purpose (See \cite[Section 1]{BuHi} for a general discussion). \bigskip

Let $X$ be a smooth projective variety over $k$ with a tilting object $T \in {\rm D}^b(X) $. The general picture is given as follows: 
\begin{thm}[{\cite[Thm. 1.8]{BuHi}}]\label{tiltingforvarieties} We have the following properties:
	\begin{enumerate}[1)]
	\item The endomorphism algebra $\Lambda = {\rm End}(T)$ is a finite dimensional $k$-algebra of finite global dimension. 
	\item There are induced equivalences of triangulated categories
\[
	{\rm RHom}(T, -): {\rm D}^b(X) \leftrightarrows {\rm D}^b(\Lambda): - \otimes^{\mathbb{L}}_{\Lambda} T.
\]
\end{enumerate}
\end{thm}
\noindent Finite dimensionality of $\Lambda$ is a consequence of the properness of $X$ over $k$, while finite global dimension actually follows from $X$ being smooth. Varieties with tilting objects include projective spaces, quadric hypersurfaces, Del Pezzo surfaces, some Toric varieties and generalised flag varieties, and furthermore any products of such $X$ or iterated projective bundles over a base $S$ with a tilting object. \bigskip

Theorem \ref{tiltingforvarieties} is a specialisation of Keller's derived Morita Theorem. 
\begin{thm}[Keller, {\cite{KellerTilting}}] Let $\T$ be an (algebraic) triangulated $k$-linear category with a tilting object $T$. Let $\Lambda = {\rm End}_\T(T)$. Then there exists an equivalence of triangulated categories
	\[
		{\rm RHom}_\T(T, -): \T \xrightarrow{\cong} {\rm D}^{\perf}(\Lambda) 
	\]
	onto the subcategory ${\rm D}^{\perf}(\Lambda) = \thick(\Lambda) \subseteq {\rm D}(\Lambda)$ of perfect complexes, meaning complexes quasi-isomorphic to bounded complexes of finitely generated projectives. 
\end{thm} 
All triangulated categories appearing in this thesis will be algebraic, see \cite{KellerTilting} or the appendix for the definition. Note that we have $ {\rm D}^{\perf}(\Lambda) = {\rm D}^b(\Lambda)$ when $\Lambda$ is Noetherian of finite global dimension. \bigskip

This thesis is concerned with tilting objects in the triangulated category $\T = \underline{{\rm MCM}}^{\mathbb{Z}}(R)$, the stable category of graded MCM modules over a graded Gorenstein $k$-algebra $R$, which is a close cousin of the triangulated categories of the form $ {\rm D}^b(X) $. The role of tilting theory for MCM modules first came to the front in the study of Kleinian singularities and the McKay correspondence, and we review some of this story. 

\subsubsection*{The McKay Correspondence}
Let $k$ here be an algebraically closed field of characteristic zero, and let $G \leq SL(2, k)$ be a finite subgroup. The possible such subgroups $G$ up to conjugacy were classified by Klein, and are in one-to-one correspondence with the simply-laced Dynkin diagrams of type ADE (where $\mu$ below denotes the number of vertices):
\begin{align*}
	&\xymatrix{ 
	A_\mu:&	\bullet \ar@{-}[r] & \bullet \ar@{-}[r] & \ar@{-}[r]  \cdots  \ar@{-}[r] & \bullet \ar@{-}[r] & \bullet   & \mu \geq 1
	}
		\\ 
	&\xymatrix@R=5pt{ 
		&&&&& \bullet \\
		D_\mu: 	&\bullet \ar@{-}[r] & \bullet \ar@{-}[r] & \ar@{-}[r]  \cdots  \ar@{-}[r] & \bullet \ar@{-}[ur] \ar@{-}[dr] &   & \mu \geq 4\\
		&&&&& \bullet 
	}\\
	&\xymatrix@R=10pt{ &&& \bullet \ar@{-}[d] \\
		E_6: &	\bullet \ar@{-}[r] & \bullet \ar@{-}[r] & \ar@{-}[r]  \bullet  \ar@{-}[r] & \bullet \ar@{-}[r] & \bullet   
	}\\
	&\xymatrix@R=10pt{ &&& \bullet \ar@{-}[d] \\
		E_7: &	\bullet \ar@{-}[r] & \bullet \ar@{-}[r] & \ar@{-}[r]  \bullet  \ar@{-}[r] & \bullet \ar@{-}[r] & \bullet  \ar@{-}[r] & \bullet   
	}\\
	&\xymatrix@R=10pt{ &&& \bullet \ar@{-}[d] \\
		E_8: &	\bullet \ar@{-}[r] & \bullet \ar@{-}[r] & \ar@{-}[r]  \bullet  \ar@{-}[r] & \bullet \ar@{-}[r] & \bullet  \ar@{-}[r] & \bullet  \ar@{-}[r] & \bullet   
	}
\end{align*}
Let $\zeta_n$ be a primitive $n$-th root of unity. The classification of finite subgroups $G \leq {\rm SL}(2, k)$ is as follows, up to conjugacy \cite[10.15]{YosBook}: 
\begin{align*}
	&A_\mu: \textup{ The cyclic group of order $\mu+1$ } &
\end{align*}
\[
	C_{\mu} := \bigg\langle \begin{pmatrix} \zeta_{\mu+1} & 0\\ 0 & \zeta_{\mu+1}^{-1} \end{pmatrix} \bigg\rangle 
\]
\begin{align*}
	&D_\mu: \textup{ The binary dihedral group of order $4(\mu-2)$ }&
\end{align*}
\[
	D_{\mu} := \bigg\langle \begin{pmatrix} 0 & \zeta_4 \\ \zeta_4 & 0 \end{pmatrix}, \ C_{2\mu+5} \bigg\rangle 
\]
\begin{align*}
	&E_6: \textup{ The binary tetrahedral group of order $24$ } &
\end{align*}
\[
	\mathbb{T} := \bigg\langle \frac{1}{\sqrt{2}}\begin{pmatrix} \zeta_8 & \zeta_8^3 \\ \zeta_8 & \zeta_8^7 \end{pmatrix}, D_4 \bigg\rangle 
\]
\begin{align*}
	&E_7: \textup{ The binary octahedral group of order $48$ } &
\end{align*}
\[
	\mathbb{O} := \bigg\langle \begin{pmatrix} \zeta_8^3 & 0 \\ 0 & \zeta_8^3 \end{pmatrix}, \mathbb{T} \bigg\rangle 
\]
\begin{align*}
	&E_8: \textup{ The binary icosahedral group of order $120$ } &
\end{align*}
\[
	\mathbb{I} := \bigg\langle \frac{1}{\sqrt{5}}\begin{pmatrix} \zeta_5^4 - \zeta_5 & \zeta_5^2-\zeta_5^3 \\ \zeta_5^2-\zeta_5^3 & \zeta_5 - \zeta_5^4 \end{pmatrix}, \frac{1}{\sqrt{5}} \begin{pmatrix} \zeta_5^2 - \zeta_5^4 & \zeta_5^4 - 1 \\ 1 - \zeta_5 & \zeta_5^3 - \zeta_5 \end{pmatrix} \bigg\rangle 
\]

Let $S = k[u,v]$ and $R = S^G$, so that $X = \spec R = \mathbb{A}^2 / G$ is the associated quotient singularity. Let $Q$ be the Dynkin graph associated to $G = G_Q$. The graph $Q$ has many natural incarnations throughout this story. In particular:
\begin{enumerate}[1)]
\item The ADE graph $Q$ is the dual graph of the exceptional divisor in the minimal resolution $\pi: \widetilde{X} \to X$. 
\item The extended ADE graph $\widetilde{Q}$ arises as the McKay graph associated to the representation theory of $G$, with extended vertex corresponding to the trivial representation of $G$. The McKay graph of $G$ has vertices $\{V_j\}$ a representative set of all irreducible representations of $G$, with multiplicity of edges between $(V_i, V_j)$ given by $\dim_k {\rm Hom}_{kG}(V_i \otimes V, V_j)$ where $V$ is the standard $2$-dimensional representation of $G \leq SL(2, k)$.\footnote{Note that the edge multiplicity is well-defined: any choice of isomorphism $\bigwedge^2 V \cong k$ with the trivial $G$-representation gives $V \cong V^*$ equivariantly, so that $\dim_k {\rm Hom}_{kG}(V_i \otimes V, V_j) = \dim_k {\rm Hom}_{kG}(V_i, V \otimes V_j) = \dim_k {\rm Hom}_{kG}(V_j \otimes V, V_i)$.}
\item Over $k = \C$, the root system of type $Q$ arises as the vanishing cohomology (or Milnor lattice) of the generic fibre $X_b$ in the semiuniversal deformation $X \to B$. 
	\end{enumerate} 

	The Kleinian singularities $X = \mathbb{A}^2//G$ form one of the most exceptional setting in all of algebra and geometry, and have many intrinsic characterisations. The coordinate ring $R = S^G$ is always generated by three fundamental invariants satisfying one relation, and so the invariant ring is isomorphic to a hypersurface ring $S^G \cong k[x,y,z]/f$, with $f$ given by:
	\begin{itemize}
		\item[$(A_\mu)$] $f = x^{2} + y^{\mu+1} + z^2,\ \mu \geq 1$
		\item[$(D_\mu)$] $f = x^2y + y^{\mu-1} + z^2, \ \mu \geq 4$
		\item[$(E_6)$] $f = x^3 + y^4 + z^2$   
		\item[$(E_7)$] $f = x^3 + xy^3 + z^2$ 
		\item[$(E_8)$] $f = x^3 + y^5 + z^2$ 
	\end{itemize}
	When working over $\C$, the above list of ADE polynomials give the local normal forms for the ``simple'' hypersurface singularities classified by the Arnold school of singularity theory, up to adding or removing squares in disjoint variables to vary dimension. \bigskip 

Taking a more commutative algebraic perspective under the Auslander school, one notes that the invariant ring $R = S^G \subseteq S$ is Gorenstein, and that $S$ is a maximal Cohen-Macaulay module over $R$. Furthermore, the natural $G$-equivariant decomposition
\begin{align*}
	S = {\rm Hom}_{kG}(kG, S) &= \bigoplus_{V_j} {\rm Hom}_{kG}(V_j, S) \\
&= \bigoplus_{V_j} {\rm Hom}_k(V_j, S)^G
\end{align*}
is a decomposition of $S$ into indecomposable MCM $R$-modules, where $\{V_j\}$ runs over the irreducible representations of $G$ with the same multiplicity as they arise in $kG$. Moreover, one can run the same construction over the completion $\widehat{S} = k[[u,v]]$ and $\widehat{R} = k[[u,v]]^G$ to obtain MCM modules over $\widehat{R}$.
 
\begin{thm}[Auslander {\cite[Cor. 10.10]{YosBook}}] The above gives a bijection between irreducible $G$-representations and indecomposable MCM $\widehat{R}$-modules. This sends the trivial representation to the free module $\widehat{R}$. 
\end{thm}
Since there are finitely many indecomposable MCM $\widehat{R}$-module, the Gorenstein ring $\widehat{R}$ is said to be of finite MCM representation type (or CM-type for short). This turns out to characterise ADE singularities. 

\begin{thm}[Buchweitz-Greuel-Kn\"orrer-Schreyer {\cite[Thm. 8.10, Cor. 12.6]{YosBook}}] Let $A = k[[x_1, \cdots, x_n]]/(g)$ be a complete hypersurface ring over $k$ with $g \in (x_1, \cdots, x_n)^2$. Then $A$ is of finite CM-type if and only if $A \cong k[[x,y, z_3, \dots, z_n]]/(f)$ with $f$ an ADE polynomial in $n \geq 2$ variables: 
	\begin{itemize}
		\item[$(A_\mu)$] $f = x^{2} + y^{\mu+1} + z_3^2 + \dots + z^2_n,\ \mu \geq 1$
		\item[$(D_\mu)$] $f = x^2y + y^{\mu-1} + z_3^2 + \dots + z^2_n,\ \mu \geq 4$
		\item[$(E_6)$] $f = x^3 + y^4 + z_3^2 + \dots + z^2_n$
		\item[$(E_7)$] $f = x^3 + xy^3 + z_3^2 + \dots + z^2_n$
		\item[$(E_8)$] $f = x^3 + y^5 + z_3^2 + \dots + z^2_n$.
	\end{itemize}
\end{thm}

\begin{thm}[Herzog {\cite[Cor. 8.16]{YosBook}}] Let $A$ be a complete Gorenstein local $k$-algebra. If $A$ is of finite CM-type, then $A$ is a hypersurface singularity, and therefore an ADE hypersurface singularity. 
\end{thm}

Let us denote by $ \underline{{\rm MCM}}(R)$ the projectively stable category of MCM $R$-modules obtained by killing morphisms factoring through projective modules. Rephrasing the above homologically, as one runs over Gorenstein complete local $k$-algebras $R$, the triangulated categories of the form $ \underline{{\rm MCM}}(R)$ with finitely many indecomposables are precisely given by the ADE hypersurface singularities. \bigskip

The search for triangulated categories of `finite representation type' is one of the central problem of abstract representation theory, where the ADE graphs feature prominently. Recall that amongst the acyclic quivers $Q$, the ones of finite representation type are precisely the quivers whose underlying graph is of type ADE. \bigskip

The representation theory of the ADE quiver $Q$ (with arbitrary choice of orientation) actually arises as part of the above storyline. To see this, consider again the Kleinian singularites $X = \mathbb{A}^2 // G$, and note that the invariant ring $R = S^G = k[u,v]^G$ inherits a grading from the standard grading on $S = k[u,v]$ with $|u| = |v| = 1$, obtained by taking $G$-invariants in each degree. For simplicity, let us restrict to the case of $|G|$ even\footnote{This covers all groups $G$ except for the odd cyclic group of order $\mu_{2n+1}$ corresponding to the $A_{2n}$ singularity. All results stated below have their analog for this case after suitable modification.}. In this case, the fundamental invariants $x,y,z \in S$ of $G$ have even degrees with $\gcd(|x|, |y|, |z|) = 2$, and $R = \bigoplus_{n \geq 0} R_{2n}$ is properly supported in even degrees. We hence regrade $R$ by halving the degrees of its homogeneous elements.\bigskip 

Since $R$ is graded, we may consider the stable category of graded MCM modules $\underline{{\rm MCM}}^{\Z}(R)$. Recall that graded modules have a degree shift operator given by $M \mapsto M(1)$ where $M(1)_n = M_{n+1}$.

\begin{thm}[Kajiura-Saito-A.Takahashi, Iyama-R.Takahashi {\cite{KST} \cite{IT}}]\label{kleiniantilting} There is a tilting object $T \in \underline{{\rm MCM}}^{\mathbb{Z}}(R)$ with endomorphism algebra ${\rm \underline{End}}_{gr R}(T) \cong kQ$, the path algebra of $Q$. Hence there is an equivalence of triangulated categories 
	\[
		\underline{{\rm MCM}}^{\mathbb{Z}}(R) \cong {\rm D}^b(kQ).
	\]
	In particular, this identifies the Grothendieck lattice $K_0(\underline{{\rm MCM}}^{\mathbb{Z}}(R))$, equipped with the symmetrized Euler pairing, with the root lattice of $Q$. Moreover, this equivalence sends the degree shift operator $M \mapsto M(1)$ to the Auslander-Reiten translate $\tau = - \otimes^{\mathbb{L}}_{kQ}D(kQ)[-1]$ on $ {\rm D}^b(kQ)$. 
\end{thm}
The indecomposable MCM $R$-modules $M_j := {\rm Hom}_k(V_j, S)^G$ are naturally graded modules. Combining the above results, we obtain bijections between isomorphisms classes of:

\begin{enumerate}[i.]
\item Irreducible representations of $G$ other than the trivial representation. 
\item Indecomposable non-free graded MCM modules up to degree shift.
\item Orbits of indecomposables under the Auslander-Reiten translate $\tau$ in ${\rm D}^b(kQ)$.
\item Vertices of $Q$.
	\end{enumerate}	
	The bijection $i \leftrightarrow ii$ is Auslander's, $ii \leftrightarrow iii$ follows from Theorem \ref{kleiniantilting}, while $iii \leftrightarrow iv$ follows since $kQ$ has finite representation type, so that each indecomposable complex is in the $\tau$-orbit of a unique indecomposable projective $kQ$-module, which are indexed by the vertices of $Q$. This gives a powerful extension of McKay's observation of the bijection $i \leftrightarrow iv$. By a theorem of Keller-Murfet-Van den Bergh \cite{KMVdB}, the equivalence above further descends to an equivalence of triangulated orbit categories (in the sense of Keller \cite{Ke})
	\[
		\underline{{\rm MCM}}(\widehat{R}) = \underline{{\rm MCM}}^{\mathbb{Z}}(R)/(1) \cong {\rm D}^b(kQ)/\tau
	\]
	and one recovers the McKay quiver (except for the vertex corresponding to the trivial representation) by taking the Auslander-Reiten quiver of either category. 

\subsection*{Existence of tilting objects} 
Fix an algebraically closed field $k$. In the first part of the thesis, we will be interested in the following central question and its applications: 
\begin{quest} Let $R = \bigoplus_{n \geq 0} R_n$ be a non-negatively graded Gorenstein $k$-algebra. When does the triangulated category $ \underline{{\rm MCM}}^{\mathbb{Z}}(R)$ admits a tilting object?
\end{quest}
Many interesting results have been obtained by Iyama and his collaborators, and we first review their work. For a more complete recent survey, see \cite{IyamaSurvey}. \bigskip

One can first extend the tilting result for Kleinian singularities to higher dimensional quotient singularities. The following is due to Iyama-Takahashi (but see also \cite{MU} for a noncommutative variant and a different interpretation). Let $S = k[x_1, \dots, x_n]$ be the standard graded polynomial algebra, and recall that we write $\Omega^i_S(-)$ to denote $i$-th syzygy. Given an $R$-module $N$ below, we let $[N]_{\mathsf{CM}}$ denote the largest MCM summand of $N$.
\begin{thm}[Iyama-Takahashi, {\cite[Thm. 1.7]{IT}}] Let $G \leq SL(n, k)$ be a finite group with $\mathsf{char}\ k \nmid |G|$, and let $R = S^G \subseteq S$ inherit the natural grading. Then $\underline{{\rm MCM}}^{\mathbb{Z}}(R)$ admits a tilting object $U$ given by
	\[
		U = \bigoplus_{i=1}^n [\Omega^i_S k(i)]_{\mathsf{CM}}.
	\]
\end{thm}
\begin{exmp}[{\cite[Ex. 7.16]{IT}}] Let $C_3 = \langle {\rm diag}(\omega, \omega, \omega) \rangle \subset SL(3, k)$ with $\omega$ a primitive third root of unity in $k$ and $\mathsf{char}\ k \neq 3$. Let $Q$ be the $3$-Kronecker quiver 
	\[
		\xymatrix{ \bullet \ar@<1ex>[r] \ar[r] \ar@<-1ex>[r] & \bullet}
	\]
	Then $ {\rm \underline{End}}_{gr R}(U) \cong kQ \times kQ \times kQ$, so that $ \underline{{\rm MCM}}^{\mathbb{Z}}(R) \cong {\rm D}^b(kQ) \times {\rm D}^b(kQ) \times {\rm D}^b(kQ)$. The $i$-th copy of $ {\rm D}^b(kQ) $ corresponds to the subcategory $ \underline{{\rm MCM}}^{i + 3\mathbb{Z}}(R)$ of graded MCM modules $M = (M_n)$ supported in degrees $n \equiv i \ (mod\ 3)$.
 \end{exmp}

\subsubsection*{Dimension zero} 
The first systematic result is the next theorem of Yamaura, which provides a complete characterisation of Gorenstein algebras with tilting objects in dimension zero. 
\begin{thm}[Yamaura, {\cite[Thm. 3.1]{Ya}}] Let $A = \bigoplus_{n \geq 0} A_n$ be a finite dimensional graded self-injective algebra $A$, so that $\underline{{\rm MCM}}^{\mathbb{Z}}(A) = {\rm \underline{\mathsf{mod}}}^{\Z} A$. The following are equivalent:
	\begin{enumerate}[i)]
	\item $\underline{\mathsf{mod}}^{\Z} A$ admits a tilting object $T$.
	\item $\gldim A_0 < \infty$.
\end{enumerate} 
When this holds, letting $a \geq 0$ be the maximal degree of $A = \bigoplus_{n = 0}^a A_n$, we may take $T = \bigoplus_{i = 0}^{a-1} A(i)_{\leq 0} := \bigoplus_{i = 0}^{a-1} A/A_{\geq i+1}(i)$. Letting $\Lambda = {\rm \underline{End}}_{gr A}(T)$, we have $\gldim \Lambda < \infty$ and we have an isomorphism of algebras
\[
	\Lambda \cong \begin{pmatrix} 
		A_0 & A_1 & \cdots & A_{a-2} & A_{a-1} \\ 
		 & A_0 & \cdots & A_{a-3} & A_{a-2} \\ 
		&  & \ddots & \vdots & \vdots \\ 
		&  & & A_{0} & A_{1} \\ 
		 & & &  & A_{0} \\ 
	\end{pmatrix}
\] 
\end{thm}

\begin{exmp} Let $A = A_0 \oplus A_1 \oplus A_2 = k \oplus V \oplus k$, with multiplication given by a perfect pairing $V \times V \to k$, extended as a unital graded algebra structure to all of $A$. It is easy to see that $A$ is self-injective, and since $A_0 = k$, Yamaura's result applies to give a tilting object with endomorphism algebra
	\[
	\Lambda = \begin{pmatrix} 
		k & V\\ 
		0 & k\\ 
	\end{pmatrix}
	\]
	Equivalently, if $n = \dim V$ with basis $\{ x_i \}$, then $\Lambda \cong kQ_n$ for $Q_n$ the $n$-Kronecker quiver
	\[
		\xymatrix{ \bullet \ar@<1ex>[r]^{ \{x_i\} } \ar[r] \ar@<-1ex>[r] & \bullet}
	\]
\end{exmp} 

\begin{exmp} Let $A = \bigwedge^*V = \bigwedge^*_k(y_1, \dots, y_n)$ be the exterior algebra on an $n$-dimensional vector space. Then $\Lambda$ is given by
\[
	\Lambda = \begin{pmatrix} 
		k & V & \cdots & \bigwedge^{n-2}V & \bigwedge^{n-1}V \\ 
		& k & \cdots & \bigwedge^{n-3}V & \bigwedge^{n-2}V \\ 
		&  & \ddots & \vdots & \vdots \\ 
		&  & & k & V \\ 
		 & & &  & k \\ 
	\end{pmatrix}
\] 
Equivalently, $\Lambda \cong kQ/I$ is given by the quiver path algebra of $Q$ 
\[
		\xymatrix{ 
		\bullet \ar@<1ex>[r]^{ \{y_i\} } \ar[r] \ar@<-1ex>[r] & \bullet \ar@<1ex>[r]^{ \{y_i\} } \ar[r] \ar@<-1ex>[r] & \ar@<1ex>[r]^{ \{y_i\} } \ar[r] \ar@<-1ex>[r] \bullet &
 \cdots \ar@<1ex>[r]^{ \{y_i\} } \ar[r] \ar@<-1ex>[r] & \bullet \ar@<1ex>[r]^{ \{y_i\} } \ar[r] \ar@<-1ex>[r] & \bullet  
	}
\]
with $n-1$ vertices, and relations on paths of length two $y_iy_j + y_jy_i = 0$ for $i \neq j$, and $y_iy_i = 0$. 
\end{exmp}

\begin{exmp}[Happel's Theorem, {\cite[Ex. 3.15]{Ya}}] Using his theorem, Yamaura gave a simple proof of a theorem of Happel. Let $\Lambda$ be a finite dimensional $k$-algebra with $\gldim \Lambda < \infty$. Define $T(\Lambda) = \Lambda \rtimes D\Lambda$ to be the trivial extension algebra of $\Lambda$ by the bimodule $D\Lambda = {\rm Hom}_k(\Lambda, k)$, with multiplication
	\[
		(x, \varphi) (x', \varphi') = (xx', x\varphi' + \varphi x').
	\]
	Grade $T(\Lambda)$ as $T(\Lambda)_0 = \Lambda$, $T(\Lambda)_1 = D\Lambda$ and $T(\Lambda)_n = 0$ otherwise. Then $T(\Lambda)$ is a graded self-injective algebra with $\gldim T(\Lambda)_0 < \infty$, and $T = \Lambda$ is a tilting object in $ \underline{\mathsf{mod}}^\Z T(\Lambda)$. It follows that we have an equivalence of triangulated categories
	\[
		{\rm D}^b(\Lambda) \cong \underline{\mathsf{mod}}^\Z \ \! T(\Lambda)
	\]
	which is the content of Happel's Theorem.
\end{exmp}

\subsubsection*{Dimension one}
The next theorem due to Buchweitz-Iyama-Yamaura \cite{BIY} establishes the existence of tilting objects in dimension one, at the cost of imposing a few conditions. We briefly review their results under the simplest hypotheses, and refer to \cite{BIY} for the general results. \bigskip

Here we let $R = \bigoplus_{n \geq 0} R_n$ be a graded connected, finitely generated, commutative reduced Gorenstein $k$-algebra of Krull dimension one. One may see \cite[Lemma 4.9(a)]{BIY} that any such algebra with $a < 0$ must be isomorphic to $R \cong k[t]$ with $|t| = -a > 0$, and so without loss of generality we impose $a \geq 0$. Let $K = Q(R)$ be the (homogeneous) total field of fraction of $R$, obtained by inverting all homogeneous non zero-divisors. We have an isomorphism of graded algebras 
\[
	K \cong \prod_{j=1}^r k[t_i, t_i^{-1}]
\]
with $|t_i| \geq 1$. The number $r$ is the number of branches of the singularity $\spec R$. Let $p \geq 1$ be the least integer degree of a non zero-divisor $x \in \m = R_{\geq 1}$. At least one non zero-divisor exists since $R$ is Gorenstein and in particular Cohen-Macaulay, and by a prime avoidance argument one can assume that $x \in \m \setminus \m^2$, and so one sees that $p = 1$ whenever $R$ is generated in degree one. 
\begin{thm}[Buchweitz-Iyama-Yamaura, \cite{BIY}]\label{BIYtilting} The category $\underline{{\rm MCM}}^{\mathbb{Z}}(R)$ admits a tilting object
	\[
		T = \bigoplus_{i = 1}^{a+p} R_{\geq i}(i) = \left(\bigoplus_{i=1}^a R_{\geq i}(i) \right) \bigoplus \left(\bigoplus_{i = a+1}^{a+p} R_{\geq i}(i)\right).
	\]
	Moreover, the objects $R_{\geq i}(i)$ for $1 \leq i \leq a$ are exceptional. The remaining module has semisimple endomorphism ring and decomposes into indecomposables as
	\[
		\bigoplus_{i=a+1}^{a+p} R_{\geq i}(i) = \bigoplus_{j=1}^r E_j^{\oplus e_j}
	\]
	with $r$ distinct indecomposables up to multiplicity. 
\end{thm}

\begin{exmp}[Buchweitz-Iyama-Yamaura, {\cite[Thm. 2.3(d)]{BIY}}] Let $f \in k[x,y]$ be a homogeneous polynomial of degree $d \geq 2$ with squarefree decomposition into irreducibles $f = f_1 f_2 \dots f_d$. Since $k$ is algebraically closed, all irreducible factors $f_i = a_ix + b_iy$ are linear forms. The ring $R = k[x,y]/(f)$ is a reduced Gorenstein ring with $a$-invariant $a = d-2$ whose spectrum gives the cone over $d$ points on $\mathbb{P}^1$, and so we have $r = d$ branches and moreover $p = 1$ since $R$ is generated by $R_1$. Writing $\m = R_{\geq 1}$, we have $T = \m(1) \oplus \m^2(2) \oplus \dots \oplus \m^{d-2}(d-2) \oplus \m^{d-1}(d-1)$, with endomorphism algebra ${\rm \underline{End}}_{gr R}(T) \cong kQ/I$ given by the quiver $Q$
	\[
\xymatrix@R=5pt{ 
	&&&&& \bullet \\
	&&&&& \bullet \\
	\bullet \ar@<.5ex>[r]^x \ar@<-.5ex>[r]_y & \bullet \ar@<.5ex>[r]^x \ar@<-.5ex>[r]_y & \cdots \ar@<.5ex>[r]^x \ar@<-.5ex>[r]_y & \bullet \ar@<.5ex>[r]^x \ar@<-.5ex>[r]_y & \bullet \ar[uur]^{p_1} \ar[ur]_{p_2} \ar[ddr]^{p_d} & \vdots \\
	&&&&& \vdots \\
	&&&&& \bullet 
}
	\]
	with a `tail' of length $a = d-2$ and $r = d$ many vertices to the right, with relations $xy = yx$ and $p_i (a_ix + b_iy) = 0$. 
 \end{exmp}

 We will make use of special cases of this theorem in Chapter 3, to investigate the construction of MCM modules over certain curve singularities of tame CM representation type. In particular, and for completeness we will prove special cases of Thm. \ref{BIYtilting}, although in cases which circumvent the main difficulties in the above theorem.\bigskip

Finally, we refer to \cite{IyamaSurvey} for more examples of graded Gorenstein algebras admitting tilting objects.
 

%% file: chapters/tilting/agv.tex
\section{Hodge theory obstruction in higher dimension}\label{hodgesection}

Let $R = k \oplus R_1 \oplus R_2 \oplus \dots $ be a graded connected, Gorenstein $k$-algebra, so that in particular $R$ is graded local with homogeneous maximal ideal $\m = R_+$. Recall that $R$ has (graded) isolated singularities if $R_\mathfrak{p}$ is regular for each homogeneous prime $\mathfrak{p} \neq \m$. Zero dimensional rings vacuously have isolated singularities, and one dimensional rings have isolated singularities if and only if they are reduced. \bigskip

We can restate the results of the previous section as follows. 
\begin{prop}[Yamaura, Buchweitz-Iyama-Yamaura] Let $R$ be a graded Gorenstein algebra satisfying $\dim R \leq 1$, with at most isolated singularities. Then $ \underline{{\rm MCM}}^{\mathbb{Z}}(R)$ has a (canonical) tilting object $T$. 
\end{prop}

Going up in dimension, one has the natural question:
\begin{quest} Let $R$ be a graded Gorenstein algebra, with $\dim R \geq 2$ and at most isolated singularities. When does $ \underline{{\rm MCM}}^{\mathbb{Z}}(R)$ have a tilting object? 
\end{quest}

The tilting problem for singularity categories is closely related to the tilting problem for $ {\rm D}^b(X) $ on an algebraic variety $X$, which is known to be heavily obstructed. Similarly, it isn't hard to find a Gorenstein algebra $R$ with no tilting object. 

\begin{exmp} Let $R = k[x,y,z]/(f)$ for $f$ a cubic polynomial with isolated singularities at the origin, and $R$ with the standard grading. Then $a = |f| - |x| - |y| - |z| = 3 - 3 = 0$ and so $\omega_R \cong R$. Therefore the Serre functor for $ \underline{{\rm MCM}}^{\mathbb{Z}}(R)$ is given by 
\[
	\mathbb{S}_R(-) = - \otimes_R \omega_R [\dim R - 1] = (-)[1]
\]
and so $ \underline{{\rm MCM}}^{\mathbb{Z}}(R)$ is $1$-Calabi-Yau. Thus there can be no equivalence $\underline{{\rm MCM}}^{\mathbb{Z}}(R) \cong {\rm D}^{\perf}(\Lambda)$ since projective modules in $\modsf \Lambda \subset {\rm D}^{\perf}(\Lambda)$ cannot exist in an $n$-Calabi-Yau category for $n \neq 0$. \qed
\end{exmp}
In the above example note that $E = \proj R$ is a plane elliptic curve, and by Orlov's theorem, since $E$ is Calabi-Yau have an equivalence of categories $\underline{{\rm MCM}}^{\mathbb{Z}}(R) \cong {\rm D}^b(E)$. The tilting problem for ${\rm D}^b(E)$ is obstructed for the same reason (and several more, see below). \newpage

\noindent For the remainder of Chapter $2$, we will impose the following simplifying assumption: 
\begin{assump} All graded $k$-algebras $R$ are assumed standard graded, meaning that $R_0 = k$ and $R = R_0[R_1]$ is generated in degree one over $R_0 = k$. 
\end{assump} 
In spite of the name, we note that this is not a standard hypothesis in tilting theory, since tilting MCM modules exist more readily over non-standard graded $k$-algebras (such as ADE hypersurface rings or invariant rings). However, letting $X = \proj R$ denote the underlying projective scheme, the above assumption will allow for an easier transfer of information between $R$ and $X$. \bigskip

In general, one may conjecture the following:\footnote{This can be stated more generally but we will stick to these hypotheses.} 
\begin{conj} Let $R$ be a standard graded Gorenstein algebra with isolated singularities. Then the following are equivalent:
	\begin{enumerate}[1)]
	\item $ \underline{{\rm MCM}}^{\mathbb{Z}}(R)$ admits a tilting object $T$.  
	\item $ {\rm D}^b(X) $ admits a tilting object $\E$.
\end{enumerate}
\end{conj}
The author has heard this conjecture from O. Iyama, who credits it to G. Stevenson. Using Orlov's semiorthogonal decompositions and making use of the recent theory of `additive invariants', one may observe that these two categories share essentially the same principal obstructions to tilting, such as $K_0$ being finitely generated free or the category not being Calabi-Yau (of dimension $n \neq 0$). In algebraic geometry, one of the most useful obstructions to tilting is the Hodge obstruction. Recall that we denote by $h^{p,q}(X) = \dim_k {\rm H}^q(X, \Omega^p_X)$ the bigraded Hodge numbers of an algebraic variety $X$. 

\begin{prop}[{\cite[Thm. 5.2, Cor. 4.2]{BuHi}}]\label{hodgeobstructionDbX} Let $k$ be a field and let $X$ be a smooth projective variety over $k$. If $ {\rm D}^b(X) $ admits a tilting object $\E$, then $h^{q,0}(X) = 0$ for $q > 0$. Moreover, if the characteristic of $k$ is zero, then $h^{p,q}(X) = 0$ for $p \neq q$.
\end{prop}

Using standard methods, it isn't hard to prove an analogous result for $ \underline{{\rm MCM}}^{\mathbb{Z}}(R)$, which should be at least expected if not well-known to experts. We will show: 
\begin{prop*}[Hodge obstruction for MCM modules, Prop. \ref{hodgeobstructionmcm}] Let $k$ be an algebraically closed field and $R$ be a standard graded Gorenstein $k$-algebra of $\dim R \geq 2$, having at most isolated singularities. Let $X = \proj R$. Assume that $ \underline{{\rm MCM}}^{\mathbb{Z}}(R)$ admits a tilting object $T$. Then: 
	\begin{enumerate}[i)]
	\item $h^{q,0}(X) = 0$ for $q > 0$. 
	\item The $a$-invariant of $R$ satisfies $a < 0$. In particular $X$ is a Fano variety. 
\end{enumerate}
Moreover, if the characteristic of $k$ is zero, then $i)$ may be strenghthened to: 
\begin{enumerate}[i')]
\item $h^{p,q}(X) = 0$ for $p \neq q$.
\end{enumerate}
\end{prop*}
The result essentially follows from Hochschild homology computations. One may deduce it from the additivity properties of Hochschild homology, which falls under the wide-reaching umbrella of the recent theory of additive invariants introduced by Tabuada \cite{TabuadaBook,TabuadaThesis, MT}.\bigskip

However, a recent embedding result of Orlov \cite{OrlovGeomReal}, which we will apply under the mild assumption that $k$ be algebraically closed\footnote{One can likely weaken this assumption, but this will be enough for us.}, allows us to avoid the use of additive invariants, relying instead on more direct arguments, at the cost of stating various facts in a weaker form that is nevertheless sufficient for our purposes. \bigskip

For $R$ standard graded Gorenstein with isolated singularities, the necessary condition $a < 0$ when $\dim R \geq 2$ is an interesting reversal from lower dimension, since we have seen that $a \geq 0$ for whenever $\dim R \leq 1$, unless $R$ is regular \cite[Lemma 4.9(a)]{BIY}. Interestingly, that last result holds also for non-standard graded $R$ with $\dim R \leq 1$. The author ignores if the condition $a < 0$ is still necessary in the case $\dim R \geq 2$, if the assumption that $R$ be generated in degree one is dropped. \bigskip

The remaining of Section \ref{hodgesection} will be devoted to the proof of Prop. \ref{hodgeobstructionmcm}. 

\subsection*{Hochschild homology of schemes}

Fix a field $k$ throughout. All $k$-varieties will be assumed quasi-projective, and we fix such an $X$ to start and set $d = \dim X$. We first review standard background following \cite{BuHi, Swan, KuzHH}. \bigskip 

Let $\Delta: X \hookrightarrow X \times X$ be the diagonal embedding and let $\OO_\Delta = \Delta_* \OO_X$ be the structure sheaf of the diagonal. Let $\mathcal{M} \in {\rm D}(X \times X)$.
\begin{defn}[{\cite[3.1]{BuHi}}] We define the Hochschild cohomology and homology of $X$ with coefficients in $\mathcal{M}$ by:

	\begin{enumerate}[1)]
	\item ${\rm HH}^*(X, \mathcal{M}) = {\rm Ext}^*_{X \times X}(\mathcal{O}_{\Delta}, \mathcal{M})$.
	\item ${\rm HH}_*(X, \mathcal{M}) = {\rm H}^{-*}(X, \mathbb{L}\Delta^* \mathcal{M})$.
\end{enumerate}
Note that when $\mathcal{M}$ is a quasicoherent sheaf we have ${\rm HH}^i(X, \mathcal{M}) = 0$ for $i < 0$, and since $\mathbb{L}\Delta^* \mathcal{M} \in {\rm D}^-(X)$ we also obtain ${\rm HH}_j(X, \mathcal{M}) = 0$ for $j < -\dim X$ by Grothendieck vanishing. When $\mathcal{M} = \OO_\Delta$, we write ${\rm HH}^*(X)$ and ${\rm HH}_*(X)$ respectively.
\end{defn}
When $X = \spec R$ is affine, the diagonal embedding $\Delta: X \to X \times X$ is induced from the multiplication map $R^{ev} := R \otimes_k R \to R$, and so Hochschild (co)homology specialises to the usual definition of ${\rm HH}^*(X, \mathcal{M}) = {\rm Ext}^*_{R^{ev}}(R, \mathcal{M})$ and ${\rm HH}_*(X, \mathcal{M}) = {\rm Tor}_*^{R^{ev}}(R, \mathcal{M})$. \bigskip

When $X$ is smooth affine over $k$ of characteristic zero, the classical Hochschild-Kostant-Rosenberg (HKR) Theorem asserts that
\[
	{\rm HH}_i(X) = {\rm H}^0(X, \Omega^i_X).
\]
R. Swan (\cite{Swan}, \cite[Rem 6.3]{Huy}) has constructed, for general $X$ smooth, two Hodge-to-Hochschild spectral sequences
\begin{align*}
	{\rm E}^{p,q}_2 &= {\rm H}^p(X, \wedge^q \T_X) \implies {\rm HH}^{p+q}(X)\\
	{\rm E}^{p,q}_2 &= {\rm H}^p(X, \Omega^{d-q}_X) \implies {\rm HH}_{p+q-d}(X)
\end{align*}
the second of which resembles the Hodge-to-DeRham spectral sequence when $k = \C$. 

\begin{thm}[{\cite[Cor. 2.6]{Swan}}]\label{HKR} When $X$ is smooth over $k$ of characteristic zero, the above spectral sequences degenerate at $E_2$. After reindexing, this gives isomorphisms: 
	\begin{align*}
	{\rm HH}^i(X) &\cong \bigoplus\limits_{p+q=i}{\rm H}^q(X, \wedge^p \T_X)\\
	{\rm HH}_i(X) &\cong \bigoplus\limits_{p-q=i}{\rm H}^q(X, \Omega^p_X).
\end{align*}
\end{thm}
The second decomposition is often referred to as the general HKR Theorem. In general characteristic, we have the weaker but useful lemma of Buchweitz-Hille, whose proof we give for completeness. 

\begin{lem}[Buchweitz-Hille, {\cite[Cor. 4.2]{BuHi}}]\label{Hodgesummand} Let $k$ be a field and $X$ a $k$-variety. Then ${\rm H}^q(X, \OO_X)$ embeds as a direct summand in ${\rm HH}_{-q}(X)$.
\end{lem}
\begin{proof}
	We have an adjoint pair $(\mathbb{L}\Delta^*, \Delta_*)$, and note that the counit map $\mathbb{L}\Delta^*(\OO_\Delta) = \mathbb{L}\Delta^*(\Delta_* \OO_X) \to \OO_X$ splits after applying $\Delta_*$ as in any adjunction. The vector space
			\begin{align*}
				{\rm H}^q(X \times X, \Delta_*\OO_X) &= {\rm Hom}_{ {\rm D}^{-}(X \times X)}(\OO_{X \times X}, \Delta_*\OO_X[q])\\
				&= {\rm Hom}_{ {\rm D}^{-}(X)}(\mathbb{L}\Delta^*(\OO_{X \times X}), \OO_X[q])\\
				&= {\rm H}^q(X, \OO_X)
			\end{align*}
			is then a summand of
	\begin{align*}
		{\rm H}^q(X \times X, \Delta_* (\mathbb{L}\Delta^*(\Delta_*\OO_X))) &= {\rm Hom}_{ {\rm D}^{-}(X \times X)}(\OO_{X \times X}, \Delta_* (\mathbb{L}\Delta^*(\Delta_* \OO_X))[q]) \\
		&= {\rm Hom}_{ {\rm D}^{-}(X)}(\mathbb{L}\Delta^*(\OO_{X \times X}), \mathbb{L}\Delta^*(\Delta_* \OO_X)[q])\\
		&= {\rm Hom}_{ {\rm D}^{-}(X)}(\OO_X, \mathbb{L}\Delta^*(\Delta_*\OO_X)[q])\\
		&= {\rm HH}_{-q}(X).
			\end{align*}
\end{proof}

The advantage of working with Hochschild homology lies in its derived invariance, giving flexibility in choosing models to compute it. Here are a few standard applications. 
\begin{exmp}[{\cite[Prop. 5.39, Rem. 6.3]{Huy}}] Let $X, Y$ be two smooth projective varieties with $ {\rm D}^b(X) \cong {\rm D}^b(Y)$. Then we have ${\rm HH}_*(X) \cong {\rm HH}_*(Y)$, and so in characteristic zero this implies
	\[
		\sum\limits_{p-q=i} h^{p,q}(X) = \sum\limits_{p-q = i} h^{p,q}(Y).
	\]\qed
\end{exmp}

\begin{exmp}[{\cite[Thm. 4.1, Cor. 4.2]{BuHi}}]\label{HHvanishing} Let $X$ be a smooth projective variety with a tilting object $T \in {\rm D}^b(X) $, and $A = {\rm End}(T)$. Then $A$ is finite dimensional over $k$, $\gldim A < \infty$ and we have an isomorphism ${\rm HH}_*(X) \cong {\rm HH}_*(A)$. In particular, ${\rm HH}_i(X) = {\rm HH}_i(A) = 0$ for $i < 0$, and when $k$ is of characteristic zero we obtain ${\rm HH}_i(X) = 0$ for all $i \neq 0$ by making use of the HKR Theorem and the Hodge symmetries $h^{p,q}(X) = h^{q,p}(X)$. This is essentially the content of Prop. \ref{hodgeobstructionDbX}. \qed
 \end{exmp} 

 The tilting hypothesis on $T$ only plays a role in showing ${\rm HH}_i(X) = 0$ for $i < 0$. Working instead with any classical generator $T$ and setting $\A = {\rm RHom}_{\OO_X}(T, T)$, the proof of Buchweitz-Hille actually shows the following, which is well-known folklore. Recall that we define the Hochschild homology of a dg $k$-algebra $\A$ by ${\rm HH}_*(\A) := {\rm HH}_*(\A/k, \A) := {\rm Tor}_*^{\A^{ev}}(\A, \A)$, where $\A^{ev} = \A^{op} \otimes_k \A$.
 \begin{prop} Let $X$ be a smooth projective variety with classical generator $T \in {\rm D}^b(X)$ and differential graded algebra $\A = {\rm RHom}_{\OO_X}(T, T)$. Then we have an isomorphism ${\rm HH}_*(X) \cong {\rm HH}_*(\A)$ of Hochschild homologies. 
 \end{prop} 
 One does not need $T$ to be tilting to have control over the groups ${\rm HH}_i(X)$. Define an object $T \in {\rm D}^b(X) $ to be \emph{silting} if it satisfies the following weaker conditions:
 \begin{enumerate}[i)]
 \item $T$ is a classical generator for $ {\rm D}^b(X) $;
 \item ${\rm Ext}^i_{\OO_X}(T, T) = 0$ for all $i > 0$. 
	 \end{enumerate}
	 It follows that ${\rm RHom}_{\OO_X}(T, T)$ may be replaced by a quasi-isomorphic dg algebra, which we denote again by $\A$ by abuse of notation, with the property that $\A^i = 0$ for $i > 0$:
	 \[
		 \cdots \xrightarrow{d} \A^{-(n+1)} \xrightarrow{d} \A^{-n} \xrightarrow{d} \cdots \xrightarrow{d} \A^{-1} \xrightarrow{d} \A^{0} \to 0
	 \]
	 Regrading $\A^{-n} = \A_{n}$, we observe that $\A = (\A_n)_{n \geq 0}$ is a non-negatively graded homological dg algebra. For such algebras, one has 
	 \begin{lem} Let $\A = \bigoplus_{n \geq 0} \A_n$ be a non-negatively graded homological dg algebra. Then ${\rm HH}_i(\A) = 0$ for $i < 0$.
	 \end{lem}
	 \begin{proof}
		 We have ${\rm HH}_i(\A) = {\rm H}_i( \A \otimes^{\mathbb{L}}_{\A^{ev}} \A) = {\rm H}_i( P \otimes_{\A^{ev}} \A)$ where $P \xrightarrow{\sim} \A$ is any h-projective resolution of $\A$ over $\A^{ev}$. Since $\A = \A_{\geq 0}$ and $\A^{ev} = \A^{ev}_{\geq 0}$, it is easy to construct $P$ with the property $P = P_{\geq 0}$ (e.g. taking $P$ to be the two-sided Bar resolution), and so ${\rm H}_i(P \otimes_{\A^{ev}} A) = 0$ for $i < 0$.
	 \end{proof}

	 \begin{cor}\label{hodgesilting} Let $X$ be a smooth projective variety over $k$ with a silting object $T \in {\rm D}^b(X) $. Then ${\rm HH}_i(X) = 0$ for $i < 0$, and moreover ${\rm HH}_i(X) = 0$ for $i \neq 0$ when the characteristic of $k$ is zero. 
	 \end{cor}

	 Before we move to establish the Hodge obstruction result for tilting objects in graded singularity categories $ \underline{{\rm MCM}}^{\mathbb{Z}}(R)$, we will need a few standard facts concerning Fourier-Mukai transforms. Let $X, Y$ denote smooth projective varieties over $k$ throughout. 
\begin{prop}[Orlov, Bondal-Van den Bergh {\cite[Thm. 5.14]{Huy}}]\label{FMT} Let $F: {\rm D}^b(X) \to {\rm D}^b(Y)$ be a fully faithful functor. Then $F$ is naturally isomorphic to a Fourier-Mukai transform $F \simeq \Phi_K$. 
\end{prop}

\begin{prop}[{\cite[Lemma 6.5]{KuzHH}}] Let $\Phi_K: {\rm D}^b(X) \to {\rm D}^b(Y)$ be a Fourier-Mukai transform. Then there is an induced linear map on Hochschild homology
\[
	\Phi_K^{ {\rm HH}}: {\rm HH}_*(X) \to {\rm HH}_*(Y).
\]
Moreover, the above satisfies the following properties:
\begin{enumerate}[i)]
\item It is natural, in that $\Phi_K^{HH} \circ \Phi_L^{HH} = \Phi^{HH}_{K \circ L}$, where $K \circ L$ is the convolution product of kernels \cite{Huy, KuzHH}. 
\item Naturally isomorphic transforms $\Phi_K \simeq \Phi_{K'}$ give rise to equal linear maps $\Phi^{HH}_K = \Phi^{HH}_{K'}$.
\item Let $K = \OO_\Delta$ so that $\Phi_{\OO_\Delta} \simeq {\rm id}_{ {\rm D}^b(X) }$. Then $\Phi_{\OO_\Delta}^{HH} = {\rm id}_{ {\rm HH}_*(X)}$. 
	\end{enumerate}
\end{prop}
\begin{proof}
The main claim along with $i)$ is \cite[Lemma 6.5]{KuzHH}, with the claims $ii)$ and $iii)$ implicit and easy to see. Note that Kuznetsov uses a different model for ${\rm HH}_*(-)$, and the equivalence with our definition is given in \cite[Prop. 8.1]{KuzHH}. 
\end{proof}
\begin{rem}
	Working over $k = \C$ and using the Hodge decomposition, one can instead work with the cohomology groups ${\rm H}^*(X; \C)$. Appropriate linear maps $\Phi_K^{H}: {\rm H}^*(X; \C) \to {\rm H}^*(Y; \C)$ were constructed in \cite[Chp. 5, Lemma 5.32, Prop. 5.33]{Huy}, and one could use these instead in all arguments below. 
\end{rem}

The following argument is due to Kiem and Lee over $k = \C$ but their argument applies equally to our situation, and we reproduce it here for convenience. 
\begin{prop}[{\cite[Prop. 4.7]{KL}}]\label{FMTHH} Assume that the Fourier-Mukai transform $\Phi_K: {\rm D}^b(X) \to {\rm D}^b(Y)$ is a fully faithful. Then the induced map $\Phi_K^{ {\rm HH}}$ is split-injective. 
\end{prop}
\begin{proof}
	By \cite[Prop. 5.9]{Huy}, the Fourier-Mukai transform $\Phi_K$ admits a right adjoint which is also a Fourier-Mukai transform, which we denote $\Phi_{K_R}$. By \cite[Cor. 1.22]{Huy}, the natural morphism ${\rm id}_{\rm D^b(X)} \to \Phi_{K_R} \circ \Phi_{K}$ is an isomorphism, and so we obtain 
	\[
		\Phi_{\OO_{\Delta}} \simeq {\rm id}_{ {\rm D}^b(X) } \simeq \Phi_{K_R} \circ \Phi_{K} \simeq \Phi_{K_R \circ K}.
	\]
	We then obtain 
	\[
		{\rm id}_{ {\rm HH}_*(X)} = \Phi^{HH}_{\OO_{\Delta}} = \Phi^{HH}_{K_R \circ K} = \Phi^{HH}_{K_R} \circ \Phi^{HH}_{K} 
	\]
	and so $\Phi^{HH}_K$ is split-injective. 
\end{proof}
Finally, we will make use of the recent `geometric realisation' theorem of Orlov. The following is a special case: 

\begin{prop}[Orlov \cite{OrlovGeomReal}]\label{Orlovembedding} Let $\Lambda$ be a finite dimensional $k$-algebra over an algebraically closed field $k$. Then there is a fully faithful exact functor 
	\[
		F: {\rm D}^{\perf}(\Lambda) \hookrightarrow {\rm D}^b(Y) 
	\]
	for some smooth projective $k$-variety $Y$. Moreover: 
	\begin{enumerate}[i)]
	\item When $\gldim \Lambda < \infty$, $F$ can be taken to have both adjoints. 
	\item $ {\rm D}^b(Y)$ contains a full strong exceptional collection of line bundles. 
\end{enumerate}
\end{prop}
\begin{rem} In particular by Ex. \ref{HHvanishing} we have ${\rm HH}_i(Y) = 0$ for $i < 0$, which improves\footnote{Of course, using Kuznetsov \cite{KuzHH} and Keller's work on additivity of Hochschild homology, one can obtain ${\rm HH}_i(Y) = 0$ for $i \neq 0$ independent of characteristic, which follows from the existence of a full exceptional collection alone. However the above will suffice.} to ${\rm HH}_i(Y) = 0$ for $i \neq 0$ in characteristic zero. 
\end{rem} 
\begin{proof}
	This is a special case of \cite[Thm. 5.2, Thm. 5.8]{OrlovGeomReal}. Let $J = \mathsf{rad}(\Lambda)$, and let $n \in \N$ be the smallest $n$ such that $J^n = 0$. Define $M = \bigoplus_{i=1}^n \Lambda/J^i$ and let $\Gamma = {\rm End}_{\Lambda}(M)$, which satisfies the following by a theorem of Auslander \cite[Thm. 5.1]{OrlovGeomReal}:
	\begin{enumerate}[1)]
	\item $\gldim \Gamma \leq n+1$;
	\item The $\Lambda-\Gamma$-bimodule $P = {\rm Hom}_{\Lambda}(M, \Lambda)$ is a finite projective module over $\Gamma$ satisfying ${\rm End}_{\Gamma}(P) \cong \Lambda$. 
\end{enumerate}
By \cite[Thm. 5.2]{OrlovGeomReal}, there is a full exceptional collection
	\[
		{\rm D}^b(\Gamma) = \big\langle E_1, \dots, E_N \big\rangle 
	\] 
	for some $N \in \N$ (here we are using a simplification afforded by $k$ being algebraically closed). Orlov then deduces the existence of a fully faithful embedding
	\[
		- \otimes_{\Lambda} P:	{\rm D}^{\perf}(\Lambda) \hookrightarrow {\rm D}^b(\Gamma).
	\]
	Lastly, by \cite[Thm. 5.8]{OrlovGeomReal}, for any (small, enhanced) triangulated category $\T$ with a full exceptional sequence, such as $\T = {\rm D}^b(\Gamma) $, there is a fully faithful embedding
	\[
		\T \hookrightarrow {\rm D}^b(Y) 
	\]
	for $Y$ a smooth projective variety, constructed by taking iterated projective bundles $Y = Y_k \to Y_{k-1} \to \dots \to Y_1 \to Y_0 = \mathbb{P}^m$. In particular $ {\rm D}^b(Y) $ also has a full strong exceptional collection of line bundles by standard results \cite{CMR}. Combining these embeddings gives the result. 
\end{proof}

After all this setup, we are now ready to establish Prop. \ref{hodgeobstructionmcm}. We restate it for convenience.
\begin{prop}[Hodge obstruction for MCM modules]\label{hodgeobstructionmcm} Let $k$ be an algebraically closed field and $R$ be a standard graded Gorenstein $k$-algebra of $\dim R \geq 2$, having at most isolated singularities. Let $X = \proj R$. Assume that $ \underline{{\rm MCM}}^{\mathbb{Z}}(R)$ admits a tilting object $T$. Then:
	\begin{enumerate}[i)]
	\item $h^{q,0}(X) = 0$ for $q > 0$.
	\item The $a$-invariant of $R$ satisfies $a < 0$. In particular $X$ is a Fano variety.
\end{enumerate}
Moreover, if $k$ has characteristic zero, then $i)$ may be strenghthened to:
\begin{enumerate}[i')]
\item $h^{p,q}(X) = 0$ for $p \neq q$.
\end{enumerate}
\end{prop}
\begin{proof}
We first prove $i)$ and $i')$. Under the hypothesis, we have $\underline{{\rm MCM}}^{\mathbb{Z}}(R) \cong {\rm D}^{\perf}(\Lambda)$ for $\Lambda = {\rm \underline{End}}_{gr R}(T)$. We split the argument according to $a > 0$, $a = 0$ or $a < 0$ and use Orlov's semiorthogonality decomposition theorem. 
	\begin{itemize}
		\item[$a>0$:] We have a fully faithful embedding $\Psi_0: {\rm D}^b(X) \hookrightarrow \underline{{\rm MCM}}^{\mathbb{Z}}(R)$, which we compose with the fully faithful functor $ \underline{{\rm MCM}}^{\mathbb{Z}}(R) = {\rm D}^{\perf}(\Lambda) \hookrightarrow {\rm D}^b(Y)$ of Prop. \ref{Orlovembedding} (since $k$ is algebraically closed) to obtain an embedding 
			\[
				F:{\rm D}^b(X) \hookrightarrow {\rm D}^b(Y).
			\]
			By Prop. \ref{FMT}, $F \simeq \Phi_K$ is naturally isomorphic to a Fourier-Mukai transform and applying Prop. \ref{FMTHH} gives an embedding 
			\[
				\Phi_K^{HH}: {\rm HH}_i(X) \hookrightarrow {\rm HH}_i(Y)
			\]
	with ${\rm HH}_i(Y) = 0$ for $i < 0$ (and $i \neq 0$ in characteristic zero), and so the same holds for ${\rm HH}_i(X)$. Making use of Lemma \ref{Hodgesummand}, the HKR Theorem and the Hodge symmetries, we obtain $i)$ and $i')$. 
		\item[$a=0$:] This case is vacuous as $ \underline{{\rm MCM}}^{\mathbb{Z}}(R)$ is $d$-Calabi-Yau for $d \geq 1$, and so admits no tilting object. 
		\item[$a<0$:] In this case we have an embedding $ \Phi_0 :\underline{{\rm MCM}}^{\mathbb{Z}}(R) \hookrightarrow {\rm D}^b(X) $ and semiorthogonal decomposition
			\[
				{\rm D}^b(X) = \big\langle E_k, E_{k-1}, \dots, E_1, E_0, \Phi_0 ( \underline{{\rm MCM}}^{\mathbb{Z}}(R)) \big\rangle
			\]
			where $E_i = \OO_X(-i)$ and $k = |a|-1$. We will show that the tilting object $\Phi_0(T) \in \Phi_0( \underline{{\rm MCM}}^{\mathbb{Z}}(R))$ extends to a silting object of $ {\rm D}^b(X) $. Semiorthogonality means that 
			\[
				{\rm Hom}_{ {\rm D}^b(X) }(\Phi_0(T), E_i[n]) = 0
			\]
			for all $i = 0, 1 \dots, k$ and all $n \in \Z$, and similarly 
			\[
				{\rm Hom}_{ {\rm D}^b(X) }(E_i, E_j[n]) = 0
			\]
			for all $0 \leq i < j \leq k$ and $n \in \Z$. Moreover, none of the objects $E_i$ and $\Phi_0(T)$ have any positive self-extensions. \bigskip
			
			Now, $ {\rm D}^b(X) $ is Ext-finite and so ${\rm Hom}_{ {\rm D}^b(X) }(\F, \G[m]) = 0$ for $m \gg 0$ for any $\F, \G$. It follows that there is an $n_0 \in \N$ such that
			\[
				{\rm Hom}_{ {\rm D}^b(X) }(E_0, \Phi_0(T)[n]) = 0 \textup{ for all } n > n_0
			\]
		or equivalently
		\[
			{\rm Ext}^n(E_0[-n_0], \Phi_0(T)) = 0 \textup{ for all } n > 0.
		\]
		Similarly, there is an $n_1 \gg n_0$ such that
		\[
			{\rm Ext}^n(E_1[-n_1], E_0[-n_0] \oplus \Phi_0(T)) = 0 \textup{ for all } n > 0.
		\]
Continuing this way, we obtain a sequence $n_0 \ll n_1 \ll \cdots \ll n_{k-1} \ll n_{k}$ such that upon setting 
		\[
		\widetilde{T} = \left(\bigoplus_{i = 0}^{k}E_i[-n_i]\right) \oplus \Phi_0(T) 
			\]
	we have ${\rm Ext}^n(\widetilde{T}, \widetilde{T}) = 0$ for all $n > 0$. Since $\widetilde{T}$ classically generated $ {\rm D}^b(X) $, $\widetilde{T}$ is a silting object. The claims $i)$ and $i')$ then follow from Cor. \ref{hodgesilting}. 
	\end{itemize} 
This proves parts $i)$ and $i')$. We claim that $i)$ implies $ii)$. To see this, using local cohomology for $R$ at $\m = R_+$ we have 
\[
	{\rm H}^q(X, \OO_X(n)) \cong {\rm H}^{q+1}_{\m}(R)_n
\]
for all $1 \leq q \leq d$ and $n \in \Z$. By graded Local duality we have ${\rm Ext}^i_R(R, \omega_R) \cong {\rm Hom}_k({\rm H}^{d+1-i}_\m(R), k)$ as graded $R$-modules, and so 
\[
	R(a) = \omega_R = {\rm Hom}_k({\rm H}^{d+1}_\m(R), k).
\]
In particular ${\rm H}^{d+1}_\m(R)_0 = {\rm H}^{d}(X, \OO_X) = 0 $ gives $R(a)_0 = R_a = 0$, and so $a < 0$. Lastly, $\OO_X(a) = \widetilde{\omega_R} \cong \omega_X$ shows that $\omega_X^{-1}$ is ample, and so $X$ is Fano. 
\end{proof}

%% file: chapters/tilting/ci.tex
\section{Cones over smooth projective complete intersections}\label{cisection} 
Complete intersections form the simplest class of (non-regular) Gorenstein algebras, and are a good starting point to investigate the converse of Prop. \ref{hodgeobstructionmcm}. Namely we are interested in the following question:
\begin{quest} Let $k$ be an algebraically closed field of characteristic zero, let $R$ be a standard graded Gorenstein $k$-algebra and $X = \proj R$. Assume that $h^{p,q}(X) = 0$ for $p \neq q$. Do the categories $ {\rm D}^b(X) $ and $ \underline{{\rm MCM}}^{\mathbb{Z}}(R)$ admit tilting objects? 
\end{quest}
Recall that we implicitly take graded to mean finitely generated in degree one. We assume that $k$ is algebraically closed of characteristic zero throughout this section. \bigskip 

The simplest class of graded Gorenstein rings on which to test this question are the complete intersections. Consider rings of the form $R = k[x_0, \dots, x_{n+c}]/(f_1, \dots, f_c)$ where $(f_1, \dots, f_c) \subseteq (x_0, \dots, x_{n+c})^2$ is a regular sequence of homogeneous polynomials, and let $\proj R = X = V(f_1, \dots, f_c) \subseteq \mathbb{P}^{n+c}$ the associated projective complete intersection. We will further assume that $R$ has isolated singularities at the origin, or equivalently that $X$ is smooth. \bigskip

A generating series for the Hodge numbers of $X$ was given by Hirzebruch \cite{HirzBook} as one of the first applications of the Hirzebruch-Riemann-Roch Theorem. Based on the above, Rapoport classified the smooth complete intersections $X$ with $h^{p,q}(X) = 0$ for $p \neq q$. 

\begin{prop}[Rapoport {\cite{Ra}}]\label{Rapoportlist} Let $X \subseteq \mathbb{P}^{n+c}$ be a smooth complete intersection with $\dim X \geq 1$ and ${\rm codim}\ \! X \geq 1$. Then $X$ satisfies $h^{p,q}(X) = 0$ for $p \neq q$ if and only if $X$ belongs to one of three families: 
	\begin{enumerate}[a)]
	\item $X = V(Q) \subseteq \mathbb{P}^{n+1}$ is a quadric hypersurface.
	\item $X = V(Q_1, Q_2) \subseteq \mathbb{P}^{2n+2}$ is a $2n$-dimensional intersection of two quadrics. 
	\item $X = X_3 \subseteq \mathbb{P}^3$ is a cubic surface. 
\end{enumerate}
\end{prop}
Using Rapoport's list, we will obtain the converse of Prop. \ref{hodgeobstructionmcm}; that is, we will deduce that the necessary vanishing of Hodge numbers is in fact sufficient to guarantee the existence of tilting objects for complete intersections. In fact, most of the following is already known from work of Kapranov (quadric hypersurfaces) and Orlov-Kuleshov (cubic surfaces), and with partial statements in the work of Kuznetsov (quadrics intersections). Our contribution will be to complete that last case, by showing the existence of a full strong exceptional collection of vector bundles on $X = V(Q_1, Q_2) \subseteq \mathbb{P}^{2n+2}$. \bigskip
 
The main theorem of this section will be: 
\begin{thm*}[Thm. \ref{hodgeconverse}] Let $k$ be algebraically closed of characteristic zero, let $R$ be a standard graded complete intersection $k$-algebra with isolated singularities, and $X = \proj R$. The following are equivalent: 
	\begin{enumerate}[1)]
	\item $ \underline{{\rm MCM}}^{\mathbb{Z}}(R)$ admits a tilting object $T$.
	\item $ {\rm D}^b(X) $ admits a tilting object $\E$. 
	\item We have $h^{p,q}(X) = 0$ for $p \neq q$.
\end{enumerate}
When $\dim R \geq 2$ and $R$ is not a polynomial algebra, this is equivalent to: 
\begin{enumerate}[1)]
\setcounter{enumi}{3}
\item $X$ belongs to one of the three familes $a)-b)-c)$. 
	\end{enumerate}
	Moreover, when these equivalent conditions hold, the tilting object $\E$ can always be assumed to come from a strong exceptional collection of vector bundles on $X$. 
\end{thm*} \bigskip

The proof of Thm. \ref{hodgeconverse} will occupy the rest of Chapter 2, and will be shown over the next ten or so pages. Let us lay out the structure of the argument: 
\begin{enumerate}[I.] 
\item We have seen the necessity implications $1) \implies 3)$ and $2) \implies 3)$; specifically, this is the content of Prop. \ref{hodgeobstructionmcm} and Prop. \ref{hodgeobstructionDbX}. 
\item The equivalence $3) \iff 4)$ under the stated assumptions is the result of Rapoport above in Prop. \ref{Rapoportlist}.
\item When $\dim R \leq 1$, $X$ is either empty or zero dimensional (a collection of points), and so the assumptions $3)$ always hold. Theorem \ref{BIYtilting} of Buchweitz-Iyama-Yamaura states the existence of a tilting object $T \in \underline{{\rm MCM}}^{\mathbb{Z}}(R)$, and so $1)$ also always holds. The statement $2)$ holds either vacuously (when $\dim R = 0$ so that $X$ is empty) or trivially as $\OO_X$ is a tilting bundle in $ {\rm D}^b(X) $ with orthogonal exceptional summands. 
\item Similarly, when $R$ is a polynomial algebra, the statement $1)$ is vacuous while statement $2)$ for $\proj R = X = \mathbb{P}^n$ follows from Beilinson's full strong exceptional collection of line bundles.  
\item It remains to show the implications $4) \implies 1)$ and $4) \implies 2)$, which we will do on a case-by-case basis. In brief: 
\begin{enumerate}[i.] 
\item The quadric hypersurface case $X = V(Q) \subseteq \mathbb{P}^{n+1}$ will follow from work of Kapranov and Buchweitz-Eisenbud-Herzog, coupled with Orlov's semiorthogonal decomposition theorem. 
\item The cubic surface case $X = X_3 \subseteq \mathbb{P}^3$ will follow from work of Kuleshov-Orlov on semiorthogonal decompositions of blow-ups, to prove $2)$, coupled with Orlov's semiorthogonal decomposition theorem to deduce $1)$. 
\item The even-dimensional complete intersection of two quadrics $X = V(Q_1, Q_2) \subseteq \mathbb{P}^{2n+2}$ with homogeneous coordinate ring $R = R_X$ will be analysed via Koszul duality by use of Buchweitz's Bernstein-Gel'fand-Gel'fand equivalence. Building on work of Buchweitz, Kuznetsov and Reiten-Van den Bergh we will deduce the existence of a full strong exceptional collection of \emph{linear} MCM modules in $ \underline{{\rm MCM}}^{\mathbb{Z}}(R)$. Said collection will then extend to a full strong exceptional collection of vector bundles on $X$ under Orlov's semiorthogonal decomposition theorem. 
	\end{enumerate}
\end{enumerate}
The points I-IV being clear, we now proceed to the case-by-case analysis of point V. 

\subsection*{Quadric hypersurfaces} 
Next, let $X = V(Q) \subseteq \mathbb{P}^{n+1}$ be a smooth quadric hypersurface for $n \geq 1$. Kapranov has shown the existence of a full strong exceptional collection of vector bundles on $X$ of the form
\begin{align*}
	{\rm D}^b(X) = \begin{cases}    \big\langle \OO_X(-(n-1)), \cdots, \OO_X(-1), \OO_X, \E \big\rangle & n \textup{ even,} \\
		\big\langle \OO_X(-(n-1)), \cdots, \OO_X(-1), \OO_X, \E_+, \E_- \big\rangle & n \textup{ odd.} \end{cases}   
\end{align*} 
Moreover, the bundles $\E$ and $\E_{\pm}$ are ACM bundles and generate a semisimple category, equivalent to the image $\Phi_0( \underline{{\rm MCM}}^{\mathbb{Z}}(R))$ under Orlov's semiorthogonal decomposition theorem. \bigskip

Let us give a full proof of the above, since it will allow us to introduce standard ideas to be reused in the case of quadrics intersections. Working somewhat anachronistically, one can establish Kapranov's Theorem from Orlov's Theorem. By \cite{BEH}, the singularity category $ \underline{{\rm MCM}}^{\mathbb{Z}}(R)$ of a quadrics hypersurface isolated singularity $R = k[x_0, \dots, x_{n+1}]/(Q)$ is semisimple, with one or two simple MCM modules (up to degree shift), say $M$ (for $n$ odd) and $M_+, M_-$ (for $n$ even). These have $2$-periodic resolutions given by matrix factorisations of $Q$ with linear entries $(A, B)$ (resp. $(A_{\pm}, B_{\pm})$) 
\[
	\dots \xrightarrow{A} R(-2)^r \xrightarrow{B} R(-1)^r \xrightarrow{A} R^r \to M \to 0
\]
and respectively
\[
	\dots \xrightarrow{A_\pm} R(-2)^{r_{\pm}} \xrightarrow{B_\pm} R(-1)^{r_\pm} \xrightarrow{A_\pm} R^{r_{\pm}} \to M_\pm \to 0 
\]
with $r = n+2$ (resp. $r_\pm = n+1$). Orlov's Theorem then gives 
\begin{align*}
	{\rm D}^b(X) = \begin{cases}    \big\langle \OO_X(-(n-1)), \dots, \OO_X(-1), \OO_X, \Phi_0(M) \big\rangle & n \textup{ even,} \\
	\big\langle \OO_X(-(n-1)), \dots, \OO_X(-1), \OO_X, \Phi_0(M_+), \Phi_0(M_-) \big\rangle & n \textup{ odd.} \end{cases}   
\end{align*} 
Let us recall how to compute the functor $\Phi_i$ for any $i \in \Z$. Given a graded MCM module $N$ with complete resolution $C$ 
\[
	\cdots \to C_{n+1} \to C_{n} \to C_{n-1} \to \cdots
\]
Let $C_{[<i]} \subseteq C$ be the subcomplex given by the summands generated in degree $< i$, and define $C_{[\geq i]} = C/C_{[<i]}$. Then $C_{[\geq i]} \in {\rm D}^b(\grmod R) $ is a lower bounded complex with bounded cohomology, and so one may sheafifify it to obtain a complex of coherent sheaves $\widetilde{C_{[\geq i]}} \in {\rm D}^b(X) $. By Lemma \ref{BuchweitzAdjoint}, we have
\[
	\Phi_i(N) \cong \widetilde{C_{[\geq i]}}.
\]
Finally, applying this to the linear resolution of $M$, $M_{\pm}$, we see that $\Phi_0(M) = \widetilde{M}$, $\Phi_0(M_{\pm}) = \widetilde{M_{\pm}}$ is simply given by sheafification. Setting $\E = \widetilde{M}$ and $\E_{\pm} = \widetilde{M_{\pm}}$, we obtain the exceptional ACM bundles of Kapranov, and it's easy to see that the relevant sheaf cohomology groups vanish and so the collection is strong. \bigskip

Summarising what we have shown:
\begin{prop}\label{quadrichypersurfacecase}
	Let $X = V(Q) \subseteq \mathbb{P}^{n+1}$ be a smooth quadric hypersurface over an algebraically closed with $k$ of characteristic zero, with homogeneous coordinate ring $R$. Then there is a full strong exceptional collection of graded MCM modules
	\[
		\underline{{\rm MCM}}^{\mathbb{Z}}(R) = \begin{cases} \big\langle M \big\rangle & n \textup{ even,} \\
		\left\langle M_{+}, M_{-} \right\rangle & n \textup{ odd,} \end{cases}  
	\]
extending under Orlov's semiorthogonal decomposition theorem to a full strong exceptional collection of vector bundles	
\begin{align*}
	{\rm D}^b(X) = \begin{cases}    \big\langle \OO_X(-(n-1)), \cdots, \OO_X(-1), \OO_X, \Phi_0(M) \big\rangle & n \textup{ even,} \\ 
	\big\langle \OO_X(-(n-1)), \cdots, \OO_X(-1), \OO_X, \Phi_0(M_+), \Phi_0(M_-) \big\rangle & n \textup{ odd.} \end{cases}   
\end{align*} \qed
\end{prop}

\subsection*{Cubic surfaces} 
Next, let $X = X_3 \subseteq \mathbb{P}^3$ be a smooth cubic surface and $R$ its coordinate ring. The next argument is due to Kuleshov-Orlov \cite{KO}, see also \cite[Thm. 2.5]{AKO}. Since $X_3$ is smooth, it is well-known that it is abstractly isomorphic to the blow-up $\pi: X \cong Bl_6 \mathbb{P}^2 \twoheadrightarrow \mathbb{P}^2$ of the projective plane in $6$ points in general position. Let $E_i = \pi^{-1}(p_i)$ be the $i$-th exceptional fibre. By \cite[Sect. 11.2]{Huy} (see also references in \cite{AKO}), the derived pullback $\mathbb{L}\pi^*: {\rm D}^b(\mathbb{P}^2) \to {\rm D}^b(X) $ is fully faithful, and there is a semiorthogonal decomposition
\begin{align*}
	{\rm D}^b(X) &= \big\langle \mathbb{L}\pi^* {\rm D}^b(\mathbb{P}^2), \OO_{E_1}, \cdots, \OO_{E_6} \big\rangle\\
	&= \big\langle \pi^*\E_1, \pi^*\E_2, \pi^*\E_3, \OO_{E_1}, \cdots, \OO_{E_6} \big\rangle 
\end{align*}
where $(\E_1, \E_2, \E_3)$ is any strong full exceptional sequence of vector bundles on $\mathbb{P}^2$. By \cite[Thm. 2.5]{AKO} (or direct calculations using the adjunction $({\rm L}\pi^*, {\rm R}\pi_*)$), the above sequence is a full strong exceptional collection of sheaves. In particular, taking $\E_1 = \OO_{\mathbb{P}^2}$, we have $\pi^*\E_1 = \OO_X$. Since $R$ has $a$-invariant $-1$, Orlov's semiorthogonal decomposition yields 
\begin{align*}
	{\rm D}^b(X) &= \big\langle \pi^*\E_1, \pi^*\E_2, \pi^*\E_3, \OO_{E_1}, \cdots, \OO_{E_6} \big\rangle \\
	&= \big\langle \OO_X , \Phi_0( \underline{{\rm MCM}}^{\mathbb{Z}}(R)) \big\rangle 
\end{align*}
and so we obtain a full strong collection of objects in $ \underline{{\rm MCM}}^{\mathbb{Z}}(R)$ by pushing the above down via the adjoint $\Psi_0$. It follows that both $ \underline{{\rm MCM}}^{\mathbb{Z}}(R)$ and $ {\rm D}^b(X) $ have tilting objects. The claim that $ {\rm D}^b(X) $ actually has a tilting object coming from a full strong exceptional collection of vector bundles (even line bundles) is due to Hille-Perling in \cite[Thm. 5.14]{LP}. \bigskip

Summarising the relevant facts from the above discussion, we have: 
\begin{prop}\label{cubicsurfacecase} Let $X = X_3 \subseteq \mathbb{P}^3$ be a smooth cubic surface over an algebraically closed field $k$ of characteristic zero, with homogeneous coordinate ring $R$. Then, with notation as above: 
	\begin{enumerate}[i)]
	\item There is a full strong exceptional collection
		\[
			\underline{{\rm MCM}}^{\mathbb{Z}}(R) = \big\langle \Psi_0(\pi^*\E_2), \Psi_0(\pi^*\E_3), \Psi_0(\OO_{E_1}), \dots, \Psi_0(\OO_{E_6}) \big\rangle. 
		\]
	\item There is a full strong exceptional collection of line bundles on $X$. 
\end{enumerate}\qed 
\end{prop}

\subsection*{Even dimensional intersections of two quadrics}
We are left with analysing the $2n$-dimensional smooth complete intersection of quadrics $X = V(Q_1, Q_2) \subseteq \mathbb{P}^{2n+2}$. Again let $R$ denote the coordinate ring. This case is far more involved, and we will need a considerable amount of setup. \bigskip 

Complete intersections of quadrics in projective space have been studied by Buchweitz, Kapranov, Bondal-Orlov and Kuznetsov. In particular Kuznetsov produced \cite[Cor. 5.7]{KuzQuadrics}, for a smooth complete intersection of two quadrics, a semiorthogonal decomposition\footnote{Kuznetsov's semiorthogonal decomposition was of the form $ {\rm D}^b(\coh X) = \langle {\rm D}^b(\coh C), \A' \rangle$, but one can exchange the two terms up to twist by the Serre functor, and we will do this here.} of $X = V(Q_1, Q_2) \subseteq \mathbb{P}^{N+2}$ 
\[
	{\rm D}^b(\coh X) = \big\langle \A, {\rm D}^b(\coh C) \big\rangle
\]
with $\A$ generated by an exceptional collection of line bundles, and $C$ is an associated hyperelliptic curve of genus $g = \frac{N+1}{2}$ (for $N$ odd) or a stacky projective line with $\Z_2$-stack structure at $N+3$ many points (for $N$ even). It was noted in \cite[Ex. 5.3, Sect. 6]{BernTabua} that in the latter case, the above semiorthogonal decomposition refines to a full exceptional sequence of objects in $ {\rm D}^b(\coh X) $ for $N$ even. However, the objects involved were given no explicit description, and it is unclear if such a collection is strong or even consists of sheaves. \bigskip

Our argument will proceed in parallel: setting $N = 2n$, using slightly different (but essentially equivalent) models, for $X = V(Q_1, Q_2) \subseteq \mathbb{P}^{2n+2}$ we will identify the non-trivial component of $ {\rm D}^b(\coh X) $ as the derived category $ {\rm D}^b(\coh \OO) $ of an hereditary order $\OO$ over $\mathbb{P}^1$, with ramification of order $2$ at $2n+3$ points. A theorem of Reiten-Van den Bergh \cite[Prop. 5.1]{ReitVdBtilting} then guarantees the existence of a full strong exceptional collection of sheaves $(\F_1, \dots, \F_{2n+5})$ in $\coh \OO \subseteq {\rm D}^b(\coh \OO) $. Our contribution will be to show that, after suitable modification, any such collection extends under Orlov's semiorthogonal decomposition to a full strong exceptional collection 
\[
{\rm D}^b(X) = \big\langle \OO_X(-2n+2)), \dots, \OO_X(-1), \OO_X, \E_1, \dots, \E_{2n+5} \big\rangle
\]
with $\E_i$ exceptional ACM vector bundles, analogous to Kapranov's decomposition in the hypersurface case. \bigskip

We begin with recalling Buchweitz's extension of the Bernstein-Gel'fand-Gel'fand correspondence to complete intersections of quadrics \cite[Appendix]{BEH}. Let $V$ be a finite dimensional vector space of dimension $2n+3$, and let $S = S(V) \cong k[x_0, \dots, x_{2n+2}]$ be the symmetric algebra on $V$. Let $W \subseteq S^2(V)$ be a $2$-dimensional subspace with basis given by a regular sequence, and let $R := S(V)/(W)$ be the complete intersection algebra. \bigskip 

Complete intersections of quadrics are Koszul, which follows from the Tate resolution of the residue field \cite{AvramovBook}. The Koszul dual $R^! = {\rm Ext}^*_R(k, k)$ takes an especially simple form. Following \cite[Appendix, Sect. A.2]{BEH}, consider the $k$-quadratic map $q: V^* \to W^*$ obtained by composition $q = \iota^* \circ \chi$ 
\[
	V^* \xrightarrow{\chi} S^2(V)^* \xrightarrow{\iota^*} W^*
\]
of the quadratic map $\chi: V^* \to S^2(V)^*$ given by $\chi(\xi)(vv') = \xi(v)\xi(v')$ with the pullback of the embedding $\iota: W \hookrightarrow S^2(V)$. Consider the symmetric algebra $S(W^*)$. We obtain an induced quadratic form 
\[
	q: V^* \to W^* \subseteq S(W^*) 	
\]
with value in the algebra $S(W^*)$, and so one can form the Clifford algebra $C = C_{S(W^*)}(V^*, q)$ given by
\[
	C = T(V^*) \otimes S(W^*)/I
\]
with $T(V^*)$ the tensor algebra on $V^*$ and $I$ the two-sided ideal generated by elements of the form 
\[
	\xi \otimes \xi - q(\xi)\  \textup{ for } \xi \in V^*. 
\] 
The Clifford algebra $C$ is $\Z$-graded by setting $|V^*| = 1$ and $|W^*| = 2$, and is finite over the central subalgebra $S(W^*) \subseteq C$. By a theorem of Sj\"odin \cite[Chp. 10]{AvramovBook}, we then have an isomorphism of graded algebras
\[
	R^! = {\rm Ext}_R^*(k, k) \cong Cl_{S(W^*)}(V^*, q) = C.
\]
Since $C$ is finite over $S(W^*)$, it is two-sided Noetherian, and Buchweitz moreover proves that it has finite injective dimension $\idim (C_C) = \idim (_C C) = {\rm codim}\ \! R = \dim W$, which in our scenario is $2$. Hence $C$ is a non-commutative Gorenstein ring, and in particular we have an equivalence $D_{sg}^{\Z}(C) \cong \underline{{\rm MCM}}^{\mathbb{Z}}(C)$ as per Buchweitz's theorem. \bigskip

The pair $(R, C)$ is then a pair of Koszul dual Noetherian Gorenstein Koszul algebras. In \cite[Appendix]{BEH}, Buchweitz proves the following extension of the Bernstein-Gel'fand-Gel'fand correspondence. Letting $A$ be either $R$ or $C$, denote by $ {\rm D}_{ \mathsf{art}}^b(A) \subseteq {\rm D}^b(\grmod A) $ and ${\rm D}^{\perf}(A) \subseteq {\rm D}^b(\grmod A) $ the full subcategories of complexes with Artinian cohomology and perfect complexes, respectively. 
\begin{thm}[Buchweitz {\cite[Appendix]{BEH}}]\label{RagnarBGG} Let $A$ be either $R$ or $C$. Then there is a functor 
	\[
		\beta: {\rm D}^b(\grmod A) \to {\rm D}^b(\grmod A^!) 
	\]
	satisfying the following properties: 
	\begin{enumerate}[1)] 
	\item $\beta$ is an equivalence of triangulated categories
		\[
	\beta:	{\rm D}^b(\grmod A) \xrightarrow{\cong} {\rm D}^b(\grmod A^!).
		\]
	\item Given $M \in {\rm D}^b(\grmod A)$, the cohomology module of $\beta(M)$ has graded components
		\[
			{\rm H}^j(\beta(M))_i = {\rm Ext}^{i+j}_{gr A}(k, M(-i)). 
		\]
	\item The equivalence $\beta$ sends $ {\rm D}_{\mathsf{art}}^b(A) $ onto $ {\rm D}^{\perf}(A^!) $ and $ {\rm D}^{\perf}(A)$ onto $ {\rm D}_{ \mathsf{art}}^b(A^!) $, and so descends to equivalences of Verdier quotients 
		\[
			{\rm D}^b(\qgr A) \cong {\rm D}_{sg}^{\Z}(A^!),
		\]
		\[
			{\rm D}_{sg}^{\Z}(A) \cong {\rm D}^b(\qgr A^!).
		\]
\end{enumerate}
\end{thm}\bigskip

We extract the following important corollary. Let $A = R$ in the theorem above, and let $\omega_R = R(a)$ and $d = \dim R$. Consider the duality $D = {\rm RHom}_R(-, \omega_R[d])$
	\[
		D: {\rm D}^b(\grmod R)^{op} \xrightarrow{\cong} {\rm D}^b(\grmod R).
	\]
	and note that $D \circ D \simeq {\rm id}$ and $D(k) \cong k$.\footnote{The complex $R(a)[d]$ is sometimes called the normalised canonical complex.} \bigskip

	Recall that a complex $F_*$ of free graded $R$-module is called linear if $F_i$ is generated in degree $i$. We will also abuse notation slightly and identify $ \underline{{\rm MCM}}^{\mathbb{Z}}(R) \cong {\rm D}^{\Z}_{sg}(R)$. 
\begin{cor} The following holds:
	\begin{enumerate}[a)]
\item Under the contravariant equivalence 
\[
	\beta \circ D:	{\rm D}^b(\grmod R)^{op} \xrightarrow{\cong} {\rm D}^b(\grmod C) 
\]
the standard t-structure on $ {\rm D}^b(\grmod C) $ pulls back to a t-structure whose heart consists of complexes $M \in {\rm D}^b(\grmod R) $ whose minimal free resolution is linear.
\item Under the induced contravariant equivalence
	\[
		\beta \circ D: \underline{{\rm MCM}}^{\mathbb{Z}}(R)^{op} \xrightarrow{\cong} {\rm D}^b(\qgr C) 
	\]
	the standard t-structure on $ {\rm D}^b(\grmod C) $ pulls back to a t-structure whose heart consists of modules $M \in \underline{{\rm MCM}}^{\mathbb{Z}}(R)$ whose minimal free resolution is {\bf eventually} linear. 
\end{enumerate}
\end{cor}
\begin{proof}
By Property 2) of Thm. \ref{RagnarBGG}, for $M \in {\rm D}^b(\grmod R) $ we have
\begin{align*}
	{\rm H}^j_i(\beta \circ D(M)) &= {\rm Ext}_{gr R}^{i+j}(k, D(M)(-i)) \\ 
	&= {\rm Hom}_{ {\rm D}^b(\grmod R) }(k, D(M)(-i)[i+j])\\
	&= {\rm Hom}_{ {\rm D}^b(\grmod R) }(k, D(M(i)[-i-j]))\\
	&= {\rm Hom}_{ {\rm D}^b(\grmod R) }(M(i)[-i-j], D(k))\\
	&= {\rm Hom}_{ {\rm D}^b(\grmod R) }(M, k(-i)[i+j])\\
	&= {\rm Ext}^{i+j}_{gr R}(M, k(-i)) 
\end{align*} 
and so $a)$ follows. To see $b)$, note that every object of $\qgr C \subseteq {\rm D}^b(\qgr C) $ arises as the ``sheafification'' of some graded module in $\grmod C \subseteq {\rm D}^b(\grmod C)$. Hence the induced heart on ${\rm D}_{sg}^{\Z}(R)^{op}$ has objects consisting of MCM approximations of complexes whose minimal resolution is linear. Those are precisely the MCM modules whose minimal resolution is eventually linear. 
\end{proof}
\begin{rem}
	There is a minor error in \cite[Appendix]{BEH}, where this corollary is misstated. In particular it is claimed that $\beta$ pulls back graded $C$-modules to complexes with linear minimal resolution, and that the induced heart on MCM modules consists precisely of linear MCM modules. The first claim is false as one needs to dualise first, and the second claim is false because MCM approximations of complexes with linear minimal resolutions can fail to be linear. Indeed this is the case for $k^{st}$ over certain Koszul algebras, where $k$ has a linear resolution but $k^{st}$ will typically admit generators in various non-zero degrees. We will see such an example in Chapter 3. 
\end{rem}\bigskip 

The induced heart on $ \underline{{\rm MCM}}^{\mathbb{Z}}(R)$ will play a central role in this thesis, and in particular in the construction of a tilting bundle for $ {\rm D}^b(X) $, and so we give it a name. We define 
\[
	\mathcal{H}^{\mathsf{lin}}(R) \subseteq \underline{{\rm MCM}}^{\mathbb{Z}}(R)
\]
the full subcategory of eventually linear stable MCM modules. In the above scenario, we have a diagram
\[
	\xymatrix{ \underline{{\rm MCM}}^{\mathbb{Z}}(R)^{op} \ar[r]^-{\beta \circ D}_{\cong} & {\rm D}^b(\qgr C) \\
		\mathcal{H}^{\mathsf{lin}}(R)^{op} \ar[u] \ar[r]^{\beta \circ D}_{\cong} & \qgr C \ar[u]
	}
\]
with vertical arrows the natural inclusions.\bigskip

Next, we summarise our approach to constructing a tilting object in $ \underline{{\rm MCM}}^{\mathbb{Z}}(R)$ and $ {\rm D}^b(\coh X) $ with $\proj R = X = V(Q_1, Q_2) \subseteq \mathbb{P}^{2n+2}$ smooth, and recall that we assume $\dim X > 0$. We first tie the two questions:
\begin{prop}\label{linearOrlovtrick} Assume that $\mathcal{H}^{\mathsf{lin}}(R)$ contains a tilting object $T$ for $ \underline{{\rm MCM}}^{\mathbb{Z}}(R)$. Then $ {\rm D}^b(\coh X) $ admits a tilting ACM vector bundle $\E$. Moreover, if the summands of $T$ form a full strong exceptional collection, then so does of $\E$. 
\end{prop} 
\begin{proof}
	Recall that $R = k[x_0, \dots, x_{2n+2}]/(Q_1, Q_2)$ has $a$-invariant $a = -2n + 1 < 0$ since we assume that $\dim X > 0$ and so $n > 0$. Consider the fully faithful embedding
	\[
		{\rm R}\Gamma_{\geq 0}: {\rm D}^b(\coh X) \hookrightarrow {\rm D}^b(\mathsf{grmod}_{\geq 0}\ \! R).
	\]
	The strong version of Orlov's Theorem (see Appendix A.3) gives a semiorthogonal decomposition 
	\[
		{\rm R}\Gamma_{\geq 0}({\rm D}^b(\coh X)) = \big\langle R(2-2n), \dots, R(-1), R, ( \underline{{\rm MCM}}^{\mathbb{Z}}(R))_{[\geq 0]}\big\rangle 
	\]
	Since $T \in \mathcal{H}^{\mathsf{lin}}(R)$, it follows that some sufficiently large syzygy $T' = \Omega^m(T)(m)$ has a linear minimal resolution, and $T'$ is also a tilting object. Since $T'$ is linear, we then have ${\rm L}_{\geq 0}(T') \cong T' \in {\rm D}^b(\grmod_{\geq 0}R) $. Since $T'$ is a module (and not a complex), we have
	\[
		{\rm Ext}^n_{gr R}(R(-i), T') = 0 \textup{ for all } n \neq 0 .
	\]
	It follows that $\bigoplus_{i=0}^{2n-2} R(-i) \oplus T'$ forms a tilting module for ${\rm R}\Gamma_{\geq 0}( {\rm D}^b(\coh X) )$ which is moreover an MCM module, and its sheafification $\E$ is a tilting ACM vector bundle for $ {\rm D}^b(\coh X) $. This proves the main claim, and the second claim is clear. 
\end{proof}

It remains to construct a tilting object $T \in \mathcal{H}^{\mathsf{lin}}(R)$. By making use of the BGG equivalence (Thm. \ref{RagnarBGG}), this reduces to constructing a tilting object in $\qgr C \subseteq {\rm D}^b(\qgr C) $. 

\subsection*{Analysing the category $\qgr C$} We review some generalities concerning $C$, which can be taken from \cite[Appendix]{BEH} or \cite[Sect. 3]{KuzQuadrics}. Recall that the graded algebra $C$ is finite projective over the central subalgebra ${\rm S}(W^*)$ with $W^*$ in degree $2$. Since $C$ is generated by $C_1$ over $C_0 = k$, passing to the $2$nd Veronese subalgebra $C_{even}$ induces an equivalence ${\rm qgr}\ C \cong {\rm qgr}\ C_{even}$, sending $M$ to $M_{even}$ (see \cite[Prop. 2.5]{Pol}). Taking central homogeneous localisations of ${\rm S}(W^*)$ is compatible with forming Clifford algebras, and so $C_{even}$ descends to a locally free sheaf of even Clifford algebras $\mathcal{O} := \widetilde{C}_{even}$ over $\mathbb{P}(W^*)$ of rank $2^{\dim V - 1} = 2^{2n+2}$. \bigskip 

Next, sheafification over $\proj S(W^*) = \mathbb{P}(W^*)$ gives an exact functor $F: \grmod C_{even} \to \coh \mathbb{P}(W^*)$ which vanishes precisely on finite length $C_{even}$-modules. It then factors as
\[
	\xymatrix@R=10pt@C=10pt{ \grmod C_{even} \ar@{->>}[dr] \ar[rr]^-{F} && \coh \mathbb{P}(W^*) \\
	& \qgr C_{even} \ar[ur]_-{\pi}} 
\] 
The functor $\pi$ is an exact functor between abelian categories, and we claim that it is faithful. Since it is exact, it suffices to show that $\pi(X) = 0$ if and only if $X = 0$. But this holds by construction since $\pi(X) = 0$ implies $F(\widetilde{X}) = 0$ for any module representative $\widetilde{X}$ of $X$, in which case $\widetilde{X}$ has finite length and $X = 0$. We have then shown that $\pi$ is faithful, and so it identifies 
\[
	\pi: \qgr C_{even} \xrightarrow{\cong} \coh \OO 
\] 
with the subcategory of coherent sheaves over $\mathbb{P}(W^*)$ admitting an action from $\OO$, or in other words coherent sheaves of modules over the sheaf of even Clifford algebras $\OO$. We are left with obtaining a good description of $\OO$. \bigskip 

For convenience, since $W$ is $2$-dimensional, picking a basis $q_0, q_\infty$ we get a parameterisation for the pencil of quadrics by $q_t = t_0q_\infty + t_1q_0$ for $t = [t_0:t_1] \in \mathbb{P}^1 = \mathbb{P}(W^*)$. The quadratic form $q_t \in W$ then defines a functional on $W^*$, which lets us define a $k$-quadratic form
\[
	V^* \xrightarrow{q} W^* \xrightarrow{q_t} k 
\]
which is by definition the quadratic form defined by $q_t \in W \subseteq S^2(V)$. The following falls out of our identifications.
\begin{lem} Under the parameterisation $\mathbb{P}^1 = \mathbb{P}(W^*)$, the fibre of $\OO$ over the point $\iota_t: \spec k(t) \hookrightarrow \mathbb{P}^1$ is isomorphic to the Clifford algebra of $q_t$, that is 
	\[
		\iota_t^* \OO \cong Cl_{k(t)}(V^*, q_t) 
	\]
	where $k(t) = k$ denotes the residue field at the point $t$. 
\end{lem}
We now use a result of Reid. Since we are in $\mathsf{char}\ \! k \neq 2$, the quadratic form $q_t$ has an associated bilinear form $B_t$, and we define the corank of $q_t$ by the dimension of the kernel of the map $V^* \to V^{**} = V$ induced from $B_t$. Recall that $q_t$ is called non-degenerate if $B_t$ defines a perfect pairing, or equivalent $q_t$ has corank zero. 
\begin{lem}[{\cite[Prop. 2.1]{ReidThesis}}]\label{ReidLemma} A complete intersection of quadrics $X = V(q_0, q_\infty) \subseteq \mathbb{P}^N$ is smooth if and only if the quadratic form $q_t$ has corank $\leq 1$ for all $t$, is generically non-degenerate and otherwise of corank $1$ at $(N+1)$-many values of $t$. 
\end{lem}
We can now prove that the sheaf of even Clifford algebras $\OO$ over $\mathbb{P}^1$ is a sheaf of hereditary orders. We begin by recalling some definitions from \cite{ReitVdBtilting}. 

\subsection*{Hereditary orders over $\mathbb{P}^1$}
Let $R = (R, \m, k)$ be a discrete valuation domain with fraction field $K$. An order $\Lambda$ over $R$ consists of a subalgebra $\Lambda \subseteq A$ of a central simple $K$-algebra $A$, such that $\Lambda$ contains $R$ as a central subring and which is finitely generated as a module over $R$. An hereditary order $\Lambda$ over $R$ is simply an order $\Lambda$ that is hereditary as an algebra. \bigskip 

Let $J = J(\Lambda)$ be the Jacobson radical. Then there is a natural number $e \in \N$ such that $J^e = \m \Lambda$, called the ramification order of $\Lambda$ over $R$, and we say that $\Lambda$ is unramified if $e = 1$.\bigskip 

Next consider $\mathbb{P}^1$ with function field $K = K(\mathbb{P}^1)$, thought of as a constant sheaf on $\mathbb{P}^1$ so that $\OO_{\mathbb{P}^1} \subseteq K$. Let $A$ be a central simple algebra over $K$. A sheaf of hereditary order $\OO$ over $\mathbb{P}^1$ is an $\OO_{\mathbb{P}^1}$-algebra $\OO \subseteq A$ which is coherent over $\mathbb{P}^1$, and which is locally an hereditary order over each DVR $\OO_{\mathbb{P}^1,t}$. \bigskip

Next, recall some facts from the structure theory of Clifford algebras, see \cite{BEH}. Let $F$ be an algebraically closed field of characteristic not $2$. Let $(U, q)$ be a finite dimensional vector space over $F$ equipped with a non-degenerate quadratic form $q$. Then the even Clifford algebra $Cl_F(U, q)_{\mathsf{even}}$ is isomorphic to a matrix algebra over $F$ when $\dim U$ is odd, and to a product of two matrix algebras when $\dim U$ is even.\bigskip

We now prove the anticipated result. This can be seen as analogous to \cite[Cor. 3.16]{KuzQuadrics}, however phrased in the language of sheaves of orders rather than stacks. 
\begin{thm} The sheaf $\mathcal{O}$ of even Clifford algebras is a sheaf of hereditary orders over $\mathbb{P}(W^*) = \mathbb{P}^1$, with ramification order over $t \in \mathbb{P}^1$ given by $e_t = 1 + corank(q_t)$, where $corank(q_t) \leq 1$. 
\end{thm}
\begin{proof}
	For $t \in \mathbb{P}^1$, set $\mathfrak{m}_t \subseteq \mathcal{O}_{\mathbb{P}^1, t}$ and $J_t \subseteq \mathcal{O}_t$ the respective Jacobson radicals. The algebra $\OO_t$ is given by the even Clifford algebra $\OO_t = Cl_{\OO_{\mathbb{P}^1, t}}(V^*, q_t)_{\mathsf{even}}$ for the induced $\OO_{\mathbb{P}^1}$-valued quadratic form $q_t$, and so is module-finite over the central DVR $\OO_{\mathbb{P}^1,t}$. We have
\[ 
\mathcal{O}_{\mathbb{P}^1, t} \subseteq \mathcal{O}_t \subseteq \mathcal{O}_t \otimes_{\mathcal{O}_{\mathbb{P}^1, t}} K = Cl_K(V^*, q)_{even} \cong M_r(K)
\]
where $r = 2^{2n+2}$. For the last isomorphism, observe that $Cl_K(V^*, q)_{even}\otimes_K \overline{K} = Cl_{\overline{K}}(V_K^*, q)_{even} \cong M_{r}(\overline{K})$ for the base change $V_K^* = V^* \otimes K$, since $q: V_K^* \to K$ is non-degenerate, and so $Cl_K(V_K^*, q)$ is a central simple algebra. However by Tsen's Theorem, there are no non-trivial central simple algebras over $K$ since $k$ is algebraically closed, hence $Cl_K(V_K^*, q)$ must be a matrix algebra. \bigskip 

Hence $\mathcal{O}_t$ is an order over $\mathcal{O}_{\mathbb{P}^1, t}$ inside a matrix algebra. To compute the ramification order $e_t$, we claim that $J_t^{1+corank(q_t)} = \m_t \mathcal{O}_t$. To see this, note that $\mathcal{O}_t / \mathcal{O}_t\m_t \cong Cl_{k}(V^*, q_t)_{even}$ and let $U_{t}$ be the kernel of the map $V^* \to V^{**} = V$ associated to the corresponding symmetric bilinear form. Then ${\rm dim}\ U_{t} = corank(q_t)$ and
$Cl_k(V^*, q_t) \cong Cl_k(\widetilde{U}, q_t) \otimes \bigwedge U_{t}$ for some complement $\widetilde{U}$ of $U_{t}$. Applying Lemma \ref{ReidLemma}, we see that the dimension of $U_t$ is zero for generic $t$, and one otherwise. When $\dim U_t = 1$, one sees that elements of the form $u \cdot u_t$ for $u \in U$ and $u_t \in U_t$ generate the Jacobson radical $J$ of $Cl_k(V^*, q_t)$, and so $J^2 = 0$ since $U_t$ is one dimensional and $\wedge^2 U_t = 0$. \bigskip

It remains to prove that $\OO_t$ is hereditary. We will do this in three steps:
\begin{enumerate}[i)]
\item Show $\OO_t$ is hereditary for $q_t$ non-degenerate, and so $t$ generic.
\item Show that the abelian category $\coh \OO$ is hereditary.
\item Deduce that the remaining ramified orders $\OO_{t_i}$ are hereditary.
	\end{enumerate}

	For the first claim, we may apply the Auslander-Goldman theorem \cite[Thm. 39.1]{ReinerBook}, which claims that an order $\Lambda$ over a DVR $R$ is hereditary if and only if the radical $J(\Lambda)$ is a projective $\Lambda$-module. When $q_t$ is non-degenerate, we have seen that $J_t = \OO_{\mathbb{P}^1,t} \m_t$, and since $\OO_{\mathbb{P}^1,t}$ is a DVR, $\m_t$ is generated by a single element, say $\pi_t$. The annihilator of $\pi_t$ is trivial in $\OO_t$ and so $J_t = \pi_t \cdot \OO_t$ is a free module of rank one. By Auslander-Goldman, $\OO_t$ is hereditary. \bigskip

	For the second claim, we will use the BGG correspondence to construct a Serre functor on ${\rm D}^b(\coh \OO)$ and deduce the result. Recall that the contravariant form of Buchweitz's BGG correspondence gave rise to an equivalence
	\[
		\underline{{\rm MCM}}^{\mathbb{Z}}(k[x_0, \cdots, x_{2n+2}]/(q_0, q_\infty))^{op} \cong {\rm D}^b(\coh \OO)
	\]
which induces identifications of corresponding hearts
\[
	\mathcal{H}^{\mathsf{lin}}(k[x_0, \cdots, x_{2n+2}]/(q_0, q_\infty))^{op} \cong \coh \OO
\]
where the left hand side is given by MCM modules with eventually linear resolutions. Now, since $k[x_0, \cdots, x_{2n+2}]/(q_0, q_\infty)$ has isolated singularities at the origin, the stable category of MCM modules has a Serre functor of the form
\begin{align*}
	\Se &= (a)[2n] \\
	 &= (1-2n)[2n] \\ 
	 &= (1-2n)[2n-1][1] \\ 
	 &= - \otimes \omega_{\mathcal{H}} [1] 
\end{align*}
where the notation $- \otimes \omega_{\mathcal{H}}$ is a placeholder for $(1-2n)[2n-1]$, but more importantly $- \otimes \omega_{\mathcal{H}}$ preserves the heart $\mathcal{H}^{\mathsf{lin}}(k[x_0, \cdots, x_{2n+2}]/(q_0, q_\infty)$ of eventually linear modules. The functor $\Se$ is sent onto the inverse Serre functor $\Se_\OO^{-1}$ on $ {\rm D}^b(\coh \OO) $, which then has the form $\Se_\OO^{-1} = - \otimes \omega^{-1}[-1]$, where $ - \otimes \omega$ is the induced autoequivalence of $\coh \OO$. From the form of Serre duality 
	\[
		{\rm Ext}_{\OO}^i(\F, \G \otimes \omega) \cong {\rm Ext}_{\OO}^{1-i}(\G, \F)^*
	\]
	one deduces that $\coh \OO$ is hereditary. \bigskip

	Lastly, we let $t_i \in \mathbb{P}^1$ be a point over which $\OO$ ramifies, and $t$ correspond to a generic (unramified) point. We first note that since $\OO_t$ is Noetherian (as it is finite over a central DVR), the full module category $\Modsf \OO_t$ is hereditary. Similarly $\OO$ is a sheaf of Noetherian $\OO_{\mathbb{P}^1}$-algebras and so $\qcoh \OO$ is hereditary. To show that $\OO_{t_i}$ is hereditary, for any $\F, \G \in \qcoh \OO$ we invoke the (first quadrant, cohomological) local-to-global spectral sequence for the ringed space $Y = (\mathbb{P}^1, \OO)$ 
\[
	{\rm H}^p(Y, \mathcal{E}xt^q(\F, \G)) \implies {\rm Ext}_{\OO}^{p+q}(\F, \G).
\]
Since $Y$ is a Noetherian topological space of dimension $1$, by the Grothendieck Vanishing Theorem we have ${\rm H}^p(Y, \A) = 0$ for all $p \geq 2$ and sheaves of abelian groups $\A$. It follows that the above spectral degenerates. \bigskip

Next, since $\OO_t$ is generically hereditary, the sheaf $\mathcal{E}xt^2(\F, \G)$ is supported at most over the ramified points $ \{ t_i \} \subseteq \mathbb{P}^1$, and we have 
\begin{align*}
	{\rm H}^0(Y, \mathcal{E}xt^2(\mathcal{F}, \mathcal{G})) \cong \prod_{t_i} {\rm Ext}_{\OO_{t_i}}^2(\F_{t_i}, \G_{t_i}). 
\end{align*}
By degeneration of the spectral sequence, the left-hand is a summand of ${\rm Ext}^2_\OO(\F, \G) = 0$ and so vanishes, and therefore so do the individual ${\rm Ext}^2_{\OO_{t_i}}(\F_{t_i}, \G_{t_i})$. \bigskip

Finally, letting $j:U \subseteq \mathbb{P}^1$ be an affine open neighbourhood of $t_i$, denote by $\OO_U$ the restriction of $\OO$ to $U$. It is easy to see that localisation $\qcoh \OO_{U} \to \Modsf \OO_{t_i}$ is essentially surjective, and any $\F \in \qcoh \OO_{U}$ extends to a quasi-coherent sheaf $j_* \F \in \qcoh \OO$ without changing its stalk at $t_i$. It follows that any $M \in \Modsf \OO_{t_i}$ arises as $\F_{t_i}$ for some $\F \in \qcoh \OO$, and we are done.
\end{proof}

Let us now put everything together. Writing $R = k[x_0, \cdots, x_{2n+2}]/(q_0, q_\infty)$, by Buchweitz's BGG correspondence we have a contravariant equivalence $\mathcal{H}^{\mathsf{lin}}(R)^{op} \cong \coh \OO$ for an hereditary order $\OO$ over $\mathbb{P}^1$. We can then apply the following theorem of Reiten and Van den Bergh.

\begin{thm}[Reiten-Van den Bergh {\cite[Prop. 5.1]{ReitVdBtilting}}] Let $\OO$ be an hereditary order on $\mathbb{P}^1$ with ramification of order $e_i$ at the point $p_i$, $i = 1, \cdots, r$. Then there is a full strong exceptional collection of sheaves \[
		{\rm D}^b(\coh \OO) = \big\langle \F_1, \dots, \F_l \big\rangle
	\]
	of length $l = 2 + \sum_{i=1}^r (e_i - 1)$. 
\end{thm}
The sheaf of even Clifford algebras $\OO$ is an hereditary order, ramified of order $2$ at $2n+3$ points. We deduce:
\begin{cor}The category $\mathcal{H}^{\mathsf{lin}}(R) \subseteq \underline{{\rm MCM}}^{\mathbb{Z}}(R)$ contains a full strong exceptional collection of length $2n+5$ for $ \underline{{\rm MCM}}^{\mathbb{Z}}(R)$.
\end{cor}

We also note the following corollary of interest. 
\begin{cor} There is a contravariant equivalence of abelian categories
	\[
		\mathcal{H}^{\mathsf{lin}}(R)^{op} \cong \coh \X
	\]
	onto the category of coherent sheaves over a weighted projective line $\X = \mathbb{P}^1(2, 2, \dots, 2)$, with weight divisor given by the ramification divisor of the hereditary order $\OO$. 
\end{cor}
\begin{proof}
	Indeed by Reiten-Van den Bergh \cite{ReitVdBtilting} we have an equivalence $\coh \OO \cong \coh \X$. 
\end{proof}

Finally, together with Prop. \ref{linearOrlovtrick}, we deduce the existence of a full strong exceptional collection of (ACM) vector bundles on $X = V(q_0, q_\infty) \subseteq \mathbb{P}^{2n+2}$  
\[
	{\rm D}^b(\coh X) = \big\langle \OO_X(2-2n), \dots, \OO_X(-1), \OO_X, \E_1, \dots, \E_{2n+5} \big\rangle
\]\bigskip

Summarising the results of this subsection, we have shown: 
\begin{prop}\label{quadriccicase}
	Let $X = V(Q_1, Q_2) \subseteq \mathbb{P}^{2n+2}$ be a smooth, even-dimensional complete intersection of quadrics over an algebraically closed field $k$ of characteristic zero, with homogeneous coordinate ring $R$. Then there is a full strong exceptional collection of MCM modules with linear minimal resolutions 
	\[
		\underline{{\rm MCM}}^{\mathbb{Z}}(R) = \big\langle M_1, M_2, \cdots, M_{2n+5} \big\rangle. 
\]
This extends under Orlov's semiorthogonal decomposition theorem to a full strong exceptional collection of vector bundles 
\[
	{\rm D}^b(\coh X) = \big\langle \OO_X(2-2n), \dots, \OO_X(-1), \OO_X, \Phi_0(M_1), \Phi_0(M_2), \dots, \Phi_0(M_{2n+5}) \big\rangle 
\]\qed
\end{prop}

\subsubsection*{Putting it all together} 
We finally obtain the proof of the main theorem of this chapter. 
\begin{thm}\label{hodgeconverse} Let $k$ be algebraically closed of characteristic zero, let $R$ be a standard graded complete intersection $k$-algebra with isolated singularities, and $X = \proj R$. The following are equivalent: 
	\begin{enumerate}[1)]
	\item $ \underline{{\rm MCM}}^{\mathbb{Z}}(R)$ admits a tilting object $T$.
	\item $ {\rm D}^b(X) $ admits a tilting object $\E$. 
	\item We have $h^{p,q}(X) = 0$ for $p \neq q$.
\end{enumerate}
When $\dim R \geq 2$ and $R$ is not a polynomial algebra, this is equivalent to: 
\begin{enumerate}[1)]
\setcounter{enumi}{3}
\item $X$ belongs to one of the three familes $a)-b)-c)$. 
	\end{enumerate}
	Moreover, when these equivalent conditions hold, the tilting object $\E$ can always be assumed to come from a strong exceptional collection of vector bundles on $X$. 
\end{thm} 
\begin{proof}
The implications $1), 2) \implies 3)$ and $3) \iff 4)$ having been previously established, we were left with $4) \implies 1), 2)$. This is the content of Prop. \ref{quadrichypersurfacecase}, Prop. \ref{cubicsurfacecase} and Prop. \ref{quadriccicase}, and so the theorem holds.
\end{proof}

%% file: chapters/classifications/classifications.tex
\chapter{Classifications of MCM modules and Betti tables over tame curve singularities}
\chaptermark{MCM modules and Betti tables over tame curve singularities}

In this section we will take up the classification of indecomposable MCM modules over certain reduced curve singularities of tame Cohen-Macaulay representation type (CM-type for short). We begin by recalling basic notions from \cite{DG}.\bigskip

Fix an algebraically closed field $k$ throughout and let $C = (C, \m, k)$ be a reduced complete local Cohen-Macaulay curve singularity over $k$, with $X = \spec C$. Any MCM $C$-module $M$ is locally free away from the singular locus and so defines a vector bundle on each irreducible component of the regular locus $X_{ \mathsf{reg}} = \coprod_{i=1}^r X_{ \mathsf{reg}, i}$, with rank vector ${\bf rk}(M) = ({\rm rk}_1, \dots, {\rm rk}_r)$ and (total) rank $N = \sum_i {\rm rk}_i$. We say that the CM-type of $C$ is of:
\begin{enumerate}[1)]
\item finite type, if $C$ has finitely many indecomposable MCM modules;
\item tame type, if $C$ has infinitely many indecomposables and the indecomposables of fixed rank $N$ can be parameterized by finitely many $1$-parameter families $\F_1, \dots, \F_{\mu(N)}$, with at most finitely many exceptions;
\item wild type, if $C$ admits $n$-parameter families of non-isomorphic indecomposables for $n$ arbitrarily large.
	\end{enumerate}
	By \cite{DG2}, the CM-type of $C$ falls precisely in one of these three cases (see also \cite{DG}). The classification of curve singularities of finite and tame CM-type is closely related to the Arnold school classification of polynomials with isolated critical points of low modality over $k = \C$. The modality of a holomorphic function germ $f: (\C^n, 0) \to (\C, 0)$ is, roughly, the minimal number $m$ for which one can obtain all isomorphism classes of deformations of $f$ by finitely many $m$-parameter families, see \cite[1.9]{AGLV} for the precise definition. The holomorphic functions of low modality were classified by Arnold and his collaborators, who produced the following list of normal forms.\bigskip

	Given a function $f = f(x_1, \dots, x_n)$ depending on variables $\{x_1, \dots, x_n\}$ and a set of disjoint variables $\{z_1, \dots, z_k\}$, we call the function 
	\[
	g(x_1, \dots, x_n, z_1, \dots, z_k) := f(x_1, \dots, x_n) + z_1^2 + \dots + z_k^2
\]
a \emph{stabilisation} of $f$. Up to stabilisation and isomorphism, the function germs of modality $m \leq 1$ are \cite[2.3]{AGLV}:

\begin{itemize}
	\item $(m=0)$: These are the simple (or ADE) singularities
\[
	\begin{tabular}{|c|c|c|c|c|}
\hline
$A_\mu,\ \mu \geq 1$ & $D_\mu,\ \mu \geq 4$ & $E_6$ & $E_7$ & $E_8$ \\
\hline
$x^{\mu + 1}$ & $x^2y + y^{\mu - 1}$ & $x^3 + y^4$ & $x^3 + xy^3$ & $x^3 + y^5$ \\
\hline
	\end{tabular}
\]
	\item $(m=1)$: There are $3$ families of parabolic singularities
		\[
	\begin{tabular}{|c|c|c|}
\hline
$P_8$ & $X_9$ & $J_{10}$ \\
\hline
$x^3 + y^3 + z^3 + axyz$ & $x^4 + y^4 + ax^2y^2$ & $x^3 + y^6 + ax^2y^2$ \\
$a^3 + 27 \neq 0$ & $a^2 \neq 4$ & $4a^3 + 27 \neq 0$ \\
\hline
	\end{tabular}
\]
as well as the hyperbolic singularities 
\[
	T_{pqr}: x^p + y^q + z^r + axyz,\ a \neq 0,\ \frac{1}{p} + \frac{1}{q} + \frac{1}{r} < 1
\]
and an additional $14$ exceptional families, whose CM-type is wild.
\end{itemize}

The above normal forms play an essential role in the classification of curves of finite and tame CM-types over a general algebraically closed field $k$. Given $C$ as above, consider its normalisation $C \subseteq \overline{C} \subseteq Q(C)$ in its total quotient ring $Q(C)$. We say that a ring $D$ birationally dominates $C$ if there are embeddings $C \subseteq D \subseteq \overline{C}$.

\begin{prop}[Greuel-Kn\"orrer {\cite[Thm. 9.2]{YosBook}}] The curve $C$ has finite CM type if and only if it birationally dominates a simple curve singularity.
\end{prop}
Note that per the results of Buchweitz-Greuel-Kn\"orrer-Schreyer-Herzog in Chapter 2, the only Gorenstein curves of finite CM type are the simple curve singularities themselves. Next, by \cite{DGK} the hyperbolic singularity of type $T_{pq2}$ are isomorphic to the stabilisation of the curves $T_{pq}$ 
\[
	T_{pq}: x^p + y^q + bx^2y^2, \textup{ with } b \neq 0, \frac{1}{p} + \frac{1}{q} < \frac{1}{2}
\]
for $b = -\frac{a^2}{4}$. Extending this family to $T_{pq}$ with $\frac{1}{p} + \frac{1}{q} = \frac{1}{2}$ enlarges it by the two parabolic singularities $T_{44} = X_9$ and $T_{36} = J_{10}$. We impose $\mathsf{char}\ k \neq 2$ for the next two propositions.

\begin{prop}[Drozd-Greuel, \cite{DG}]The curve $C$ is of tame CM-type if and only if it birationally dominates a curve of type $T_{pq}$ with $\frac{1}{p} + \frac{1}{q} \leq \frac{1}{2}$. 
\end{prop}

Given a curve $C$ of tame type, we say that $C$ is tame of \emph{domestic} representation type if there is a uniform bound $\mu(n) \leq N$ on the number of $1$-parameter families of indecomposables of fixed rank $n$. 

\begin{prop}[Drozd-Greuel, \cite{DG}]Let $C$ be a curve of tame CM-type. Then $C$ is of domestic type if and only if it properly birationally dominates a curve of type $T_{44}$ or $T_{36}$.
\end{prop} 

Tameness of the $T_{pq}$ was originally established in the hyperbolic case by indirect methods via deformation theory \cite{DG}, with the indecomposables later described by use of the minimal resolution of the surface singularity $T_{pq2}$ \cite{DGK}. The parabolic case $T_{44}$, $T_{36}$ was first studied by Dieterich \cite{Diet}, who proved tameness by somewhat indirect methods; in particular this left open the description of the indecomposables. The question of obtaining explicit presentations of the indecomposables, or equivalently of writing down the indecomposable matrix factorizations, was raised by Drozd and Tovpyha in \cite{DT}, where they produced some of the indecomposable matrix factorizations of $T_{44}$.

\subsection*{Results}
We now outline the results of this chapter. Let $k$ be an algebraically closed field. We will classify the indecomposable graded MCM modules over the following graded algebras:
\begin{itemize}
	\item The homogeneous coordinate ring $R_{Y_2}$ of $4$ points $Y_2 \subseteq \mathbb{P}^2$ in general linear position, meaning that $Y_2$ arises as the complete intersection of two conics. 
	\item The homogeneous coordinate ring $R_{Y_1}$ of $4$ points $Y_1 \subseteq \mathbb{P}^1$, which can be written in normalised form as the hypersurface ring $R_{Y_1} = k[x,y]/(f_\lambda)$ with $f_\lambda = l_1l_2l_3l_4$ a product of linear forms 
		\[
			f_\lambda = xy(x-y)(x-\lambda y),\ \lambda \neq 0,1.
		\]
	\end{itemize}
	To do this, we will prove and use the following equivalences of categories. Both are special cases of the general theorem of Buchweitz-Iyama-Yamaura \cite{BIY} see in Chapter 2, but we will prove these directly. 
	\begin{thm}\label{D4tilt} There is an equivalence of triangulated categories
	\[
		\underline{{\rm MCM}}^{\mathbb{Z}}(R_{Y_2}) \cong {\rm D}^b(kQ)
	\]
	with $Q$ the $\widetilde{D}_4$ quiver
	{\tiny{
			\[
	\xymatrix@R=5pt@C=20pt{&\bullet\\
		&\bullet\\
		\bullet \ar[uur]^{p_1} \ar[ur]_{\ p_2} \ar[dr]^{\ p_3} \ar[ddr]_{p_4}& \\
		&\bullet\\
		&\bullet
	}
	\]
}}
\end{thm}

\begin{thm}\label{squidtilt} There is an equivalence of triangulated categories
	\[
		\underline{{\rm MCM}}^{\mathbb{Z}}(R_{Y_1}) \cong {\rm D}^b(Sq(2,2,2,2;\lambda))
	\]
	where the ``Squid'' algebra ${\rm Sq}(2,2,2,2;\lambda)$ is the path algebra of the quiver
{\tiny{
		\[
	\xymatrix@R=5pt@C=20pt{&&\bullet\\
		&&\bullet\\
		\bullet \ar@<.5ex>[r]^y \ar@<-.5ex>[r]_x &\bullet \ar[uur]^{p_1} \ar[ur]_{\ p_2} \ar[dr]^{\ p_3} \ar[ddr]_{p_4}& \\
		&&\bullet\\
		&&\bullet
	}
\]
}}
with relations $p_i l_i(x,y) = 0$ for $i = 1,2,3,4$.
\end{thm}

Both path algebras are derived tame and have a well-studied representation theory. \bigskip

In the case of $R_{Y_2}$, we will construct all indecomposable graded MCM modules and compare them to the known classification of indecomposable representations of $kQ$. We will see that one then recovers the regular components of $kQ$ out of the classical pencil of conics construction. We will also write down the Betti tables of complete resolutions of indecomposables. \bigskip

In the case $R_{Y_1}$, we will make use of the derived equivalence
\[
	{\rm D}^b(Sq(2,2,2,2;\lambda)) \cong {\rm D}^b(\X)
\]
with $\X = \mathbb{P}^1(2,2,2,2;\lambda)$ a weighted projective line of genus one in the sense of Geigle-Lenzing. We will produce the MCM modules corresponding to the simple torsion sheaves on $\X$. By results of Lenzing and Meltzer \cite{LM}, one can obtain all indecomposable sheaves by iterated applications of two `twists' autoequivalences
\[
	T_1, T_2: {\rm D}^b(\X) \xrightarrow{\cong} {\rm D}^b(\X)
\]
applied to indecomposable torsion sheaves, which send the simple torsion sheaves to the stable sheaves, and so one reduces the classification to understanding the action of $T_1, T_2$ on the MCM modules above. Finally, we will completely classify the Betti tables of indecomposable graded MCM modules. \bigskip

For either algebra $R_{Y_i}$, $i = 1,2$, we will see that all indecomposable MCM $\widehat{R_{Y_i}}$-modules arise as the completion of some graded MCM $R_{Y_i}$-modules, and so the classification results extend to $ \underline{{\rm MCM}}(\widehat{R_{Y_i}})$.

\input{chapters/classifications/tiltingmodules}

\input{chapters/classifications/conics}
\input{chapters/classifications/elliptic}

%% file: chapters/classifications/tiltingmodules.tex
\section{The tilting modules}

Let us now prove the above theorems by exhibiting an appropriate tilting MCM module. The proofs make use of Orlov's semiorthogonal decomposition theorem. Recall that $k$ is algebraically closed throughout this chapter, and we write $R = S/I$ for $S$ the ambient polynomial ring and $R = R_{Y_i}$, $i = 1,2$. To simplify calculations, we shall make use of the Orlov-Buchweitz embedding (see Appendix) 
\[
	\underline{{\rm MCM}}^{\mathbb{Z}}(R) \hookrightarrow {\rm D}^b(\mathsf{grmod}_{\geq 0} R)
\]
\[
	M \mapsto M_{[\geq 0]}
\]
where $M_{[\geq 0]}$ is the complex with bounded cohomology obtained by taking a complete resolution $C$ of $M$ and killing generators of degree $< 0$. This is most useful in the following situation:
\begin{lem}\label{cutoff} Let $R$ be a graded connected Gorenstein algebra, and let $M$ be a graded MCM $R$-module generated in degree zero with no free summand. Then $M_{[\geq 0]} \cong M$ in $ {\rm D}^b(\mathsf{grmod}_{\geq 0} R) $. In particular, if $M, N$ are both generated in degree zero without free summands, we have
	\[
		\underline{{\rm Ext}}_{gr R}^n(M, N) = \begin{cases} {\rm Ext}^n_{gr R}(M, N) & n \geq 0 \\ 
		0 & n < 0.\end{cases}
	\]
\end{lem}
\begin{proof}
	Since $R$ is graded connected, minimal complete resolutions exist for MCM modules without free summands\footnote{Note that minimality of $C$ also \emph{requires} $M$ to have no free summands, since any summand of $M$ isomorphic to $R$ would produce a summand of the form $ 0 \to R \xrightarrow{id} R \to 0$ of $C$ in degree $0, -1$.}. Let $C$ be a minimal complete resolution of $M$:
	\[
	\xymatrix@C=20pt@R=5pt{\dots \ar[r] & C_1 \ar[r] & C_0 \ar[dr] \ar[rr] && C_{-1} \ar[r]& \dots \\
		& & & M \ar[ur] & }
	\]
	Since $M$ is generated in degree zero, so is $C_0$ by minimality the graded free modules $C_{n}$ are generated in positive degrees for $n > 0$ and in negative degrees for $n < 0$. Killing negative degree generators returns the minimal free resolution $M_{[\geq 0]} \xrightarrow{\cong} M$ of $M$. For the second claim, since $(-)_{[\geq 0]}$ is fully faithful, we have $ \underline{{\rm Ext}}^*_{gr R}(M, N) \cong {\rm Hom}_{ {\rm D}^b(\grmod R) }(M_{[\geq 0]}, N_{[\geq 0]}) \cong {\rm Ext}_{gr R}^*(M, N)$ and the result follows. 
\end{proof}
\subsubsection*{The algebra $R_{Y_2}$}

Now let $R = R_{Y_2}$. Under the assumption that $Y_2 \subseteq \mathbb{P}^2$ be in general position, meaning no three points of $Y_2$ lie on a line, then $Y_2 = V(Q, Q')$ is a complete intersection of conics, and so we have $R_{Y_2} \cong k[x,y,z]/(Q, Q')$. In particular, $R_{Y_2}$ is Gorenstein with $a$-invariant $a = 1$.

\begin{thm}\label{D4tiltproof} Let $R = R_{Y_2}$, so that $X = \proj R = Y_2 \subseteq \mathbb{P}^2$ is the set of $4$ points $\{p_i\}$. Let $L_i = R/I(p_i)$ be the homogeneous coordinate ring of $p_i$, thought of as an $R$-module. Then there is a full strong exceptional collection 
	\[
		\underline{{\rm MCM}}^{\mathbb{Z}}(R) = \langle \m(1), L_1, L_2, L_3, L_4 \rangle
	\]
	with the endomorphism algebra of $T = \m(1) \oplus \left( \bigoplus_{i=1}^4 L_i\right)$ given by $kQ$ as above.
\end{thm}
\begin{proof}
	We have $a = a(R) = 1$. Applying Orlov's theorem with cutoff $i = -a = -1$ gives a semiorthonal decomposition
	\begin{align*}
		\underline{{\rm MCM}}^{\mathbb{Z}}(R) &= \langle k^{st}(1), {\rm st} \circ {\rm R}\Gamma_{\geq i+a} {\rm D}^b(X) \rangle \\
&= \langle k^{st}(1), L_1, L_2, L_3, L_4 \rangle
\end{align*}
since ${\rm R}\Gamma_{\geq i+a}(X, \mathcal{O}_{p_i})^{st} = \Gamma_{\geq 0}(X, \mathcal{O}_{p_i})^{st} = L_i^{st} = L_i$ since $L_i = R/I(p_i) \cong k[z_i]$ has depth $1$ and is already MCM. This exceptional sequence is not strong, but we claim that it becomes strong upon replacing $k^{st}(1)$ by $\m(1) = k^{st}(1)[-1]$. \bigskip

Since ${\rm st}\circ {\rm R}\Gamma_{\geq i+a}(X, -)$ is fully faithful the $L_i$ are pairwise orthogonal. To verify that the exceptional sequence is strong, we calculate the remaining extension groups by Serre duality:
\begin{align*}
	\underline{{\rm Ext}}_{gr R}^n(\m(1), L_i) &= {\rm D}\underline{{\rm Ext}}_{gr R}^n(L_i, \m(2)) \\
	&= {\rm D} \underline{{\rm Ext}}_{gr R}^{n-1}(L_i, k^{st}(2))\\
	&= {\rm D} \underline{{\rm Ext}}_{gr R}^{n-1}(L_i, k(2)).
\end{align*}
We claim that
\[
	\dim_k \underline{{\rm Ext}}_{gr R}^n(\m(1), L_i) = \dim_k \underline{{\rm Ext}}_{gr R}^{n-1}(L_i, k(2)) = \begin{cases} 1, & n = 0 \\ 0, & n \neq 0. \end{cases}
\]
Equivalently, the minimal complete resolution of $L_i$ over $R$ looks like
	\[
		\xymatrix@C=20pt@R=5pt{\cdots \ar[r] & C_1 \ar[r] & R \ar[dr] \ar[rr] && R(2) \ar[r]& C_{-2} \ar[r] & \cdots \\
		& & & L_i \ar[ur] & }
	\]
The right hand tail of the minimal complete resolution dualises to a minimal projective resolution of the $R$-dual
\[
	\dots \to C_{-2}^* \to C_{-1}^* \to L_i^* \to 0
\]
and so the above claim $C_{-1} \cong R(2)$ is equivalent to $C_{-1}^* \cong R(-2)$, meaning that $L_i^*$ is singly generated in degree $2$. Now, writing $R = S/(Q, Q')$ for $S = k[x,y,z]$ and $(Q, Q')$ a regular sequence of quadrics, by a well-known change-of-rings result (see \cite[Lemma 2.5]{BEH}) we have isomorphisms of graded $R$-module
\[
	{\rm Ext}^n_{S}(L_i, S) \cong {\rm Ext}^{n-2}_{R}(L_i, R).
\]
Setting $n = 2$ gives ${\rm Ext}^2_S(L_i, S) \cong {\rm Hom}_R(L_i, R) = L_i^*$. Writing $L_i = S/(l, l')$ for a regular sequence of linear forms $(l, l')$ in $S$, self-duality of the Koszul complex $\mathbb{K}_S(l, l')$ gives 
\[
	{\rm Ext}^2_S(L_i, S) \cong L_i(-2)
\]
and so $L_i^* \cong L_i(-2)$ is singly generated in degree $2$, as we wanted, and so we obtain that
\[
	\dim_k \underline{{\rm Ext}}_{gr R}^n(\m(1), L_i) = \dim_k \underline{{\rm Ext}}_{gr R}^{n-1}(L_i, k(2)) = \begin{cases} 1, & n = 0 \\ 0, & n \neq 0. \end{cases}
\]
This proves that the exceptional sequence is strong. \bigskip

Finally, for dimension reasons we have ${\rm \underline{End}}_{gr R}(T) \cong kQ$, where $Q$ is the quiver
	{\tiny{
			\[
	\xymatrix@R=5pt@C=20pt{&\bullet\\
		&\bullet\\
		\bullet \ar[uur]^{p_1} \ar[ur]_{\ p_2} \ar[dr]^{\ p_3} \ar[ddr]_{p_4}& \\
		&\bullet\\
		&\bullet
	}
	\]
}}
\end{proof}

\subsubsection*{The algebra $R_{Y_1}$}
Let $R = R_{Y_1}$ for $4$ points $Y_1 \subseteq \mathbb{P}^1$, which is isomorphic to the hypersurface ring $R_{Y_1} = k[x,y]/(f_\lambda)$, with $f_\lambda = l_1l_2l_3l_4 = xy(x-y)(x-\lambda y)$, $\lambda \neq 0,1$, a product of $4$ linear forms in normalised form. Using the same approach, we now prove:
\begin{thm}\label{squidtiltproof} Let $R = R_{Y_1}$, so that $X = \proj R = Y_1 \subseteq \mathbb{P}^1$ is a set of $4$ points $\{p_i = V(l_i)\}$. Let $L_i = R/l_i$. Then there is a full strong exceptional collection
	\[
	\underline{{\rm MCM}}^{\mathbb{Z}}(R) = \langle \m(1), \m^2(2), L_1, L_2, L_3, L_4 \rangle
	\]
	with the endomorphism algebra of $T = \m(1) \oplus \m^2(2) \oplus \left( \bigoplus_{i=1}^4 L_i \right)$ given by the squid algebra $Sq(2,2,2,2;\lambda)$ presented above.
\end{thm}

As before, the calculations of the Tate cohomology groups $ \underline{{\rm Ext}}_{gr R}^n(\m^j(j), L_i)$ are the most involved part of the proof, and we collect them in the next lemma.
\begin{lem}\label{squidstrongcalc} We have
	\begin{align*}
		\dim_k \underline{{\rm Ext}}_{gr R}^n(\m(1), L_i) &= \begin{cases} 1, & n = 0 \\ 0, & n \neq 0 \end{cases} \\
		\dim_k \underline{{\rm Ext}}_{gr R}^n(\m^2(2), L_i) &= \begin{cases} 1, & n = 0 \\ 0, & n \neq 0 \end{cases} \\
		\dim_k \underline{{\rm Ext}}_{gr R}^n(\m^1(1), \m^2(2)) &= \begin{cases} 2, & n = 0 \\ 0, & n \neq 0. \end{cases}
	\end{align*}
\end{lem}
Before we prove this, we begin by writing down some relevant complete resolutions for calculations. Note that the complete resolution of $L_i$ is easy to obtain, since $L_i = R/l_i$ and $l_i$ is part of the matrix factorization $(l_i, f_\lambda/l_i)$ in $S = k[x,y]$. Its complete resolution is then
\[
	C(L_i):	\cdots \to R(-4) \xrightarrow{f_\lambda/l_i} R(-1) \xrightarrow{l_i} R \xrightarrow{f_\lambda/l_i} R(3) \xrightarrow{l_i} R(4) \to \cdots
\]
The complete resolution of $k^{st}$ is known and can be obtained by a classical method of Eisenbud, see \cite{Eis}. In this case, decompose $f_\lambda = x\cdot f_x + y \cdot f_y$ in $S$ for two cubic polynomials $f_x, f_y$. Note that in $\mathsf{char}\ \! k \neq 2$, we can take $f_x = \frac{1}{4} \partial_x f_\lambda$ and $f_y = \frac{1}{4} \partial_y f_\lambda$ by the Euler identity. We then have a matrix factorization $(A,B)$ of $f_\lambda$
\begin{align*}
A=\begin{pmatrix} x & y \\ -f_y & f_x \end{pmatrix} & B=\begin{pmatrix} f_x & -y \\ f_y & x \end{pmatrix}
\end{align*}
giving rise to the minimal complete resolution of $k^{st} = {\rm coker}(A)$ below
\[
	C(k^{st}):	\cdots \to R(-4) \oplus R(-2) \xrightarrow{B} R(-1)^{\oplus 2} \xrightarrow{A} R \oplus R(2) \xrightarrow{B} R(3)^{\oplus 2} \xrightarrow{A} R(4)\oplus R(6) \to \dots
\]
We can now prove the lemma.
\begin{proof}
	First, letting $T = \m(1) \oplus \m^2(2) \oplus \left( \bigoplus_{i=1}^4 L_i \right)$, note that Lemma \ref{cutoff} gives
	\[
		\underline{{\rm Ext}}^n_{gr R}(T, T) = {\rm Ext}_{gr R}^n(T, T)
	\]
	which vanishes for $n < 0$, taking care of all negative groups. Now, we have $a = a(R) = 2$, hence Serre duality gives
\begin{align*}
	\underline{{\rm Ext}}_{gr R}^n(\m(1), L_i) &= {\rm D}\underline{{\rm Ext}}_{gr R}^n(L_i, \m(3)) \\
	&= {\rm D} \underline{{\rm Ext}}_{gr R}^{n-1}(L_i, k^{st}(3))\\
	&= {\rm D} \underline{{\rm Ext}}_{gr R}^{n-1}(L_i, k(3))
\end{align*}
and the first equality follows from $C(L_i)$. For the next equality, consider the extension
\[
	\xi: 0 \to \m/\m^2 \to R/\m^2 \to k \to 0
\]
where $\m/\m^2 \cong k(-1)^{\oplus 2}$. This gives rise to an exact triangle in $ \underline{{\rm MCM}}^{\mathbb{Z}}(R)$
\[
	\xi^{st}:	(\m/\m^2)^{st} \to (R/\m^2)^{st} \to k^{st} \to (\m/\m^2)^{st}[1]
\]
which can be rewritten as
\[
	\xi^{st}: (k^{st}(-1))^{\oplus 2} \to \m^2[1] \to k^{st} \to (k^{st}(-1))^{\oplus 2}[1].
\]
Serre duality gives
\begin{align*}
	\underline{{\rm Ext}}_{gr R}^n(\m^2(2), L_i) &= {\rm D} \underline{{\rm Ext}}_{gr R}^n(L_i, \m^2(4)) \\
	&= {\rm D}\underline{{\rm Ext}}^{n}_{gr R}(L_i, \m^2[2]) \\
	&= {\rm D}\underline{{\rm Ext}}^{n+1}_{gr R}(L_i, \m^2[1])
\end{align*}
Applying the long exact sequence of Tate cohomology to $\xi^{st}$ and reading off from $C(L_i)$, while remembering that the above groups vanish for $n < 0$, we obtain $ \underline{{\rm Ext}}^{n+1}_{gr R}(L_i, \m^2[1]) = 0$ for $n \neq 0$ and the long exact sequence amounts to
\[
	0 \to \underline{{\rm Ext}}^{0}_{gr R}(L_i, k^{st}) \to \underline{{\rm Ext}}^{1}_{gr R}(L_i, k^{st}(-1))^{\oplus 2} \to \underline{{\rm Ext}}^{1}_{gr R}(L_i, \m^2[1]) \to 0.
\]
A dimension count then gives $\dim_k \underline{{\rm Ext}}^1_{gr R}(L_i, \m^2[1]) = \dim_k \underline{{\rm Ext}}_{gr R}^0(\m^2(2), L_i) = 1.$

Finally we compute
\begin{align*}
\underline{{\rm Ext}}_{gr R}^n(\m(1), \m^2(2)) &= {\rm D} \underline{{\rm Ext}}_{gr R}^n(\m^2, \m(1))\\
&= {\rm D}\underline{{\rm Ext}}_{gr R}^n(\m^2[1], \m[1](1))\\
&= {\rm D}\underline{{\rm Ext}}_{gr R}^n(\m^2[1], k^{st}(1))\\
&= {\rm D}\underline{{\rm Ext}}_{gr R}^n(\m^2[1], k(1))\\
\end{align*}
Applying the long exact sequence from $\xi^{st}$, note that $ \underline{{\rm Ext}}^*_{gr R}(k^{st}, k(1)) = 0$ from the structure of $C(k^{st})$, giving $ \underline{{\rm Ext}}^n_{gr R}(\m^2[1], k(1)) \cong \underline{{\rm Ext}}^n_{gr R}(k^{st}(-1), k(1))^{\oplus 2}$ whose dimension is as stated.
\end{proof}

We can now prove the theorem.
\begin{proof}
	We use Orlov's theorem with cutoff $i = -a = -2$ and $X = \proj R$ to obtain a semiorthogonal decomposition (living inside $ {\rm D}^b(\mathsf{grmod}_{\geq -2}R)$)
	\begin{align*}
		(-)_{[\geq -2]}\left( \underline{{\rm MCM}}^{\mathbb{Z}}(R)\right) &= \langle k(2), k(1), {\rm R}\Gamma_{\geq 0} {\rm D}^b(X) \rangle \\
		&= \langle k(2), k(1), L_1, \dots, L_4 \rangle.
	\end{align*}
	This full exceptional collection is not strong, and we will apply a right mutation to obtain
	\[
		R: \left\langle k(2), k(1) \right\rangle \mapsto \left\langle k(1), R/\m^2(2) \right\rangle 
	\]
	To see this, note that the extension
	\[
		\xi: \mathfrak{m}/\mathfrak{m}^2 \to R/\mathfrak{m}^2 \to k \to \mathfrak{m}/\mathfrak{m}^2[1]
	\] is isomorphic to the universal extension
	\[
		{\rm Ext}^1(k, k(-1))^* \otimes_k k(-1) \to R_{k}(k(-1)) \to k \xrightarrow{coev} {\rm Ext}^1(k, k(-1))^{*} \otimes_k k(-1)[1].
	\]
	This calculates the right mutation $R_{k}(k(-1)) = R/\mathfrak{m}^2$, and similarly $R_{k(2)}\left(k(1)\right) = \left(R/\mathfrak{m}^2\right)(2)$. After mutating and desuspending the first two terms, we obtain the resulting exceptional collection
\[
	\left\langle k(1)[-1], \left(R/\mathfrak{m}^2\right)(2)[-1], L_1, \dots, L_4\right\rangle
\]
Upon stabilising, this is sent to the full exceptional collection
\[
	\underline{{\rm MCM}}^{\mathbb{Z}}(R) = \langle \m(1), \m^2(2), L_1, L_2, L_3, L_4 \rangle
\]
which is strong by Lemma \ref{squidstrongcalc}. It remains to calculate the endomorphism algebra. Letting $T = \m(1) \oplus \m^2(2) \oplus \left( \bigoplus_{i=1}^4 L_i\right)$, we have seen that ${\rm \underline{End}}_{gr R}(T) \cong {\rm End}_{gr R}(T)$. Consider the following morphisms:
	\[
		\xymatrix@C2pc@R0.3pc{	&& L_1 \\
				&& L_2 \\
				\mathfrak{m}(1) \ar@<.5ex>[r]^{x} \ar@<-.5ex>[r]_{y} & \mathfrak{m}^2(2)  \ar[ruu]^{q_1} \ar[ru]_{q_2} \ar[rd]^{q_3} \ar[rdd]_{q_4} \\
				&& L_3 \\
				&& L_4 \\
}
	\]
	with $q_i: \mathfrak{m}^2(2) \to L_i$ induced by $x,y \mapsto \overline{x}, \overline{y} \in R/l_i = L_i$. These satisfy the relations of the squid algebra $Sq(2,2,2,2;\lambda)$ and so there is an algebra morphism
	\[
	\varphi: Sq(2,2,2,2;\lambda) \to {\rm End}_{gr R}(T). 
	\]
	We have computed above the dimensions of the Hom spaces between the indecomposable summands of $T$, and so one sees that the morphisms considered generate ${\rm End}_{gr R}(T)$, so that $\varphi$ is surjective. Comparing the dimensions of ${\rm End}_{gr R}(T)$ and $Sq(2,2,2,2;\lambda)$ shows that they are equal, and so $\varphi$ is an isomorphism. 
\end{proof}

Finally, to end this section and as alluded to above, let us mention a basic fact relating graded MCM modules over the previous rings $R$ to MCM modules over the completion $\widehat{R}$ at $\m$. It is well-known (see \cite{YosBook}) that the completion functor $\widehat{M} \cong M \otimes_R \widehat{R}$ preserves indecomposables MCM modules, and that two graded modules satisfy $\widehat{M} \cong \widehat{N}$ if and only if $M \cong N(n)$ for some $n \in \Z$. It isn't however always true that indecomposable MCM $\widehat{R}$-modules arise as the completion of a graded module, but this does hold in special circumstances. For instance:
\begin{cor} Let $\widehat{R}$ be the completion of $R = R_{Y_i}$ for $i = 1,2$ at $\mathfrak{m} = R_{\geq 1}$. Then every indecomposable MCM $\widehat{R}$-module is the completion $\widehat{M}$ of an indecomposable graded MCM $R$-module $M$. 
\end{cor}
\begin{proof} By \cite[Prop. 1.5]{KMVdB}, the completion functor $\widehat{(-)}: \underline{{\rm MCM}}^{\Z}(R) \to \underline{{\rm MCM}}(\widehat{R})$ identifies with the universal functor to the triangulated hull of the orbit category $\underline{{\rm MCM}}^{\Z}(R)/(1)$. By Keller's theorem \cite{Ke}, this functor is essentially surjective whenever $ \underline{{\rm MCM}}^{\mathbb{Z}}(R) \cong {\rm D}^b(\mathcal{H})$ for $\mathcal{H}$ an hereditary category, and $(1)$ moves away from the heart $\mathcal{H}$ in that $X(n) \notin \mathcal{H}$ for all $n \gg 0$ and $X \in \mathcal{H}$. We do this in each case as follows:
	\begin{enumerate}[i)]
	\item $(R = R_{Y_2})$: Take $\mathcal{H} = \modsf kQ$, then $(1) = \Se = \tau [-1]$ moves away from the heart $\mathcal{H}$.
	\item $(R = R_{Y_1})$: Take $\mathcal{H} = \coh \X$ where $\X$ is the weighted projective line of type $(2,2,2,2;\lambda)$ derived equivalent to $Sq(2,2,2,2;\lambda)$, then $\Se_R = (2)$ corresponds to $\Se_\X = - \otimes \omega_\X[1]$ on $ {\rm D}^b(\X)$, and one easily sees that $(1)$ moves away from the heart $\mathcal{H}$. 
\end{enumerate}
\end{proof}

%% file: chapters/classifications/conics.tex
\section{Graded MCM modules over the cone of $4$ points on $\mathbb{P}^2$ in general position}
Let $R = R_{Y_2} = k[x,y,z]/(Q, Q')$. We now investigate the structure of $ \underline{{\rm MCM}}^{\mathbb{Z}}(R)$. Thinking ahead, we will modify the equivalence of the previous section to simplify calculations down the line. Recall that $\Se_R(-) = - \otimes_R \omega_R[\dim R - 1] = (1)$ is a Serre functor for $ \underline{{\rm MCM}}^{\mathbb{Z}}(R)$ since $R$ has isolated singularities. Writing $k^{st} = \Se^{-1}_R(k^{st}(1))$, the full strong exceptional sequence of the previous section
\[
	\langle k^{st}(1)[-1], L_1, \dots, L_4 \rangle
\]
can be exchanged for the exceptional sequence
\[
	\langle L_1, \dots, L_4, k^{st} \rangle.
\]
This sequence is also full by Serre duality, since it has trivial right orthogonal category, and is also strong as one immediately verifies that $ \underline{{\rm Ext}}_{gr R}^n(L_i, k^{st}) = \underline{{\rm Ext}}_{gr R}^n(L_i, k) = 0$ for $n \neq 0$. Letting $T = (\bigoplus_{i=1}^4 L_i) \oplus k^{st}$, we have ${\rm \underline{End}}_{gr R}(T) = kQ$ with $Q$ the ``four subspace'' quiver
			\[
				\xymatrix@R=5pt@C=20pt{1 \ar[ddr]^{p_1}&\\
					2 \ar[dr]_{p_2}&\\
		& 0 \\
		3 \ar[ur]^{p_3}&\\
		4 \ar[uur]_{p_4}&
	}
\]
The module $T$ is also a tilting object in the opposite category $ \underline{{\rm MCM}}^{\mathbb{Z}}(R)^{op}$ with endomorphism algebra $kQ^{op}$, and $\modsf kQ^{op} = kQ\ \! \mathsf{mod}$. From Tilting theory we obtain:
\begin{prop}\label{oppD4tilt} There is a contravariant equivalence of triangulated categories
\[
	F: \underline{{\rm MCM}}^{\mathbb{Z}}(R)^{op} \xrightarrow{\cong} {\rm D}^b(kQ\ \! \mathsf{mod}) 
\]
onto the bounded derived category of left $kQ$-modules, or equivalently covariant quiver representations of $Q$.
\end{prop}
We will need an explicit description of the functor $F$. Given $M, N \in \underline{{\rm MCM}}^{\mathbb{Z}}(R)$, recall that we write
\[
	{\rm R\underline{Hom}}_{gr R}(M , N) = {\rm Hom}_{gr R}(C(M), C(N))
\]
where $C(M), C(N)$ are complete resolutions of $M, N$, so that 
\[
	{\rm H}^n{\rm R\underline{Hom}}_{gr R}(M , N) \cong \underline{{\rm Ext}}_{gr R}^n(M, N).
\]
Let $\mathcal{R} = {\rm R\underline{Hom}}_{gr R}(T , T) $ be the derived endomorphism algebra of $T$, quasi-isomorphic to the usual endomorphism algebra ${\rm \underline{End}}_{gr R}(T) = kQ$ via the standard truncation zigzag
\[
	\varphi: \mathcal{R} \xleftarrow{\sim} \tau^{\leq 0}\mathcal{R} \xrightarrow{\sim} {\rm H}^0\mathcal{R} = kQ.
\]
The zigzag of quasi-isomorphisms $\varphi$ induces an equivalence of derived categories
\[
	\varphi_*: {\rm D}^{\perf}(\mathcal{R}\ \! \mathsf{Mod}) \xrightarrow{\cong} {\rm D}^{\perf}(kQ\ \! \mathsf{Mod}) = {\rm D}^b(kQ\ \! \mathsf{mod}) 
\]
given by
\[
	X \mapsto kQ \otimes^{\mathbb{L}}_{\tau^{\leq 0}\mathcal{R}} X.
\]
We then define $F$ as the composite equivalence $F = \varphi_* \circ {\rm R\underline{Hom}}_{gr R}(- , T) $
\[
	\xymatrix@C6pc{ F:\underline{{\rm MCM}}^{\mathbb{Z}}(R)^{op} \ar[r]^{\cong}_-{{\rm R\underline{Hom}}_{gr R}(- , T)} & {\rm D}^{\perf}(\mathcal{R}\ \! \mathsf{Mod}) \ar[r]^{\cong}_-{\varphi_*} & {\rm D}^b(kQ\ \! \mathsf{mod}). }
\]

Working with the contravariant equivalence $F$ will turn out to be easier in practice. Note that since $kQ$ is hereditary, any complex in $ {\rm D}^b(kQ\ \! \mathsf{mod}) $ is formal, that is to say
\[
	X \cong \bigoplus_{n \in \Z}{\rm H}^n(X)[-n]
\]
and it suffices to understand the cohomology modules ${\rm H}^n(X)$.

\begin{lem}\label{modulestructure} Let $M \in \underline{{\rm MCM}}^{\mathbb{Z}}(R)$. Then for each $n \in \Z$, we have an isomorphism of left $kQ$-module
	\[
		{\rm H}^n(F(M)) \cong \underline{{\rm Ext}}_{gr R}^n(M, T).
	\]
\end{lem}
\begin{proof}
	The action of $Z^{0}\mathcal{R} \subset \tau^{\leq 0}\mathcal{R} \subset \mathcal{R}$ on $ {\rm R\underline{Hom}}_{gr R}(M , T) $ by post-composition descends to the action of ${\rm H}^0\mathcal{R} = kQ$ on $\varphi_* ({\rm R\underline{Hom}}_{gr R}(M, T) )$ by post-composition. Taking $n$-th cohomology gives the above module structure.
\end{proof}
Equivalently, the quiver representation ${\rm H}^n(F(M))$ is given by

\[
	\xymatrix@R=10pt@C=30pt{ 
		\underline{{\rm Ext}}_{gr R}^n(M, L_1) \ar[ddr] & \\
		\underline{{\rm Ext}}_{gr R}^n(M, L_2) \ar[dr] & \\
		& \underline{{\rm Ext}}_{gr R}^n(M, k^{st}) \\
		\underline{{\rm Ext}}_{gr R}^n(M, L_3) \ar[ur] & \\
		\underline{{\rm Ext}}_{gr R}^n(M, L_4) \ar[uur] &	}
\]
with linear maps induced from the morphism $L_i \to k^{st}$. To classify the indecomposable graded MCM $R$-module, it then suffices to classify the indecomposable representations of $Q$ and find the indecomposables $M$ for which ${\rm H}^0(F(M))$ exhaust this list. One obtains the remaining indecomposables by suspension $M \mapsto M[n]$.

\subsection{Four Subspace Problem}
In this section we review the known classification of representations of the $\widetilde{D}_4$ quiver $Q$ with ``four subspace'' orientation. A good reference for this classification is \cite[XIII.3]{SimSkow} and Happel's monograph \cite{HappelBook} for the derived category aspects, and we refer to Appendix \ref{artinappendix} for standard definitions and generalities on quiver representations.

Since $Q$ is an extended Dynkin quiver, the structure of the module category $kQ\ \! \modsf$ contains three types of Auslander-Reiten components, namely the preprojective, preinjective components and regular components, and the first two become attached in the derived category. Let 
\[
	\Se_{kQ}(-) = - \otimes^{\mathbb{L}}_{kQ} {\rm D}(kQ)
\]
be the Serre functor on $ {\rm D}^b(kQ\ \! \mathsf{mod})$, and $\tau = \Se_{kQ} \circ [-1]$ the Auslander-Reiten translate\footnote{This is somewhat anachronistic. See Appendix \ref{artinappendix}.}. Letting $e_i$ be the idempotent at the $i$-th vertex, the indecomposable projectives $P(i) = kQ e_i$ are given for $i = 0,1,2,3,4$ by
\[
\xymatrix@R=3pt@C=15pt{ 
	0\ar[ddr] &  && k\ar[ddr]^1 &  &&  0\ar[ddr] &  &&  0\ar[ddr] &  &&  0\ar[ddr] &   \\
	0\ar[dr] &  && 0\ar[dr] & &&  k \ar[dr]_1 &  &&  0 \ar[dr] &  &&  0 \ar[dr] &   \\
	  & k &&  & k &&   & k &&   & k &&   & k  \\
	  0\ar[ur] &  && 0\ar[ur]  &  &&  0\ar[ur]  &  &&  k\ar[ur]^1   &  &&  0\ar[ur] &   \\
	0\ar[uur] &  && 0\ar[uur] &  &&  0\ar[uur] &  &&  0\ar[uur] &  &&  k\ar[uur]_1 & 
}
\]
and the indecomposable injectives $I(i) = {\rm D}(P(i)_{kQ}) = {\rm D}(e_ikQ)$ are given by 
\[
\xymatrix@R=3pt@C=15pt{ 
	k\ar[ddr]^1 &  && k\ar[ddr] &  &&  0\ar[ddr] &  &&  0\ar[ddr] &  &&  0\ar[ddr] &   \\
	k\ar[dr]_1 &  && 0\ar[dr] & &&  k \ar[dr] &  &&  0 \ar[dr] &  &&  0 \ar[dr] &   \\
	  & k &&  & 0 &&   & 0 &&   & 0 &&   & 0  \\
	  k\ar[ur]^1 &  && 0\ar[ur]  &  &&  0\ar[ur]  &  &&  k\ar[ur]   &  &&  0\ar[ur] &   \\
	k\ar[uur]_1 &  && 0\ar[uur] &  &&  0\ar[uur] &  &&  0\ar[uur] &  &&  k\ar[uur] & 
}
\]

The Auslander-Reiten translate exchanges those modules up to suspension, via $\tau P(i) = I(i)[-1]$ and $\tau^{-1}I(i) = P(i)[1]$.
The family of indecomposable complexes $\{\tau^m P(i)\}$ for $i = 0,1,2,3,4$ and $m \in \Z$ form a connected component in the Auslander-Reiten quiver $\Gamma ({\rm D}^b(kQ\ \! \mathsf{mod}))$ of the form (see \cite[I.5]{HappelBook})
{\tiny{
\[
\xymatrix@R.8pc@C1.3pc{
	\cdots \ar[ddr] &&I(1)[-1] \ar[ddr] && P(1) \ar[ddr] && \tau^{-1}P(1) \ar[ddr]\\
	\cdots \ar[dr] &&I(2)[-1] \ar[dr] && P(2)  \ar[dr] &&  \tau^{-1}P(2)  \ar[dr]\\
&I(0)[-1] \ar[uur] \ar[ur] \ar[dr] \ar[ddr] && P(0)\ar[uur] \ar[ur] \ar[dr] \ar[ddr] && \tau^{-1}P(0)\ar[uur] \ar[ur] \ar[dr] \ar[ddr] && \cdots \\
\cdots \ar[ur]&&I(3)[-1] \ar[ur] && P(3) \ar[ur] && \tau^{-1}P(3) \ar[ur]\\
\cdots \ar[uur] &&I(4)[-1] \ar[uur] && P(4) \ar[uur] && \tau^{-1}P(4) \ar[uur]\\
}
\]
}}
called a \emph{transjective} component. The indecomposables complexes $\tau^{m}P(i)$ for $m \leq 0$ are modules called preprojective modules, while the indecomposables complexes of the form $\tau^m I(i)$ for $m \geq 0$ are modules called preinjective modules. Note that the transjective component consists of the preprojective component in $kQ\ \! \mathsf{mod}$ attached to the suspended preinjective component in $(kQ\ \! \mathsf{mod})[-1]$. We denote the transjective component by $\mathcal{PI}$ and its $n$-th suspension by $\mathcal{PI}[n]$.

One can write down the preprojective and preinjective modules explicitly but we will be satisfied with the description 
\[
	\bigcup_{n \in \Z} \mathcal{PI}[n] = \{ \tau^{m} P(i)[n] \ | \ i = 0,1,2,3,4, \textup{ and } m,n \in \Z \}.
\]

\subsubsection*{Regular components}
A module whose indecomposable summands are neither preprojective nor preinjective is called \emph{regular}. Let $\mathcal{R}(Q)$ denote the category of regular modules. Its structure is as follows.

\begin{prop} The category $\mathcal{R}(Q)$ is a full abelian subcategory closed under extension and under the Auslander-Reiten $\tau$. Moreover, $\mathcal{R}(Q)$ is serial in that every object has a unique finite composition series with simple regular factors. 
\end{prop}

Since $Q$ is extended Dynkin, the Auslander-Reiten components in $\mathcal{R}(Q)$ break down into a $\mathbb{P}^1$ family of disjoint tubes  of finite ranks $\{\T_{\lambda} \}_{\lambda \in \mathbb{P}^1}$. The additive closures $\mathsf{add}\ \! \T_{\lambda}$ are abelian subcategories closed under extension and under $\tau$, with finitely many simples (whose number is the rank of the tube). Moreover, the categories $\mathsf{add}\ \! \T_\lambda$ are pairwise Hom and Ext orthogonal.

We begin with a description of the simple regular modules in $\mathcal{R}(Q)$. Consider quiver representations of the form
\[
	\xymatrix@R=3pt@C=20pt{
		k  \ar[ddr]^{\varphi_1}& \\
		k  \ar[dr]_{\varphi_2} &\\
		& k^2 \\
		k  \ar[ur]^{\varphi_3}& \\
		k  \ar[uur]_{\varphi_4}&
	}
\]
with all $\varphi_i \neq 0$. One can identify this with the set $\{ {\rm im}\ \! \varphi_i \}$ of $4$ lines through the origin in $k^2$. Assume that the first three lines are distinct, so that up to change of basis the above representation is given by 
\[
	(\varphi_1, \varphi_2, \varphi_3, \varphi_4) = (\begin{bsmallmatrix} 0 \\1 \end{bsmallmatrix}, \begin{bsmallmatrix} 1 \\1 \end{bsmallmatrix}, \begin{bsmallmatrix} 1 \\0 \end{bsmallmatrix}, \begin{bsmallmatrix} \lambda_0 \\ \lambda_1 \end{bsmallmatrix} ).
	\]
	Writing $\lambda = (\lambda_0, \lambda_1)$, we denote the above representation by $R_\lambda$. It is clear that the isomorphism class of $R_\lambda$ only depends on the point $\lambda = [\lambda_0 : \lambda_1] \in \mathbb{P}^1$.
	
	Again considering the $4$-tuples of lines given by $\big( \varphi_1, \varphi_2, \varphi_3, \varphi_4 \big)$, now assume that exactly two of the lines collide. This yields a partition of $\{\varphi_1, \varphi_2, \varphi_3, \varphi_4 \} = \{\varphi_i, \varphi_j \} \coprod \{\varphi_p, \varphi_q \}$ where ${\rm im}\ \varphi_i = {\rm im}\ \varphi_j$ and ${\rm im}\ \varphi_p \neq {\rm im}\ \varphi_q$. Keeping track of ordering, there are $\binom{4}{2} = 6$ such partitions. Let us define corresponding representations as
	\begin{align*}
R_0^+: (\begin{bsmallmatrix} 0 \\1 \end{bsmallmatrix}, \begin{bsmallmatrix} 1 \\1 \end{bsmallmatrix}, \begin{bsmallmatrix} 1 \\0 \end{bsmallmatrix}, \begin{bsmallmatrix} 0 \\ 1 \end{bsmallmatrix} ) && R_0^-: (\begin{bsmallmatrix} 1 \\1 \end{bsmallmatrix}, \begin{bsmallmatrix} 0 \\1 \end{bsmallmatrix}, \begin{bsmallmatrix} 0 \\1\end{bsmallmatrix}, \begin{bsmallmatrix} 1 \\ 0 \end{bsmallmatrix} ) \\
R_1^+: (\begin{bsmallmatrix} 0 \\1 \end{bsmallmatrix}, \begin{bsmallmatrix} 1 \\1 \end{bsmallmatrix}, \begin{bsmallmatrix} 1 \\0 \end{bsmallmatrix}, \begin{bsmallmatrix} 1 \\ 1 \end{bsmallmatrix} ) && R_1^-: (\begin{bsmallmatrix} 1 \\1 \end{bsmallmatrix}, \begin{bsmallmatrix} 0 \\1 \end{bsmallmatrix}, \begin{bsmallmatrix} 1 \\1 \end{bsmallmatrix}, \begin{bsmallmatrix} 1 \\ 0 \end{bsmallmatrix} ) \\
R_\infty^+: (\begin{bsmallmatrix} 0 \\1 \end{bsmallmatrix}, \begin{bsmallmatrix} 1 \\1 \end{bsmallmatrix}, \begin{bsmallmatrix} 1 \\0 \end{bsmallmatrix}, \begin{bsmallmatrix} 1 \\ 0 \end{bsmallmatrix} ) && R_\infty^-: (\begin{bsmallmatrix} 1 \\0 \end{bsmallmatrix}, \begin{bsmallmatrix} 1 \\0 \end{bsmallmatrix}, \begin{bsmallmatrix} 0 \\1 \end{bsmallmatrix}, \begin{bsmallmatrix} 1 \\ 1 \end{bsmallmatrix} )
	\end{align*}
We have chosen this normalization with the following properties in mind:
\begin{enumerate}[i.]
\item For $\lambda = 0, 1, \infty$, we have $R_{\lambda} = R_{\lambda}^{+}$.
\item The involution $R_{\lambda}^{+} \leftrightarrow R_{\lambda}^{-}$ corresponds to interchanging $\{\varphi_i, \varphi_j\}$ and $\{\varphi_p, \varphi_q \}$.
\end{enumerate}

\begin{figure}\label{Spm}
\begin{align*}
	\xymatrix@R=1pt@C=15pt{ 
		& k \ar[ddr]^1 &\\
		& 0 \ar[dr] &\\
S_0^+:		&  & k\\
		& 0 \ar[ur] &\\
		& k \ar[uur]_1 &
}
&&
	\xymatrix@R=1pt@C=15pt{ 
		& 0 \ar[ddr] &\\
		& k \ar[dr]_1 &\\
S_0^-:		&  & k\\
		& k \ar[ur]^1 &\\
		& 0 \ar[uur] &
}
\\
\\
	\xymatrix@R=1pt@C=15pt{ 
		& 0 \ar[ddr] &\\
		& k \ar[dr]_1 &\\
S_1^+:		&  & k\\
		& 0 \ar[ur] &\\
		& k \ar[uur]_1 &
}
&&
	\xymatrix@R=1pt@C=15pt{ 
		& k \ar[ddr]^1 &\\
		& 0 \ar[dr] &\\
S_1^-:		&  & k\\
		& k \ar[ur]^1 &\\
		& 0 \ar[uur] &
}
\\
\\
	\xymatrix@R=1pt@C=15pt{ 
		& 0 \ar[ddr] &\\
		& 0 \ar[dr] &\\
S_\infty^+:		&  & k\\
		& k \ar[ur]^1 &\\
		& k \ar[uur]_1 &
}
&&
	\xymatrix@R=1pt@C=15pt{ 
		& k \ar[ddr]^1 &\\
		& k \ar[dr]_1 &\\
S_\infty^-:		&  & k\\
		& 0 \ar[ur] &\\
		& 0 \ar[uur] &
}
\end{align*}
\caption{The simple regular modules $S_{0,1,\infty}^{\pm}$.}
\end{figure}
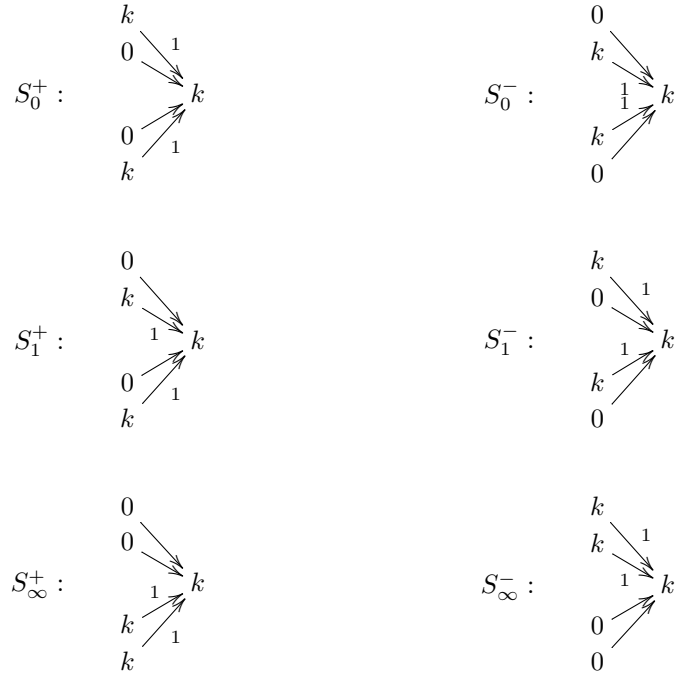

Finally, let us introduce additional representations $\{S_0^{\pm}, S_1^{\pm}, S_\infty^{\pm} \}$ as in Figure \ref{Spm}. We now have a complete set of simple regular modules for $kQ$.
 
\begin{prop}\label{regsimples} The following properties hold:
	\begin{enumerate}[1)]
	\item The set $\mathcal{S} = \{R_\lambda\}_{\lambda \in \mathbb{P}^1\setminus \{0, 1, \infty\}} \cup \{S_\lambda^{\pm}\}_{\lambda = 0, 1, \infty}$ is a complete set of isomorphism classes of simples in $\mathcal{R}(Q)$. In particular, each $S \in \mathcal{S}$ is indecomposable with ${\rm End}_{kQ}(S) = k$.
	\item The Auslander-Reiten translate acts by
		\begin{align*}
			&\tau R_\lambda \cong R_\lambda,\ \lambda \neq 0, 1, \infty\\
	&\tau S^{\pm}_\lambda = S^{\mp}_\lambda,\ \lambda = 0,1,\infty.
\end{align*}
\item For each $\lambda = 0, 1, \infty$, we have non-trivial short exact sequences
	\[
		0 \to S_\lambda^{\pm} \to R_\lambda^{\pm} \to S_\lambda^{\mp} \to 0.
	\]
\end{enumerate}
\end{prop}

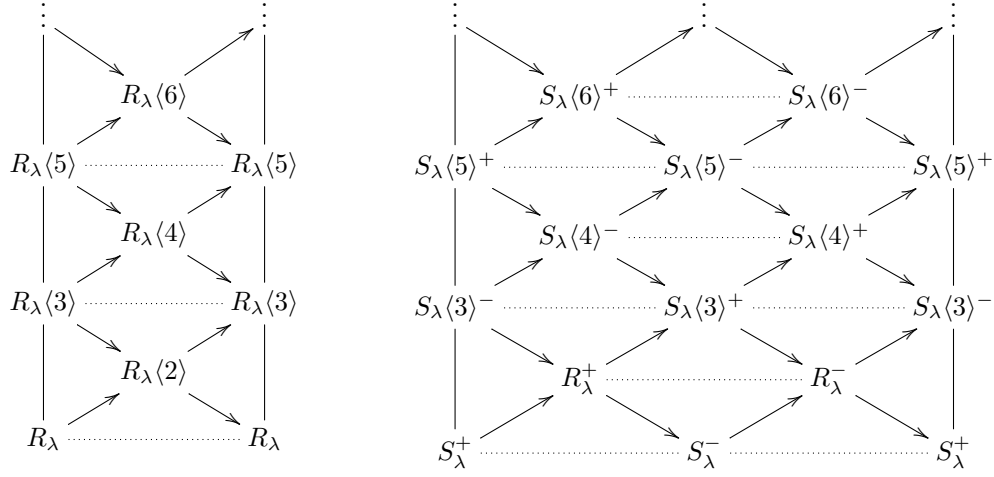
\begin{figure}\label{tubes}
\begin{align*}
\xymatrix@R=10pt@C=10pt{
	\vdots \ar@{-}[dd] \ar[dr] && \vdots \ar@{-}[dd] \\
	& R_\lambda \langle 6 \rangle \ar[ur] \ar[dr] & \\
	R_\lambda \langle 5 \rangle \ar[ur] \ar[dr] \ar@{.}[rr] \ar@{-}[dd] && R_\lambda \langle 5 \rangle \ar@{-}[dd]\\
& R_{\lambda}\langle 4 \rangle \ar[dr] \ar[ur] & \\
R_{\lambda}\langle 3 \rangle \ar[dr] \ar[ur] \ar@{.}[rr] \ar@{-}[dd] && R_{\lambda}\langle 3 \rangle \ar@{-}[dd] \\
& R_{\lambda} \langle 2 \rangle \ar[dr] \ar[ur] & \\
R_{\lambda} \ar[ur] \ar@{.}[rr] && R_{\lambda}
}
&&
\xymatrix@R=10pt@C=10pt{
\vdots \ar[dr] \ar@{-}[dd] && \vdots \ar[dr] && \vdots \ar@{-}[dd] \\
& S_{\lambda}\langle 6 \rangle^{+} \ar[dr] \ar[ur] \ar@{.}[rr] && S_\lambda \langle 6 \rangle^{-} \ar[dr] \ar[ur] & \\
S_\lambda \langle 5 \rangle^{+} \ar[dr] \ar[ur] \ar@{-}[dd] \ar@{.}[rr] && S_\lambda \langle 5 \rangle^{-} \ar[dr] \ar[ur] \ar@{.}[rr] && S_\lambda \langle 5 \rangle^{+} \ar@{-}[dd] \\
& S_\lambda \langle 4 \rangle^{-} \ar[dr] \ar@{.}[rr] \ar[ur] && S_\lambda \langle 4 \rangle^{+} \ar[dr] \ar[ur] \\
S_\lambda \langle 3 \rangle^{-} \ar@{-}[dd] \ar@{.}[rr] \ar[dr] \ar[ur] && S_\lambda \langle 3 \rangle^{+} \ar[dr] \ar[ur] \ar@{.}[rr] && S_\lambda \langle 3 \rangle^{-} \ar@{-}[dd] \\
& R_\lambda^{+} \ar[dr] \ar[ur] \ar@{.}[rr] && R_\lambda^{-} \ar[dr] \ar[ur] \\
S_\lambda^{+} \ar[ur] \ar@{.}[rr] && S_\lambda^{-} \ar[ur] \ar@{.}[rr] && S_\lambda^{+}
}
\end{align*}
\caption{Tubes $\T_\lambda$ of rank one and two.}
\end{figure}

The simples $\{R_\lambda\}_{\lambda \in \mathbb{P}^1 \setminus \{0, 1, \infty\}}$ each generate a rank one tube $\T_\lambda$, and $\{S_{0,1,\infty}^{\pm}\}$ generate tubes $\T_{0,1,\infty}$ of rank two, as in Figure \ref{tubes} with the edges attached. Since each category $\mathsf{add}\ \! \T_\lambda$ is a serial abelian category with finitely many simples, each indecomposable in $\T_\lambda$ is uniquely determined by its simple socle and its length (or equivalently its height in the tube). 
\begin{defn} Let $S \in \mathcal{S}$ be a simple regular module. We denote by $S\langle r \rangle$ the unique indecomposable regular module of length $r$ with socle $S$.
\end{defn}

Note that by Prop. \ref{regsimples} we have $R^{\pm}_\lambda = S^{\pm}\langle 2 \rangle$ for $\lambda = 0,1,\infty$. The remaining indecomposables are constructed by iterated extensions as follows. Let $S\langle r \rangle$ be an indecomposable regular with simple socle $S$. From Auslander-Reiten duality (or Serre duality) we have
\[
	{\rm Ext}^1_{kQ}(S\langle r \rangle, \tau S) \cong {\rm D}{\rm Hom}_{kQ}(S, S \langle r \rangle) = k
\]
and so there is a unique non-split extension
\[
	\xi: 0 \to \tau S \to (\tau S)\langle r+1 \rangle \to S\langle r \rangle \to 0.
\]
Its middle term must be indecomposable, since a decomposable module would have summands of length $\leq r$ and, as $S \langle r \rangle$ is uniserial, any surjection $(\tau S)\langle r +1 \rangle \twoheadrightarrow S \langle r \rangle$ would have to restrict to a surjection on some indecomposable summand, which would be an isomorphism for length reason and hence create a splitting of $\xi$. Applying this to our set of simple regulars gives unique short exact sequences 
\[
	0 \to R_\lambda \to R_\lambda \langle r+1 \rangle \to R_\lambda \langle r \rangle\to 0,\quad \quad \lambda \neq 0,1,\infty
\]
\[
	0 \to S_\lambda^{\pm} \to S_\lambda^{\pm} \langle r+1 \rangle \to S^{\mp}_\lambda \langle r \rangle\to 0,\quad \quad \lambda = 0,1,\infty
\]
and we have a description of the indecomposables in each tube $\{\T_\lambda\}_{\lambda \in \mathbb{P}^1}$, as in figure \ref{tubes}. As before, denote by $\T_\lambda[n]$ the $n$-th suspension of a tube in $ {\rm D}^b(kQ\ \! \mathsf{mod}) $.

\begin{thm}\label{ARquiverQ} The Auslander-Reiten quiver of ${\rm D}^b(kQ\ \! \mathsf{mod})$ is given by the union of disjoint components
\[
	\Gamma({\rm D}^b(kQ\ \! \mathsf{mod})) = \bigg(\bigcup_{n \in \Z} \mathcal{PI}[n] \bigg) \cup \bigg(\bigcup_{\lambda \in \mathbb{P}^1, n \in \Z} \T_{\lambda}[n]\bigg).
\]
In particular this yields a full classification of indecomposables in ${\rm D}^b(kQ\ \! \mathsf{mod})$.
\end{thm}

\subsection{Indecomposable graded MCM modules}

We want to describe the classification of indecomposables on the other side via $F: \underline{{\rm MCM}}^{\mathbb{Z}}(R)^{op} \xrightarrow{\cong} {\rm D}^b(kQ\ \! \mathsf{mod})$. First recall that the Serre functor on $ \underline{{\rm MCM}}^{\mathbb{Z}}(R)$ is given by $\Se_R = (1)$ and the Auslander-Reiten translate by $\tau = (1)[-1] = {\rm syz}_R^1(-)(1).$

Since $F$ is a contravariant equivalence, we have $F \circ \Se_R = \Se_{kQ}^{-1} \circ F$ and $F \circ \tau = \tau^{-1} \circ F$ by uniqueness of Serre functors\footnote{The functor $\Se^{-1}$ is a Serre functor for the opposite triangulated category.}. Since the indecomposable summands of a tilting object are sent onto the indecomposable projectives, we immediately have:
\begin{prop} We have isomorphisms
	\begin{align*}
		F\left({\rm syz}_R^{n-m}k^{st}(-m)\right) &\cong \tau^{m}P(0)[n]\\
		F\left( {\rm syz}_R^{n-m}L_i(-m) \right) &\cong \tau^{m} P(i)[n]
\end{align*}
for all $m,n \in \Z$ and $i=1,2,3,4$.
\end{prop}

The indecomposable regular modules in $kQ\ \! \mathsf{mod} \subset {\rm D}^b(kQ\ \! \mathsf{mod})$ are characterised amongst indecomposables by being $\tau$-periodic, of period $1$ or $2$. In $ \underline{{\rm MCM}}^{\mathbb{Z}}(R)$ this corresponds to indecomposables with a periodic minimal free resolution (then of period $1$ or $2$), and we have seen that they vary in families. 

Before we describe the periodic modules, let us apply some normalisations. For the remainder of this section, assume that $\mathsf{char}\ k \neq 2$. The complete intersection $R = k[x,y,z]/(Q, Q')$ is the homogeneous coordinate ring of a set $X$ of $4$ distinct points in $\mathbb{P}^2$ in general position\footnote{A set of $n+2$ points in $\mathbb{P}^n$ is in general position if no $n+1$ of them lie in a proper linear subspace. Here this means that no $3$ points in $X$ are collinear.}, and conversely any set of $4$ points in general position arise as the complete intersection of two conics. Moreover, there is a unique pencil of conics through $X$ and this pencil contains $3$ singular conics. 

\begin{figure}\label{fourlines}
\[
\xymatrix{
&  &  & \\
& p_1 \ar@{-}[u]^{l_1} \ar@{-}[l] \ar@{-}[d] \ar@{-}[r] & p_2 \ar@{-}[u] \ar@{-}[d] \ar@{-}[r]^{l_2} & \\
& p_4 \ar@{-}[d] \ar@{-}[l]^{l_4} \ar@{-}[r] & p_3\ar@{-}[d]^{l_3} \ar@{-}[r] & \\
&  &  & \\
}
\]
\caption{The four points $X = V(Q_0) \cap V(Q_\infty)$ in $\mathbb{P}^2$.}
\end{figure}
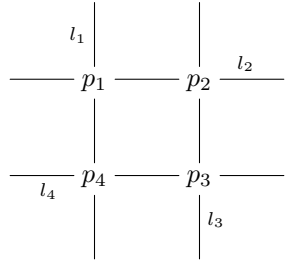

It is well-known that any sets of $4$ points in general position are related by a projective transformation, and so up to coordinate change we can assume that the points are given by $[\pm1: \pm1: \pm1] \in \mathbb{P}^2$, and so that the singular conics are given by equations
\begin{align*}
&Q_0 = x^2 - y^2\\
&Q_1 = x^2 - z^2\\
&Q_\infty = y^2 - z^2
\end{align*}
with general conic in the pencil given by $\{Q_\lambda = 0\}$ with $Q_\lambda = \lambda_0Q_\infty + \lambda_1Q_0$, $\lambda = [\lambda_0 : \lambda_1] \in \mathbb{P}^1$. With this normalisation, we have an isomorphism $R \cong k[x,y,z]/(Q_0, Q_\infty)$. We may picture $X = V(Q_0) \cap V(Q_\infty)$ as in figure \ref{fourlines}, where we let $l_i$ stand for both the line and the corresponding linear form in the factorizations
\begin{align*}
& Q_0 = (x-y)(x+y) = l_1l_3\\
& Q_\infty = (y-z)(y+z) = l_2l_4
\end{align*}
and from the figure \ref{fourlines} we have $I(p_i) = (l_i, l_{i+1})$, using cyclic indexing. Recall that for each point $p_i \in X$, we defined $L_i = R/I(p_i) = S/(l_i, l_{i+1})$ as its homogeneous coordinate ring thought of as an $R$-module, with $S = k[x,y,z]$.

Now for $\lambda = (\lambda_0, \lambda_1) \in k^2\setminus \{0\}$, let $(\Phi_{\lambda}^{+}, \Phi_{\lambda}^{-})$ be the pair of matrices over $S$ given by
\begin{align*}
\Phi_{\lambda}^{+} = \begin{bmatrix} \lambda_1 (x+y) & y+z \\
										 \lambda_0(z-y) & x-y \end{bmatrix}
										 &&
\Phi_{\lambda}^{-} = \begin{bmatrix} x-y & -y-z \\
                                        \lambda_0(y-z)  & \lambda_1(x+y)  \end{bmatrix}
\end{align*}
We have $\Phi_{\lambda}^{+} \Phi_{\lambda}^{-} = Q_{\lambda} \cdot {\rm I}_2 = \Phi_{\lambda}^{-} \Phi_{\lambda}^{+}$. Letting $\mu = (\mu_0, \mu_1)$ correspond to a different point in $\mathbb{P}^1$, the sequence $(Q_{\mu}, Q_{\lambda})$ is regular in $S$ and the pair $(\Phi^{+}_{\lambda}, \Phi^{-}_{\lambda})$ defines a matrix factorization of $Q_{\lambda}$ over $S' = S/(Q_{\mu})$. The module
\[
N_{\lambda} = {\rm coker}\big(R(-1)^{\oplus 2} \xrightarrow{\overline{\Phi}_{\lambda}^{+}} R^{\oplus 2} \big)
\]
is then an MCM module over $R \cong S'/(Q_\lambda)$, with isomorphism class independent of representative of $[\lambda_0: \lambda_1] \in \mathbb{P}^1$. Its complete resolution is then given by
\[
	C(N_\lambda):	\cdots \to R(-2)^{\oplus 2} \xrightarrow{ \overline{\Phi}_\lambda^-} R(-1)^{\oplus 2} \xrightarrow{ \overline{\Phi}_\lambda^+} R^{\oplus 2} \xrightarrow{ \overline{\Phi}_\lambda^-} R(1)^{\oplus 2} \to \cdots
\]

\begin{prop} The following holds for all $\lambda \in \mathbb{P}^1$:
	\begin{enumerate}[i)]
	\item We have ${\rm End}_{gr R}(N_{\lambda}) = {\rm \underline{End}}_{gr R}(N_\lambda) = k$, hence the modules $N_\lambda$ are indecomposable.
	\item We have $\tau^2 N_\lambda \cong N_\lambda$.
\end{enumerate}
\end{prop}
\begin{proof}
For part $i)$, direct calculations show that the only scalar matrices $A, B$ fitting in a commutative diagram
\[
\xymatrix@R1pc@C2pc{
R(-1)^{\oplus 2} \ar[d]^A \ar[r]^{\Phi_{\lambda}^{+}} & R^{\oplus 2} \ar[d]^B \\
R(-1)^{\oplus 2} \ar[r]^{\Phi_{\lambda}^{+}} & R^{\oplus 2}
}
\]
are given by $A = B = c \cdot{\rm I}_2$ for some $c \in k$, for any $\lambda \in \mathbb{P}^1$. Thus ${\rm End}_{gr R}(N_{\lambda}) = \underline{{\rm End}}_{gr R}(N_{\lambda}) = k$. Part $ii)$ follows from the definition of $\tau = {\rm syz}_R^1(-)(1)$ and the minimal complete resolution of $N_\lambda$ is $2$-periodic. 
\end{proof}

Next, consider $\lambda = 0, 1, \infty$ corresponding to the $3$ singular conics listed above. Then the factorizations
\begin{align*}
    & Q_0 = x^2 - y^2 = (x-y)(x+y)\\
    & Q_1 = x^2 - z^2 = (x-z)(x+z)\\
    & Q_{\infty} = y^2 - z^2 = (y-z)(y+z)
\end{align*}
are size one matrix factorizations of $Q_{\lambda}$, netting us additional MCM modules\footnote{The choice of sign will become clear in Theorem \ref{D4maincalc}.}
\begin{align*}
    D_0^{+} = {\rm coker}\big( R(-1) \xrightarrow{x-y} R\big) &&  D_0^{-} = {\rm coker}\big( R(-1) \xrightarrow{x+y} R\big)\\
    D_1^{+} = {\rm coker}\big( R(-1) \xrightarrow{x+z} R\big) &&  D_1^{-} = {\rm coker}\big( R(-1) \xrightarrow{x-z} R\big)\\
    D_{\infty}^{+} = {\rm coker}\big( R(-1) \xrightarrow{y+z} R\big) &&  D_{\infty}^{-} = {\rm coker}\big( R(-1) \xrightarrow{y-z} R\big)
\end{align*}

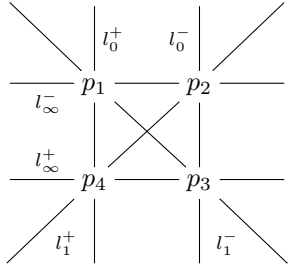
\begin{figure}\label{pairsoflines}
\[
\xymatrix{
	\ar@{-}[dr]&  &  & \ar@{-}[dl]\\
& p_1 \ar@{-}[dr] \ar@{-}[u]_{l_0^+} \ar@{-}[l]^{l_\infty^-} \ar@{-}[d] \ar@{-}[r] & p_2 \ar@{-}[dl] \ar@{-}[u]^{l_0^-} \ar@{-}[d] \ar@{-}[r] & \\	
& p_4 \ar@{-}[dl]^{l_1^+} \ar@{-}[d] \ar@{-}[l]_{l_\infty^+} \ar@{-}[r] & p_3 \ar@{-}[dr]_{l_1^-} \ar@{-}[d] \ar@{-}[r] & \\
&  &  & \\
}
\]
\caption{The three pairs of lines corresponding to singular conics.}
\end{figure} 

We will write $D_\lambda^{\pm} = R/l_\lambda^{\pm}$ for $l_\lambda^{\pm}$ the corresponding linear form\footnote{There is some ambiguity in that $(l_1, l_2, l_3, l_4) = (l_0^+, l_\infty^-, l_0^-, l_\infty^+)$. In practice we will use $\{l_i\}$ solely to refer to $L_i = R/(l_i, l_{i+1})$ and the $\{l_\lambda^{\pm}\}$ to refer to $D_\lambda^{\pm} = R/l_\lambda^{\pm}$. }. Each pair $(l_\lambda^{+}, l_\lambda^{-})$ corresponds to a pair of lines in Figure \ref{pairsoflines} forming the singular conic $V(Q_\lambda)$. We then have minimal complete resolutions
\[
	C(D_\lambda^{\pm}): \cdots \to R(-2)\xrightarrow{l_\lambda^{\mp}} R(-1) \xrightarrow{l_\lambda^{\pm}} R \xrightarrow{l_\lambda^{\mp}} R(1) \xrightarrow{l_\lambda^{\pm}} R(2) \to  \cdots
\]
\begin{prop} The following holds for each $\lambda = 0,1,\infty$:
	\begin{enumerate}[i)]
	\item We have ${\rm End}_{gr R}(D_{\lambda}^{\pm}) = {\rm \underline{End}}_{gr R}(D_\lambda^{\pm}) = k$, hence the modules $D_{\lambda}^{\pm}$ are indecomposable.
\item We have $\tau D_{\lambda}^{\pm} \cong D_{\lambda}^{\mp}$.
\end{enumerate}
\end{prop}
\begin{proof}
We have ${\rm End}_{gr R}(D_{\lambda}^{\pm}) \cong (R/l_\lambda^{\pm})_0 = k$ hence $i)$ follows, and $ii)$ follows from the definition $\tau = {\rm syz}_R^1(-)(1)$.
\end{proof}

We now have enough indecomposable MCM modules to produce all simple regular modules over $kQ$. The following is the main calculation of this section.

\begin{thm}\label{D4maincalc} We have the following isomorphisms in ${\rm D}^b(kQ\ \! \mathsf{mod})$:
	\begin{enumerate}[i)]
\item $F(N_{\lambda}) \cong R_{\lambda}$ for all $\lambda \in \mathbb{P}^1$.
\item $F(D_{\lambda}^{\pm}) \cong S_{\lambda}^{\pm}$ for $\lambda = 0, 1, \infty$.
\end{enumerate}
\end{thm}
\begin{proof}
	For ease of calculations we will use Orlov's fully faithful embedding
	\[
		(-)_{[\geq 0]}: \underline{{\rm MCM}}^{\mathbb{Z}}(R) \hookrightarrow {\rm D}^b(\mathsf{grmod}_{\geq 0}\ \! R).
	\]
	Note that we have $(L_i)_{[\geq 0]} = L_i$, $(N_\lambda)_{[\geq 0]} = N_\lambda$, $(D_\lambda^{\pm})_{[\geq 0]} = D_\lambda^{\pm}$ by Lemma \ref{cutoff} since they are generated in degree zero, and $(k^{st})_{[\geq 0]} = k$ since $a > 0$ by Cor. \ref{a>0} to Orlov's Theorem.
	
We first show $i)$. Let $U = T_{[\geq 0]} = \big(\bigoplus_{i = 1}^4 L_i \big) \oplus k$. We have
\[
\underline{{\rm Hom}}_{gr R}(N_{\lambda}, T) \cong {\rm Hom}_{gr R}(N_{\lambda}, U) \twoheadrightarrow {\rm Hom}_{gr R}(N_{\lambda}, k) \neq 0
\]
since $N_{\lambda}$ is generated in degree $0$. This means that ${\rm H}^0(F(N_\lambda)) = {\rm \underline{Hom}}_{gr R}(N_\lambda, T) \neq 0$ and so ${\rm H}^i(F(N_\lambda)) = 0$ for $i \neq 0$ since $N_{\lambda}$ is indecomposable and $F(N_\lambda)$ is formal in $ {\rm D}^b(kQ\ \! \mathsf{mod}) $. Hence $F(N_\lambda)$ is a $kQ$-module. 

By Lemma \ref{modulestructure} we have to compute the module structure on $\underline{{\rm Hom}}_{gr R}(N_{\lambda}, T) \cong {\rm Hom}_{gr R}(N_{\lambda}, U)$ under endomorphisms of $U$, or equivalently the maps in the diagram
\[
	\xymatrix@R=10pt@C=30pt{ 
		{\rm Hom}_{gr R}(N_\lambda, L_1) \ar[ddr] & \\
		{\rm Hom}_{gr R}(N_\lambda, L_2) \ar[dr] & \\
		& {\rm Hom}_{gr R}(N_\lambda, k) \\
		{\rm Hom}_{gr R}(N_\lambda, L_3) \ar[ur] & \\
		{\rm Hom}_{gr R}(N_\lambda, L_4) \ar[uur] &	}
\]
induced by the canonical quotient $L_i \twoheadrightarrow k$. Let $\underline{\dim} (N_\lambda) := \underline{\dim} (F(N_\lambda)) = (d_0, d_1, d_2, d_3, d_4)$ be the dimension vector of the above quiver representation. We have $d_0 = \dim_k {\rm Hom}_{gr R}(N_\lambda, k) = 2$ with the obvious basis, and we claim that $d_i = \dim_k {\rm Hom}_{gr R}(N_\lambda, L_i) = 1$ for $i = 1,2,3,4$. We prove this by constructing explicit bases. Present $L_i = R/I(p_i) = R/(l_i, l_{i+1})$ as a quotient of two linear forms as above. Consider the morphisms
\[
\xymatrix@C2.8pc{ R(-1)^{\oplus 2} \ar[r]^{\Phi_{\lambda}^{+}} \ar[d]^{\begin{bsmallmatrix} 0 & 1\\ -\lambda_0 & 0 \end{bsmallmatrix}} & R^{\oplus 2} \ar[d]^{\begin{bsmallmatrix} 0 & 1 \end{bsmallmatrix}} \ar@{->>}[r] & N_{\lambda} \ar@{.>}[d] & R(-1)^{\oplus 2} \ar[r]^{\Phi_{\lambda}^{+}} \ar[d]^{\begin{bsmallmatrix} \lambda_1 & 1\\ -\lambda_0 & -1 \end{bsmallmatrix}} & R^{\oplus 2} \ar[d]^{\begin{bsmallmatrix} 1 & 1 \end{bsmallmatrix}} \ar@{->>}[r] & N_{\lambda} \ar@{.>}[d]\\
R(-1)^{\oplus 2} \ar[r]_{\begin{bsmallmatrix} x-y & y-z \end{bsmallmatrix}} & R \ar@{->>}[r] & L_1 & R(-1)^{\oplus 2} \ar[r]_{\begin{bsmallmatrix} x+y & y-z \end{bsmallmatrix}} & R \ar@{->>}[r] & L_2 \\
}
\]
\[
\xymatrix@C2.8pc{ R(-1)^{\oplus 2} \ar[r]^{\Phi_{\lambda}^{+}} \ar[d]^{\begin{bsmallmatrix} \lambda_1 & 1\\ 0 & 0 \end{bsmallmatrix}} & R^{\oplus 2} \ar[d]^{\begin{bsmallmatrix} 1 & 0 \end{bsmallmatrix}} \ar@{->>}[r] & N_{\lambda} \ar@{.>}[d] & R(-1)^{\oplus 2} \ar[r]^{\Phi_{\lambda}^{+}} \ar[d]^{\begin{bsmallmatrix} \lambda_0 \lambda_1 & \lambda_1 \\ \lambda_0 \lambda_1 & \lambda_0 \end{bsmallmatrix}} & R^{\oplus 2} \ar[d]^{\begin{bsmallmatrix} \lambda_0 & \lambda_1 \end{bsmallmatrix}} \ar@{->>}[r] & N_{\lambda} \ar@{.>}[d]\\
R(-1)^{\oplus 2} \ar[r]_{\begin{bsmallmatrix} x+y & y+z \end{bsmallmatrix}} & R \ar@{->>}[r] & L_3 & R(-1)^{\oplus 2} \ar[r]_{\begin{bsmallmatrix} x-y & y+z \end{bsmallmatrix}} & R \ar@{->>}[r] & L_4 \\
}
\]
where we recall that $\Phi_{\lambda}^{+} = \begin{bsmallmatrix} \lambda_1 (x+y) & y+z \\
							\lambda_0(z-y) & x-y
							\end{bsmallmatrix}$.
							These chain-maps cannot be nullhomotopic for degree reasons as there are no non-zero morphism $h:R \to R(-1)$, and so $d_i \geq 1$ for $i = 1,2,3,4$. Since $\tau^2 N_\lambda \cong N_\lambda$, we know that $F(N_\lambda)$ must be a regular indecomposable with $d_0 = 2$, and an appeal to the classification of regular indecomposables then shows that $d_i = 1$ for $i=1,2,3,4$. With the bases constructed above, it is now clear that $F(N_\lambda) = R_\lambda$. 

						We now show $ii)$. By the same argument as above, $F(D_\lambda^{\pm})$ must be a regular indecomposable $kQ$-module and we are left with computing the module structure on ${\rm Hom}_{gr R}(D_{\lambda}^{\pm}, U)$. Let $ \underline{\dim}(D_\lambda^{\pm}) := \underline{\dim}(F(D^{\pm}_\lambda) = (d_0, d_1, d_2, d_3, d_4)$ be its dimension vector. We have $ d_0 = \dim_k {\rm Hom}_{gr R}(D_\lambda^{\pm}, k) = 1$ with generator given by the canonical quotient $D_\lambda^{+} \twoheadrightarrow k$. Writing $D_\lambda^{\pm} = R/l_\lambda^{\pm}$ for $l_\lambda^{\pm}$ the corresponding linear form, we have ${\rm Hom}_{gr R}(R/l_\lambda^{\pm}, L_i) \cong {\rm Hom}_{gr S}(S/l_\lambda^{\pm}, L_i)$ and one calculates its dimension from the incidence relations
\[
d_i = \dim_k {\rm Hom}_{gr S}(S/l_\lambda^{\pm}, L_i) = \begin{cases} 1, & p_i \in V(l_\lambda^{\pm}) \\ 0, & p_i \notin V(l_\lambda^{\pm}) \end{cases} 
\]
which can be read from Figure \ref{pairsoflines}. Since $F(D_\lambda^{\pm})$ is an indecomposable regular module, by comparing dimension vectors we obtain $F(D_{\lambda}^{\pm}) = S_\lambda^{\pm}$ as claimed.

\end{proof}
\begin{rem} We can pick explicit bases for ${\rm Hom}_{gr R}(D^+_\lambda, L_i)$ as below:
\[
	\xymatrix@C0.6pc{ R(-1) \ar@/^0.1pc/[r]^{x-y} \ar[d]^{\begin{bsmallmatrix}  1 \\ 0 \end{bsmallmatrix}} & R \ar[d]^{\begin{bsmallmatrix} 1 \end{bsmallmatrix}} \ar[r] & D_0^+ \ar@{.>}[d] & R(-1) \ar@/^0.1pc/[r]^{x-y} \ar[d]^{\begin{bsmallmatrix}  0 \\ 0  \end{bsmallmatrix}} & R \ar[d]^{\begin{bsmallmatrix} 0 \end{bsmallmatrix}} \ar[r] & D_0^+ \ar@{.>}[d] & R(-1) \ar@/^0.1pc/[r]^{x-y} \ar[d]^{\begin{bsmallmatrix}  0 \\ 0  \end{bsmallmatrix}} & R \ar[d]^{\begin{bsmallmatrix} 0 \end{bsmallmatrix}} \ar[r] & D_0^+ \ar@{.>}[d] & R(-1) \ar@/^0.1pc/[r]^{x-y} \ar[d]^{\begin{bsmallmatrix} 1 \\ 0  \end{bsmallmatrix}} & R \ar[d]^{\begin{bsmallmatrix} 1 \end{bsmallmatrix}} \ar@{->>}[r] & D_0^+ \ar@{.>}[d]\\
R(-1)^{\oplus 2} \ar@/_0.1pc/[r]_{\begin{bsmallmatrix} x-y & y-z \end{bsmallmatrix}} & R \ar[r] & L_1 & R(-1)^{\oplus 2} \ar@/_0.1pc/[r]_{\begin{bsmallmatrix} x+y & y-z \end{bsmallmatrix}} & R \ar[r] & L_2 & R(-1)^{\oplus 2} \ar@/_0.1pc/[r]_{\begin{bsmallmatrix} x+y & y+z \end{bsmallmatrix}} & R \ar[r] & L_3 & R(-1)^{\oplus 2} \ar@/_0.1pc/[r]_{\begin{bsmallmatrix} x-y & y+z \end{bsmallmatrix}} & R \ar[r] & L_4 \\
}
\]
\[
	\xymatrix@C0.6pc{ R(-1) \ar@/^0.1pc/[r]^{x+z} \ar[d]^{\begin{bsmallmatrix}  0 \\ 0 \end{bsmallmatrix}} & R \ar[d]^{\begin{bsmallmatrix} 0 \end{bsmallmatrix}} \ar[r] & D_1^+ \ar@{.>}[d] & R(-1) \ar@/^0.1pc/[r]^{x+z} \ar[d]^{\begin{bsmallmatrix}  1 \\ -1  \end{bsmallmatrix}} & R \ar[d]^{\begin{bsmallmatrix} 1 \end{bsmallmatrix}} \ar[r] & D_1^+ \ar@{.>}[d] & R(-1) \ar@/^0.1pc/[r]^{x+z} \ar[d]^{\begin{bsmallmatrix}  0 \\ 0  \end{bsmallmatrix}} & R \ar[d]^{\begin{bsmallmatrix} 0 \end{bsmallmatrix}} \ar[r] & D_1^+ \ar@{.>}[d] & R(-1) \ar@/^0.1pc/[r]^{x+z} \ar[d]^{\begin{bsmallmatrix} 1 \\ 1 \end{bsmallmatrix}} & R \ar[d]^{\begin{bsmallmatrix} 1 \end{bsmallmatrix}} \ar[r] & D_1^+ \ar@{.>}[d]\\
R(-1)^{\oplus 2} \ar@/_0.1pc/[r]_{\begin{bsmallmatrix} x-y & y-z \end{bsmallmatrix}} & R \ar[r] & L_1 & R(-1)^{\oplus 2} \ar@/_0.1pc/[r]_{\begin{bsmallmatrix} x+y & y-z \end{bsmallmatrix}} & R \ar[r] & L_2 & R(-1)^{\oplus 2} \ar@/_0.1pc/[r]_{\begin{bsmallmatrix} x+y & y+z \end{bsmallmatrix}} & R \ar[r] & L_3 & R(-1)^{\oplus 2} \ar@/_0.1pc/[r]_{\begin{bsmallmatrix} x-y & y+z \end{bsmallmatrix}} & R \ar[r] & L_4 \\
}
\]
\[
	\xymatrix@C0.6pc{ R(-1) \ar@/^0.1pc/[r]^{y+z} \ar[d]^{\begin{bsmallmatrix}  0 \\ 0 \end{bsmallmatrix}} & R \ar[d]^{\begin{bsmallmatrix} 0 \end{bsmallmatrix}} \ar[r] & D_\infty^+ \ar@{.>}[d] & R(-1) \ar@/^0.1pc/[r]^{y+z} \ar[d]^{\begin{bsmallmatrix}  0 \\ 0  \end{bsmallmatrix}} & R \ar[d]^{\begin{bsmallmatrix} 0 \end{bsmallmatrix}} \ar[r] & D_\infty^+ \ar@{.>}[d] & R(-1) \ar@/^0.1pc/[r]^{y+z} \ar[d]^{\begin{bsmallmatrix}  0 \\ 1  \end{bsmallmatrix}} & R \ar[d]^{\begin{bsmallmatrix} 1 \end{bsmallmatrix}} \ar[r] & D_\infty^+ \ar@{.>}[d] & R(-1) \ar@/^0.1pc/[r]^{y+z} \ar[d]^{\begin{bsmallmatrix} 0 \\ 1  \end{bsmallmatrix}} & R \ar[d]^{\begin{bsmallmatrix} 1 \end{bsmallmatrix}} \ar[r] & D_\infty^+ \ar@{.>}[d]\\
R(-1)^{\oplus 2} \ar@/_0.1pc/[r]_{\begin{bsmallmatrix} x-y & y-z \end{bsmallmatrix}} & R \ar[r] & L_1 & R(-1)^{\oplus 2} \ar@/_0.1pc/[r]_{\begin{bsmallmatrix} x+y & y-z \end{bsmallmatrix}} & R \ar[r] & L_2 & R(-1)^{\oplus 2} \ar@/_0.1pc/[r]_{\begin{bsmallmatrix} x+y & y+z \end{bsmallmatrix}} & R \ar[r] & L_3 & R(-1)^{\oplus 2} \ar@/_0.1pc/[r]_{\begin{bsmallmatrix} x-y & y+z \end{bsmallmatrix}} & R \ar[r] & L_4 \\
}
\]
This clearly shows that $F(D^+_\lambda) = S^+_\lambda$, and so that $F(D^-_\lambda) = S^-_\lambda$ by using $\tau$.
\end{rem}
\begin{rem} As a result, we see that the pencil of conics $\{V(Q_\lambda) \}_{\lambda \in \mathbb{P}^1}$ serves as natural parameter space for the family of tubes $\{\T_\lambda\}_{\lambda \in \mathbb{P}^1}$, with the rank of $\T_\lambda$ given by the number of branches of $V(Q_\lambda)$.
\end{rem}

Finally we obtain the remaining indecomposables by taking extensions. Let $\mathcal{S}' = \{N_\lambda\}_{\lambda \neq 0, 1, \infty} \cup \{D_{\lambda}^{\pm}\}_{\lambda = 0,1,\infty}$, which is sent contravariantly onto the set of simple regular modules $\mathcal{S} \subset {\rm D}^b(kQ\ \! \mathsf{mod}) $ by $F$. For any MCM modules $M, N$, we have $ \underline{{\rm Ext}}_{gr R}^1(M, N) = {\rm Ext}^1_{gr R}(M, N)$. We can define MCM modules $N_\lambda \langle r \rangle$, $D_\lambda^{\pm} \langle r \rangle$ for $r \geq 1$ iteratively as the unique modules fitting inside non-trivial short exact sequences
\[
	0 \to N_{\lambda}\langle r \rangle \to N_{\lambda} \langle r+1 \rangle \to N_\lambda \to 0,\quad \quad \lambda \neq 0,1,\infty
\]
\[
	0 \to D^{\mp}_{\lambda}\langle r \rangle \to D^{\pm}_{\lambda} \langle r+1 \rangle \to D^{\pm}_\lambda \to 0,\quad \quad \lambda = 0,1,\infty
\]
corresponding to the short exact sequences
\[
	0 \to R_\lambda \to R_\lambda \langle r+1 \rangle \to R_\lambda \langle r \rangle\to 0,\quad \quad \lambda \neq 0,1,\infty
\]
\[
	0 \to S_\lambda^{\pm} \to S_\lambda^{\pm} \langle r+1 \rangle \to S^{\mp}_\lambda \langle r \rangle\to 0,\quad \quad \lambda = 0,1,\infty
\]
under the isomorphisms 
\[
	{\rm Ext}_{gr R}^1(N_\lambda, N_\lambda\langle r \rangle) \cong \underline{{\rm Ext}}_{gr R}^1(N_\lambda, N_\lambda \langle r \rangle) \cong {\rm Ext}^1_{kQ}(R_\lambda \langle r \rangle, R_\lambda) = k
\]
and respectively for $D_\lambda^{\pm}\langle r \rangle$. Note that under this notation we have $N_\lambda = D_\lambda^+\langle 2 \rangle$ and ${\rm syz}_R^1(N_\lambda)(1) = \tau N_\lambda = D_{\lambda}^-\langle 2 \rangle$ for $\lambda = 0,1,\infty$, while $\tau N_\lambda \cong N_\lambda$ for $\lambda \neq 0,1,\infty$ by Prop. \ref{regsimples}.

A priori, the modules $N_\lambda \langle r \rangle$, $D^{\pm}_\lambda \langle r \rangle$ in ${\rm MCM}^{\Z}(R)$ only have indecomposable images in $ \underline{{\rm MCM}}^{\mathbb{Z}}(R)$, and so might contain a free summand. We show that this is not so.
\begin{lem} The modules $N_\lambda \langle r \rangle$, $D^{\pm}_\lambda \langle r \rangle$ are indecomposable in the module category. 
\end{lem}
\begin{proof}
	Let $M$ be such a module, and we can write $M = F \oplus [M]$ for $F$ a maximal free summand and $[M]$ indecomposable. From the above short exact sequences one sees that $M$ is generated in degree zero by $\beta_{0,0}(M)$ generators, and Betti numbers are subadditive under short exact sequences. We have $\beta_{0,0}([M]) \leq \beta_{0,0}(M)$. We will show that the numbers $\beta_{0,0}([M])$ are additive under the above short exact sequences, thus reversing the inequality $\beta_{0,0}([M]) \geq \beta_{0,0}(M)$ and proving $M = [M]$.

	We have $\beta_{0,0}([M]) = \dim_k {\rm Hom}_{gr R}([M], k)$, and we claim that ${\rm Hom}_{gr R}([M], k) = {\rm \underline{Hom}}_{gr R}([M], k)$. To see this, note that a non-zero map $f:[M] \to k$ factoring through a free module $G$ must surject onto one of its summand to reach $k$, thus splitting off a free summand of $[M]$, a contradiction. We then obtain
\begin{align*}
		\beta_{0,0}([M]) &= \dim_k {\rm Hom}_{gr R}([M], k)\\
		&= \dim_k {\rm \underline{Hom}}_{gr R}([M], k)\\
		&= \dim_k {\rm \underline{Hom}}_{gr R}(M, k)\\
		&= \dim_k {\rm \underline{Hom}}_{gr R}(M, k^{st}) \\
		&= \dim_k {\rm Hom}_{kQ}(P(0), F(M))
	\end{align*}
	and $\dim_k {\rm Hom}_{kQ}(P(0), -)$ is additive on short exact sequence of $kQ$-modules, from which the above extensions come from. This shows $\beta_{0,0}(M) \leq \beta_{0,0}([M])$ and so $M = [M]$.
\end{proof}

This completes the classification of graded MCM $R$-modules. To state the result in full, recall that the complexity of a module $M$ is the least integer $c = cx(M)$ such that the minimal free resolution $F_*$ of $M$ has ranks $rk (F_n)$ with growth of order $O(n^{c-1})$. 

Let $\CC$ be the Auslander-Reiten component of $\Gamma( \underline{{\rm MCM}}^{\mathbb{Z}}(R))$ containing $L_1, \dots, L_4, k^{st}$, and let $\mathcal{Q}_\lambda$ be the component contaning $N_\lambda$ for $\lambda \in \mathbb{P}^1$. All indecomposables in $\mathcal{Q}_{\lambda}$ have periodic minimal free resolution period $1$ or $2$ according to whether $\lambda \neq 0, 1, \infty$ or $\lambda = 0, 1, \infty$, corresponding to the $\tau$-period. Hence these modules have complexity one. In particular they satisfy $M[2] \cong M(2)$. In the next section, we will construct the minimal complete resolutions of $L_1, \dots, L_4, k^{st}$ which will be of complexity two. 

Summarising this section, we have shown:

\begin{thm} The indecomposable (non-free) graded MCM $R$-modules are listed in the following table, up to degree shift:
\[
\centerline{
\begin{tabular}{|c|c|c|}
  \hline
   \multicolumn{3}{|c|}{Indecomposable objects (up to degree shift)}\\
   \hline
 Complexity $2$ & \multicolumn{2}{c|}{Complexity $1$}\\
 \hline
 & $\lambda \in \mathbb{P}^1 \setminus \{0, 1, \infty \}$ & $\lambda = 0, 1, \infty$ \\
 \cline{2-3}
 & & \\
 $k^{st}[n]$, $n \in \Z$ & $N_{\lambda}\langle r \rangle$, $r \geq 1$ & $D_{\lambda}\langle r \rangle^{+}$, $r \geq 1$ \\
 $L_i[n]$, $n \in \Z$ &  & $D_{\lambda}\langle r \rangle^{-}$, $r \geq 1$ \\
& & \\
 \hline
\end{tabular}}
\]
\end{thm}

\begin{thm}
The Auslander-Reiten quiver of $\underline{{\rm MCM}}^\Z(R)$ is given by the union of disjoint components
\[
\Gamma\big(\underline{{\rm MCM}}^\Z(R)\big) = \bigg(\bigcup_{n \in \Z} \mathcal{C}[n] \bigg) \cup \bigg(\bigcup_{\lambda \in \mathbb{P}^1, n \in \Z} \mathcal{Q}_{\lambda}[n]\bigg).
\]
The components $\CC$ can be drawn as

{\tiny{
\[
\xymatrix@R.8pc@C1.6pc{
	\cdots\ar[ddr] &&  L_1(1)[-1] \ar[ddr] && L_1 \ar[ddr] && L_1(-1)[1] \ar[ddr] \\
	\cdots \ar[dr] && L_2(1)[-1] \ar[dr] && L_2 \ar[dr] && L_2(-1)[1]  \ar[dr] \\
	& k^{st}(2)[-2] \ar[uur] \ar[ur] \ar[dr] \ar[ddr] && k^{st}(1)[-1] \ar[uur] \ar[ur] \ar[dr] \ar[ddr] && k^{st}\ar[uur] \ar[ur] \ar[dr] \ar[ddr] && \cdots  \\
\cdots \ar[ur]&& L_3(1)[-1] \ar[ur] && L_3 \ar[ur] && L_3(-1)[1] \ar[ur]\\
\cdots \ar[uur] && L_4(1)[-1] \ar[uur] && L_4 \ar[uur] && L_4(-1)[1] \ar[uur] \\
}
\]
}}
and the components $\mathcal{Q}_\lambda$ as below, with the edges identified
\begin{align*}
	\xymatrix@R.8pc@C.2pc{\vdots \ar[dr] \ar@{-}[dd] && \vdots \ar@{-}[dd] &&  \vdots \ar[dr] \ar@{-}[dd] && \vdots \ar[dr] && \vdots \ar@{-}[dd] \\
 & N_\lambda \langle 6 \rangle \ar[ur] \ar[dr] & && & D_\lambda^-\langle 6 \rangle \ar[dr] \ar[ur] \ar@{.}[rr] && D_\lambda^+\langle 6 \rangle \ar[dr] \ar[ur] & \\
 N_\lambda\langle 5 \rangle \ar[dr] \ar[ur] \ar@{.}[rr] \ar@{-}[dd] && N_\lambda \langle 5 \rangle \ar@{-}[dd] && D_\lambda^+\langle 5 \rangle \ar[dr] \ar[ur] \ar@{-}[dd] \ar@{.}[rr] && D_\lambda^-\langle 5 \rangle \ar[dr] \ar[ur] \ar@{.}[rr] && D_\lambda^+\langle 5 \rangle \ar@{-}[dd] \\
 & N_\lambda\langle 4 \rangle \ar[dr] \ar[ur] & && & D_\lambda^+\langle 4 \rangle \ar[dr] \ar@{.}[rr] \ar[ur] && D_\lambda^-\langle 4 \rangle \ar[dr] \ar[ur] \\
 N_\lambda\langle 3 \rangle \ar[dr] \ar[ur] \ar@{.}[rr] \ar@{-}[dd] && N_\lambda\langle 3 \rangle \ar@{-}[dd] && D_\lambda^-\langle 3 \rangle \ar@{-}[dd] \ar@{.}[rr] \ar[dr] \ar[ur] && D_\lambda^+\langle 3 \rangle \ar[dr] \ar[ur] \ar@{.}[rr] && D_\lambda^-\langle 3 \rangle \ar@{-}[dd] \\
 & N_\lambda \langle 2 \rangle \ar[dr] \ar[ur] & && & D_\lambda^{-}\langle 2 \rangle \ar[dr] \ar[ur] \ar@{.}[rr] && D_{\lambda}^{+}\langle 2 \rangle \ar[dr] \ar[ur] \\
 N_\lambda \ar[ur] \ar@{.}[rr] && N_\lambda && D_\lambda^{+} \ar[ur] \ar@{.}[rr] && D_\lambda^{-} \ar[ur] \ar@{.}[rr] && D_\lambda^{+}\\
 & \lambda \in \mathbb{P}^1 \setminus \{0, 1, \infty \} & && && \lambda = 0, 1, \infty. &&
 }
\end{align*}
\end{thm}

\subsection{Betti tables of indecomposables}
Let $M$ be a graded MCM $R$-module with minimal complete free resolution
\[
\cdots \to \bigoplus_{j \in \Z} R(-j)^{\oplus \beta_{i, j}} \to \cdots \to \bigoplus_{j \in \Z} R(-j)^{\oplus \beta_{1, j}} \to \bigoplus_{j \in \Z} R(-j)^{\oplus \beta_{0, j}} \to \bigoplus_{j \in \Z} R(-j)^{\oplus \beta_{-1, j}} \to \cdots.
\]
The Betti table of $M$ is the table whose entry in the $i$-th column and $j$-th row is given by $\beta_{i, i+j}$.
\[
		\begin{tabular}{c|c c c c c c c c}
		& $\cdots$ & 0 & 1 & $\cdots$ & a-2 & a-1 & a & $\cdots$ \\
		\hline
		\vdots & & \vdots & \vdots & & \vdots & \vdots & \vdots\\ 
		-2 & $\cdots$ &$\beta_{0, -2}$&$\beta_{1,-1}$ & $\cdots$ & $\beta_{a-2, a-4}$ & $\beta_{a-1, a-3}$ & $\beta_{a, a-2}$ & $\cdots$ \\ 
		-1 &$\cdots$ & $\beta_{0, -1}$&$\beta_{1,0}$ & $\cdots$ & $\beta_{a-2, a-3}$ & $\beta_{a-1, a-2}$ & $\beta_{a, a-1}$ & $\cdots$ \\ 
		0 &$\cdots$ & $\beta_{0, 0}$&$\beta_{1,1}$ & $\cdots$ & $\beta_{a-2, a-2}$ & $\beta_{a-1, a-1}$ & $\beta_{a, a}$ & $\cdots$ \\ 
		1 &$\cdots$ & $\beta_{0, 1}$&$\beta_{1,2}$ & $\cdots$ & $\beta_{a-2, a-1}$ & $\beta_{a-1, a}$ & $\beta_{a, a+1}$ & $\cdots$ \\ 
		\vdots & & \vdots & \vdots & & \vdots & \vdots & \vdots\\ 
		\end{tabular}
\]
We list here all Betti tables of indecomposables. Note that $\tau^{\pm1} = {\rm syz}^{\pm 1}_R(-)(\pm1)$ and $(\pm1)$ correspond to horizontal and vertical shifts respectively.

It is simple to obtain the Betti tables from our classifications. Let $F: \underline{{\rm MCM}}^{\mathbb{Z}}(R)^{op} \xrightarrow{\cong} {\rm D}^b(kQ\ \! \modsf)$ be the previous equivalence. For any quiver representation $X$ of $Q$, we write $X_0$ for the vector space sitting over the $0$ vertex. We will use the following lemma for calculations. 
\begin{lem} For any $M \in \underline{{\rm MCM}}^{\mathbb{Z}}(R)$, we have
	\[
		\beta_{i,i+j}(M) = \dim_k {\rm H}^{-j}(\tau^{-i-j}F(M))_0. 
	\]
\end{lem}
\begin{proof}
	We have $\beta_{i,i+j}(M) = \dim_k \underline{{\rm Ext}}_{gr R}^i(M, k(-i-j)) = \dim_k \underline{{\rm Ext}}_{gr R}^{i}(M, k^{st}(-i-j))$. We also have $\tau M = M(1)[-1]$, and so we obtain
	\begin{align*}
		\beta_{i,i+j}(M) &= \dim_k \underline{{\rm Ext}}^{i}_{gr R}(M, k^{st}(-i-j)) \\
		&= \dim_k \underline{{\rm Hom}}_{gr R}(M, k^{st}(-i-j)[i]) \\
		&= \dim_k \underline{{\rm Hom}}_{gr R}(M, \tau^{-i-j}k^{st}[-j]) \\ 
		&= \dim_k \underline{{\rm Hom}}_{gr R}(\tau^{i+j}M, k^{st}[-j]) \\ 
		&= \dim_k {\rm Hom}_{ {\rm D}^b(kQ^{op}) }(F(k^{st}[-j]), F(\tau^{i+j}M))\\
		&= \dim_k {\rm Hom}_{ {\rm D}^b(kQ^{op}) }(P(0)[j], \tau^{-i-j}F(M))\\
		&= \dim_k {\rm H}^{-j}(\tau^{-i-j}F(M))_0.
	\end{align*}
\end{proof} 

It follows that the Betti table of any indecomposable $M$ can be computed from the dimension vectors of the $\tau$-orbit of quiver representation $\{ \tau^n F(M) \}_{n \in \Z}$. We deduce that, if $X \in kQ\ \modsf \subseteq {\rm D}^b(kQ \ \! \modsf) $ is a $kQ$-module such that $\tau^{-n}X$ is also a $kQ$-module for all $n \geq 0$, then any MCM $R$-module corresponding to $X$ has a linear resolution in that $\beta_{i,i+j}(M) = 0$ for $j \neq 0$ for all $i \geq 0$. Moreover, the regular $kQ$-modules correspond precisely to the MCM modules which are completely linear, meaning that $\beta_{i, i+j}(M) = 0$ for any $j \neq 0$ and all $i \in \Z$. \bigskip

The dimension vectors of indecomposable representations of $Q$ are written down in \cite[XIII.3]{SimSkow}, and from this it is easy to obtain the corresponding Betti tables of indecomposables MCM modules. We record this in the next proposition.

\begin{prop}\label{betti} The Betti tables of indecomposable graded MCM modules are given up to syzygy and degree shifts by:
\begin{align*}
N_{t}\langle r \rangle: && D_{t}\langle r \rangle^{\pm}: \\
& \begin{tabular}{ c | c  c  c  c  c  c }
 & \ldots & -1 & 0 & 1 & 2 & \ldots \\
\hline
 -1 & - & - & - & - & - & - \\
  0 & \ldots & 2r & 2r & 2r & 2r & \ldots \\
  1 & - & - & - & -  & - & -\\
\end{tabular}
&& \begin{tabular}{ c | c  c  c  c  c  c }
 & \ldots & -1 & 0 & 1 & 2 & \ldots \\
\hline
 -1 & - & - & - & - & - & - \\
  0 & \ldots & r & r & r & r & \ldots \\
  1 & - & - & - & -  & - & -\\
\end{tabular}
\end{align*}

\begin{align*}
 k^{st}: &&  L_i: & \\
 & \begin{tabular}{ c | c  c  c  c  c  c  c  c }
 & \ldots & -2 & -1 & 0 & 1 & 2 & 3 & \ldots \\
\hline
 --2 & - & - & - & - & - & - & - & -\\
  -1 & \ldots & 5 & 3 & 1 & - & - & - & - \\
 0 & - & - & - & 1 & 3 & 5 & 7 & \ldots \\
  1 & - & - & - & - & -  & - & - & - \\
\end{tabular}
&& \begin{tabular}{ c | c  c  c  c  c  c  c  c }
 & \ldots & -2 & -1 & 0 & 1 & 2 & 3 & \ldots \\
\hline
 -2 & - & - & - & - & - & - & - & -\\
 -1 & \ldots & 2 & 1 & - & - & - & - & - \\
  0 & - & - & - & 1 & 2 & 3 & 4 & \ldots \\
  1 & - & - & - & - & -  & - & - & - \\
\end{tabular}
\end{align*}
\qed
\end{prop}
\newpage

%% file: chapters/classifications/elliptic.tex
\section{Graded MCM modules over the cone of $4$ points on $\mathbb{P}^1$}
The algebra $R_{Y_2}$ had tame domestic CM-type and almost all of the isomorphism classes of indecomposable MCM modules were exhausted by a $1$-parameter family of tubes $\{ \T_\lambda \}_{\lambda \in \mathbb{P}^1}$. In contrast, the algebra $R_{Y_1}$ studied in this section is of tame but non-domestic CM-type, and we will see that all indecomposables live in $1$-parameter familes of tubes 
\[
\{ \T_{\lambda, q} \}_{\lambda \in \mathbb{P}^1} \quad 	q \in \Q \cup \{\infty\}
\]
indexed by a choice of slope $q$ for an appropriate stability condition. We have shown the existence of an equivalence of triangulated categories 
\[
	\underline{{\rm MCM}}^{\mathbb{Z}}(R_{Y_1}) \cong {\rm D}^b(Sq(2,2,2,2;\lambda))
\]
onto the derived category of the squid algebra ${\rm Sq}(2,2,2,2;\lambda)$ with quiver
{\tiny{
		\[
	\xymatrix@R=5pt@C=20pt{&&\bullet\\
		&&\bullet\\
		\bullet \ar@<.5ex>[r]^y \ar@<-.5ex>[r]_x &\bullet \ar[uur]^{p_1} \ar[ur]_{\ p_2} \ar[dr]^{\ p_3} \ar[ddr]_{p_4}& \\
		&&\bullet\\
		&&\bullet
	}
\]
}}
and relations $p_i l_i(x,y) = 0$ for $i = 1,2,3,4$. Squid algebras arise as endomorphism algebras of certain tilting sheaves on Geigle-Lenzing weighted projective line of the corresponding weight type, in this case giving 
\[
	{\rm D}^b(Sq(2,2,2,2;\lambda)) \cong {\rm D}^b(\X) 
\]
for $\X = \mathbb{P}^1(2,2,2,2;\lambda)$ a weighted projective line of genus one. The above description of the set of indecomposables was obtained by Geigle, Lenzing and Meltzer and parallels the Atiyah classification of vector bundles on the elliptic curve. The main aspect of the classification is the existence of two `twist' autoequivalences
\[
	T_1, T_2 : {\rm D}^b(\X) \xrightarrow{\cong} {\rm D}^b(\X) 
\]
generating a braid group on three strands $\mathcal{B}_3 = \langle T_1, T_2 \rangle$, from which one obtains all indecomposables (up to suspension) by successive applications starting from the category of torsion sheaves $\mathsf{coh}_0\ \! \X$.

Let $R = R_{Y_1} = S/(f_\lambda)$ for $S = k[x,y]$ and $f_\lambda = xy(x-y)(x-\lambda y)$, $\lambda \neq 0,1$. The aim of this section is to understand the parallel classification in the category $ \underline{{\rm MCM}}^{\mathbb{Z}}(R)$, or equivalently in the homotopy category of matrix factorisations ${\rm \underline{MF}}(S, f_\lambda) \cong \underline{{\rm MCM}}^{\mathbb{Z}}(R)$. In this chapter, we will do the following things:

\begin{enumerate}[1)]
\item We will write down the matrix factorisations corresponding to the simple torsion sheaves, from which the remaining indecomposables can be produced by taking extensions and applying $\B_3$.
\item We will give formulas for the rank and degree of a (complex of) sheaves $\F_M \in {\rm D}^b(\X) \cong \underline{{\rm MCM}}^{\mathbb{Z}}(R)$ in terms of the Betti table of the corresponding MCM module $M$, allowing an intrinsic description of the `charge' of a graded MCM module $\ZZ : K_0( \underline{{\rm MCM}}^{\mathbb{Z}}(R)) \to \Z^{\oplus 2}$ defined as
\[
\ZZ(M) = \begin{pmatrix} rk\ \! M \\ deg\ \! M \end{pmatrix} := \begin{pmatrix} rk\ \! \F_M \\ deg\ \! \F_M \end{pmatrix}
\]
Along with the t-structure giving rise to $\coh \X \subset \underline{{\rm MCM}}^{\mathbb{Z}}(R)$, this gives the data of a stability condition on $ \underline{{\rm MCM}}^{\mathbb{Z}}(R)$ in the sense of Bridgeland. 
\item We will write down the action of the operations $M \mapsto M(1)$ on the charge $\ZZ(M)$, giving rise to an action of the cyclic group $C_4$ of order $4$ on the lattice $\Z^{\oplus 2}$. We then give a fundamental domain for the action of $C_4$, so that it suffices to consider indecomposables with $\ZZ(M)$ in that fundamental domain. 
\item Finally as main result, making use of the above fundamental domain, we will completely classify the Betti tables of indecomposable graded MCM modules $M$ up to degree shift and syzygy. 
	\end{enumerate}

	\subsubsection*{Relations to previous work}
The parabolic surface singularities 
		\[
	\begin{tabular}{|c|c|c|}
\hline
$P_8$ & $X_9$ & $J_{10}$ \\
\hline
$x^3 + y^3 + z^3 + axyz$ & $x^4 + y^4 + ax^2y^2 + z^2$ & $x^3 + y^6 + ax^2y^2 + z^2$ \\
$a^3 + 27 \neq 0$ & $a^2 \neq 4$ & $4a^3 + 27 \neq 0$ \\
\hline
	\end{tabular}
\]
give the cone over the embedding of an elliptic curve $E_a$ inside $\mathbb{P}^2$, $\mathbb{P}(1,1,2)$ and $\mathbb{P}(1, 2, 3)$ respectively. In his thesis, A. Pavlov classified the Betti tables of graded MCM modules over a hypersurface ring $A = k[x,y,z]/(f)$ for $f$ in the above table by making use of Orlov's equivalence
\[
	\underline{{\rm MCM}}^{\mathbb{Z}}(A) \cong {\rm D}^b(E_a).
\]
to reduce calculations of Betti tables to sheaf cohomology calculations. Our methods are directly inspired from his, although the presence of exceptional objects (sitting in tubes of rank two) adds a layer of complexity for which further ideas are required. \bigskip

Lastly, the above picture suggests that similar results hold for the curve singularity of type $T_{36}$. This requires more involved (but similar) calculations, and we will not go through this here. 
\subsubsection*{Weighted projective lines} 
We now recall standard notation, definitions and background results which will be used implicitly throughout this section. The reader is referred to \cite{GL, LM, MeltBook} for a more in-depth view of the topic. Note that we will only use the weighted projective line of genus one $\X = \mathbb{P}^1({\bf p}, \pmb{\lambda}) = \mathbb{P}^1(2,2,2,2; \lambda)$ with $\pmb{\lambda} = (0,\infty,1, \lambda)$. The derived categories of weighted projective lines of genus one were thoroughly investigated by Lenzing and Meltzer in \cite{LM, MeltBook}.\bigskip

Let us fix notation. Given a set of weights ${ \bf p} = (p_1, \dots, p_n)$ and points $\pmb{\lambda} = (\lambda_1, \dots, \lambda_n)$ on $\mathbb{P}^1$, a weighted projective line $\X = \mathbb{P}^1({\bf p}, \pmb{\lambda})$ is constructed from $\mathbb{P}^1$ from an appropriate `root construction' of order $p_i$ at the point $\lambda_i$, with the resulting geometric object living in one of many categories according to one's taste (as a Deligne-Mumford stack, as an orbifold $\mathbb{P}^1$, as a `noncommutative' projective variety, \dots ). The resulting category of coherent sheaves $\coh \X$ being the main object of interest, we will use the original (and most tractable) model for it introduced by Geigle and Lenzing \cite{GL}. 

Given a set $( {\bf p}, \pmb{\lambda})$ of weighted points on $\mathbb{P}^1$, we can and will assume that $p_i \geq 2$ and that $\pmb{\lambda} = (0,\infty,1,  \lambda_4, \dots \lambda_n)$. Let $S = k[u,v]$ be the homogeneous coordinate ring of $\mathbb{P}^1$. Introduce the 'homogeneous coordinate ring' of $\X$ by adding $p_i$-roots at the points $\lambda_i$ 
\[
	S({\bf p}, \pmb{\lambda}) = \frac{S[x_1, \dots, x_n]}{\left(l_i(u,v) - x_i^{p_i} \right)_{i = 1, \dots, n}}
\]
where $l_i \in k[u,v]$ is the linear form cutting out $\lambda_i \in \mathbb{P}^1$. Since we assume that $\lambda_1 = [0:1]$ and $\lambda_2 = [1:0]$ the first two relations become $u = x_1^{p_1}, v = x_2^{p_2}$, hence it is customary to write
\[
	S({\bf p}, \pmb{\lambda}) = \frac{k[x_1, \dots, x_n]}{\left( l_i(x_1^{p_1}, x_2^{p_2}) - x_i^{p_i} \right)_{i=3, \dots, n}}.
\]
To the set of weights ${\bf p}$, one associates a rank one abelian group $\mathbb{L} = \mathbb{L}({\bf p})$ with presentation
\[
	\mathbb{L}({\bf p}) = \langle \vec{x}_1, \dots, \vec{x}_n, \vec{c} \ | \ p_1 \vec{x}_1 = \dots = p_n \vec{x}_n = \vec{c} \rangle.
\]
Setting $p = {\rm lcm}(p_1, \dots, p_n)$, there is a group homomorphism $\delta: \mathbb{L} \to \Z$ given by $\delta(\vec{x}_i) = \frac{p}{p_i}$, with finite kernel. We denote by $\vec{\omega} = (n-2)\vec{c} - \sum_{i=1}^n \vec{x}_i$ the canonical element in $\mathbb{L}$. The coordinate ring $S({\bf p}, \pmb{\lambda})$ admits a grading by the group $\mathbb{L}$ by setting $|x_i| := \vec{x}_i$. Taking a cue from the Serre's Theorem, we define the abelian category of coherent sheaves $\coh \X$ as the Serre quotient
\[
	\coh \X = \frac{\grmod^{\mathbb{L}} S(\bf{p}, \pmb{\lambda})}{\mathsf{grmod}_0^{\mathbb{L}} S(\bf{p}, \pmb{\lambda})}
\]
of the finitely generated graded $S(\bf{p}, \pmb{\lambda})$-modules by the subcategory of finite length modules, with $\qcoh \X$ defined similarly. Alternatively, these categories have models as actual (quasi-)coherent sheaves on a ringed spaced, see \cite{GL}. \bigskip

Writing $B = S({\bf p}, \pmb{\lambda})$ to alleviate notation, each $\vec{x} \in \mathbb{L}$ gives rise to a corresponding line bundle $\OO_\X(\vec{x}) = \widetilde{B(\vec{x})}$, and we denote by
\[
	\widetilde{M} \otimes \OO_\X(\vec{x}) := \widetilde{M(\vec{x})}
\]
the twisting operator on sheaves. Note that by \cite[Appendix]{LenzingKTheory}, $\coh \X$ has a symmetric closed monoidal structure with unit $\OO_\X$ for which the above line bundles are the invertible objects, compatible with sheafification, so that $\vec{x} \mapsto \OO_\X(\vec{x})$ gives an isomorphism of abelian groups $\mathbb{L} \cong {\rm Pic}(\X)$. We denote $\omega_\X = \OO_\X(\vec{\omega})$ the canonical line bundle. 

The categories $\coh \X$ and $\qcoh \X$ share the same formal properties as those of smooth projective curves. In particular $\coh \X$ is an Ext-finite hereditary category with Serre duality 
\[
	{\rm Ext}^i(\F, \G \otimes \omega_\X) \cong {\rm D}{\rm Ext}^{1-i}(\G, \F)
\]
and we set $\Se_\X(-) = -\otimes \omega_\X [1]$ the Serre functor on $ {\rm D}^b(\X) $ with $\tau = - \otimes \omega_\X$ the Auslander-Reiten translate. 

We can compute the cohomology of line bundles as ${\rm H}^0(\X, \OO_\X(\vec{x})) \cong B_{\vec{x}}$, and ${\rm H}^1(\X, \OO_\X(\vec{x}))$ by Serre duality. Every coherent sheaf on $\X$ is the direct sum of its torsion subsheaf and quotient torsion-free sheaf, which is then a vector bundle, and vector bundles admit finite filtrations with line bundle successive quotients. One is lead to understand torsion sheaves and vector bundles separately. 

The indecomposable torsion sheaves are supported over a single point $x \in \mathbb{P}^1$. We say that $x$ is \emph{ordinary} if it lies outside of the set $\pmb{\lambda}$, and \emph{exceptional} otherwise. Torsion sheaves supported over $x$ form a serial abelian subcategory, with unique simple sheaf over $x$ ordinary and $p_i$-many simple sheaves $\{S_{i, j} \}_{j \in \Z/p_i\Z}$ over $x = \lambda_i$ exceptional. These have presentations 
\[
	0 \to \OO_\X((j-1)\vec{x}_i) \xrightarrow{x_i} \OO_\X(j\vec{x}_i) \to S_{i, j} \to 0. 
\]
In particular we single out $S_{i, 0}$ as the unique simple sheaf with a non-zero section\footnote{This agrees with the notation in \cite{CB} but disagrees with \cite{MeltBook}.}, and we have ${\rm Hom}(\OO_\X, S_{i, 0}) = k$ and $S_{i, j} \otimes \omega_\X = S_{i, j+1}$. There is a family of indecomposable ``ordinary'' torsion sheaves $S_x$ for any $x$, with presentations
\[
0 \to \OO_\X(-\vec{c}) \xrightarrow{l(u,v)} \OO_\X \to S_x \to 0
\]
where $l(u,v) \in S = k[u,v] = k[x_1^{p_1}, x_2^{p_2}] \subseteq B$ is the linear form cutting down the point $x \in \mathbb{P}^1$. The sheaf $S_x$ has length one when $x$ is ordinary and length $p_i$ over $x = \lambda_i$, with the $S_{i, j}$ as simple composition factors. The Auslander-Reiten quiver of the subcategory of torsion sheaves $\mathsf{coh}_0\ \! \X$ then forms a $\mathbb{P}^1$ family of tubes $\{ \T_\lambda \}_{\lambda \in \mathbb{P}^1}$, of rank $p_i$ over $\lambda_i$ and rank one elsewhere. 

In contrast, the classification of indecomposable vector bundles on $\X$ differs greatly in complexity according to whether the \emph{virtual genus} $g_\X \in \Q$ satisfies $g_\X < 1$, $g_\X = 1$ or $g_\X > 1$ as in the case of algebraic curves, and a complete classification is only attainable for $g_\X \leq 1$.  See \cite{GL, MeltBook} for more details. 

The rank and degree of a sheaf define maps on $K_0(\X) := K_0(\coh \X)$, uniquely determined by additivity from their value on line bundles as
\begin{align*}
	& rk(\OO_\X(\vec{x})) = 1 \\
	& deg(\OO_\X(\vec{x})) = \delta(\vec{x}). 
\end{align*}
In particular $deg(S_x) = deg(\OO_\X(\vec{c})) = p$ and $deg(S_{i,j}) = 1$. 

We make note of two particular tilting sheaves in $ {\rm D}^b(\X) $. To fix notation, note that $\mathbb{L}$ is an ordered abelian group with positive cone $\mathbb{L}_+ = \N \cdot \{\vec{x}_1, \dots, \vec{x}_n \}$, and let $S_{i,j}\langle r \rangle$ be the unique indecomposable torsion sheaf of length $r$ supported over $\lambda_i$ with simple socle $S_{i,j}$. Then we have a tilting bundle 
\[
	T_{can} = \bigoplus\limits_{\vec{0} \leq \vec{x} \leq \vec{c}} \OO_\X(\vec{x})
\]
whose endomorphism algebra is given by Ringel's canonical algebra $C(\bf{p}, \pmb{\lambda})$ $= kQ/I$ with quiver
\[
	\xymatrix@R=10pt@C=30pt{ & \OO_\X(\vec{x_1}) \ar[r]^-{x_1} & \OO_\X(2\vec{x_1}) \ar[r]^-{x_1} & \cdots \ar[r]^-{x_1} & \OO_\X((p_1-1)\vec{x_1}) \ar[rdd]^-{x_1} & \\
		& \OO_\X(\vec{x_2}) \ar[r]^-{x_2} & \OO_\X(2\vec{x_2}) \ar[r]^-{x_2} & \cdots \ar[r]^-{x_2} & \OO_\X((p_2-1)\vec{x_2}) \ar[rd]^-{x_2} & \\
		\OO_\X \ar[ruu]^-{x_1} \ar[ru]^-{\quad x_2} \ar[rd]^-{x_n} &\vdots & \vdots & \vdots & \vdots & \OO_\X(\vec{c}) \\
		& \OO_\X(\vec{x}_n) \ar[r]^-{x_n} & \OO_\X(2\vec{x}_n) \ar[r]^-{x_n} & \cdots \ar[r]^-{x_n} & \OO_\X((p_n -1)\vec{x}_n) \ar[ru]^-{x_n} & \\
	}
\]
and relations $l_i(x_1^{p_1}, x_2^{p_2}) - x_i^{p_i} = 0$, as well as a tilting sheaf \cite[Ex. 4.4]{LM2} 
\[
	T_{squid} = \OO_\X \oplus \OO_\X(\vec{c}) \oplus \left( \bigoplus\limits_{i=1}^n \bigoplus\limits_{r = 1}^{p_i-1} S_{i,0}\langle r \rangle \right)
\]
with endomorphism algebra a Squid algebra $Sq(\bf{p}, \pmb{\lambda})$ $= kQ/I$ with quiver
\[
	\xymatrix@R=10pt@C=30pt{ 
		&& S_{1,0} \ar[r] & S_{1,0}\langle 2 \rangle \ar[r] & \cdots \ar[r] & S_{1,0}\langle p_1 - 1\rangle \\
		&& S_{2,0} \ar[r] & S_{2,0}\langle 2 \rangle \ar[r] & \cdots \ar[r] & S_{2,0}\langle p_2 - 1\rangle \\
		\OO_\X \ar@<.5ex>[r]^-{u} \ar@<-.5ex>[r]_-{v} & \OO_\X(\vec{c}) \ar[ruu]^-{p_1} \ar[ru]_-{p_2} \ar[rd]^-{p_{n-1}}\ar[rdd]_-{p_n} &\vdots & \vdots & \vdots & \vdots  \\
		&& S_{n-1,0} \ar[r] & S_{n-1,0}\langle 2 \rangle \ar[r] & \cdots \ar[r] & S_{n-1,0}\langle p_{n-1} - 1\rangle \\
		&& S_{n,0} \ar[r] & S_{n,0}\langle 2 \rangle \ar[r] & \cdots \ar[r] & S_{n,0}\langle p_n - 1\rangle \\
	}
\] 
and relations $p_i l_i(u,v) = 0$. Note that the indecomposable summands of $T_{can}$ and $T_{squid}$ naturally form full strong exceptional collections since the above quivers have no cycles. Since the head of the Squid algebra is a Kronecker quiver, we can interpret that second full exceptional collection as
\[
	{\rm D}^b(\X) = \left\langle {\rm D}^b(\mathbb{P}^1), S_{1,0}, \cdots, S_{1, 0}\langle p_1 - 1 \rangle, \cdots, S_{n,0}, \cdots, S_{n,0} \langle p_n - 1 \rangle \right\rangle.
\]
Moreover, this admissible embedding $ {\rm D}^b(\mathbb{P}^1) \hookrightarrow {\rm D}^b(\X) $ sends
\begin{align*}
	\OO_{\mathbb{P}^1}(n) &\mapsto \OO_\X(n\vec{c})\\
 	k(x) &\mapsto S_x
\end{align*}
for $n \in \Z$ and any $x \in \mathbb{P}^1$. In this way, the category $\coh \X$ can be thought of as an enlargement of $\coh \mathbb{P}^1$, see \cite{LenzingKTheory}.

\subsubsection*{Setup}
We can now go ahead with the results of this section. Thinking ahead, we will normalise our calculations by replacing the tilting object $T = \m(1) \oplus \m^2(2) \oplus \left( \bigoplus_{i = 1}^4 L_i \right)$ by 
\[
	U = T(-3)[1] = k^{st}(-2) \oplus (R/\m^2)^{st}(-1) \oplus \left( \bigoplus_{i=1}^4 L_i(-3)[1] \right)
\]
with same endomorphism algebra ${\rm \underline{End}}_{gr R}(U) \cong Sq(2,2,2,2;\lambda)$. From previous results, we obtain:
\begin{cor}\label{wpltilt} For $\X = \mathbb{P}^1(2,2,2,2; \lambda)$, we have equivalences of triangulated categories
	\[
		\underline{\rm MCM}^\Z(R) \cong {\rm D}^b(Sq(2,2,2,2;\lambda)) \cong {\rm D}^b(\X).
	\]
The composed equivalence sends the full strong exceptional collection 
	\[
		\left(k^{st}(-2), (R/\mathfrak{m}^2)^{st}(-1), L_1(-3)[1], \dots, L_4(-3)[1]\right)
	\]
	to 
	\[
		\left( \OO_\X, \OO_\X(\vec{c}), S_{1,0}, \dots, S_{4, 0} \right).
	\]
\end{cor}
Next, we review the structure of $ {\rm D}^b(\X) $ for $\X = \mathbb{P}^1(2,2,2,2;\lambda)$.

\subsection{Weighted projective lines of genus one and braid group actions}

In this subsection we review the classification of indecomposable coherent sheaves over a weighted projective line of genus one due to Lenzing and Meltzer \cite{LM}, which closely mirrors Atiyah's classification of sheaves on an elliptic curve. Everything here is due to them, and we follow \cite{LM, MeltBook}. Let
\[
	\ZZ: K_0(\X) \to \Z^{\oplus 2}
\]
be the `charge'\footnote{This follows Bridgeland's terminology in \cite{BridgelandStab}.}, sending the class $[\F]$ of a coherent sheaf to 
\[
	\ZZ(\F) = \begin{pmatrix} rk(\F) \\ deg(\F) \end{pmatrix}
\]
Let $\mu(\F) = \frac{deg(\F)}{rk(\F)}$ be the slope of $\F$. We say that $\F$ is \emph{semistable} (resp. \emph{stable}) if for each proper subsheaf $0 \neq \F' \subsetneq \F$ we have $\mu(\F') \leq \mu(\F)$ (resp. $\mu(\F') < \mu(\F)$). Denote by $\CC_q$ the category of semistable sheaves of slope $q \in \Q \cup \{ \infty \}$. Note that $\CC_\infty = \mathsf{coh}_0\ \! \X$ is the subcategory of torsion sheaves. The category $\CC_q$ is a full abelian subcategory of $\coh \X$ closed under extension for any $q$, with simple objects given by the stable sheaves. Weighted projective lines of genus one are characterised by $deg(\omega_\X) = 0$, or equivalently $\omega_\X$ has finite order in ${\rm Pic}(\X)$, and so $\CC_q$ is closed under the Auslander-Reiten translate $\tau = - \otimes \omega_\X$. In the genus one case, indecomposable sheaves are semistable \cite[Thm. 5.6]{GL} and we are lead to describe the categories $\CC_q$ for each $q$. We will do this by means of the Telescopic functors of Lenzing and Meltzer.

Let $\U$ be the $\tau$-orbit of a stable sheaf, e.g. $\U = \{\OO_\X, \omega_\X, \dots, \omega_\X^{p-1}\},\ \U = \{S_{i,0}, S_{i,1}, \dots, S_{i, p_i-1} \}$, or $\U = \{S_x\}$. We define the left mutation $L_\U$ (respectively right mutation $R_\U$) by $\U$ acting on $\F \in {\rm D}^b(\X) $ via the distinguished triangles
\[
	L_\U(\F)[-1] \to \bigoplus_{\E \in \U} {\rm Hom}^{\bullet}(\E, \F) \otimes_k \E \xrightarrow{ev} \F \to L_\U(\F)
\]
\[
	R_\U(\F) \to \F \xrightarrow{coev} \bigoplus_{\E \in \U}{\rm Hom}^{\bullet}(\F, \E)^* \otimes_k \E \to R_\U(\F)[1]
\]
The $L_\U,\ R_\U$ are called \emph{tubular mutations} and recover the notion of spherical twists when $\U = \{\E\}$.
\begin{prop}[{\cite[Thm. 5.1.3]{MeltBook}}] The constructions $L_\U,\ R_\U$ are functorial in $\F$. Moreover, they give inverse autoequivalences $L_\U,\ R_\U: {\rm D}^b(\X) \xrightarrow{\cong} {\rm D}^b(\X) $.
\end{prop}
For simplicity and to fix notation, let us restrict to the case of interest $\X = \mathbb{P}^1(2,2,2,2;\lambda)$. In this case, $\omega_\X$ has order $2$ in ${\rm Pic}(\X)$. In \cite{LM}, \cite[Chp. 5]{MeltBook}, Lenzing and Meltzer consider two autoequivalences $T_1, T_2: {\rm D}^b(\X) \xrightarrow{\cong} {\rm D}^b(\X) $ given respectively by right mutation by the orbit of the structure sheaf $\U = \{\OO_\X, \omega_\X\}$, and left mutation by the orbit of a simple torsion sheaf $\U = \{S_{1, 0}, S_{1, 1}\}$
\[
	T_1(\F) \to \F \xrightarrow{coev} \bigoplus_{j \in \Z_2}{\rm Hom}^{\bullet}(\F, \omega_{\X}^{\otimes j})^* \otimes_k \omega_{\X}^{\otimes j} \to T_1(\F)[1]
\]
\[
	T_2(\F)[-1] \to \bigoplus_{j \in \Z_2} {\rm Hom}^{\bullet}(S_{1, j}, \F) \otimes_k S_{1, j} \xrightarrow{ev} \F \to T_2(\F)
\]
These act on rank and degree by
\[
	\begin{pmatrix} rk(T_1\F) \\ deg(T_1\F) \end{pmatrix} = \begin{pmatrix} 1 & 1 \\ 0 & 1 \end{pmatrix} \begin{pmatrix} rk(\F) \\ deg(\F) \end{pmatrix}
\]
\[
	\begin{pmatrix} rk(T_2\F) \\ deg(T_2\F) \end{pmatrix} = \begin{pmatrix} 1 & 0 \\ 1 & 1 \end{pmatrix} \begin{pmatrix} rk(\F) \\ deg(\F) \end{pmatrix}
\]
and restrict to equivalences
\begin{align*}
	&T_1: \mathcal{C}_q \xrightarrow{\cong} \mathcal{C}_{\frac{q}{q+1}}\quad \quad \quad \textup{ for all } 0 \leq q \leq \infty\\
	&T_2: \mathcal{C}_q \xrightarrow{\cong} \mathcal{C}_{q+1}\quad \quad \quad \textup{ for all } q \in \Q \cup \{\infty \}\\
	&T_1: \CC_\infty \xrightarrow{\cong} \CC_1.
\end{align*}
Let $\mathcal{B}_3 = \langle \sigma_1, \sigma_2 \ | \ \sigma_1 \sigma_2 \sigma_1 = \sigma_2 \sigma_1 \sigma_2 \rangle$ be the braid group on $3$ strands.
\begin{prop}[{\cite[Sect. 5.3.6]{MeltBook}}] The functors $T_1, T_2^{-1}$ satisfy the Braid relations
	\[
		T_1T_2^{-1}T_1 \cong T_2^{-1} T_1 T_2^{-1}.
	\]
	Moreover, the homomorphism $\mathcal{B}_3 \to Aut({\rm D}^b(\X))$ given by $\sigma_1 \mapsto T_1$, $\sigma_2 \mapsto T_2^{-1}$ is fully faithful.
\end{prop}
The induced action on rank and degree then gives rise to the well-known homomorphism $\mathcal{B}_3 \twoheadrightarrow SL(2, \Z)$. It is also well-known that the transformations $T_1: q \mapsto \frac{q}{q+1}$ and $T_2: q \mapsto q+1$ defines an action of the free monoid on two words $F\{T_1,T_2\}$ on the positive rationals $\mathbb{Q}_{+}$, which is free and transitive with single generator $1$. This gives the positive rationals the structure of an infinite binary tree called the \emph{Calkins-Wilf tree}. Writing $q \in \mathbb{Q}_{+}$ as $w_q(T_1, T_2)\cdot 1$ for some unique word $w_q(T_1, T_2)$, one deduces the existence of an autoequivalence
\[
	\Phi_{q, \infty}:= w_q(T_1, T_2) \circ T_1: {\rm D}^b(\X) \xrightarrow{\cong} {\rm D}^b(\X)
\]
restricting to $\mathcal{C}_{\infty} \xrightarrow{\cong} \mathcal{C}_q$ for any $q > 0$. One then extends $\Phi_{q, \infty}$ with the same property to any $q \leq 0$ by $\Phi_{q, \infty} := T_2^{-n} \circ \Phi_{q+n, \infty}$ for $n \gg 0$, with result independent of $n$. The functors $\Phi_{q, \infty}$ are called \emph{Telescopic Functors}.

\begin{figure}\label{freesemigroup}
\[
\xymatrix@C=5pt{ 
	& & & & & & & &\ar@{-}[dllll]_{T_1}1\ar@{-}[drrrr]^{T_2}& & & & & & & &	 \\
	&&&&\ar@{-}[dll]_{T_1}\frac{1}{2}\ar@{-}[drr]^{T_2}&&&& &&&&\ar@{-}[dll]_{T_1}2\ar@{-}[drr]^{T_1}&&&&	 \\
	&&\ar@{-}[dl]_{T_1}\frac{1}{3}\ar@{-}[dr]^{T_2}&& &&\ar@{-}[dl]_{T_1}\frac{3}{2}\ar@{-}[dr]^{T_2}&& &&\ar@{-}[dl]_{T_1}\frac{2}{3}\ar@{-}[dr]^{T_2}&& &&\ar@{-}[dl]_{T_1}3\ar@{-}[dr]^{T_2}&&	 \\
&\cdots&&\cdots& &\cdots&&\cdots& &\cdots&&\cdots& &\cdots&&\cdots&	 
}
\]
\caption{The Calkins-Wilf tree.}
\end{figure}
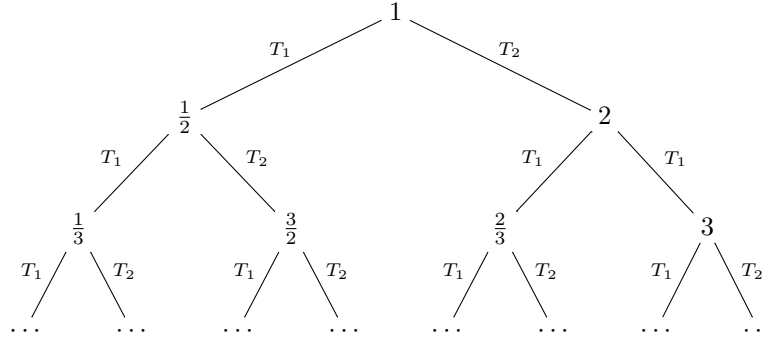

Since the category $\mathcal{C}_{\infty}$ consists of all skyscraper sheaves, the category $\mathcal{C}_q$ is serial for any $q$, with simples given by the stable sheaves of slope $q$. As the type of $\X$ is $(2,2,2,2)$, the Auslander-Reiten quiver of $\mathcal{C}_q$ breaks down into tubes of rank one indexed by the ordinary points $x \in \mathbb{P}^1$ and tubes of rank two indexed by the exceptional points $x_i \in \mathbb{P}^1$. These correspond to indecomposables for which $\F \otimes \omega_{\X} \cong \F$ and $\F \otimes \omega_{\X} \ncong \F$, respectively. Moreover, the exceptional sheaves are precisely the stable sheaves living in rank two tubes. 

For computing morphism spaces in $ {\rm D}^b(\X) $, we have the following well-known results.
\begin{lem}[Lemma 4.1, \cite{LM}]\label{slopevanishing} Let $\F, \G$ be semistable sheaves of slopes $q, q'$.
	\begin{enumerate}
		\item If $q > q'$, then ${\rm Hom}(\F, \G) = 0$.
		\item If $q < q'$, then ${\rm Ext}^1(\F, \G) = 0$.
	\end{enumerate}
\end{lem}

\begin{prop}[Weighted Riemann-Roch, \cite{LM}]\label{RR}We have
	\[
		\chi(\F, \G) + \chi(\F, \G \otimes \omega_{\X}) = \begin{vmatrix} rk(\F) & rk(\G) \\ deg(\F) & deg(\G) \end{vmatrix}.
	\]
	In particular, we get
	\[
		\chi(\F) + \chi(\F \otimes \omega_{\X}) = deg(\F).
	\]
\end{prop}

\subsection{Matrix factorisations corresponding to simple torsion sheaves}
We see that to understand the classification of indecomposables in $ \underline{{\rm MCM}}^{\mathbb{Z}}(R) \cong {\rm D}^b(\X) $, we must: 

\begin{enumerate}[1)]
\item Understand the rank and degree maps 
	\[
		rk, deg: K_0( \underline{{\rm MCM}}^{\mathbb{Z}}(R)) \to \Z
	\]
	and thus define the slope of an (indecomposable) graded MCM module.
\item Recover the category $\CC_\infty$ inside $ \underline{{\rm MCM}}^{\mathbb{Z}}(R)$. 
\item Understand the action of $T_1, T_2$ on $\CC_\infty \subseteq \underline{{\rm MCM}}^{\mathbb{Z}}(R)$. 
	\end{enumerate}

	In order to situate ourselves, let us first calculate the images of $\{k^{st}(-j)\}_{j \in \Z}$ under the above equivalence. By $2$-periodicity we have $k^{st}(-j-4) \cong k^{st}(-j)[-2]$, and so it suffices to compute the image of
\[
	k^{st}, k^{st}(-1), k^{st}(-2), k^{st}(-3).
\]
The Serre functor $\mathbb{S}_R(M) = M(2)$ is sent to the Serre functor $\mathbb{S}_\X(\F) = \F \otimes \omega_{\X}[1]$. The periodicity identity $(4) = [2]$ corresponds to the fact that $\omega_{\X}$ is $2$-torsion. Keeping this in mind, we will prove the following:

\begin{thm}\label{residuefields} Under the equivalence of Cor. \ref{wpltilt}, we have
	\begin{align*}
		 k^{st}	 &\mapsto \omega_{\X}[1] \\
		 k^{st}(-1) &\mapsto  \OO_\X(-\vec{c})[1] \\
		 k^{st}(-2) &\mapsto  \OO_\X\\
		 k^{st}(-3) &\mapsto  \OO_\X(-\vec{c}) \otimes \omega_{\X}.
	\end{align*}
\end{thm}
\begin{proof} By Corollary \ref{wpltilt} we already have $k^{st}(-2) \mapsto \OO_\X$, and so $k^{st} = \mathbb{S}_R\left(k^{st}(-2)\right) \mapsto \mathbb{S}_{\X}\left(\OO_\X\right) = \omega_{\X}[1]$. Let $\F_{k^{st}(-1)}$ correspond to $k^{st}(-1)$. We previously obtained the exceptional pair $(k^{st}(-2), \left(R/\mathfrak{m}^2\right)^{st}(-1))$ as the right mutation
	\[
		R: (k^{st}(-1), k^{st}(-2)) \mapsto \left(k^{st}(-2), R_{k^{st}(-1)}(k^{st}(-2)) \right)
	\]
	and so we can obtain $(\OO_\X, \OO_\X(\vec{c}))$ as the right mutation
	\[
	R: (\F_{k^{st}(-1)}, \OO_\X) \mapsto \left(\OO_\X, R_{\F_{k^{st}(-1)}}(\OO_\X) \right)
	\]
	Since left and right mutations are inverses (Prop. \ref{braidaction}), we can recover $\F_{k^{st}(-1)}$ as the left mutation 
	\[
		\F_{k^{st}(-1)} \cong L_{\OO_\X} \left( R_{\F_{k^{st}(-1)}}(\OO_\X)\right) \cong L_{\OO_\X}(\OO_\X(\vec{c}))
	\]
	calculated by the distinguished triangle
	\[
		L_{\OO_\X}(\OO_\X(\vec{c}))[-1] \to {\rm RHom}(\OO_\X, \OO_\X(\vec{c})) \otimes_k \OO_\X \xrightarrow{ev} \OO_\X(\vec{c}) \to L_{\OO_\X}(\OO_\X(\vec{c})).
	\]
	This calculation can be done inside $ {\rm D}^b(\mathbb{P}^1) \subseteq {\rm D}^b(\X) $, where the above corresponds to the Euler sequence (up to twist)
	\[
		0 \to \OO_{\mathbb{P}^1}(-1) \to {\rm H}^0(\mathbb{P}^1, \OO_{\mathbb{P}^1}(1)) \otimes_k \OO_{\mathbb{P}^1} \xrightarrow{ev} \OO_{\mathbb{P}^1}(1) \to 0.
	\]
	This gives $\F_{k^{st}(-1)} = L_{\OO_\X}(\OO_\X(\vec{c})) \cong  \OO_\X(-\vec{c})[1]$, and finally $k^{st}(-3) \mapsto \OO_\X(-\vec{c}) \otimes \omega^{-1}_\X = \OO_\X(-\vec{c}) \otimes \omega_\X$.
\end{proof}
\begin{rem} Our choice of tilting object $U = T(-3)[1]$ in Cor. \ref{wpltilt} was made in anticipation of this result, and this normalisation will be helpful in calculations of Betti tables in further sections.
\end{rem}\bigskip

We can produce the indecomposable matrix factorisations corresponding to the simple objects in $\mathcal{C}_\infty$. For simplicity we assume $\mathsf{char}\ k \neq 2$ for the calculation. We already have the bases of the rank two tubes, given by
\[
	L_i(-3)[1] = {\rm coker}(R(-3) \xrightarrow{f_\lambda/l_i} R) = R/(\frac{f_\lambda}{l_i})
\]
which are sent to $S_{i,0}$ under Cor. \ref{wpltilt}, and $\tau L_i(-3)[1] = L_i(-1) = R/l_i(-1)$ sent to $S_{i,1}$. By the previous theorem, the pair $\left( k^{st}(-1)[-1], k^{st}(-2)\right)$ corresponds to $\left(\OO_\X(-\vec{c}), \OO_\X\right)$, with $2$-dimensional morphism space. Taking a basis\footnote{This may not correspond a priori to the basis $\{u,v\}$ of ${\rm H}^0(\X, \OO_\X(\vec{c}))$ on the other side. However, the choices we will make in Prop. \ref{morphismbasis} will turn out to correspond to $\{u,v\}$ up to rescaling, see Lemma \ref{SESexceptionaltubes} and following remark.} $\phi_0, \phi_\infty$ for $ {\rm \underline{Hom}}_{gr R}(k^{st}(-1)[-1], k^{st}(-2))$, for any $p = [p_0 : p_1] \in \mathbb{P}^1$ we define the MCM module $M_p$ as the cone of $\phi_p = p_1 \phi_0 + p_0 \phi_\infty$
\[
	k^{st}(-1)[-1] \xrightarrow{\phi_p} k^{st}(-2) \to M_p \to k^{st}(-1)[2]
\]
whose distinguished triangle is sent to
\[
	\OO_\X(-\vec{c}) \xrightarrow{s_p} \OO_\X \to S_p \to \OO_\X(-\vec{c})[1].
\]
for the corresponding cosection $s_p$, with cokernel an ordinary skyscraper sheaf. We can produce the associated matrix factorisations. We have $f_\lambda = xy(x-y)(x- \lambda y) = x^3y - (1+\lambda)x^2y^2 + \lambda xy^3$. Write $f_\lambda = x f_x + y f_y$ for $f_x = \frac{1}{4} \frac{\partial f_\lambda}{\partial x}$ and $f_y = \frac{1}{4}\frac{\partial f_\lambda}{\partial y}$, and note that $x | f_y$ and $y | f_x$. We have already seen part $1)$ of the next result, which we restate for convenience. 
\begin{prop}\label{morphismbasis} Assume $\mathsf{char}\ k \neq 2$. We have the following explicit presentations.
	\begin{enumerate}[1)]
		\item $k^{st}$ corresponds to the matrix factorisation
			\[
			\xymatrix@C2pc{S(-4) \oplus S(-2)  \ar[r]^-{B} &S(-1) \oplus S(-1)  \ar[r]^-{A} & S \oplus S(2)}
			\]
			with
			\begin{align*}
				A = \begin{pmatrix} x & y \\ -f_y & f_x \end{pmatrix} && B = \begin{pmatrix} f_x & -y \\ f_y & x \end{pmatrix}.
			\end{align*}
		\item A basis of morphisms $\phi_0, \phi_\infty: k^{st}(-1)[-1] \to k^{st}(-2)$ can be taken as
			\[
				\xymatrix@C2pc{S(-6) \oplus S(-6) \ar@<.5ex>[d]^{\varphi_\infty} \ar@<-.5ex>[d]_{\varphi_0}     \ar[r]^-{-A}& S(-5) \oplus S(-3) \ar@<.5ex>[d]^{\psi_\infty} \ar@<-.5ex>[d]_{\psi_0}   \ar[r]^-{-B}& S(-2) \oplus S(-2) \ar@<.5ex>[d]^{\varphi_\infty} \ar@<-.5ex>[d]_{\varphi_0}  \\
				S(-6) \oplus S(-4)  \ar[r]^-{B}& S(-3) \oplus S(-3)  \ar[r]^-{A}& S(-2) \oplus S}
			\]
			with matrices given by
			\begin{align*}
				&\varphi_0 = \begin{pmatrix} 0 & 1 \\ 0 & -\frac{f_y}{x}  \end{pmatrix} && \psi_0 = \begin{pmatrix} -\frac{f_y}{x} & -1 \\ 0 & 0\end{pmatrix}\\
				&\varphi_\infty = \begin{pmatrix} 1 & 0 \\ \frac{f_x}{y} & 0  \end{pmatrix} && \psi_\infty = \begin{pmatrix} 0 & 0\\ -\frac{f_x}{y} & 1\end{pmatrix}.
			\end{align*}
	\end{enumerate}
\end{prop}
\begin{proof} Part $1)$ follows from the Tate resolution, see \cite{Eis}. For part $2)$, note that the MCM approximation $k^{st}(-2) \to k(-2)$ corresponds to the natural projection and induces natural isomorphisms on Tate cohomology
\[
	\underline{{\rm Ext}}^0_{gr R}\left(k^{st}(-1)[-1], k^{st}(-2)\right) \xrightarrow{\cong} \underline{{\rm Ext}}^0_{gr R}\left(k^{st}(-1)[-1], k(-2)\right).
\]
The morphisms $\phi_0, \phi_\infty$ descend to the natural basis on the latter. \end{proof}
Now let $\phi_p = p_1 \phi_0 + p_0 \phi_\infty$. Taking ${\rm Cone}(\phi_p)$ yields a $4 \times 4$ matrix factorisation
\[
	\xymatrix@R0.01pc@C5pc{S(-6) & S(-3) & S(-2)\\
	S(-4)  &S(-3)   &S\\
		S(-5)\ar@<2ex>[r]^{\begin{pmatrix} B & \psi_p \\ 0 & B \end{pmatrix}} & S(-2) \ar@<2ex>[r]^{\begin{pmatrix} A & \varphi_p \\ 0 & A \end{pmatrix}} & S(-1)\\
		S(-3) & S(-2) & S(1)
	}
\]
with $\varphi_p = p_1\varphi_0 + p_0\varphi_\infty$, $\psi_p = p_1\psi_0 + p_0\psi_\infty$ and matrices given by
\begin{align*}
	\begin{pmatrix} A & \varphi_p \\ 0 & A \end{pmatrix} &=
	\begin{pmatrix} 
	x & y & p_0 & p_1 \\
	-f_y & f_x & p_0\frac{f_x}{y} & -p_1\frac{f_y}{x} \\
	0 & 0 & x & y \\ 
	0 & 0 & -f_y & f_x 
\end{pmatrix}
		\\
		\\
	\begin{pmatrix} B & \psi_p  \\0 & B \end{pmatrix} &=
	\begin{pmatrix} 
		f_x & -y & -p_1\frac{f_y}{x} & -p_1 \\
		f_y & x & -p_0\frac{f_x}{y} & p_0 \\
	0& 0 & f_x & -y \\ 
	0& 0 & f_y & x 
\end{pmatrix}
\end{align*}
The matrices $\varphi_p,\ \psi_p$ have two scalar entries, and so this matrix factorisation is stably equivalent to a $2 \times 2$ matrix factorisation. Direct calculations show the following:
\begin{prop}\label{ordinaryskyscraper} The module $M_p$ is given by the reduced matrix factorisation
			\[
				\xymatrix@C2pc{S(-5) \oplus S(-4)  \ar[r]^-{B_p}& S(-2) \oplus S(-3)  \ar[r]^-{A_p}& S(-1) \oplus S}
			\]
			where $(A_p, B_p)$ for $p_1 \neq 0$ are given by
			\begin{align*}
				A_p = \begin{pmatrix} x - \frac{p_0}{p_1} y & \frac{1}{p_1}y^2 \\ -p_0 \frac{f_\lambda}{xy} & \frac{f_\lambda}{x} \end{pmatrix} && B_p = \begin{pmatrix} \frac{f_\lambda}{x} & -\frac{1}{p_1}y^2 \\ p_0 \frac{f_\lambda}{xy} & x - \frac{p_0}{p_1} y \end{pmatrix}
			\end{align*}
			and for $p_1 = 0$, $p_0 \neq 0$ by
			\begin{align*}
				A_p = \begin{pmatrix} y & x^2 \\ 0 & \frac{f_\lambda}{y} \end{pmatrix} && B_p = \begin{pmatrix} \frac{f_\lambda}{y} & -x^2 \\ 0 & y \end{pmatrix}
			\end{align*}
\end{prop}

Recall that $L_i(-3)[1] = R/(\frac{f_\lambda}{l_i})$ and $L_i(-1) = R/l_i(-1)$ correspond to the simple sheaves $S_{i, 0},\ S_{i,1}$. From the above presentation, one verifies:
\begin{lem}\label{SESexceptionaltubes} For each exceptional point $p_i = V(l_i) \in \mathbb{P}^1$, there are short exact sequences of MCM modules
	\[
		0 \to	R/l_i(-1) \to M_{p_i} \to R/(\frac{f_\lambda}{l_i}) \to 0.
	\]
\end{lem}
Hence the $M_{p_i}$ corresponds to $S_{p_i}$, the `ordinary' skyscraper sheaf over the exceptional point $p_i$ for $p_i = 0, \infty, 1, \lambda$.
\begin{rem} Let $\{u', v'\}$ be a basis of ${\rm Hom}(\OO_\X, \OO_\X(-\vec{c})) \cong {\rm H}^0(\X, \OO_\X(\vec{c}))$ corresponding to $\{\phi_0, \phi_\infty\}$. There is an invertible matrix taking $\{u', v'\}$ to $\{u, v\}$. However, the induced transformation on $\mathbb{P}^1$ fixes the $4$ exceptional points and so must be trivial, hence the basis $\{u', v'\}$ is given by $\{u,v\}$ up to rescaling. 
\end{rem}

Summarising, we have shown:
\begin{prop}\label{simpleinfty} Assume $\mathsf{char}\ k \neq 2$. The set of indecomposable MCM modules
	\[
			\{ M_p \}_{p \neq 0, \infty,1, \lambda} \cup \{R/(\frac{f_\lambda}{l_i}), R/l_i(-1) \}_{i=1,2,3,4}
	\]
	correspond under the equivalence of Corollary \ref{wpltilt} to the set of simple torsion sheaves in $\mathcal{C}_\infty$
	\[
		 \{S_p\}_{p \neq 0,\infty,1, \lambda} \cup \{S_{i,0}, S_{i,1}\}_{i = 1,2,3,4}. 
	\]
\end{prop}

We now investigate the shape of the Betti table of indecomposable graded MCM modules in a more systematic way. This will occupy the remaining sections. 

\subsection{Betti tables from cohomology tables}

Now given $M$, write $\F_M$ for the corresponding complex of coherent sheaves. In Thm. \ref{residuefields} we saw that
	\begin{align*}
		\F_{k^{st}}	 &= \omega_{\X}[1] \\
		\F_{k^{st}(-1)} &=  \OO_\X(-\vec{c})[1] \\
		\F_{k^{st}(-2)} &=  \OO_\X\\ 
		\F_{k^{st}(-3)} &=  \OO_\X(-\vec{c}) \otimes \omega_{\X}.
	\end{align*}
	with $\F_{k^{st}(-j-4)} = \F_{k^{st}(-j)}[-2]$. Let $C$ be the minimal complete resolution of $M$, which looks like \[
		C: \cdots \to \bigoplus_{j \in \Z} R(-j)^{\oplus \beta_{i+1,j}} \to \bigoplus_{j \in \Z} R(-j)^{\oplus \beta_{i,j}} \to \bigoplus_{j \in \Z} R(-j)^{\oplus \beta_{i-1,j}} \to  \cdots
	\]
	We can calculate the graded Betti numbers $\beta_{i,j}$ by
	\[
		\beta_{i,j} = \dim_k {\rm Hom}_{gr R}(C, k[i](-j)) = \dim_k \underline{{\rm Ext}}_{gr R}^i(M, k(-j)) = \dim_k \underline{{\rm Ext}}_{gr R}^i(M, k^{st}(-j)) 
	\]
	and so by the dimension of the corresponding morphism space in $ {\rm D}^b(\X) $. This idea was used by A. Pavlov in his thesis to produce classifications of Betti tables of graded MCM modules over the cone of various embeddings of an elliptic curve \cite{Pavlov}. In our context, Thm. \ref{residuefields} implies:

\begin{cor} We can calculate Betti numbers $\beta_{i,j} = \beta_{i,j}(M)$ as follows:
	\begin{align*}
	& \beta_{i,0} = \dim_k {\rm Ext}^i(\F_{M}, \omega_{\X}[1]) = h^{-i}(\F_{M}) \\
	& \beta_{i,1} = \dim_k {\rm Ext}^i(\F_{M}, \OO_\X(-\vec{c})[1]) = h^{-i}(\F_M(\vec{c}) \otimes \omega_{\X})\\
	& \beta_{i,2} = \dim_k {\rm Ext}^i(\F_{M}, \OO_\X) = h^{1-i}(\F_{M} \otimes \omega_{\X})\\
	& \beta_{i,3} = \dim_k {\rm Ext}^i(\F_{M}, \OO_\X(-\vec{c}) \otimes \omega_{\X}) = h^{1-i}(\F_M(\vec{c})).
	\end{align*}
When $\F_M$ is a coherent sheaf, collecting terms via the periodicity $\beta_{i, j} = \beta_{i+2, j+4}$, the only possible non-trivial Betti numbers for $M$ of the form $\beta_{0, *}, \beta_{1, *}$ are
\[
	\beta_{0, 0}, \beta_{0, 1}, \beta_{0, 2}, \beta_{0,3}, \beta_{1,2}, \beta_{1,3}, \beta_{1,4}, \beta_{1,5}.
\]
\end{cor}

Since $\coh \X$ is hereditary, indecomposables in ${\rm D}^b(\X)$ are of the form $\F[n]$ for $\F$ an indecomposable coherent sheaf and $n \in \mathbb{Z}$, and it suffices to work out Betti tables corresponding to coherent sheaves. In this case, the data is best expressed in the following table:
\begin{align*}
	\beta(M) = \begin{pmatrix*}[c]\beta_{0,0} & \beta_{1,2} \\ \beta_{0,1} & \beta_{1,3} \\ \beta_{0,2} & \beta_{1,4} \\ \beta_{0,3} & \beta_{1,5} \\ \end{pmatrix*}
	= \begin{pmatrix*}[l] h^0(\F) & h^0(\F \otimes \omega_{\X}) \\ h^0(\F(\vec{c}) \otimes \omega_{\X}) & h^0(\F(\vec{c})) \\ h^1(\F\otimes \omega_{\X}) & h^1(\F) \\ h^1(\F(\vec{c})) & h^1(\F(\vec{c}) \otimes \omega_{\X}) \\   \end{pmatrix*}
\end{align*}
where $\F = \F_M$. We will refer to the latter table as the cohomology table $\beta(\F)$.

\begin{exmp} Since $\F_{k^{st}(-2)} = \OO_\X$, we can calculate
\begin{align*}
	\beta(k^{st}(-2)) = \begin{pmatrix*}[c] 1 & 0 \\ 0 & 2 \\ 1 & 0 \\ 0 & 0 \\ \end{pmatrix*}
	= \begin{pmatrix*}[l] h^0(\OO_\X) & h^0(\omega_{\X}) \\ h^0(\OO_\X(\vec{c}) \otimes \omega_{\X}) & h^0(\OO_\X(\vec{c})) \\ h^1(\omega_{\X}) & h^1(\OO_\X) \\ h^1(\OO_\X(\vec{c})) & h^1(\OO_\X(\vec{c}) \otimes \omega_{\X}) \\   \end{pmatrix*}
\end{align*}
and we recover the Betti table of Prop. \ref{morphismbasis}. 
\end{exmp}
\begin{exmp} Since $\F_{M_p} = S_p$, we have 
\begin{align*}
	\beta(M_p) = \begin{pmatrix*}[c]1 & 1 \\ 1 & 1 \\ 0 & 0 \\ 0 & 0 \\ \end{pmatrix*}
	= \begin{pmatrix*}[l] h^0(S_p) & h^0(S_p \otimes \omega_{\X}) \\ h^0(S_p(\vec{c}) \otimes \omega_{\X}) & h^0(S_p(\vec{c})) \\ h^1(S_p\otimes \omega_{\X}) & h^1(S_p) \\ h^1(S_p(\vec{c})) & h^1(S_p(\vec{c}) \otimes \omega_{\X}) \\   \end{pmatrix*}
\end{align*}
as $S_p \otimes \mathcal{L} \cong S_p$ for any line bundle $\mathcal{L}$ and `ordinary' skyscraper sheaf $S_p$. This recovers the Betti table of Prop. \ref{ordinaryskyscraper}.
\end{exmp}

\subsubsection*{The $C_4$ group action} 

The autoequivalence $M \mapsto M(1)$ on $ \underline{{\rm MCM}}^{\mathbb{Z}}(R)$ induces an autoequivalence $\F_M \mapsto \F_M\{1\}$ on ${\rm D}^b(\X)$. Since shifting the grading acts by translation on Betti tables, it suffices to compute one Betti table in the orbit $\{M(n)\}_{n \in \Z}$, or equivalently one cohomology table in the orbit $\{\F_M\{n\} \}_{n \in \Z}$. First, we calculate the effect of $(1)$ on rank and degree. We begin by finding an intrinsic description of the maps
\[
	rk, deg: K_0( \underline{{\rm MCM}}^{\mathbb{Z}}(R)) \to \Z
\]
defined by $deg(M) := deg(\F_M)$ and $rk(M) := rk(\F_M)$.
\begin{lem}\label{rkdeg} For any graded MCM module M, we have 
	\begin{align*}
		deg(M)  &= \chi(M, k^{st}) - \chi(M, k^{st}(-2)) \\
			&= \sum_{i \in \Z}(-1)^{i} \beta_{i, 0} - \sum_{i \in \Z}(-1)^{i} \beta_{i, 2},	\\
			rk(M) 	&= \frac{1}{2} \chi(M, k^{st}(-1)\oplus k^{st}(-2)) - \frac{1}{2} \chi( M, k^{st} \oplus k^{st}(-3)) \\
			&= \frac{1}{2} \sum_{i \in \Z} (-1)^{i}(\beta_{i, 1} + \beta_{i, 2}) - \frac{1}{2} \sum_{i \in \Z} (-1)^{i}(\beta_{i, 0} + \beta_{i,3}).
	\end{align*}
When M corresponds to a coherent sheaf, this simplifies to
\begin{align*}
	deg(M) &= (\beta_{0, 0} + \beta_{1,2}) - (\beta_{0,2} + \beta_{1,4}), \\
	rk(M) &= \frac{1}{2}(\beta_{0,1} + \beta_{0,2} + \beta_{1,3} + \beta_{1,4}) - \frac{1}{2}(\beta_{0,0} + \beta_{0,3} + \beta_{1,2} + \beta_{1,5}).
\end{align*}
\end{lem}
\begin{proof} The formula for $deg(M)$ falls out of the weighted Riemann-Roch theorem via Thm. \ref{residuefields}. To deduce the formula for $rk(M)$, we use $deg(\F(\vec{c})) = deg(\F) + 2 rk(\F)$, so that $rk(\F) = \frac{1}{2}(deg(\F(\vec{c})) - deg(\F))$, then collect terms via Thm. \ref{residuefields}.
\end{proof}

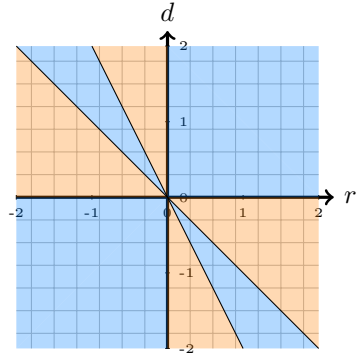
\begin{figure}
\begin{align*}
   \begin{minipage}[r]{4cm}
       \begin{tikzpicture}
	    \draw[gray!50, thin, step=0.3] (-2,-2) grid (2,2);
	    \draw[very thick,->] (-2,0) -- (2.2,0) node[right] {$r$};
	    \draw[very thick,->] (0,-2) -- (0,2.2) node[above] {$d$};
	    \foreach \x in {-2,...,2} \draw (\x,0.03) -- (\x,-0.03) node[below] {\tiny\x};
	    \foreach \y in {-2,...,2} \draw (-0.03,\y) -- (0.03,\y) node[right] {\tiny\y};
	    \fill[blue!50!cyan,opacity=0.3] (0,0) -- (0,2) -- (2,0);
	    \fill[blue!50!cyan,opacity=0.3] (2,0) -- (0,2) -- (2,2);
	    \fill[blue!50!cyan,opacity=0.3] (0,0) -- (-2,0) -- (-2,-2);
	    \fill[blue!50!cyan,opacity=0.3] (0,0) -- (0,-2) -- (-2,-2);
	    \fill[blue!50!cyan,opacity=0.3] (0,0) -- (-1,2) -- (-2,2);
	    \fill[blue!50!cyan,opacity=0.3] (0,0) -- (2,-2) -- (1,-2);
	    \fill[orange,opacity=0.3] (0,0) -- (-1,2) -- (0,2);
	    \fill[orange,opacity=0.3] (0,0) -- (1,-2) -- (0,-2);
	    \fill[orange,opacity=0.3] (0,0) -- (-2,0) -- (-2,2);
	    \fill[orange,opacity=0.3] (0,0) -- (2,0) -- (2,-2);
	    \draw (-1,2) -- node[below,sloped] {} (1,-2);
	    \draw (-2,2) -- node[below,sloped] {} (2,-2);
       \end{tikzpicture}
   \end{minipage}
\end{align*}
\caption{Orbits of $\langle r \rangle = C_4$.}\label{orbits}
\end{figure}

We now calculate the effect of $M \mapsto M(1)$ on the vector 
\[
	\ZZ(M) = \begin{pmatrix} rk(M) \\ deg(M) \end{pmatrix}
\]
\begin{prop}\label{matrix} For any graded MCM module $M$, we have 
	\begin{align*}
		\begin{pmatrix} rk(M(1)) \\ deg(M(1)) \end{pmatrix} &= \begin{pmatrix} -1 & -1 \\ 2 & 1 \end{pmatrix} \begin{pmatrix} rk(M) \\ deg(M) \end{pmatrix}. \\
	\end{align*}
\end{prop}
\begin{proof} Writing $K_0 := K_0( \underline{{\rm MCM}}^{\mathbb{Z}}(R))$, we will check that the following commutes:
	\[
		\xymatrix@C5pc{
			K_0 \ar[d]^{(1)}\ar[r]^{\ZZ} & \Z^2 \ar[d]^{\begin{psmallmatrix} -1 & -1 \\ 2 & 1 \end{psmallmatrix}} \\
		K_0 \ar[r]^{\ZZ} & \Z^2 \\
	}
\]
The full exceptional collection of Orlov's Theorem $\langle k^{st}(-1), k^{st}(-2), L_1(-3)[1], \dots, L_4(-3)[1] \rangle$ gives a $\Z$-basis of $K_0$ on which this can be checked. We have
\begin{align*}
	&\ZZ\left(k^{st}(-2)\right) = \ZZ\left(\OO_\X\right) = \begin{pmatrix} 1 \\ 0 \end{pmatrix} \\ 
	&\ZZ\left(k^{st}(-1)\right) = \ZZ\left(\OO_\X(-\vec{c})[1]\right) =\begin{pmatrix} -1 \\ 2 \end{pmatrix}\\
	&\ZZ\left(k^{st}\right) = \ZZ\left(\omega_{\X}[1]\right) = \begin{pmatrix} -1 \\ 0 \end{pmatrix}\\
	\\
	&\ZZ\left(L_i(-3)[1]\right) = \ZZ\left(S_{i, 0} \right) = \begin{pmatrix} 0 \\ 1\end{pmatrix} \\
	&\ZZ\left(L_i(-2)[1]\right) = \begin{pmatrix} -1 \\ 1 \end{pmatrix}.
\end{align*}
The last line is calculated by Lemma \ref{rkdeg}, and so the above diagram commutes. 
\end{proof} 
Let $r \in {\rm GL}_2(\Z)$ be the matrix
\begin{align*}
	r = \begin{pmatrix} -1 & -1 \\ 2 & 1 \end{pmatrix}.
\end{align*}
Note that $r^2 = -{\rm I}_2$, which also follows formally from $(2) = \Se_R$ corresponding to $\Se_\X = - \otimes \omega_\X[1]$ in $ {\rm D}^b(\X)$, since $deg(\omega_\X) = 0$. Let $C_4 = \langle r \rangle$. Our goal will be to describe the possible $\beta(M)$ for $M$ indecomposable with $(rk(M), deg(M)) = (r,d)$ fixed, and it is sufficient to do this on a fundamental domain for $C_4$. The action of $r$ preserves convex cones and integrality and induces a chamber decomposition of $\Z^{\oplus 2}$, which simple calculations show can be pictured as in Figure \ref{orbits}.
We will implicitly ignore $(0,0) \in \Z^{\oplus 2}$ in statements to come as there are no indecomposable sheaves of this type, and consider only the action on $\Z^{\oplus 2} \setminus \{(0,0)\}$. From Figure \ref{orbits}, one sees:
\begin{prop} The union of the three regions $\mathcal{R}_1 \cup \mathcal{R}_2 \cup \mathcal{R}_3$ below is a fundamental domain for the action of $C_4$ on $\Z^{\oplus 2} \setminus \{ (0,0) \}$: 
		\begin{align*}
			&\mathcal{R}_1 = \{ (r,d) \ | \ r \geq 0,\ d > 0 \}\\
			&\mathcal{R}_2 = \{ (r,d) \ | \ r > 0,\ d = 0 \}\\
			&\mathcal{R}_3 = \{ (r,d) \ | \ r > 0,\ d < -2r \}.
		\end{align*}
\end{prop}
\begin{figure}
\centerline{
       \begin{tikzpicture}
	    \draw[gray!50, thin, step=0.3] (-2,-2) grid (2,2);
	    \draw[very thick,->] (-2,0) -- (2.2,0) node[right] {$r$};
	    \draw[very thick,->] (0,-2) -- (0,2.2) node[above] {$d$};
	    \foreach \x in {-2,...,2} \draw (\x,0.03) -- (\x,-0.03) node[below] {\tiny\x};
	    \foreach \y in {-2,...,2} \draw (-0.03,\y) -- (0.03,\y) node[right] {\tiny\y};
	    \fill[blue!50!cyan,opacity=0.3] (0,0) -- (0,2) -- (2,0);
	    \fill[blue!50!cyan,opacity=0.3] (2,0) -- (0,2) -- (2,2);
	    \fill[blue!50!cyan,opacity=0.3] (0,0) -- (1,-2) -- (0,-2);
	    \draw (-1,2) -- node[below,sloped] {} (1,-2);
	    \draw (-2,2) -- node[below,sloped] {} (2,-2);
       \end{tikzpicture}}
       \caption{Fundamental domain for $C_4$.}
	\label{FundamentaldomainC4}
\end{figure}
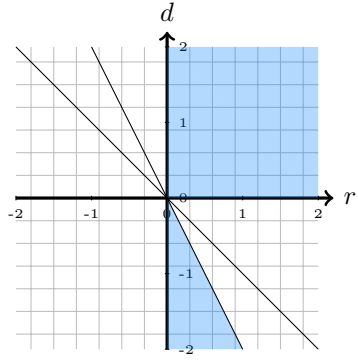

The choice of domain for $C_4$ may appear odd, but note that since $r \geq 0$ throughout, each pair $(r,d)$ is realized by a coherent sheaf and we need not consider complexes. Another reason for this choice is to maximize vanishing patterns in the cohomology table
\[
	\beta(\F) = \begin{pmatrix*}[l] h^0(\F) & h^0(\F \otimes \omega_{\X}) \\ h^0(\F(\vec{c}) \otimes \omega_{\X}) & h^0(\F(\vec{c})) \\ h^1(\F\otimes \omega_{\X}) & h^1(\F) \\ h^1(\F(\vec{c})) & h^1(\F(\vec{c}) \otimes \omega_{\X}) \\   \end{pmatrix*}
\]

\begin{lem}\label{h-vanishing} Let $\F$ be an indecomposable coherent sheaf.
	\begin{enumerate}[1)]
	\item For $\F$ in region $\mathcal{R}_1$, we have $h^1(\F) = h^1(\F \otimes \omega_{\X}) = h^1(\F(\vec{c})) = h^1(\F(\vec{c})\otimes \omega_{\X}) = 0.$
\item For $\F$ in region $\mathcal{R}_3$, we have $h^0(\F) = h^0(\F \otimes \omega_{\X}) = h^0(\F(\vec{c})) = h^0(\F(\vec{c})\otimes \omega_{\X}) = 0.$	\end{enumerate}
\end{lem}
\begin{proof} These follow from Lemma \ref{slopevanishing} by slope arguments, using the formula $$\mu(\F \otimes \mathcal{L}) = \mu(\F) + deg(\mathcal{L})$$ for a line bundle $\mathcal{L}$ and recalling that $deg(\omega_\X) = 0$ and $deg(\OO_\X(\vec{c})) = 2$. \end{proof}
This lemma reduces calculations in regions $\mathcal{R}_1, \mathcal{R}_3$ to computing Euler characteristics, and region $\mathcal{R}_2$ can be dealt with by hand. In the upcoming sections we will completely classify the Betti tables $\beta(M)$ of indecomposables with $\ZZ(M)$ in the fundamental domain of $C_4$. The smaller domain $\Z_{\geq 0}^{\oplus 2}$ for $G$ will play a role in a later section, where we will discuss its role in explicit constructions of matrix factorisations.

\subsection{Cohomology tables of indecomposable coherent sheaves}

\subsubsection*{Cohomology tables for rank one tubes}

We are now in a position to compute the cohomology tables of indecomposable sheaves. We will list the corresponding possible Betti tables of matrix factorisations in a later section, under a different normalisation. We begin with indecomposables living in rank one tubes, or equivalently which satisfy $\F \otimes \omega_{\X} \cong \F$. Recall that we denote by $\Phi_{q, \infty}$ the Telescopic autoequivalence of $ {\rm D}^b(\X) $ which restricts to $\Phi_{q, \infty}: \CC_\infty \xrightarrow{\cong} \CC_q$.

\begin{thm}\label{rankone} Let $(r,d)$ be in the fundamental domain with $q = \frac{d}{r}$. Consider $\F$ with $\F \otimes \omega_{\X} \cong \F$ and $(rk(\F), deg(\F)) = (r,d)$. An indecomposable such $\F$ exists if and only if $gcd(r,d)$ is even, in which case $\frac{|gcd(r,d)|}{2}$ gives the length of $\F$ in $\mathcal{C}_q$. The cohomology table $\beta(\F)$ is then given as:
\[
	\begin{tabular}{|c|c|c|}
	\hline
	$r \geq 0,\ d > 0$ & $r > 0,\ d = 0$ & $r > 0,\ d < -2r$ \\
	\hline
	&& \\
	$\begin{pmatrix} \frac{d}{2} & \frac{d}{2} \\ \frac{d}{2} + r & \frac{d}{2} + r \\ 0 & 0 \\ 0 & 0 \\ \end{pmatrix}$
	& $\begin{pmatrix} 0 & 0 \\ r & r \\ 0 & 0 \\ 0 & 0 \\ \end{pmatrix}$
	& $\begin{pmatrix} 0 & 0 \\ 0 & 0 \\ -\frac{d}{2} & -\frac{d}{2} \\ -\frac{d}{2} - r & -\frac{d}{2} - r \end{pmatrix}$ \\
	&& \\
	\hline
\end{tabular}
\]
\end{thm}
\begin{proof} An indecomposable $\F$ is in the image of $\Phi_{q, \infty}: \CC_\infty \xrightarrow{\cong} \CC_q$, and this functor acts on $(r,d)$ by ${\rm SL}(2, \Z)$ transformation and therefore preserves $\gcd(r,d)$. In $\CC_\infty$, the lowest value of $\gcd(r,d) = \gcd(0, d) = d$ possible for indecomposables in rank one tubes is $2$, realised by the `ordinary sheaves' $S_x$, with higher values $2n$ realised by $S_x \langle n \rangle$. This proves the claim except for the shape of $\beta(\F)$.

	Now, $\F \mapsto \F \otimes \omega_{\X}$ acts by column change on cohomology tables, and so $\beta(\F)$ is symmetrical. By Riemann-Roch we have $2 \cdot \chi(\F) = d$ and $2 \cdot \chi(\F(\vec{c})) = d + 2r$. Combining this with Lemma $\ref{h-vanishing}$ determines tables in region $\mathcal{R}_1, \mathcal{R}_3$, and we now consider the region $\mathcal{R}_2$ given by sheaves of degree zero. Since $\F$ lives in a rank one tube, it lives in a disjoint component from $\OO_\X$. Applying $\Phi_{0, \infty}^{-1}$ sends them to torsion sheaves with disjoint supports, and therefore ${\rm Ext}^*(\OO_\X, \F) = 0$. An application of Lemma \ref{h-vanishing} and Riemann-Roch as above determines $\beta(\F)$.
\end{proof}

\subsubsection*{Cohomology tables for rank two tubes}
We now study indecomposable sheaves $\F$ with $\F \otimes \omega_{\X} \ncong \F$. We begin with some generalities, most of which is well-known.

\begin{prop}\label{ranktwoprop} Let $(r,d)$ be in the fundamental domain and $q = \frac{d}{r}$. The following hold:
	\begin{enumerate}[i)]
		\item For any $(r,d)$, there is an indecomposable $\F$ with $(rk(\F), deg(\F)) = (r,d)$ and $\F \otimes \omega_\X \ncong \F$.
		\item Any such indecomposable $\F$ has length $|gcd(r,d)|$ in $\mathcal{C}_q$.
		\item There are finitely many indecomposable sheaves of type $(r,d)$ if and only if $gcd(r,d)$ is odd, in which case there are exactly eight.
		\item There is an exceptional sheaf of type $(r,d)$ if and only if $|gcd(r,d)| = 1$, in which case all such indecomposables are exceptional.
		\item When $gcd(r,d)$ is even and $d \neq 0$, $\beta(\F) = \beta(\widetilde{\F})$ where $\widetilde{\F}$ is indecomposable of same rank and degree, and $\widetilde{\F} \otimes \omega_{\X} \cong \widetilde{\F}$.
	\end{enumerate}
\end{prop}
\begin{proof} The first four points follow from the autoequivalence $\Phi_{q, \infty}$ as in the proof of Theorem \ref{rankone}. For $v)$, similarly reduce to skyscraper sheaves. Let $S\langle 2n \rangle$ be an indecomposable torsion sheaf supported over the exceptional point $x_i$ of degree $2n$. Then $[S\langle 2n \rangle]$ has height $2n$ in its tube of rank two, and computing Grothendieck classes gives $[S\langle 2n \rangle] = n [S_{i, 0}] + n[S_{i, 1}]$. In particular the ``ordinary'' torsion sheaf $S_{x}$ for $x = x_i$ has degree $2$, and we have $[S_x] = [S_{i,0}] + [S_{i,1}]$. From the presentation 
\[
	0 \to \OO_\X(-\vec{c}) \to \OO_\X \to S_x \to 0
\]
we see that $[S_x] = [S_{x'}]$ for any ordinary point $x'$, and so $[S\langle 2n \rangle] = n[S_x] = n[S_{x'}] = [S_{x'} \langle n \rangle]$ where $S_{x'}\langle n \rangle$ is a length $n$ indecomposable sheaf supported at $x'$.

Applying $\Phi_{q, \infty}$, we deduce that for any indecomposable $\F$ of type $(r,d)$ with $gcd(r,d)$ even, there is another indecomposable $\widetilde{\F}$ of type $(r,d)$ with $[\widetilde{\F}] = [\F]$ and $\widetilde{\F}\otimes \omega_{\X} \cong \widetilde{\F}$. Since $[\widetilde{\F}] = [\F]$, we have $\chi(\widetilde{\F} \otimes \mathcal{L}) = \chi(\F \otimes \mathcal{L})$ for any line bundle $\mathcal{L}$, and outside of the case $d = 0$, those values determine the cohomology table, hence $\beta(\F) = \beta(\widetilde{\F})$.
\end{proof}

The remainder of the section will be aimed at proving the next theorem.
\begin{thm}\label{ranktwothm} Let $(r,d)$ be in the fundamental domain. The $\beta(\F)$ of indecomposables of type $(r,d)$ satisfying $\F \otimes \omega_{\X} \ncong \F$ are listed as follows:
\[
\resizebox{\columnwidth}{!}{%
\begin{tabular}{|c|c|c|}
	\hline
	$(r,d)$ & $r \geq 0,\ d > 0$ & $r > 0,\ d < -2r$ \\
	\hline
	&& \\
$d$ odd & $\begin{pmatrix} \frac{d \pm 1}{2} & \frac{d\mp 1}{2} \\ \frac{d\mp1}{2} + r & \frac{d\pm1}{2} + r \\ 0 & 0 \\ 0 & 0 \\ \end{pmatrix}_4$ & $\begin{pmatrix} 0 & 0 \\ 0 & 0 \\ -\frac{d\pm1}{2} & -\frac{d\mp1}{2} \\ -\frac{d\mp1}{2} - r & -\frac{d\pm1}{2} - r  \\ \end{pmatrix}_4$ \\
	&& \\
	\hline
	&& \\
	(odd, even) & $\begin{pmatrix} \frac{d}{2} \pm 1 & \frac{d}{2} \mp 1 \\ \frac{d}{2} \mp 1 + r  & \frac{d}{2} \pm 1 + r \\ 0 & 0 \\ 0 & 0 \\ \end{pmatrix}_1  \begin{pmatrix}\frac{d}{2} & \frac{d}{2} \\ \frac{d}{2} + r & \frac{d}{2} + r \\ 0 & 0 \\ 0 & 0 \\ \end{pmatrix}_6$ & $\begin{pmatrix} 0 & 0 \\ 0 & 0 \\ -\frac{d}{2} \pm 1 & -\frac{d}{2} \mp 1 \\ -\frac{d}{2} \mp 1 - r & -\frac{d}{2} \pm1 - r \\ \end{pmatrix}_1 \begin{pmatrix} 0 & 0 \\ 0 & 0 \\ -\frac{d}{2} & -\frac{d}{2} \\ -\frac{d}{2} - r & -\frac{d}{2} - r \\ \end{pmatrix}_6$ \\
	&& \\
	\hline
	&& \\
	(even, even) & $\begin{pmatrix} \frac{d}{2} & \frac{d}{2} \\ \frac{d}{2} + r & \frac{d}{2} + r  \\ 0 & 0 \\ 0 & 0 \\ \end{pmatrix}_8$ & $\begin{pmatrix} 0 & 0 \\ 0 & 0 \\ -\frac{d}{2} & -\frac{d}{2} \\ -\frac{d}{2} - r & -\frac{d}{2} - r \\ \end{pmatrix}_8$  \\
	&& \\
	\hline
\end{tabular}%
}
\]
\[
\begin{tabular}{|l|c|c|}
	\hline
	$r > 0,\ d = 0$ & Tube contains $\OO_{\X}$ & Tube does not contain $\OO_{\X}$ \\
	\hline
	&& \\
	r odd & $\begin{pmatrix} 1 & 0 \\ r-1 & r+1 \\ 1 & 0 \\ 0 & 0 \\ \end{pmatrix}_1 \begin{pmatrix} 0 & 1 \\ r+1 & r-1 \\ 0 & 1 \\ 0 & 0 \\ \end{pmatrix}_1$  & $\begin{pmatrix} 0 & 0 \\ r & r \\ 0 & 0 \\ 0 & 0 \\ \end{pmatrix}_6$ \\
	&& \\
	\hline
	&& \\
	r even & $\begin{pmatrix} 1 & 0 \\ r & r \\ 0 & 1 \\ 0 & 0 \\ \end{pmatrix}_1$ $\begin{pmatrix} 0 & 1 \\ r & r \\ 1 & 0 \\ 0 & 0 \\ \end{pmatrix}_1$ & $\begin{pmatrix} 0 & 0 \\ r & r \\ 0 & 0 \\ 0 & 0 \\ \end{pmatrix}_6$ \\
	&& \\
	\hline
\end{tabular}
\]
The subscript counts the number of indecomposables satisfying $\F \otimes \omega_{\X} \ncong \F$  with given cohomology table.
\end{thm}
{\bf The proof strategy}. We will make use of Crawley-Boevey's generalisation of Kac's Theorem for weighted projective lines. By the previous proposition, indecomposables with $gcd(r,d)$ odd correspond to real roots of the associated root system, which are enumerated in a standard basis for $K_0$. Going through the list, one tabulates all triples $(rk(\F), deg(\F), \chi(\F))$ coming from real roots $[\F]$, and this triple completely determines $\beta(\F)$ in regions $\mathcal{R}_1, \mathcal{R}_3$. The region $\mathcal{R}_2$ is then dealt with by hand. We first recall the needed notions, following \cite{CB, Schiff}.

\subsection*{Kac's Theorem, after Schiffmann-Crawley-Boevey} We follow \cite{CB, Schiff} for the material that follows. Let $\X = \mathbb{P}^1({\bf p}, \pmb{\lambda})$ be a general weighted projective line for now, and let $T_{can} = \bigoplus_{\vec{0} \leq \vec{x} \leq \vec{c}} \OO_\X(\vec{x})$ be the canonical tilting object with endomorphism algebra $C(\bf{p}, \pmb{\lambda})$ $= kQ/I$ with quiver $Q$:
\[
	\xymatrix@R=10pt{ & \vec{x_1} \ar[r] & 2\vec{x_1} \ar[r] & \dots \ar[r] & (p_1-1)\vec{x_1} \ar[rdd] & \\
		& \vec{x_2} \ar[r] & 2\vec{x_2} \ar[r] & \dots \ar[r] & (p_2-1)\vec{x_2} \ar[rd] & \\
	 \vec{0} \ar[ruu] \ar[ru] \ar[rd] &\vdots & \vdots & \vdots & \vdots & \vec{c} \\
	 & \vec{x}_n \ar[r] & 2\vec{x}_n \ar[r] & \dots \ar[r] & (p_n -1)\vec{x}_n \ar[ru] & \\
	}
\]
Let $Q'$ be the tree subquiver corresponding to $T' = \bigoplus_{\vec{0} \leq \vec{x} < \vec{c}} \OO_\X(\vec{x})$, which we label differently as
\[
	\xymatrix@R=5pt{ & 1,1 \ar[r] & 1,2 \ar[r] & \dots \ar[r] & 1, p_1-1 \\
		& 2,1 \ar[r] & 2,2 \ar[r] & \dots \ar[r] & 2, p_2-1 \\
	 0 \ar[ruu] \ar[ru] \ar[rd] &\vdots&\vdots&\vdots&\vdots \\
	 & n,1 \ar[r] & n,2 \ar[r] & \dots \ar[r] & n, p_n -1 \\
	}
\]
Let $\mathfrak{g}$ be its associated Kac-Moody algebra with root system $\Gamma$, with simple roots $\varepsilon_0$, $\varepsilon_{i,j}$, and let $L\mathfrak{g} = \mathfrak{g}[t, t^{-1}]$ its loop algebra with root system $\widehat{\Gamma} = \Z \delta \oplus \Gamma$, with symmetric form extended by $(\delta, -) = 0$. The derived equivalence
	\[
		{\rm RHom}(T_{can}, -):	{\rm D}^b(\X) \xrightarrow{\cong} {\rm D}^b(C(\bf{p}, \pmb{\lambda}))
	\]
	sends $\{\OO_\X, S_{i, j}\}_{j \neq 0}$ to the simple modules $S(0), S(i,j)$ supported over $Q'$. This identifies the summand $\Z[\OO_\X] \oplus (\bigoplus_{i, j\neq 0} \Z[S_{i,j}])$ of $K_0(\X)$ with $\Gamma$, sending the symmetrised Euler form with the Weyl-invariant symmetric bilinear form on $\Gamma$. As Schiffmann then shows \cite{Schiff}, this extends to a full isomorphism $K_0(\X) \xrightarrow{\cong} \widehat{\Gamma}$ sending $[S_x]$ to $\delta$. The induced positive cone given by classes of coherent sheaves on $\widehat{\Gamma}$ is given by nonnegative combinations of
	\[
		\varepsilon_0,\ \varepsilon_0 + n\delta,\ \varepsilon_{i,j},\ \delta - \sum_{j \neq 0} \varepsilon_{i,j}, \quad n \in \Z
	\]
	with $[\OO_\X(n\vec{c})] \mapsto \varepsilon_0 + n\delta$ and $[S_{i,0}] \mapsto \delta - \sum_{j \neq 0} \varepsilon_{i,j}$. A version of Kac's Theorem then holds for coherent sheaves on $\X$, which we only state in a weak form:
	\begin{prop}[Crawley-Boevey, \cite{CB}] The isomorphism $K_0(\X) \xrightarrow{\cong} \widehat{\Gamma}$ induces a bijection between Grothendieck classes of indecomposable coherent sheaves and the positive roots of $L\mathfrak{g}$.
		\begin{enumerate}[1)]
			\item When $\beta$ is a positive real root, then there is a unique indecomposable $\F$ such that $[\F] \mapsto \beta$.
			\item When $\beta$ is a positive imaginary root, then there are infinitely many indecomposables $\F$ for which $[\F] \mapsto \beta$.
		\end{enumerate}
	\end{prop}
	The suspension $[1]$ acts on $K_0(\X) \cong K_0({\rm D}^b(\X))$ by $[\F[1]] = -[\F]$, and so this extends to a bijection between all (positive and negative) roots of $L\mathfrak{g}$ and orbits of indecomposables in $ {\rm D}^b(\X) $ under $\F \mapsto \F[2]$. Letting $\Delta$ be the set of roots of $\mathfrak{g}$, the roots of $L\mathfrak{g}$ are given by $\{\alpha + n\delta\ |\ \alpha \in \Delta,\ n \in \Z\}$, and the real roots are those of the form $\alpha + n\delta$ with $\alpha \in \Delta^{re}$.\bigskip

	Coming back to $\X = \mathbb{P}^1(2,2,2,2;\lambda)$, $\mathfrak{g}$ is of affine type $\widetilde{D}_4$ and we write $\varepsilon_{i}$ for $\varepsilon_{i,1}$. The positive real roots of $\widetilde{D}_4$ are given by solutions $\alpha = \sum_{i=0}^4 \alpha_i \varepsilon_i$ to
	\[
		q(\alpha) = \left(\alpha_0^2 + \alpha_1^2 + \dots + \alpha_4^2 \right) - \left(\alpha_0\alpha_1 + \dots + \alpha_0\alpha_4 \right) = 1
	\]
	with $\alpha_i \in \Z_{\geq 0}$. Writing $q(\alpha) = \sum_{i = 1}^4 \left(\alpha_i - \frac{1}{2}\alpha_0 \right)^2 \geq 0$, one sees that the solutions are as listed in Figure \ref{realrootsD4}. Alternatively, the positive real roots are given by the dimension vectors of indecomposables over the quiver $Q$ of type $\widetilde{D}_4$, except for those which are multiples of $(\alpha_0, \alpha_1, \alpha_2, \alpha_3, \alpha_4) = (2,1,1,1,1)$. The dimension vectors are written down in \cite[XIII.3]{SimSkow}. 
	\begin{figure}
\[
	\centerline{
\small{\begin{tabular}{|c|c|c|c|c||c|c|c|}
	\hline
	$\alpha_0$ & $\alpha_1$ & $\alpha_2$ & $\alpha_3$ & $\alpha_4$ &  $r$ & $d$ & $\chi$ \\
	\hline
	2m & m+1 & m & m & m &  2m & 4m+1+2n & 2m+n \\
	2m & m & m+1 & m & m &  2m & 4m+1+2n & 2m+n \\
	2m & m & m & m+1 & m &  2m & 4m+1+2n  & 2m+n \\
	2m & m & m & m & m+1 &  2m & 4m+1+2n  & 2m+n \\
	&&&&&&& \\
	2m+1 & m & m & m & m &  2m+1 & 4m+2n  & 2m+1+n  \\
	&&&&&&& \\
	2m+1 & m+1 & m & m & m  & 2m+1 & 4m+1+2n  & 2m+1+n  \\
	2m+1 & m & m+1 & m & m  & 2m+1 & 4m+1+2n  & 2m+1+n  \\
	2m+1 & m & m & m+1 & m  & 2m+1 & 4m+1+2n  & 2m+1+n  \\
	2m+1 & m & m & m & m+1  & 2m+1 & 4m+1+2n  & 2m+1+n  \\
	&&&&&&& \\
	2m+1 & m+1 & m+1 & m & m &  2m+1 & 4m+2+2n  & 2m+1+n  \\
	2m+1 & m & m & m+1 & m+1 &  2m+1 & 4m+2+2n  & 2m+1+n  \\
	2m+1 & m+1 & m & m+1 & m &  2m+1 & 4m+2+2n  & 2m+1+n  \\
	2m+1 & m & m+1 & m & m+1 &  2m+1 & 4m+2+2n  & 2m+1+n  \\
	2m+1 & m+1 & m & m & m+1 &  2m+1 & 4m+2+2n  & 2m+1+n  \\
	2m+1 & m & m+1 & m+1 & m &  2m+1 & 4m+2+2n  & 2m+1+n  \\
	&&&&&&& \\
	2m+1 & m & m+1 & m+1 & m+1  & 2m+1 & 4m+3+2n  & 2m+1+n  \\
	2m+1 & m+1 & m & m+1 & m+1  & 2m+1 & 4m+3+2n  & 2m+1+n  \\
	2m+1 & m+1 & m+1 & m & m+1  & 2m+1 & 4m+3+2n  & 2m+1+n  \\
	2m+1 & m+1 & m+1 & m+1 & m  & 2m+1 & 4m+3+2n  & 2m+1+n  \\
	&&&&&&& \\
	2m+1 & m+1 & m+1 & m+1 & m+1 &  2m+1 & 4m+4+2n  & 2m+1+n  \\
	&&&&&&& \\
	2m+2 & m & m+1 & m+1 & m+1 &  2m+2 & 4m+3+2n  & 2m+2+n   \\
	2m+2 & m+1 & m & m+1 & m+1 &  2m+2 & 4m+3+2n  & 2m+2+n   \\
	2m+2 & m+1 & m+1 & m & m+1 &  2m+2 & 4m+3+2n  & 2m+2+n  \\
	2m+2 & m+1 & m+1 & m+1 & m &  2m+2 & 4m+3+2n  & 2m+2+n \\
	\hline
	\end{tabular}}}
\]
\caption{Positive real roots of $\widetilde{D}_4$ and triples $(r, d, \chi)$, where $m \geq 0,\ n \in \Z$.}\label{realrootsD4}
\end{figure}

The last three columns record $(rk(\F), deg(\F), \chi(\F))$ where $\F$ is any indecomposable corresponding to $\beta = \alpha + n\delta$. To see this, note that we must have $[\F] = \alpha_0 [\OO_\X] + \sum_{i=1}^4 \alpha_i [S_{i,1}] + n[S_x]$ and that $\chi(\OO_\X) = 1$, $\chi(S_{i, 1}) = 0$ and $\chi(S_x) = 1$. A general real root then has the form $\beta = \pm \alpha + n \delta$, with $\alpha$ in the above table and $n \in \Z$. Finally, let us record a lemma before moving on to the proof of Theorem \ref{ranktwothm}. 
\begin{lem}\label{chi} Let $\F \in {\rm D}^b(\X)$ be an indecomposable complex with $[\F] \mapsto \beta$ real, with $r \geq 0,\ d \in \Z$. Then the possible values of $\chi(\F)$ only depend on $d$ and are listed below:
	\[
		\begin{tabular}{|c|c|c|}
			\hline
			$d$ & $\chi$ \\
			\hline
			& \\
			$d$ odd & $\left(\frac{d-1}{2}\right)_4,\ \left(\frac{d+1}{2}\right)_4$ \\
			& \\
			\hline
			& \\
			$d$ even & $\left(\frac{d}{2}-1\right)_1,\ \left(\frac{d}{2}\right)_6,\ \left(\frac{d}{2}+1\right)_1$ \\
			& \\
			\hline
		\end{tabular}
	\]
	The subscript indicates how many times $\chi$ appears for fixed $(r,d)$.
\end{lem}
\begin{proof}  This follows by inspection of Figure \ref{realrootsD4}, where the case $r > 0$ corresponds to $\beta = \alpha + n\delta$ for $\alpha$ in the table, and $r = 0$ uses $\pm \alpha$ for $\alpha$ in the first four rows.
\end{proof}
We can now compute cohomology tables of indecomposables in rank two tubes
\[
	\beta(\F) = \begin{pmatrix*}[l] h^0(\F) & h^0(\F \otimes \omega_{\X}) \\ h^0(\F(\vec{c}) \otimes \omega_{\X}) & h^0(\F(\vec{c})) \\ h^1(\F\otimes \omega_{\X}) & h^1(\F) \\ h^1(\F(\vec{c})) & h^1(\F(\vec{c}) \otimes \omega_{\X}) \\   \end{pmatrix*}
\]
\begin{proof} Let $\F$ be indecomposable with $\F \otimes \omega_{\X} \ncong \F$, of type $(r,d)$. First assume that $(r,d)$ is in region $\mathcal{R}_1$, the case $\mathcal{R}_3$ being similar. As previously, the bottom half of $\beta(\F)$ vanishes for slope reasons (Lemma \ref{h-vanishing}). When $\gcd(r,d)$ is even then by Proposition \ref{ranktwoprop} the table $\beta(\F)$ is symmetrical under exchanging columns, and therefore is forced to be as in Theorem \ref{rankone}. When $\gcd(r,d)$ is odd, then by \ref{ranktwoprop} and Kac's Theorem the class $[\F]$ must correspond to a real root, and so the possible values of $\chi$ are listed in Lemma \ref{chi}, which determines the possible values of $h^0(\F)$. By Riemann-Roch we have
	\[
		h^0(\F) + h^0(\F \otimes \omega_{\X}) = d
	\]
	\[
		h^0(\F(\vec{c})) + h^0(\F(\vec{c}) \otimes \omega_{\X}) = d+2r
	\]
	Now, keeping in mind that $\beta(\F) = \beta(M)$ where $M$ is presented by a matrix factorisation, the sum of each column must be equal. From this one sees that the tables must be of the form
\[
\centerline{
	\begin{tabular}{|c|c|}
	\hline
	$(r,d)$ & $r \geq 0,\ d > 0$ \\
	\hline
	& \\
$d$ odd & $\begin{pmatrix} \frac{d \pm 1}{2} & \frac{d \mp 1}{2} \\ \frac{d \mp 1}{2} + r & \frac{d \pm 1}{2} + r \\ 0 & 0 \\ 0 & 0 \\ \end{pmatrix}_4$ \\
	& \\
	\hline
	& \\
	$d$ even & $\begin{pmatrix} \frac{d}{2} \pm 1 & \frac{d}{2} \mp 1 \\ \frac{d}{2} \mp 1 + r  & \frac{d}{2} \pm 1 + r \\ 0 & 0 \\ 0 & 0 \\ \end{pmatrix}_{1}  \begin{pmatrix}\frac{d}{2} & \frac{d}{2} \\ \frac{d}{2} + r & \frac{d}{2} + r \\ 0 & 0 \\ 0 & 0 \\ \end{pmatrix}_6$ \\
	& \\
	\hline
\end{tabular}}
\]

The same argument determines tables in region $\mathcal{R}_3$, so we are left with region $\mathcal{R}_2$ where $r > 0,\ d = 0$. We first note that the base of rank two tubes in $\mathcal{C}_{0}$ consists of line bundles of degree zero, with one distinguished tube having $\{\OO_\X, \omega_{\X}\}$ as base. For the other tubes, an appeal to the autoequivalence $\Phi_{0, \infty}: \mathcal{C}_0 \xrightarrow{\cong} \mathcal{C}_\infty$ shows that ${\rm Ext}^*(\OO_\X, \F) = 0 = {\rm Ext}^*(\OO_\X, \F \otimes \omega_{\X})$. The first row and third row of $\beta(\F)$ vanishes, and slope considerations show vanishing of the fourth row. The table must then be
\[
	\begin{pmatrix*}[l] h^0(\F) & h^0(\F \otimes \omega_{\X}) \\ h^0(\F(\vec{c}) \otimes \omega_{\X}) & h^0(\F(\vec{c})) \\ h^1(\F\otimes \omega_{\X}) & h^1(\F) \\ h^1(\F(\vec{c})) & h^1(\F(\vec{c}) \otimes \omega_{\X}) \\   \end{pmatrix*} = \begin{pmatrix} 0 & 0 \\ r & r \\ 0 & 0 \\ 0 & 0 \end{pmatrix}
\]
We are down to $\F$ in the Auslander-Reiten component of $\{\OO_\X, \omega_{\X}\}$. Now, $\F$ is uniserial with socle either $\OO_\X$ or $\omega_\X$. Assume the first. From the structure of a rank two tube, the simple top of $\F$ is $\omega_{\X}$ when $r$ is even, and $\OO_\X$ for $r$ odd. This determines the dimensions of ${\rm Hom}(\OO_\X, \F),\ {\rm Hom}(\omega_{\X}, \F),\ {\rm Hom}(\F, \OO_\X),\ {\rm Hom}(\F, \omega_{\X})$ as
\begin{align*}
	&{\rm dim}_k \ {\rm Hom}(\OO_\X, \F) = 1\\
	&{\rm dim}_k \ {\rm Hom}(\omega_\X, \F) = 0\\
	&{\rm dim}_k \ {\rm Hom}(\F, \OO_\X) = \begin{cases} 1 & \textup{r odd} \\ 0 & \textup{r even} \end{cases} \\
	&{\rm dim}_k \ {\rm Hom}(\F, \omega_\X) = \begin{cases} 0 & \textup{r odd} \\ 1 & \textup{r even} \end{cases}
\end{align*}
and from Serre duality one deduces the shape of the first and third rows. This is enough to determine the tables as
\[
\centerline{
	\begin{tabular}{|c|c|}
	\hline
	r odd & r even\\
	\hline
	& \\
	$\begin{pmatrix} 1 & 0 \\ r-1 & r+1 \\ 1 & 0 \\ 0 & 0 \\ \end{pmatrix}$ &
	$\begin{pmatrix} 1 & 0 \\ r & r \\ 0 & 1 \\ 0 & 0 \\ \end{pmatrix}$ \\
	& \\
	\hline
\end{tabular}}
\]
The case of $\F$ with socle $\omega_{\X}$ is then given by the mirrored table. \end{proof}

\subsection{Betti tables of indecomposables}
Finally, we collect and list the Betti tables of indecomposable MCM modules in the standard format
\[
	\begin{tabular}{c|c c }
		  & 0 & 1 \\
		\hline
		0 & $\beta_{0,0}$&$\beta_{1,1}$  \\
		1 & $\beta_{0,1}$&$\beta_{1,2}$  \\
		2 & $\beta_{0,2}$&$\beta_{1,3}$  \\
		3 & $\beta_{0,3}$&$\beta_{1,4}$  \\
		4 & $\beta_{0,4}$&$\beta_{1,5}$  \\
		5 & \vdots & \vdots \\
	\end{tabular}
\]
One can then obtain the complete Betti table by extending by $2$-periodicity via $\beta_{i,j} = \beta_{i+2, j+4}$.

In the previous sections, we produced the cohomology tables $\beta(F_M)$ of coherent sheaves with $(r,d)$ in a chosen appropriate fundamental domain for the action of $M \mapsto M(1)$. In what follows, we will use a slightly different fundamental domain, better adapted to displaying the Betti tables $\beta(M)$. Call an indecomposable $M$ of the first kind if it belongs to the same Auslander-Reiten component as some $\Omega^{n}k^{st}(m)$, and of the second kind otherwise. 
\begin{cor}\label{bettifirstkind} The indecomposables of the first kind are uniquely determined by their Betti table. Up to translation, these are all tables of the form
\begin{align*}
&	\begin{tabular}{c|c c }
		  & 0 & 1 \\
		\hline
		0 & 1&-  \\
		1 & r-1&-  \\
		2 & 1&r+1  \\
		3 & -&-  \\
	\end{tabular}
\quad
	\begin{tabular}{c|c c }
		  & 0 & 1 \\
		\hline
		0 & r+1&1  \\
		1 & -&r-1 \\
		2 & -&1  \\
		3 & -&-  \\
	\end{tabular}
	\quad \textup{for $r > 0$ odd}, \\
	& \\
&	\begin{tabular}{c|c c }
		  & 0 & 1 \\
		\hline
		0 & 1&-  \\
		1 & r&- \\
		2 & -&r  \\
		3 & -&1  \\
	\end{tabular}
\quad\quad\quad
	\begin{tabular}{c|c c }
		  & 0 & 1 \\
		\hline
		0 & r&-  \\
		1 & 1&- \\
		2 & -&1  \\
		3 & -&r  \\
	\end{tabular}
	\quad\quad \ \  \textup{for $r > 0$ even}.
\end{align*}
The indecomposables with the above tables have degree $0$ and rank $r$.
\end{cor}
\begin{proof} Assuming the hypothesis, $\Omega^{-n}M(-m-2)$ is in the same Auslander-Reiten component as $k^{st}(-2)$ which corresponds to $\OO_{\X}$, then apply theorem \ref{ranktwothm} and translate the resulting tables in the above form. 
\end{proof}

Most indecomposables are of the second kind.
\begin{cor}\label{bettisecondkind} Up to translation, the Betti tables of indecomposables of the second kind are all tables of type ${\rm I-V}$:
\[
\resizebox{\columnwidth}{!}{%
\begin{tabular}{c|c c }
		I  & 0 & 1 \\
		\hline
		0 & a&-  \\
		1 & b&a  \\
		2 & -&b  \\
		3 & -&-  \\
	\end{tabular}
\quad
\quad
	\begin{tabular}{c|c c }
		II  & 0 & 1 \\
		\hline
		0 & a+1&0  \\
		1 & b&a  \\
		2 & 0&b+1  \\
		3 & 0&0  \\
	\end{tabular}
\quad
\quad
	\begin{tabular}{c|c c }
		III  & 0 & 1 \\
		\hline
		0 & a&0  \\
		1 & b+1&a+1  \\
		2 & 0&b  \\
		3 & 0&0  \\
	\end{tabular}
\quad
\quad
	\begin{tabular}{c|c c }
		IV  & 0 & 1 \\
		\hline
		0 & a+2&0  \\
		1 & b&a  \\
		2 & 0&b+2  \\
		3 & 0&0  \\
	\end{tabular}
\quad
\quad
	\begin{tabular}{c|c c }
		V  & 0 & 1 \\
		\hline
		0 & a&0  \\
		1 & b+2&a+2  \\
		2 & 0&b  \\
		3 & 0&0  \\
\end{tabular}%
}
\]
with $a,b \geq 0$, where we have $b \neq 0$ for tables of type ${\rm I}$ and $b-a$ odd for tables of type ${\rm IV-V}$. Here the degree is given by $d = \beta_{0,0} + \beta_{1,2} = 2a,\ 2a+1,\ 2a+2$ and the rank by $r = b-a$.
\end{cor}
\begin{proof} Let $M$ be indecomposable of the second kind. We claim that, up to translation, $\beta(M)$ can be put in the form
	\[
		\beta(M) = \begin{pmatrix} \beta_{0,0} & \beta_{1,2} \\ \beta_{0,1} & \beta_{1,3} \\ \beta_{0,2} & \beta_{1,4} \\ \beta_{0,3} & \beta_{1,5} \\ \end{pmatrix} = \begin{pmatrix} \alpha & \beta \\ \beta + r & \alpha + r \\ 0 & 0 \\ 0 & 0 \end{pmatrix}
	\]
	for some $\alpha, \beta$ and $r \in \Z$. Running over possibilities in Thm. \ref{rankone}, \ref{ranktwothm}, this is already the case in regions $\mathcal{R}_1, \mathcal{R}_2$, where $d, r \geq 0$. For $(rk(M),\ deg(M))$ belonging to $\mathcal{R}_3$, applying $M \mapsto M(2)$ will put $\beta(M)$ in the above form. Note that this sends $(r,d) \mapsto (-r, -d)$, and so the above tables coming from region $(3)$ will have $r < 0$. This will change our fundamental domain to $r > -\frac{d}{2},\ d \geq 0$:
\[
       \begin{tikzpicture}
	    \draw[gray!50, thin, step=0.3] (-2,-2) grid (2,2);
	    \draw[very thick,->] (-2,0) -- (2.2,0) node[right] {$r$};
	    \draw[very thick,->] (0,-2) -- (0,2.2) node[above] {$d$};

	    \foreach \x in {-2,...,2} \draw (\x,0.03) -- (\x,-0.03) node[below] {\tiny\x};
	    \foreach \y in {-2,...,2} \draw (-0.03,\y) -- (0.03,\y) node[right] {\tiny\y};

	    \fill[blue!50!cyan,opacity=0.3] (2,0) -- (0,2) -- (2,2);
	    \fill[blue!50!cyan,opacity=0.3] (0,0) -- (0,2) -- (2,0);
	    \fill[blue!50!cyan,opacity=0.3] (0,0) -- (-1,2) -- (0,2);

	    \draw (-1,2) -- node[below,sloped] {} (1,-2);
       \end{tikzpicture}
\]
Now, the case $\alpha = \beta = 0$ corresponds to $d = 0$, or tables in region $(2)$. The other two regions run over the same pairs $(\alpha, \beta)$, with only difference whether $r \geq 0$ or $r < 0$. Next, set $a = \alpha$, $b = \alpha + r$. Running over the possible $\alpha$ in Thm. \ref{rankone}, \ref{ranktwothm} shows that tables must have shapes ${\rm I-V}$. In particular in type ${\rm I}$, note that $b = 0$ implies $r = -\frac{d}{2}$ which falls outside of our domain. Lastly, fixing the type ${\rm I-V}$ of a table, note that $(a,b)$ and $(r,d)$ uniquely determine each other via $r = b-a$ and $d = 2a,\ 2a+1,\ 2a+2$, and so the classification is complete.
\end{proof}

%% file: chapters/fano/fano.tex
\chapter{Fano algebras, higher preprojective algebras and Artin-Schelter regular algebras} 
\chaptermark{Fano, higher preprojective and Artin-Schelter regular algebras}
\input{chapters/fano/preproj}
\input{chapters/fano/fanofromkoszul}

%% file: chapters/fano/preproj.tex
\section{Finite dimensional Fano algebras and noncommutative projective geometry}

This part of the thesis concerns an application of singularity categories to the representation theory of finite dimensional algebras and to noncommutative projective geometry in the sense of Artin-Zhang, which leads naturally to homological results on Koszul duality for Gorenstein algebras. We begin with some motivation and by reviewing basic definitions and examples. Let $k$ be a field throughout, which we assume algebraically closed in the first section for simplicity. 

Given a smooth projective $k$-variety $X$ with tilting complex $\E$, the endomorphism algebra $\Lambda = {\rm End}(\E)$ is finite dimensional of finite global dimension, and one has a derived equivalence 
\[
	{\rm D}^b(X) \cong {\rm D}^b(\Lambda).
\]The finer properties of the (anti-)canonical bundle of $X$ exert some control over the representation theory of $\Lambda$, through the identification of Serre functors 
\[
	\Se_X = - \otimes_{\mathcal{O}_X} \omega_{X}[{\rm dim}\ X] \mapsto \Se_{\Lambda} = - \otimes_{\Lambda}^{\mathbb{L}} \omega_{\Lambda}[{\rm gldim}\ \Lambda]
\]
for a suitable canonical complex $\omega_{\Lambda} = D\Lambda [-{\rm gldim}\ \Lambda] \in {\rm D}^{b}(\Lambda^{ev})$, where $D\Lambda = {\rm Hom}_k(\Lambda, k)$. For a general finite dimensional algebra $\Lambda$ of finite global dimension, the notion of ampleness of the anti-canonical complex $$\omega_{\Lambda}^{-1} = {\rm RHom}_{\Lambda}(D\Lambda, \Lambda)[{\rm gldim}\ \Lambda]$$ was introduced and studied in depth by Minamoto \cite{Mina}. This leads to the notion of a Fano algebra $\Lambda$, as well as a suitable ``anticanonical ring'' 
\[
	\Pi(\Lambda) = \bigoplus_{m \geq 0}{\rm Hom}_{ {\rm D}^b(\Lambda)}(\Lambda, \omega_{\Lambda}^{-m})
\]
where $\omega^{-m}_{\Lambda} = (\omega_{\Lambda}^{-1})^{\otimes_{\Lambda}^{\mathbb{L}} m}$, mimicking the classical anticanonical ring
\[
	R(X) = \bigoplus\limits_{m \geq 0}{\rm H}^0(X, \omega_X^{-m}) = \bigoplus\limits_{m \geq 0}{\rm Hom}_{ {\rm D}^b(X) }(\OO_X, \omega_X^{-m}).
\]

Now fix $\Lambda$ a finite dimensional algebra of finite global dimension $\gldim \Lambda = d < \infty$. Let ${\rm Pic}(\Lambda)$ be the group of isomorphism classes of invertible complexes of bimodules $L \in {\rm D}^b(\Lambda^{ev}) $ under the derived tensor product $- \otimes^{\mathbb{L}}_{\Lambda}-$, thought of as the group ${\rm Pic}(X)$ of line bundles $\mathcal{L}$ on a projective variety $X$. Recall that the Serre criterion for ampleness states that $\mathcal{L}$ is ample if and only if for all coherent sheaves $\F$
\[
	{\rm H}^s(X, \F \otimes \mathcal{L}^{m}) = 0 \textup{ for } s > 0 \textup{ for all } m \gg 0.
\]
Let $({\rm D}^{\leq 0}, {\rm D}^{\geq 0})$ be the standard t-structure on $ {\rm D}^b(X) $, and consider the pair of full subcategories
\begin{align*}
	{\rm D}^{\leq 0, \mathcal{L}} &= \{ \F \in {\rm D}^b(X) \ | \ {\rm H}^s(X, \F \otimes \mathcal{L}^{m}) = 0 \textup{ for } s > 0 \textup{ for all } m \gg 0\}\\
	{\rm D}^{\geq 0, \mathcal{L}} &= \{ \F \in {\rm D}^b(X) \ | \ {\rm H}^s(X, \F \otimes \mathcal{L}^{m}) = 0 \textup{ for } s < 0 \textup{ for all } m \gg 0\}
\end{align*}
The hypercohomology spectral sequence shows that $\mathcal{L}$ is ample if and only if the above gives rise to the standard t-structure on $ {\rm D}^b(X) $. By analogy, given $ L \in {\rm Pic}(\Lambda) $, writing $L^m := L^{\otimes_\Lambda^{\mathbb{L}} m}$ we consider the pair of full subcategories ${\rm D}^L := ({\rm D}^{\leq 0, L}, {\rm D}^{\geq 0, L})$ 
\begin{align*}
	{\rm D}^{\leq 0, L} &= \{ M \in {\rm D}^b(\Lambda) \ | \ {\rm H}^s(M \otimes^{\mathbb{L}}_{\Lambda}L^{m}) = 0 \textup{ for } s > 0 \textup{ for all } m \gg 0\}\\
	{\rm D}^{\geq 0, L} &= \{ M \in {\rm D}^b(\Lambda) \ | \ {\rm H}^s(M \otimes^{\mathbb{L}}_{\Lambda}L^{m}) = 0 \textup{ for } s < 0 \textup{ for all } m \gg 0\}
\end{align*}
obtained by formally substituting $\OO_X$ for $\Lambda$, since ${\rm H}^s(X, \F) = {\rm Hom}_{{\rm D}^b(X) }(\OO_X, \F[s])$ and ${\rm H}^s(M) = {\rm Hom}_{{\rm D}^b(\Lambda) }(\Lambda, M[s])$. Next, we say that a complex $M \in {\rm D}^b(\Lambda)$ is \emph{pure} if it is in the essential image of $\modsf \Lambda \hookrightarrow {\rm D}^b(\Lambda)$, or equivalently if ${\rm H}^s(M) = 0$ for $s \neq 0$.
\begin{defn}[Minamoto \cite{Mina}] Let $L \in {\rm Pic}(\Lambda) $ be an invertible complex of bimodules. 
	\begin{enumerate}[i)]
	\item We say that $L$ is \emph{almost} ample (\emph{almost} very ample) if $L^m$ is pure for all $m \gg 0$ ($m \geq 0$).
	\item We say that $L$ is ample (very ample) if $L$ is almost ample (almost very ample) and furthermore ${\rm D}^L$ forms a t-structure for $ {\rm D}^b(\Lambda) $.
\end{enumerate}
\begin{rem}
Our terminology differs slightly from \cite{Mina}, where only ii) is used. Our notion of `very ample' also corresponds to `extremely ample' in \cite[Defn. 3.4]{Mina}, since Minamoto reserves the adjective `very' for an intermediate condition. The adjective `almost' also replaces `quasi' in \cite{MM}. As we will not make extended usage of these notions, the change in terminology should not prove confusing. 
\end{rem}
\end{defn}

When $L$ is ample, note that the t-structure ${\rm D}^L$ on $ {\rm D}^b(\Lambda) $ can be very different than the standard one, and the heart $\mathcal{H}^L$ is of considerable interest (see \cite[Sect. 5]{Mina}). The terminology should feel natural by the following lemma: 
\begin{lem}[{\cite[Lemma 3.11]{Mina}}] Let $L \in {\rm Pic}(\Lambda)$ be (almost) ample. Then $L^m$ is (almost) very ample for all $m \gg 0$.
\end{lem}

\begin{defn}[{\cite[Def. 4.1]{Mina}}] Let $\Lambda$ be a finite dimensional algebra of finite global dimension $\gldim \Lambda = d < \infty$, and consider the anti-canonical complex $\omega_\Lambda^{-1} \in {\rm Pic}(\Lambda)$
	\[
		\omega_\Lambda^{-1} = {\rm RHom}_\Lambda(D\Lambda, \Lambda)[d].
	\]
	We say that $\Lambda$ is (almost) Fano if $\omega^{-1}_\Lambda$ is (almost) ample.
\end{defn}
\begin{defn} Let $\Lambda$ be an almost Fano algebra. The $\N$-graded algebra 
\[
		\Pi(\Lambda) := \bigoplus\limits_{m \geq 0}{\rm Hom}_{{\rm D}^b(\Lambda)}(\Lambda, \omega^{-m}_\Lambda) = \bigoplus\limits_{m \geq 0}{\rm H}^0(\omega^{-m}_\Lambda)
\]
is called the \emph{higher preprojective algebra} of $\Lambda$.
\end{defn}

Consider the $d$-shifted Serre functor $\Se_d := \Se_\Lambda \circ [-d]$ on $ {\rm D}^b(\Lambda) $ and its inverse
\begin{align*}
\Se_d(M) &= M\otimes^{\mathbb{L}}_\Lambda \omega_\Lambda = M \otimes^{\mathbb{L}}_\Lambda D\Lambda[-d]\\
\Se_d^{-1}(M) &= M \otimes^{\mathbb{L}}_\Lambda \omega_\Lambda^{-1} = M \otimes^{\mathbb{L}}_\Lambda {\rm RHom}_\Lambda(D\Lambda, \Lambda)[d] \cong {\rm RHom}_\Lambda(D\Lambda, M)[d].
\end{align*}
Then $\Lambda$ is almost Fano if and only if $\Se_d^{-m}(\Lambda)$ is pure for all $m \gg 0$, that is 
\[
	\Se_d^{-m}(\Lambda) \in \modsf \Lambda \subseteq {\rm D}^b(\Lambda) \textup{ for all } m \gg 0.
\]
The stronger condition that $\omega_\Lambda^{-1}$ be almost very ample is equivalent to 
\[
	\Se_d^{-m}(\Lambda) \in \modsf \Lambda \subseteq {\rm D}^b(\Lambda) \textup{ for all } m \geq 0.
\]
One can picture this as the infinite sequence of conditions
\begin{align*}
	&{\rm Ext}^i_{\Lambda}(D\Lambda, \Lambda) = 0 &\textup{ for all } i < d \\
	&{\rm Ext}^i_{\Lambda}(D\Lambda, {\rm Ext}^d_\Lambda(D\Lambda, \Lambda)) = 0 &\textup{ for all } i < d \\
	&{\rm Ext}^i_{\Lambda}(D\Lambda, {\rm Ext}^{d}_\Lambda(D\Lambda, {\rm Ext}^d_\Lambda(D\Lambda, \Lambda))) = 0 &\textup{ for all } i < d
\end{align*}
\[
\vdots
\]
Since ${\rm RHom}_\Lambda(D\Lambda, M) \cong M \otimes^{\mathbb{L}}_\Lambda {\rm RHom}_\Lambda(D\Lambda, \Lambda)$ for all $M \in {\rm D}^b(\Lambda)$, these can be equivalently phrased as the Tor-vanishing conditions
\begin{align*}
	&{\rm Ext}^i_{\Lambda}(D\Lambda, \Lambda) = 0 &\textup{ for all } i < d \\
	&{\rm Tor}_{d-i}^{\Lambda}({\rm Ext}^d_\Lambda(D\Lambda, \Lambda), {\rm Ext}^d_\Lambda(D\Lambda, \Lambda)) = 0 &\textup{ for all } i < d \\
	&{\rm Tor}_{d-i}^{\Lambda}({\rm Ext}^d_\Lambda(D\Lambda, \Lambda)^{\otimes_\Lambda 2}, {\rm Ext}^d_\Lambda(D\Lambda, \Lambda)) = 0 &\textup{ for all } i < d 
\end{align*}
\[
\vdots
\]

\begin{exmp}[$d=1$] Let $\Lambda$ be a finite dimensional hereditary algebra. Let $M$ be an indecomposable $\Lambda$-module. Then 
	\[
		\Se_1^{-1}(M)=\tau^{-1}(M) = {\rm RHom}_\Lambda(D\Lambda, M)[1] \cong \begin{cases} {\rm Ext}^1_{\Lambda}(D\Lambda, M) & M \textup{ is not injective} \\ {\rm Hom}_\Lambda(D\Lambda, M)[1] & M \textup{ is injective}  \end{cases} 
	\]	
	since ${\rm Hom}_\Lambda(D\Lambda, M) = 0$ for $M$ with no injective summands. It follows that $\omega_\Lambda^{-1}$ is almost very ample if and only if for every indecomposable projective $P$, the modules 
	\[
		P, \tau^{-1}P, \tau^{-2}P, \tau^{-3}P, \dots, \tau^{-m}P, \dots 
	\]
	are never injective. Equivalently, $\Lambda$ is of infinite representation type. We will see later that $\omega^{-1}_\Lambda$ is actually very ample in this case \cite[Prop. 5.1]{Mina}. The algebra 
	\[
		\Pi(\Lambda) = \bigoplus\limits_{m \geq 0} {\rm Hom}_{\Lambda}(\Lambda, \tau^{-m}\Lambda)
	\]
	is the classical preprojective algebra of Gel'fand-Ponomarev and Baer-Geigle-Lenzing \cite{BGL}.\qed
\end{exmp}

On the other hand, when $\Lambda$ is hereditary of finite representation type the Serre functor $\Se_\Lambda$ satisfies the \emph{fractionally Calabi-Yau} property $\Se_\Lambda^n \cong [m]$ for some $n, m \in \Z$ (see \cite[Thm. 5.1]{Mina}, \cite[Ex. 8.3(2)]{Ke}), and so $\omega_\Lambda$ behaves like a torsion line bundle in this case. \bigskip

As part of a larger program to study algebras of higher global dimension with `hereditary behavior', Herschend, Iyama and Oppermann \cite{HIO} have introduced the class of `$d$-hereditary` algebras, which break down into `$d$-representation finite' and `$d$-representation infinite' algebras. We will be interested mainly in the latter. 
\begin{defn}[Herschend-Iyama-Oppermann \cite{HIO}] A finite dimensional algebra $\Lambda$ with $\gldim \Lambda = d < \infty$ is called $d$-representation infinite if $\omega_\Lambda^{-1}$ is almost very ample. We will also call $\Lambda$ higher representation infinite when $d$ is implicit. 
\end{defn}
These notions are closely related to the study of geometric helices by Bondal-Polishchuk \cite{BP} and Bridgeland-Stern \cite{BS}. 

\begin{exmp}\label{helixexample} Let $X$ be a smooth Fano variety with a geometric helix of sheaves $(E_i)_{i \in \Z}$ of period $n$, meaning that for all $i \in \Z$:
\begin{enumerate}[i)]
\item $(E_{i}, \dots, E_{i+n-1})$ is a full exceptional collection;
\item $E_{i+n} = E_{i} \otimes \omega_X^{-1}$;
\item for all $i \leq j$ we have ${\rm Hom}(E_i, E_j[s]) = 0 \textup{ whenever } s \neq 0$.
	\end{enumerate}	
That last condition implies that $(E_{i+1}, \dots, E_{i+n})$ is strong, and let $\E = \bigoplus_{j=1}^n E_{i+j}$ be the tilting sheaf with $\Lambda = {\rm End}(\E)$. If $\gldim \Lambda = \dim X$, then $\Lambda$ is higher representation infinite by condition $ii)$ and $iii)$. Its higher preprojective algebra $\Pi(\Lambda)$ is the rolled-up helix algebra of \cite{BS}. Moreover, letting $\pi: Y \to X$ be the total space of the canonical bundle $\omega_X$, the sheaf $\pi^*\E$ pulls back to a tilting sheaf on $Y$ and we additionally have $\Pi(\Lambda) = {\rm End}(\pi^* \E)$ (see \cite{BS}). \qed
\end{exmp} 

Now let $X$ be a Fano variety. Recall that when $\omega_X^{-1}$ is very ample, the anticanonical ring $R(X)$ is a Noetherian graded algebra, generated in degree one over $ R(X)_0 = {\rm H}^0(X, \OO_X)$ and of Krull dimension $\dim R(X) = \dim X + 1$. Moreover one recovers $X = \proj R(X)$, as well as the category of (quasi-)coherent sheaves $\coh X = \qgr R(X)$ ($\qcoh X = \QGr R(X)$) as the Serre quotients
\begin{align*}
	\qgr R(X) &= \frac{\grmod R(X)}{\mathsf{grmod}_0 R(X)} &
	\QGr R(X) &= \frac{\Grmod R(X)}{\mathsf{Grmod}_0 R(X)}
\end{align*}
where $\mathsf{grmod}_0 R(X)$ and $\mathsf{Grmod}_0 R(X)$ denote the subcategories of right bounded (or Artinian) modules $\{ M \ | \ M_{\geq n} = 0,\ n \gg 0\}$. The situation is more complicated for the higher preprojective algebra $\Pi(\Lambda)$ of an almost Fano algebra, but some analogous results hold.
\begin{lem} Assume that $\omega_\Lambda^{-1}$ is almost very ample. Then
\begin{align*}
		\Pi(\Lambda) = \bigoplus\limits_{m \geq 0} {\rm H}^0(\omega^{-m}_\Lambda) &\cong \bigoplus\limits_{m \geq 0} {\rm H}^0(\omega^{-1}_\Lambda)^{\otimes_\Lambda m}
= {\rm T}_{\Lambda}({\rm Ext}^d_\Lambda(D\Lambda, \Lambda)).
	\end{align*}
	Hence $\Pi(\Lambda)$ is finitely generated in degree one over $\Pi(\Lambda)_0 = \Lambda$. Moreover we have $\gldim \Pi(\Lambda) = \gldim \Lambda + 1$.
\end{lem}
\begin{proof}
	The first statement is a consequence of the Tor-vanishing conditions and the K\"unneth Theorem, see \cite{Mina}. For the global dimension claim see \cite{HIO}.
\end{proof}
In contrast, higher preprojective algebras are rarely Noetherian. In the classical case of an hereditary algebra of infinite representation type, the preprojective algebra $\Pi(\Lambda)$ is Noetherian if and only if $\Lambda$ is of tame representation type, if and only if $\Lambda$ is Morita equivalent to $kQ$ with $Q$ a quiver of extended Dynkin type. This is in keeping with expectations from noncommutative projective geometry where the graded algebras arising through natural constructions are rarely Noetherian. According to Polishchuk's \cite{Pol}, the more natural condition is that of \emph{coherence}. We follow \cite{Pol} and \cite[Sect. 2.1]{Mina} from here. Let $B = B_0 \oplus B_1 \oplus \dots $ be a locally finite graded $k$-algebra, and recall that modules are implicitly taken to be right modules. 
\begin{defn} Let $M$ a graded $B$-module. We say that $M$ is coherent if:
	\begin{enumerate}[a)]
	\item $M$ is a finitely generated graded $B$-module. 
	\item Every map $f: P \to M$ from a finitely generated projective graded $B$-module has finitely generated kernel. 
\end{enumerate}
We say that $B$ is graded  coherent if $B$ and $B/B_{\geq 1}$ are coherent graded modules.
\end{defn} 
\begin{rem}
	There is the analogous notion of coherence for the underlying ungraded algebra of $B$, and some authors sometimes denote the above notion by `graded coherence'. For us, coherence of a graded algebra will implicitly refer to the graded notion, and so there should be no ambiguity. 
\end{rem}
Denote by $\coh B \subseteq \grmod B \subsetneq \Grmod B$ the full subcategory of coherent and finitely presented graded modules, respectively. The next lemma is originally a result of Serre for coherent sheaves over ringed spaces \cite[Sect. II.13]{FAC}, whose proof goes through unchanged as pointed out in \cite[Prop. 1.1]{Pol} and \cite[Sect. 2]{Mina}. Note that one can equivalently recast condition $b)$ in terms of maps from finitely generated graded free modules. \bigskip 

Recall that a full subcategory $\CC \subseteq \Grmod B$ is called a Serre subcategory if for every short exact sequence in $\Grmod B$
	\[
		0 \to F \to G \to K \to 0
	\]
	if two terms of the sequence are in $\CC$, so is the third.
	\begin{lem}The subcategory $\coh B \subseteq \Grmod B$ is a Serre subcategory closed under finitely generated submodules. In particular $\coh B \subseteq \Grmod B$ is an abelian subcategory closed under extension. Moreover, when $B$ is coherent we have $\coh B = \grmod B$, and the latter then contains all graded modules of finite length. 
\end{lem}

The full subcategory $\mathsf{Grmod}_0\ \! B \subsetneq \Grmod B$ of right bounded modules also forms a Serre subcategory, and we may form the Serre quotient to obtain an abelian category
\[
	\QGr B := \frac{\Grmod B}{\mathsf{Grmod}_0 B}
\]
thought of as the category of quasicoherent sheaves on some noncommutative variety. When $B$ is coherent, finite length modules are finitely presented and so the full subcategory $\mathsf{grmod}_0\ \! B \subseteq \grmod B$ of finite length modules is a Serre subcategory with which to form the Serre quotient
\[
	\qgr B := \frac{\grmod B}{\mathsf{grmod}_0 B}
\]
to obtain the corresponding category of coherent sheaves.\bigskip

The importance of coherence for graded algebras is in the following beautiful result of Minamoto. A large supply of $d$-representation infinite algebras arise as endomorphism algebras of (good) tilting object on a Fano variety $X$ so that
\[
	{\rm D}^b(\Lambda) \cong {\rm D}^b(\coh X)
\]
and they are in general expected to `come from geometry'. In general, let $\Lambda$ be a higher representation infinite algebra $\Lambda$ with anticanonical complex $\omega^{-1}_\Lambda$, and recall that ${\rm D}^{\omega^{-1}_\Lambda} = ({\rm D}^{\leq 0, \omega^{-1}_\Lambda}, {\rm D}^{\geq 0, \omega^{-1}_\Lambda})$ denotes the pair of full subcategories 
\begin{align*}
	{\rm D}^{\leq 0, \omega^{-1}_\Lambda} &= \{ M \in {\rm D}^b(\Lambda) \ | \ {\rm H}^s(M \otimes^{\mathbb{L}}_{\Lambda}\omega_\Lambda^{-m}) = 0 \textup{ for } s > 0 \textup{ for all } m \gg 0\}\\
	{\rm D}^{\geq 0, \omega^{-1}_\Lambda} &= \{ M \in {\rm D}^b(\Lambda) \ | \ {\rm H}^s(M \otimes^{\mathbb{L}}_{\Lambda}\omega_\Lambda^{-m}) = 0 \textup{ for } s < 0 \textup{ for all } m \gg 0\}.
\end{align*}
Additionally, we write $\mathcal{H}^{\omega_{\Lambda}^{-1}}(\Lambda) = {\rm D}^{\leq 0, \omega^{-1}_{\Lambda}} \cap\ \! {\rm D}^{\geq 0, \omega^{-1}_{\Lambda}}$ for the full subcategory of `eventual modules', or simply $\mathcal{H}(\Lambda)$ if $\omega^{-1}_{\Lambda}$ is understood. 

\begin{thm}[{\cite[Thm. 3.7, Cor. 3.12]{Mina}}]\label{Minamototheorem} Let $\Lambda$ be a higher representation infinite algebra with higher preprojective algebra $\Pi(\Lambda)$. Then:
	\begin{enumerate}[1)]
	\item There are equivalences of derived categories
		\[
			-\otimes_\Lambda^{\mathbb{L}} \Pi(\Lambda):			{\rm D}(\Modsf \Lambda) \xrightarrow{\cong} {\rm D}(\QGr \Pi(\Lambda)).
		\]
	\item The following are equivalent:
		\begin{enumerate}[i)]
		\item $\Pi(\Lambda)$ is coherent. 
		\item ${\rm D}^{\omega_\Lambda^{-1}}$ forms a t-structure on $ {\rm D}^b(\Lambda) $.
	\end{enumerate}
	When either of these holds, the above equivalence descends to
	\[
		{\rm D}^b(\Lambda) \cong {\rm D}^b(\qgr \Pi(\Lambda))
	\]
	and ${\rm D}^{\omega_\Lambda^{-1}}$ is the pullback of the standard t-structure on $ {\rm D}^b(\qgr \Pi(\Lambda)) $.
\end{enumerate}
\end{thm}
\begin{rem} A similar statement appears in \cite{Pol} in the context of a triangulated category containing a geometric helix, where a certain pair of subcategories form a t-structure if and only if the graded algebra associated to this helix is coherent. \end{rem}

\begin{exmp}\cite[Prop. 5.1, Cor. 3.6]{Mina} Let $\Lambda = kQ$ be an hereditary algebra of infinite representation type. Then $\Pi(\Lambda) = \Pi(Q)$ is always coherent. \qed
\end{exmp}

\begin{exmp} Let $X$ be a smooth Fano variety with a geometric helix $\{E_i\}_{i \in \Z}$ of period $n$ as in Example \ref{helixexample}. Assume furthermore that the $E_i$ consist of sheaves and we let $\E = \bigoplus_{j=1}^{n} E_{i+j}$ and $\Lambda = {\rm End}(\E)$. Then by \cite[Thm 3.6]{BS}, $\Pi(\Lambda)$ is Noetherian and finite over its centre. Moreover when $\E$ is a vector bundle it isn't hard to see that the t-structure ${\rm D}^{\omega_\Lambda^{-1}}$ is the pushforward of the standard t-structure under the equivalence
	\[
		{\rm RHom}(\E, -): {\rm D}^b(X) \xrightarrow{\cong} {\rm D}^b(\Lambda)
	\]
	hence we have
	\[
		\coh X \cong \mathcal{H}(\Lambda) \cong \qgr \Pi(\Lambda).
	\]\qed
\end{exmp}
	
One may ask whether coherence holds in general \cite[Question 4.37]{HIO}. The following was conjectured by Minamoto in 2012 in Banff. 
\begin{conj}[Minamoto]\label{Minamotoconjecture} The higher preprojective algebra $\Pi(\Lambda)$ of a higher representation infinite algebra $\Lambda$ is always coherent. 
\end{conj}

It is an important and difficult task in general to determine when coherence holds for algebras of interest. A related class of algebras for which this problem is unresolved are the Artin-Schelter regular algebras. Let $A = k \oplus A_1 \oplus A_2 \oplus \dots $ be a graded connected $k$-algebra, finitely generated in degree one. We say that $A$ is Artin-Schelter regular if
\begin{enumerate}[1)]
\item $\gldim A = d < \infty$;
\item $A$ satisfies the Gorenstein condition
	\[
		{\rm Ext}_A^i(k, A) = \begin{cases} 0 & i \neq d \\ k(-a) & i = d \end{cases} 
	\]
	for some $a \in \Z$.
	\end{enumerate}
	Note that one sometimes requires additonal finiteness conditions, but we do not do this here. In this generality, one has:
	\begin{conj}[Bondal, \cite{Pio}]\label{Bondalconjecture} Artin-Schelter regular algebras are always coherent.
\end{conj}

We finally arrive at the results of this chapter. Let $n$ denote the global dimension of an Artin-Schelter regular algebra or of a higher preprojective algebra. Both conjectures have affirmative answers when $n \leq 2$. A main result of this thesis is a negative answer to both conjectures in all dimensions $n \geq 4$.

\begin{thm*}[Thm. \ref{abscounterexamples}]\label{counterexamples} There are higher preprojective algebras $\{\Pi_n\}_{n \geq 4}$ and Artin-Schelter regular algebras $\{E_n\}_{n \geq 4}$ of global dimension $n \geq 4$, all of which fail to be coherent. 
\end{thm*}

The bound $n \geq 4$ seems to be the lowest we can go with our methods, and so the following may be warranted.
\begin{conj*} Let $\Pi_3$ be a higher preprojective algebra and $E_3$ be an Artin-Schelter regular algebra, both of global dimension $3$. Then $\Pi_3$ and $E_3$ are coherent graded algebras. 
\end{conj*}

We will actually characterise coherence for a restricted class of algebras. Let $A$ be a Koszul Frobenius algebra, graded connected over $k$. To $A$ we will associate a higher representation infinite algebra $B$, along with a contravariant equivalence of triangulated categories 
\[
	{\rm D}_{sg}^{\Z}(A)^{op} \cong {\rm D}^b(\modsf B).
\]
We then associate to $A$ a higher preprojective algebra $\Pi = \Pi(B)$ and an Artin-Schelter regular algebra $E = {\rm Ext}^*_A(k, k)^{op}$, and recast both conjectures in terms of ${\rm D}^{\Z}_{sg}(A)$. We will prove: 

\begin{thm*}[Thm. \ref{summary6.2}]\label{characterisationcoherence} For $A$ as above, the following are equivalent:
	\begin{enumerate}[1)]
	\item $\Pi$ is coherent. 
	\item $E$ is coherent.
	\item $A$ is absolutely Koszul in the sense of Herzog-Iyengar \cite{HI}. 
\end{enumerate}
\end{thm*}
This last condition has been heavily studied in commutative algebra this last decade, and holds fairly generally. In spite of this, we will obtain counterexamples by constructing pathological commutative Frobenius Koszul algebras, building on examples due to J.-E. Roos. \bigskip 

In the above setting, we are lead to study a natural pair of full subcategories ${\rm D}^{ \mathsf{lin}} = ({\rm D}^{\leq 0}_{sg}(A), {\rm D}^{\geq 0}_{sg}(A))$ of ${\rm D}_{sg}^{\Z}(A)$, whose intersection 
\[
	\mathcal{H}^{ \mathsf{lin}}(A) = {\rm D}^{\leq 0}_{sg}(A) \cap {\rm D}^{\geq 0}_{sg}(A)
\]
consists of objects whose minimal graded free resolution is eventually linear. We will study this pair of subcategories for a general Koszul Gorenstein algebra $A$. Our principal result, from which all above stated results will follow, will be the following pair of theorems. These should be considered as sharp generalisations of the Bernstein-Gel'fand-Gel'fand correspondence for Koszul Gorenstein algebras. 

\begin{thmx}\label{ThmA} The following are equivalent for a Gorenstein Koszul algebra $A$: 
	\begin{enumerate}[i)]
	\item ${\rm D}^{\mathsf{lin}}$ forms a bounded t-structure on ${\rm D}^{\Z}_{sg}(A)$. 
	\item $A$ is absolutely Koszul in the sense of Herzog-Iyengar. 
\end{enumerate}
When either of these equivalent conditions hold, the natural realisation functor 
\[
	{\rm real}: {\rm D}^b(\mathcal{H}^{\mathsf{lin}}(A)) \xrightarrow{\cong} {\rm D}^\Z_{sg}(A).
\]
is an equivalence of triangulated categories.
\end{thmx}

\begin{thmx}\label{ThmB} Let $A$ be an absolutely Koszul Gorenstein algebra. Then $E = (A^!)^{op} = {\rm Ext}_A^*(k, k)^{op}$ is coherent, and we have a contravariant equivalence of abelian categories
	\[
		\mathcal{H}^{\mathsf{lin}}(A)^{op} \xrightarrow{\cong} \qgr E
	\]
	sending $M$ to ${\rm Ext}^*_A(M, k)$. \bigskip

	\noindent The converse holds in the Artinian case: if $A$ is Artinian Koszul Gorenstein with ${\rm Ext}^*_A(k, k)^{op}$ coherent, then $A$ is absolutely Koszul.
\end{thmx} 

Combining them, we obtain:
\begin{thmx}\label{ThmC} Let $A$ be an absolutely Koszul Gorenstein algebra, with $E = (A^!)^{op} = {\rm Ext}^*_A(k, k)^{op}$. Then we have equivalences of triangulated categories
	\[
	{\rm D}_{sg}^{\Z}(A)^{op} \cong {\rm D}^b(\qgr E)
	\]
	such that the t-structure ${\rm D}^{\mathsf{lin}}$ is the pullback of the standard t-structure on the right hand side.
\end{thmx} 
\noindent We will prove Theorems A, B and C in Chapter 5 and give the aforementioned applications in Chapter 6.

%% file: chapters/fano/fanofromkoszul.tex
\section{Fano algebras from Koszul Frobenius algebras}
%

We now set conventions for the remaining of the chapter. In order to encompass natural examples coming from quiver path algebras we will work in slightly greater generality than the previous section, but the reader will not lose out on any of the essential ideas by taking all algebras to be graded connected over a field.\bigskip

Let $k$ be a field throughout and $\Bbbk$ a finite dimensional semisimple $k$-algebra. A graded algebra $A = \Bbbk \oplus A_1 \oplus A_2 \oplus \dots $ will mean a locally finite graded $\Bbbk$-algebra, finitely generated by $A_1$ over $A_0 = \Bbbk$. By $\Bbbk$-algebra we mean that the product $A \otimes_{\Bbbk} A \to A$ is bilinear over $\Bbbk$, e.g. $A = kQ/I$ is a graded path algebra over $\Bbbk = kQ_0$. Recall that we write $D = {\rm Hom}_k(-, k)$ for $k$-duality.
\begin{defn}
	A graded Frobenius $\Bbbk$-algebra $A$ is a finite $\Bbbk$-algebra such that $D(A) \cong A(a)$ independently both as left and right $A$-modules, for some $a \in \Z$. A graded symmetric algebra is a graded Frobenius algebra such that $_A D(A)_A \cong$ $_A A(a)_A$ as bimodules.
\end{defn}

Given a graded Frobenius $A$, the invariant $a$ is then the top degree of $A$, and we will always assume that $a \geq 1$ since $a = 0$ gives $A = \Bbbk$. By \cite[Lemma 2.9]{MM}, there is a graded $k$-algebra automorphism $\nu: A \to A$, uniquely defined up to inner automorphisms, such that $_A{\rm D}(A)_A \cong$ $_{1}A(a)_{\nu}$ as bimodules, using subscripts to denote the $\nu$-twisted module structure on the right. This is the Nakayama automorphism $\nu = \nu_A$, and graded symmetric algebras are those for which we can take $\nu = {\rm id}$. Note that pulling back through $\nu$ gives an isomorphism of bimodules $_1 A(a)_{\nu} \cong$ $_{\nu^{-1}} A(a)_1$. Moreover, since $\nu$ preserves grading, it descends to an automorphism of $\Bbbk$ and similarly $_1 \Bbbk_{\nu} \cong$ $_{\nu^{-1}} \Bbbk_1$ as $A$-bimodules. In particular $\Bbbk_\nu \cong \Bbbk$ as right $A$-modules, but note that $\nu$ typically permutes the simple summands. 

Since $\Bbbk$ is semisimple, it is itself symmetric as a $k$-algebra \cite[Cor. 5.17]{FrobeniusBook} and we have natural isomorphisms between the various duality functors 
\[
	D ={\rm Hom}_k(-, k) \cong {\rm Hom}_{\Bbbk}(-, {\rm Hom}_k(\Bbbk, k)) \cong {\rm Hom}_{\Bbbk}(-, \Bbbk).
\]

Assume for now that $A$ is a graded Frobenius $\Bbbk$-algebra. Then $A$ is self-injective and so we have a natural equivalence ${\rm D}_{sg}^{\Z}(A) = \underline{\mathsf{mod}}^{\Z} A$. We will need the following. 
\begin{prop}[Auslander-Reiten Duality \cite{ARS}] Let $A$ be graded Frobenius. Writing $\omega_A = D(A) \cong$ $_1 A(a)_\nu$, the category $ \underline{\mathsf{mod}}^\Z A$ has a Serre functor given by $\Se_A = - \otimes_A \omega_A[-1]$.
\end{prop}

Next, note that every module is vacuously MCM over a self-injective algebra, and in particular every $M \in \underline{\mathsf{mod}}^\Z A$ admits a complete resolution. The minimal graded projective resolution $P_* \xrightarrow{\sim} \Bbbk = A/A_{\geq 1}$ of $\Bbbk$ dualises to a minimal projective coresolution $\Bbbk \xrightarrow{\sim} D(P_*)$, and so the minimal complete resolution $C(\Bbbk)$ of $\Bbbk$ looks like
\[
	\xymatrix@R=5pt@C=15pt{ \cdots \ar[r] & A(-1)^{\oplus \dim A_1} \ar[r] & A \ar[dr] \ar[rr] && _{\nu^{-1}}A(a)_1 \ar[r] & _{\nu^{-1}}A(a+1)_1^{\oplus \dim A_1} \ar[r] & \cdots \\
		& & &  \Bbbk \ar[ur] 
	}
\]
from which we can read the dimension of $ \underline{{\rm Ext}}_{gr A}^i(\Bbbk, \Bbbk(-j))$ and so we have $\underline{{\rm Ext}}_{gr A}^n(\Bbbk, \Bbbk) = 0$ for $n \neq 0$, and ${\rm \underline{Hom}}_{gr A}(\Bbbk, \Bbbk)$ is a quotient algebra of $\Bbbk$ and therefore semisimple. A basic set of non-projective indecomposable summands of $\Bbbk$ are then orthogonal w-exceptional\footnote{Meaning weakly exceptional, in that we allow for division rings as endomorphism algebras instead of $k$.} objects in $ \underline{\mathsf{mod}}^\Z A$, and $\thick(\Bbbk)$ is a semisimple subcategory.\bigskip

We call a thick subcategory generated by finitely many orthogonal exceptional (or w-exceptional) objects a \emph{block}. From the previous paragraph, we see that $\thick(\Bbbk) = \mathsf{add} \{ \Bbbk[i] \ | \ i \in \Z \}$ forms a block in $\underline{\mathsf{mod}}^\Z A$. One extends the usual notions of exceptional sequences and geometric helices to blocks in the natural way. Recall that $A$ is Koszul if the minimal projective resolution $P_* \xrightarrow{\sim} \Bbbk$ is linear, that is $P_i$ is generated in degree $i$, and that we denote by $\Se_n = \Se_A \circ [-n]$ the desuspended Serre functor. The following is well-known, but we will give a complete proof.

\begin{prop}\label{helicity} Let $A$ be a Koszul Frobenius $\Bbbk$-algebra with $a \geq 1$, and let $E_i = \Bbbk(-i)[i]$. Then the sequence $(E_i)_{i \in \Z}$ forms a block helix for $ \underline{\mathsf{mod}}^{\Z} A$ of period $a$. That is:
	\begin{enumerate}[i)]
	\item There is a full block exceptional collection $\underline{\mathsf{mod}}^\Z A = \left\langle E_i, E_{i+1}, \dots, E_{i+a-1} \right\rangle$.
	\item We have $E_{i+a} = \Se^{-1}_{a-1} E_i$.
\item $(E_i)_{i \in \Z}$ forms a block geometric helix: for every pair $(i,j)$ with $i \leq j$, we have $\underline{{\rm Ext}}^s_{gr A}(E_i, E_j) = 0 \textup{ for all } s \neq 0$. 
	\end{enumerate}
	\begin{proof} $i)$. The conditions $\underline{{\rm Ext}}_{gr A}^s(E_k, E_l) = \underline{{\rm Ext}}^{s+l-k}_{gr A}(\Bbbk, \Bbbk(k-l)) = 0 \textup{ for } i \leq l < k \leq i+a-1 \textup{ and all } s \in \Z$
and 
\[
	\underline{{\rm Ext}}_{gr A}^s(E_k, E_k) = \begin{cases} 0 & s \neq 0 \\
	\textup{semisimple} & s = 0\end{cases} 	
\]
follow from the structure of the minimal complete resolution $C(\Bbbk)$, hence the above forms a block exceptional decomposition. To see that it is full, note that $soc(A) \cong \Bbbk(-a)$ gives a short exact sequence 
\[
\xi:	0 \to \Bbbk(-a) \to A \to A/soc(A) \to 0
\]
which shows that $\Bbbk(-a) \cong A/soc(A)[-1]$ in $\underline{\mathsf{mod}}^\Z A$. Taking a Jordan-H\"older filtration of $A/soc(A)$ shows that $\Bbbk(-a) = A/soc(A)[-1] \in \thick(\Bbbk, \Bbbk(-1), \dots, \Bbbk(-a+1))$, and iteratively grade shifting $\xi$ by $(-1)$ shows that $\Bbbk(-a-j) \in \thick(\Bbbk, \Bbbk(-1), \dots, \Bbbk(-a+1))$ for all $j \geq 0$. Dualising $\xi$ and applying the same argument shows that $\Bbbk(-j) \in \thick(\Bbbk, \Bbbk(-1), \dots, \Bbbk(-a+1))$ for all $j \leq 0$, and applying Jordan-H\"older filtrations to a general finite dimensional module shows that $(\Bbbk(-i), \Bbbk(-i-1), \dots, \Bbbk(-i-a+1))$ is full for any $i \in \Z$. \bigskip 

\noindent $ii)$ This follows since 
\[
	\Se_{a-1}^{-1}(\Bbbk(-i)[i]) = \Bbbk(-i)[i] \otimes_A \omega_A^{-1}[a] = \Bbbk_{\nu^{-1}}(-i-a)[i+a] \cong \Bbbk(-i-a)[i+a].
\]

Lastly, $ \underline{{\rm Ext}}_{gr A}^{s}(E_i, E_j) = \underline{{\rm Ext}}_{gr A}^{s+j-i}(\Bbbk, \Bbbk(i-j)) = 0$ for $s \neq 0$ whenever $i-j \leq 0$ is a condition which only involves the nonnegative part $C_{\geq 0}(\Bbbk)$ of the complete resolution $C(\Bbbk)$, where it reduces to ${\rm Ext}_{gr A}^{s+l}(\Bbbk, \Bbbk(-l)) = 0$ for $s \neq 0$ and all $l \geq 0$ which is the definition of Koszul. 
\end{proof}
	
\end{prop}
When $A$ is Koszul Frobenius, it is immediate that the opposite category $(\underline{\mathsf{mod}}^\Z A)^{op}$ inherits a block geometric helix by setting ${\rm E}_i^{op} = E_{-i}$. Let $T_i = \bigoplus_{j=0}^{a-1} E_{i+j}$ (resp. $T_i^{op} = \bigoplus_{j=0}^{a-1} E^{op}_{i+j}$) be the associated tilting objects in $ \underline{\mathsf{mod}}^\Z A$ (resp. $(\underline{\mathsf{mod}}^\Z A)^{op}$). Let $\Lambda = {\rm \underline{End}}_{gr A}(T_i)$, which is independent of choice of $i \in \Z$, and note that $\Lambda^{op} = {\rm \underline{End}}_{gr A}(T_i^{op})$. From the geometric helix condition, we obtain: 

\begin{prop}\label{tiltingequivalence4.2} Let $A$ be a Koszul Frobenius $\Bbbk$-algebra, and let $d = a-1$. The algebras $\Lambda$, $\Lambda^{op}$ are both $d$-representation infinite algebras, and the tilting objects $T_i$, $T_i^{op}$ induce equivalences of triangulated categories 
	\begin{align*}
		F_i&: \underline{\mathsf{mod}}^\Z A  \xrightarrow{\cong} {\rm D}^b(\Lambda) \\
		F_i^{op}&: (\underline{\mathsf{mod}}^\Z A)^{op} \xrightarrow{\cong} {\rm D}^b(\Lambda^{op})
	\end{align*} 
	sending $T_i$ to $\Lambda$ (resp. $T_i^{op}$ to $\Lambda^{op}$).
\end{prop} 
\begin{proof}
That $\Lambda$ and $\Lambda^{op}$ are $d$-representation infinite is a rephrasing of condition $ii)$ and $b)$ in Prop. \ref{helicity}. The equivalence of triangulated categories is standard for any tilting object. 
\end{proof} 

Let us look at some example applications of the above proposition. 
\begin{exmp} The exterior algebra $A = \bigwedge_k^*(y_0, \dots, y_n)$ is Koszul Frobenius of top degree $a = n+1$. The algebra $\Lambda = kQ/I$ is the Beilinson algebra, with quiver given by
\[
		\xymatrix{ 
		\bullet \ar@<1ex>[r]^{ \{x_i\} } \ar[r] \ar@<-1ex>[r] & \bullet \ar@<1ex>[r]^{ \{x_i\} } \ar[r] \ar@<-1ex>[r] & \ar@<1ex>[r]^{ \{x_i\} } \ar[r] \ar@<-1ex>[r] \bullet &
 \cdots \ar@<1ex>[r]^{ \{x_i\} } \ar[r] \ar@<-1ex>[r] & \bullet \ar@<1ex>[r]^{ \{x_i\} } \ar[r] \ar@<-1ex>[r] & \bullet  
	}
\]
with $n+1$ vertices and relations $I = (x_ix_j - x_jx_i)$. 
\end{exmp}
\begin{exmp} Let $V$ be a vector space over $k$ of dimension $n \geq 2$ and $V \times V \to k$ a perfect pairing, extended to a graded algebra structure on $A = A_0 \oplus A_1 \oplus A_2 = k \oplus V \oplus k$. Then $A$ is Frobenius of top degree $a=2$, and Koszul\footnote{This reference covers the case when $A$ is commutative, or equivalently when the pairing is symmetric, but the general case follows from Prop. \ref{helicity}.} by \cite{AL}. The algebra $\Lambda \cong kQ$ is the path algebra of the $n$-Kronecker quiver
	\[
		\xymatrix{ \bullet \ar@<1ex>[r]^{ \{x_i\} } \ar[r] \ar@<-1ex>[r] & \bullet}
	\]
	where $\{ x_i \}$ is a basis for $V$. \qed
\end{exmp}
\begin{exmp} More generally, let $B = k \oplus B_1 \oplus \cdots$ be a Koszul $k$-algebra of finite global dimension $n$. Then $A = B^! = {\rm Ext}_B^*(k, k)$ is a Koszul Frobenius $k$-algebra if and only if $B$ is an Artin-Schelter regular $k$-algebra \cite[Thm. 5.10]{SMI}. The top degree of $A$ is then $a = n$. Picking a basis $\{ y_i \}$ for $B_1 \cong {\rm Ext}_{gr A}^1(k, k(-1))$, we have $\Lambda = kQ/I$ for the quiver $Q$
\[
		\xymatrix{ 
		\bullet \ar@<1ex>[r]^{ \{y_i\} } \ar[r] \ar@<-1ex>[r] & \bullet \ar@<1ex>[r]^{ \{y_i\} } \ar[r] \ar@<-1ex>[r] & \ar@<1ex>[r]^{ \{y_i\} } \ar[r] \ar@<-1ex>[r] \bullet &
 \cdots \ar@<1ex>[r]^{ \{y_i\} } \ar[r] \ar@<-1ex>[r] & \bullet \ar@<1ex>[r]^{ \{y_i\} } \ar[r] \ar@<-1ex>[r] & \bullet  
	}
\]
with $n$ vertices, and quadratic relations amongst $\{ y_i \}$ inherited from the quadratic algebra $B$. \qed
\end{exmp}
Interestingly, this construction attaches an $(n-1)$-representation infinite algebra $\Lambda$, and therefore a $n$-preprojective algebra $\Pi = \Pi(\Lambda)$, to any Koszul Artin-Schelter regular algebra $B$ with $\gldim B = n$. \bigskip

Next, let us see an example over a semisimple base $\Bbbk$ which is not a field. 
\begin{exmp} Let $A = kQ/(Q_2)$ be the radical square zero algebra over $\Bbbk = kQ_0$ with quiver $Q$ an oriented cycle 
	{\tiny{
	\begin{align*}
	\xymatrix@R=15pt@C=15pt{
		& \bullet \ar[r] & \bullet \ar[dr] \\
		\bullet \ar[ur]	& && \bullet \ar[d] \\
		\bullet \ar@{.}[u] &&& \bullet \ar[dl] \\
		& \bullet \ar[ul] & \ar[l] \bullet
	}
	\end{align*}
}}
	Computing the indecomposable projectives $P(i)$ and injectives $I(i)$, one sees that $A$ is a basic graded self-injective algebra, and so is graded Frobenius, with top degree is $a = 1$. Moreover, radical square zero algebras are always Koszul, and so the previous proposition applies. The algebra $\Lambda \cong kQ_0$ is semisimple. It follows that $ \underline{\mathsf{mod}}^\Z A$ is a semisimple category. \qed
\end{exmp}\bigskip

Next, for the rest of this subsection, we fix a Koszul Frobenius algebra $A$ of top degree $a$, with associated $d$-representation infinite algebra $\Lambda$ with $d = a-1$. For our purposes we will single out the equivalence $G := F^{op}_{-a+1}$ and spell out its properties. The tilting object is given by $T_{-a+1}^{op} = \bigoplus_{j=0}^{a-1} E_{j-a+1}^{op} = \bigoplus_{j=0}^{a-1} E_{a-1-j} = \bigoplus_{j=0}^{a-1} E_{j}$. Let us write for the record: 
\begin{cor}\label{KoszulFanocontravariantequivalence} There is a contravariant equivalence of categories 
	\[
		G: (\underline{\mathsf{mod}}^\Z A)^{op} \xrightarrow{\cong} {\rm D}^b(\Lambda^{op})
	\]
	sending $T_{-a+1} = \bigoplus_{j=0}^{a-1} E_j = \bigoplus_{j=0}^{a-1} \Bbbk(-j)[j]$ to $\Lambda^{op}$. 
\end{cor}
Our aim for the remainder of this subsection is to characterise coherence of the higher preprojective algebra $\Pi(\Lambda^{op})$ in terms of $A$. \bigskip

Let us write $\beta_{i,j} = \beta_{i,j}(M) = \dim_k \underline{{\rm Ext}}_{gr A}^i(M, \Bbbk(-j))$ for the graded Betti numbers of $M$. Of course when $\Bbbk$ decomposes we can refine the numbers $\beta_{i,j}$ further, but we will not do this. Note that $\beta_{i,j}(M_\nu) = \dim_k \underline{{\rm Ext}}_{gr A}^i(M_\nu, \Bbbk(-j)) = \dim_k \underline{{\rm Ext}}_{gr A}^i(M, \Bbbk_{\nu^{-1}}(-j)) = \dim_k \underline{{\rm Ext}}_{gr A}^i(M, \Bbbk(-j)) = \beta_{i,j}(M)$. 

\begin{lem}\label{bettichp4} Let $M \in \underline{\mathsf{mod}}^\Z A$. Then: 
	\begin{enumerate}[i)]
	\item $\dim_k {\rm H}^s(G M) = \sum\limits_{l = 0}^{a-1} \beta_{s+l, l}$
	\item $\dim_k {\rm H}^s\left((G M) \otimes_{\Lambda^{op}}^{\mathbb{L}} \omega_{\Lambda^{op}}^{-m} \right) = \sum\limits_{l = 0}^{a-1} \beta_{ma + s+l, ma+l}$.
	\end{enumerate}	
\end{lem}
\begin{proof} We have 
	\begin{align*}
		{\rm H}^s(G M) 	&= {\rm Hom}_{ {\rm D}^b(\Lambda^{op})}(\Lambda^{op}, G(M) [s]) \\
		&={\rm Hom}_{ {\rm D}^b(\Lambda^{op})}( G(T^{op}_{-a+1}), G(M) [s]) \\  
		&={\rm Hom}_{ {\rm D}^b(\Lambda^{op})}( G(T^{op}_{-a+1}), G(M[-s])) \\  
		&= {\rm \underline{Hom}}_{gr A}(M, T^{op}_{-a+1}[s])  \\
		&= \bigoplus_{l = 0}^{a-1}\underline{ {\rm Ext}}_{gr R}^{s+l}(M, \Bbbk(-l)).
	\end{align*}
This proves $i)$. For $ii)$, note that $(G M) \otimes_{\Lambda}^{\mathbb{L}} \omega_{\Lambda}^{-m} = \Se_d^{-m} \left( G M\right) = G (\Se_{a-1}^m M) = G( M_{\nu^m}(ma)[-ma])$, then apply $i)$.
\end{proof}

We can picture these cohomology groups using the Betti table $\beta(M)$ whose entries are given by $\beta_{i, i+j}$ in the $i$-th column and $j$-th row for $i, j \in \Z$:
\[
		\begin{tabular}{c|c c c c c c c c}
		& $\cdots$ & 0 & 1 & $\cdots$ & a-2 & a-1 & a & $\cdots$ \\
		\hline
		\vdots & & \vdots & \vdots & & \vdots & \vdots & \vdots\\ 
		-2 & $\cdots$ &$\beta_{0, -2}$&$\beta_{1,-1}$ & $\cdots$ & $\beta_{a-2, a-4}$ & $\beta_{a-1, a-3}$ & $\beta_{a, a-2}$ & $\cdots$ \\ 
		-1 &$\cdots$ & $\beta_{0, -1}$&$\cellcolor{gray!30}\beta_{1,0}$ & $\cellcolor{gray!30}\cdots$ & $\cellcolor{gray!30}\beta_{a-2, a-3}$ & $\cellcolor{gray!30}\beta_{a-1, a-2}$ & $\cellcolor{gray!30}\beta_{a, a-1}$ & $\cdots$ \\ 
		0 &$\cdots$ & $\cellcolor{gray!30}\beta_{0, 0}$&$\cellcolor{gray!30}\beta_{1,1}$ & $\cellcolor{gray!30}\cdots$ & $\cellcolor{gray!30}\beta_{a-2, a-2}$ & $\cellcolor{gray!30}\beta_{a-1, a-1}$ & $\beta_{a, a}$ & $\cdots$ \\ 
		1 &$\cdots$ & $\beta_{0, 1}$&$\beta_{1,2}$ & $\cdots$ & $\beta_{a-2, a-1}$ & $\beta_{a-1, a}$ & $\beta_{a, a+1}$ & $\cdots$ \\ 
		\vdots & & \vdots & \vdots & & \vdots & \vdots & \vdots\\ 
		\end{tabular}
\]
The sum of the zeroth row highlights gives the dimension of ${\rm H}^0(G M)$, and the sum of the $(-1)$-th row highlights gives that of ${\rm H}^{1}(G M)$. More generally the dimension vector of ${\rm H}^s\left((G M) \otimes_{\Lambda^{op}}^{\mathbb{L}} \omega^{-m}_{\Lambda^{op}} \right)$ represents $a$-many consecutive entries on row $-s$ and these entries move to the right as $m$ increases. This is because the operation 
\[
	\beta(M) \mapsto \beta(\Se_{a-1} M) = \beta(M_\nu(a)[-a]) = \beta(M(a)[-a]) 
\]
translates the table to the left by $a$-many columns.\bigskip 

Now let us focus on the equivalence $G: ( \underline{\mathsf{mod}}^{\Z}\ \! A)^{op} \xrightarrow{\cong} {\rm D}^b(\Lambda^{op})$. By Minamoto's Theorem, coherence of $\Pi(\Lambda^{op})$ is equivalent to the pair of full subcategories ${\rm D}^{\omega^{-1}} = ({\rm D}^{\leq 0, \omega^{-1}_{\Lambda^{op}}}, {\rm D}^{\geq 0, \omega^{-1}_{\Lambda^{op}}})$
\begin{align*}
	{\rm D}^{\leq 0, \omega^{-1}_{\Lambda^{op}}} &= \{ X \in {\rm D}^b(\Lambda^{op}) \ | \ {\rm H}^s(X \otimes^{\mathbb{L}}_{\Lambda^{op}}\omega_{\Lambda^{op}}^{-m}) = 0 \textup{ for } s > 0 \textup{ for all } m \gg 0\}\\
	{\rm D}^{\geq 0, \omega^{-1}_{\Lambda^{op}}} &= \{ X \in {\rm D}^b(\Lambda^{op}) \ | \ {\rm H}^s(X \otimes^{\mathbb{L}}_{\Lambda^{op}}\omega_{\Lambda^{op}}^{-m}) = 0 \textup{ for } s < 0 \textup{ for all } m \gg 0\}
\end{align*}
forming a t-structure on $ {\rm D}^b(\Lambda^{op}) $, which is then bounded. We can translate this into a statement about $\underline{\mathsf{mod}}^\Z A$. Before moving further, let us first recall basic notions concerning t-structures (see \cite[Chp. 8]{HTT}). \bigskip

Let $\T$ be a triangulated category with a pair of full subcategories $t = (\T^{\leq 0}, \T^{\geq 0})$, and denote $\T^{\leq n} = \T^{\leq 0} \circ [-n]$ and $\T^{\geq n} = \T^{\geq 0} \circ [-n]$. The pair $t$ forms a t-structure if and only if these satisfy the three defining properties:
\begin{enumerate}[T1.]
	\item $\T^{\leq 0} \subseteq \T^{\leq 1}$ and $\T^{\geq 1} \subseteq \T^{\geq 0}$.
	\item ${\rm Hom}_\T(\T^{\leq 0}, \T^{\geq 1}) = 0$.
	\item For all $X$ in $\T$ there is a distinguished triangle $X^{\leq 0} \to X \to X^{\geq 1} \to X^{\leq 0}[1]$ with $X^{\leq 0}$ in $\T^{\leq 0}$ and $X^{\geq 1}$ in $\T^{\geq 1}$.
\end{enumerate}
When this holds, the embedding $i_n: \T^{\leq n} \hookrightarrow \T$ and $j_n: \T^{\geq n} \hookrightarrow \T$ admit right and left adjoints 
\begin{align*}
	 &i_n: \T^{\leq n} \leftrightarrows \T: \tau^{\leq n} \\
	 &\tau^{\geq n}: \T \leftrightarrows \T^{\geq n}: j_n.
\end{align*}
The unit and counit maps of the adjunctions gives rise to maps $X \to \tau^{\geq n}X$ and  $\tau^{\leq n}X \to X$, which fit into distinguished triangles
\[
	\tau^{\leq n}X \to X \to \tau^{\geq n+1}X \to (\tau^{\leq n}X)[1].
\]
Letting $\CC = \T^{\leq 0} \cap \T^{\geq 0}$ stand for the heart of this t-structure, which is naturally an abelian category, one defines the t-cohomology objects $\mathcal{H}^0X = \tau^{\geq 0} \tau^{\leq 0}X$ and $\mathcal{H}^nX = \mathcal{H}^0(X[n]) = (\tau^{\geq n}\tau^{\leq n}X)[n]$ living in $\CC \subseteq \T$. The functor $\mathcal{H}^0: \T \to \CC$ is then a cohomological functor. A t-structure is bounded if every object $X \in \T$ has finitely many non-zero cohomology objects $\mathcal{H}^nX$. \bigskip

The axioms T1, T2, T3 are self-dual, and so a t-structure $t = (\T^{\leq 0}, \T^{\geq 0})$ on $\T$ induces an opposite t-structure $t^{op} = (\T^{op, \leq 0}, \T^{op, \geq 0})$ on $\T^{op}$ by setting
	\begin{align*}
		\T^{op, \leq 0} &= (\T^{\geq 0})^{op} \\
		\T^{op, \geq 0} &= (\T^{\leq 0})^{op}.
	\end{align*}
 Note that the suspension on $\T^{op}$ is given by $[-1]$, and so we have
	\begin{align*}
		\T^{op, \leq n} &= (\T^{\geq -n})^{op} \\
		\T^{op, \geq n} &= (\T^{\leq -n})^{op}.
	\end{align*}
Now consider the pair $({\rm D}^{\leq 0, \omega_{\Lambda^{op}}^{-1}},\ {\rm D}^{\geq 0, \omega_{\Lambda^{op}}^{-1}})$ on ${\rm D}^b(\Lambda^{op})$. We define the pair $(\underline{\mathsf{mod}}^{\leq 0}A,\ \underline{\mathsf{mod}}^{\geq 0}A)$ on $\underline{\mathsf{mod}}^\Z A$ by pulling back the above pair through the contravariant equivalence $G$, so that for every $n \in \Z$ we have
	\begin{align*}
		\underline{\mathsf{mod}}^{\leq n} A &= G^{-1}\left({\rm D}^{\geq -n, \omega_{\Lambda^{op}}^{-1}} \right)\\
		\underline{\mathsf{mod}}^{\geq n} A &= G^{-1}\left({\rm D}^{\leq -n, \omega_{\Lambda^{op}}^{-1}} \right).
	\end{align*}

	\begin{prop}\label{pullbacktstructure} The subcategories $ \underline{\mathsf{mod}}^{\leq n} A$, $ \underline{\mathsf{mod}}^{\geq n}A $, and $\underline{\mathsf{mod}}^{\leq 0}A\ \! \cap\ \! \underline{\mathsf{mod}}^{\geq 0}A $ of $ \underline{\mathsf{mod}}^\Z A$ are given by
	\begin{align*}
		\underline{\mathsf{mod}}^{\leq n} A &= \{ M \in \underline{\mathsf{mod}}^{\Z} A \ | \ {\rm H}^s\left( (G M) \otimes_{\Lambda^{op}}^{\mathbb{L}} \omega_{\Lambda^{op}}^{-m} \right) = 0 \textup{ for } s < -n \textup{ whenever } m \gg 0 \}\\
		&= \{ M \in \underline{\mathsf{mod}}^{\Z} A \ | \ \beta_{i,j}(M) = 0 \textup{ for } j-i > n \textup{ whenever } j \gg 0 \}\\
		&= \{ M  \in \underline{\mathsf{mod}}^{\Z} A\ | \ \beta_{i,j}(M) = 0 \textup{ for } j-i > n \textup{ whenever } i \gg 0 \}\\
		 \\
		 \underline{\mathsf{mod}}^{\geq n}A &= \{ M  \in \underline{\mathsf{mod}}^{\Z} A \ | \ {\rm H}^s\left( (G M) \otimes_{\Lambda^{op}}^{\mathbb{L}} \omega_{\Lambda^{op}}^{-m} \right) = 0 \textup{ for } s > -n \textup{ whenever } m \gg 0 \}\\
		 &= \{ M \in  \underline{\mathsf{mod}}^{\Z} A\ | \ \beta_{i,j}(M) = 0 \textup{ for } j-i < n \textup{ whenever } j \gg 0 \}\\
		 &= \{ M \in  \underline{\mathsf{mod}}^{\Z} A\ | \ \beta_{i,j}(M) = 0 \textup{ for } j-i < n \textup{ whenever } i \gg 0 \}\\
		 \\
		 \underline{\mathsf{mod}}^{\leq 0}A \cap \underline{\mathsf{mod}}^{\geq 0}A  &=	\{ M \in \underline{\mathsf{mod}}^{\Z} A \ | \ \beta_{i,j}(M) = 0 \textup{ for } i \neq j \textup{ whenever } j \gg 0 \} \\
		 &= \{ M \in \underline{\mathsf{mod}}^{\Z} A \ | \ \beta_{i,j}(M) = 0 \textup{ for } i \neq j \textup{ whenever } i \gg 0 \}.
	\end{align*}
\end{prop}
\begin{proof} By Lemma \ref{bettichp4} we have $\dim_k {\rm H}^s(G(M) \otimes^{\mathbb{L}}_{\Lambda^{op}} \omega_{\Lambda^{op}}^{-m}) = \sum_{l=0}^{a-1}\beta_{ma + s+ l, ma + l}$. The first equality follows by setting $s = i-j$, and the second by noting that the indices $i, j$ in $\{(i,j) \ | \ \beta_{i,j} \neq 0 \}$ go to $+\infty$ together.
\end{proof}
	\begin{rem} We can picture these categories in terms of the shape of the Betti table $\beta(M)$. The category $\underline{\mathsf{mod}}^{\leq n} A$ consists of modules whose Betti table is eventually supported on the $n$-th row or above (resp. $\underline{\mathsf{mod}}^{\geq n} A$ consists of modules with Betti table eventually supported on the $n$-th row or below). 
\end{rem}\bigskip

Now, pulling back the pair of subcategories $(\underline{\mathsf{mod}}^{\leq 0}A,\ \underline{\mathsf{mod}}^{\geq 0}A)$ of Prop. \ref{pullbacktstructure} onto Minamoto's pair $({\rm D}^{\leq 0, \omega_{\Lambda^{op}}^{-1}},\ {\rm D}^{\geq 0, \omega_{\Lambda^{op}}^{-1}})$ through the contravariant equivalence $G$ of Prop. \ref{KoszulFanocontravariantequivalence}, we may rephrase Minamoto's Theorem \ref{Minamototheorem} as follows: 
\begin{prop}\label{summary4.2} The following are equivalent.
	\begin{enumerate}[i)]
	\item The higher preprojective algebra $\Pi(\Lambda^{op})$ is coherent. 
	\item The pair ${\rm D}^{\omega^{-1}} = ( {\rm D}^{\leq 0, \omega_{\Lambda^{op}}^{-1}}, {\rm D}^{\geq 0, \omega_{\Lambda^{op}}^{-1}})$ forms a t-structure on ${\rm D}^b(\Lambda^{op})$, which is then bounded. 
		\item The pair $( \underline{\mathsf{mod}}^{\leq 0}A, \underline{\mathsf{mod}}^{\geq 0}A)$ forms a t-structure on $ \underline{\mathsf{mod}}^{\Z}\ \! A$, which is then bounded. 
\end{enumerate}
\end{prop}

In the next section we will investigate when a natural generalisation of the above forms a t-structure for the graded singularity category over an arbitrary Koszul Gorenstein $\Bbbk$-algebra.

%% file: chapters/tstructure/tstructure.tex
\chapter{Absolutely Koszul algebras and t-structures of Koszul type}

Throughout this chapter, $A$ will denote more generally a Koszul Gorenstein algebra satisfying some mild finiteness hypotheses to be set-out shortly; certainly two-sided Noetherian and graded connected over a field suffices, and the arguments in the general case will not differ substantially. Over such an algebra $A$, define the following full subcategories of $ \underline{{\rm MCM}}^{\mathbb{Z}}(A)$:
	\begin{align*}
		\underline{\rm{MCM}}^{\leq n}(A) &= \{ M  \in \underline{{\rm MCM}}^{\mathbb{Z}}(A) \ | \ \beta_{i,j}(M) = 0 \textup{ for } j-i > n \textup{ whenever } i \gg 0 \}\\
		\underline{\rm{MCM}}^{\geq n}(A) &= \{ M  \in \underline{{\rm MCM}}^{\mathbb{Z}}(A) \ | \ \beta_{i,j}(M) = 0 \textup{ for } j-i < n \textup{ whenever } i \gg 0 \}.
	\end{align*}
	Moreover, we will consider the pair $t^{\mathsf{lin}} = (\underline{{\rm MCM}}^{\leq 0}(A),\ \underline{{\rm MCM}}^{\geq 0}(A))$ as a candidate pair for a t-structure, with intersection
	\[
		\mathcal{H}^{\mathsf{lin}}(A) := \underline{{\rm MCM}}^{\leq 0}(A) \ \! \cap \ \! \underline{{\rm MCM}}^{\geq 0}(A)
	\]
	the category of eventually linear stable MCM modules. The goal of this chapter is to give a proof of the following theorems.
	\begin{thmy}\label{thmA} Let $A$ be a Koszul Gorenstein algebra. The following are equivalent:
	\begin{enumerate}[i)]
	\item $A$ is absolutely Koszul.
	\item $t^{\mathsf{lin}}$ forms a bounded t-structure.
\end{enumerate}
When either of these equivalent conditions hold, the natural realisation functor
\[
	\mathsf{real}: {\rm D}^b(\mathcal{H}^{\mathsf{lin}}(A)) \xrightarrow{\cong} \underline{{\rm MCM}}^{\mathbb{Z}}(A) 
\]
is an equivalence of triangulated categories.
\end{thmy}

\begin{thmy}\label{thmB} Let $A$ be an absolutely Koszul, Gorenstein $\Bbbk$-algebra, which we assume homologically homogeneous. Then $E= (A^!)^{op} = {\rm Ext}_A^*(\Bbbk, \Bbbk)^{op}$ is coherent and we have a contravariant equivalence of abelian categories
	\[
		\mathcal{H}^{\mathsf{lin}}(A)^{op} \xrightarrow{\cong} \qgr E 
	\]
	sending $M$ to ${\rm Ext}^*_A(M, \Bbbk)$. \bigskip 
	
	\noindent Conversely, if $A$ is Artinian Koszul Gorenstein with $E =(A^!)^{op} = {\rm Ext}^*_A(\Bbbk, \Bbbk)^{op}$ coherent, then $A$ is absolutely Koszul. 
\end{thmy} \bigskip

Putting the two theorems together, we then obtain: 
\begin{thmy}\label{thmC} Let $A$ be an absolutely Koszul, Gorenstein $\Bbbk$-algebra, which is assumed homologically homogeneous, and let $E = (A^!)^{op} = {\rm Ext}^*_A(\Bbbk, \Bbbk)^{op}$. Then there exists an equivalence of triangulated categories 
	\[
		\underline{{\rm MCM}}^{\mathbb{Z}}(A)^{op} \cong {\rm D}^b(\qgr E)
	\]
	such that $t^{\mathsf{lin}}$ arises as the pullback of the standard t-structure on the right hand side. 
\end{thmy}

\input{chapters/tstructure/thmA}

\input{chapters/tstructure/thmB}

\input{chapters/tstructure/singdim}

%% file: chapters/tstructure/thmA.tex
\section{Linearity defect and Theorem A}

\subsection*{Standing hypotheses}
We now impose hypotheses for the remainder of this chapter, which are slightly more permissive than those used in Chapter 4. As before $k$ denotes a fixed field and $\Bbbk$ will denote a fixed finite-dimensional semisimple $k$-algebra. A graded $\Bbbk$-algebra will always mean a graded algebra $S = \Bbbk \oplus S_1 \oplus S_2 \oplus \dots$, with multiplication $S \otimes_{\Bbbk} S \to S$ bilinear over $\Bbbk$, and we further assume that $\dim_k S_i < \infty$ for all $i \geq 0$. \bigskip 

The weakened assumptions are meant to capture the following natural examples:
\begin{enumerate}[a)]
\item $S = k[x_0, \dots, x_n]/I$ is a commutative graded $k$-algebra of finite type over $k = \Bbbk$;
\item $S = kQ/I$ is a graded path algebra for a finite bound quiver $(Q, I)$ with path-length homogeneous relations $I$, thought of as an algebra over $\Bbbk := kQ_0$;
\item skew group algebras $S * G$ over $\Bbbk = kG$, with underlying graded vector space $(S*G)_i := (S \otimes_k kG)_i = S_i \otimes_k kG$ and twisted multiplication
	\[
		(s, g)* (s', g') = (sg(s'), gg')
	\]
where $S = k \oplus A_1 \oplus \cdots $ a standard graded algebra and $G \leq {\rm Aut}_0(S)$ a finite subgroup of homogeneous automorphisms of $S$, over $k$ of characteristic not dividing the order of the group; 
\item more generally the Koszul duals $S^! = {\rm Ext}^*_S(\Bbbk, \Bbbk)$ of Koszul algebras of the above form.
	\end{enumerate} 
	Additionally, working over a semisimple base $\Bbbk$ instead of the field $k$, starting in Section 5.2 we will need to assume that our algebras are homologically homogeneous, meaning that the simple summands of $\Bbbk$ have the same projective dimension. The relevance is only to the proof of Theorem B and we will impose it from Section 5.2 onwards.

\subsection*{Basic properties} 
Recall that modules over a graded algebra $S$ will by default refer to graded right $S$-modules. Under the assumption that $S$ is (right) coherent, the category of finitely presented modules $\grmod S \subsetneq \Grmod S$ is an abelian subcategory closed under extensions and under taking finitely generated submodules, and contains all finite length modules. Moreover, each $M \in \grmod S$ admits a projective resolution 
\[
\cdots \to P_2 \to P_1 \to P_0 \to M \to 0
\]
by finitely generated projectives. Recall that a complex of projectives $P_* = (P_*, \partial)$ is minimal if $P_* \otimes_S \Bbbk$ has trivial differential. Writing $\mathsf{r}_S := \rad(S) = S_{\geq 1}$ for the homogeneous Jacobson radical of $S$, this is equivalently to the statement $\partial(P_{n+1}) \subseteq P_{n}\mathsf{r}_S$ for all $n$. \bigskip

Recall that a general graded algebra $T$ is semiperfect if $T/J(T)$ is semisimple and the idempotents of $T/J(T)$ lift to $T$. This holds in our setting as $S/\mathsf{r}_S = S/S_{\geq 1} = \Bbbk$ is semisimple, and idempotent lifting holds for degree reasons. By \cite{Krause}, each $M \in \grmod S$ admits a projective cover $\pi: P \twoheadrightarrow M$, meaning that $P$ is finite projective with ${\rm ker}(\pi) \subseteq P\mathsf{r}_S$, and so each $M \in \grmod S$ admits a minimal projective resolution $P_* \xrightarrow{\sim} M$, unique up to (non-canonical) isomorphism. Moreover, by \cite{Krause} the category $\grmod S$ has the Krull-Schmidt property. \bigskip 

The map $S \twoheadrightarrow S/\mathsf{r}_S = \Bbbk$ for $S$ semiperfect always induces a bijection on indecomposable finitely generated (graded) projectives \cite{Krause}, with $\Bbbk$ considered as a graded algebra sitting in degree zero. Since we have a splitting $S \leftrightarrows \Bbbk$, here this bijection takes the form
\begin{align*}
P_S = P &\mapsto P \otimes_S \Bbbk \\ 
V_\Bbbk = V &\mapsto V\otimes_\Bbbk S.
\end{align*}
In particular every finitely generated projective $S$-module can be written as $P = V \otimes_\Bbbk S$ for some unique graded $\Bbbk$-module $V$.\bigskip

Throughout this chapter, $A$ will always denote a graded Gorenstein $\Bbbk$-algebra, which is assumed two-sided coherent and of Gorenstein dimension $d = \idim(_A A) = \idim(A_A) < \infty$. Given $M \in \grmod A$ with finite presentation 
\[
P_1 \to P_0 \to M \to 0
\]
the dual (left) module $M^* = {\rm Hom}_A(M, A)$ is the kernel of a map between finitely generated projectives over $A^{op}$ 
\[
0 \to M^* \to P_0^* \to P_1^*
\]
and so must be finitely presented by two-sided coherence of $A$. Taking dual modules then restricts to a functor between finitely presented modules 
\[
	(-)^*: \grmod A \to \grmod A^{op}.
\]
By \cite[Lemma 5.3]{BurkeStev}, whose proof in the Noetherian case applies in the coherent setting (see also \cite[Prop. 3.5]{Yeku}), this extends to a duality on the bounded derived category
\[
	(-)^{\vee} := {\rm RHom}_A(-, A): {\rm D}^b(\grmod A)^{op} \xrightarrow{\cong} {\rm D}^b(\grmod A^{op})
\]
in that $(-)^{\vee \vee} \simeq {\rm id}$. We form in the usual ways the categories $\underline{{\rm MCM}}^{\mathbb{Z}}(A) \subseteq \underline{\mathsf{mod}}^\Z A$ and ${\rm D}_{sg}^{\Z}(A)$. The proofs of Buchweitz's theorems \cite{BuchweitzManuscript} apply verbatim in our situation to give equivalences
\[
	\mathcal{K}_{ac}({\rm proj}^\Z \ \! A) \cong \underline{{\rm MCM}}^{\mathbb{Z}}(A) \cong {\rm D}_{sg}^{\Z}(A)
\]
with the middle category inheriting the structure of a triangulated category through either equivalence. The stable category $ \underline{{\rm MCM}}^{\mathbb{Z}}(A)$ also inherits the Krull-Schmidt property from ${\rm MCM}^{\Z}(A)$ since endomorphism algebras of graded modules have finite length over $\Bbbk$ and are then Artinian, and so have the idempotent lifting property (see \cite{Krause} for more details). \bigskip

Given $N \in \grmod A$, by the Krull-Schmidt property we have an essentially unique decomposition $N = [N] \oplus F$ with $F$ the largest projective summand of $N$. Note that $N$ and $N'$ are stably isomorphic if and only if $[N]$ and $[N']$ are isomorphic. We will need the following lemma, which was used previously in Chapter 3 but which is worth spelling out in full here. 
\begin{lem}\label{unstablelemma} Let $N \in \grmod A$. Then ${\rm \underline{Hom}}_{gr A}\left(N, \Bbbk(j)\right) \cong {\rm Hom}_{gr A}\left([N], \Bbbk(j)\right)$ for all $j \in \Z$. 
\end{lem}
\begin{proof} We have ${\rm \underline{Hom}}_{gr A}\left(N, \Bbbk(j)\right) \cong {\rm \underline{Hom}}_{gr A}\left([N], \Bbbk(j)\right)$ and so it is sufficient to check that any morphism $f: [N] \to \Bbbk(j)$ factoring through a projective 
	\[
		\xymatrix@R=10pt{[N] \ar[dr]^-{\alpha} \ar[rr]^f && \Bbbk(j) \\
		& P \ar[ur]^-{\beta}&}
	\]	
	must be zero. Assume it is not. Since $\Bbbk(j)$ is semisimple, the image ${\rm im}(f) = {\rm im}(\beta)$ is a semisimple summand supported in degree $-j$ and without loss of generality we can assume that $f$ and $\beta$ are surjective. The map $\beta$ then further factors through the projective cover $\pi: A(j) \twoheadrightarrow \Bbbk(j)$ as depicted 
	\[
		\xymatrix@R=10pt{[N] \ar[dr]^-{\alpha} \ar[ddr]_{\widetilde{\alpha}} \ar[rr]^f && \Bbbk(j) \\
			& P \ar@{..>}[d] \ar[ur]^-{\beta}&\\
		& A(j) \ar[uur]_{\pi}& 
	}
	\]
	That $f$ is onto is implies that $\widetilde{\alpha}$ is also onto for degree reasons, in which case $[N]$ contains a projective summand, a contradiction.
\end{proof}

\subsection*{Linearity defect and t-structures} 
We are finally ready to investigate when the pair of full subcategories $t^{\mathsf{lin}} = ( \underline{{\rm MCM}}^{\leq 0}(A),\ \underline{{\rm MCM}}^{\geq 0}(A))$ forms a t-structure. Recall that we define the (stable)\footnote{It would be more appropriate to call these the stable graded Betti numbers as these differ slightly from the usual definition, in that ours satisfy $\beta_{0,*}(P) = 0$ for any finite projective $P$. Since we will only use the notation $\beta_{i,j}$ as defined here, the author hopes that there will be no ambiguity.} graded Betti numbers of $M \in \underline{{\rm MCM}}^{\mathbb{Z}}(A)$ as the dimensions 
\[
	\beta_{i,j}(M) = \dim_k \underline{{\rm Ext}}^i_{gr A}(M, \Bbbk(-j))	 
\]
and that for $i \geq 0$ we have agreements
\[ 
	\beta_{i,j}(M) = \beta_{i,j}([M]) = \dim_k {\rm Ext}^i_{gr A}([M], \Bbbk(-j)) = \dim_k {\rm Tor}^A_i([M], \Bbbk)_j. 
\]
We have defined full subcategories $ \underline{{\rm MCM}}^{\leq n}(A)$, $ \underline{{\rm MCM}}^{\geq n}(A)$ by
	\begin{align*}
		\underline{\rm{MCM}}^{\leq n}(A) &= \{ M  \in \underline{{\rm MCM}}^{\mathbb{Z}}(A) \ | \ \beta_{i,j}(M) = 0 \textup{ for } j-i > n \textup{ whenever } i \gg 0 \}\\
		\underline{\rm{MCM}}^{\geq n}(A) &= \{ M  \in \underline{{\rm MCM}}^{\mathbb{Z}}(A) \ | \ \beta_{i,j}(M) = 0 \textup{ for } j-i < n \textup{ whenever } i \gg 0 \}.
	\end{align*}
	Note that $ \underline{{\rm MCM}}^{\leq 0}(A)[-n] = \underline{{\rm MCM}}^{\leq n}(A)$ and $ \underline{{\rm MCM}}^{\geq 0}(A)[-n] = \underline{{\rm MCM}}^{\geq n}(A)$. We say that $M$ is $n$-linear if ${\rm Tor}^S_i(M, \Bbbk)_j = 0$ for $j-i \neq n$, and that $M$ is linear if it is $0$-linear. The intersection
	\[
		\mathcal{H}^{\mathsf{lin}}(A) := \underline{{\rm MCM}}^{\leq 0}(A) \ \! \cap \ \! \underline{{\rm MCM}}^{\geq 0}(A)
	\]
	is the subcategory of eventually linear stable MCM modules. The above pair typically occurs as the t-structure dual to the standard t-structure under (contravariant) Koszul duality equivalences. Let us begin with a classical example. 

	\begin{exmp}[BGG correspondence for exterior and symmetric algebras] Let $A = \bigwedge_k^*(y_0, \cdots, y_n)$ be an exterior algebra over $k$ with Koszul dual symmetric algebra $B := {\rm Ext}_A^*(k, k) \cong k[x_0, \cdots, x_n]$. Following \cite{EFS}, there are natural functors $L, R$ giving rise to an equivalence of triangulated categories
		\[
			L: {\rm D}^b(\grmod A) \leftrightarrows {\rm D}^b(\grmod B): R
		\]
		We have a natural isomorphism of complexes $L(M) \simeq {\rm RHom}_A(k, M)$ up to a natural regrading, which gives \cite[Prop. 2.3]{EFS}
		\[
			{\rm H}^j(L(M))_i = {\rm Ext}^{i-j}_{gr A}(k, M(-i))
		\]
		Consider the duality 
		\[
			{\rm D}: {\rm D}^b(\grmod A)^{op} \xrightarrow{\cong} {\rm D}^b(\grmod A) 
		\]
		given by $D = {\rm RHom}_A(-, \omega_A)$, where $\omega_A = A(a) = A(n+1)$, and $A^{op} \cong A$ naturally identified via the antipode map. Then $D$ satisfies $D^2 \simeq {\rm id}$ and $D(k) \cong k$. Combining with the above, we obtain a contravariant equivalence
		\[
			L D: {\rm D}^b(\grmod A)^{op} \xrightarrow{\cong} {\rm D}^b(\grmod B)
		\]
		with a natural isomorphism $L D \simeq {\rm RHom}_A(k, D(-)) \simeq {\rm RHom}_A(-, k)$ up to regrading, and we have the more natural 
		\[
			{\rm H}^j(LD(M))_i = {\rm Ext}^{i-j}_{gr A}(M, k(-i)).
		\]
		The functor $LD$ sends perfect complexes to complexes with Artinian (or torsion) cohomology, and so descends to an equivalence
		\[
		LD: (\underline{\mathsf{mod}}^\Z A)^{op} \xrightarrow{\cong} {\rm D}^b(\coh \mathbb{P}^n). 
		\]
		From the above cohomology calculations, one sees that $ \underline{\mathsf{mod}}^{\geq 0} A = LD^{-1}({\rm D}^{\leq 0})$ and $ \underline{\mathsf{mod}}^{\leq 0} A = LD^{-1}({\rm D}^{\geq 0})$, where $({\rm D}^{\leq 0}, {\rm D}^{\geq 0})$ is the standard t-structure on $ {\rm D}^b(\coh \mathbb{P}^n) $. \qed
\end{exmp}
	 
The above example generalises easily, and one can streamline the proof using Morita theory (see Appendix A.2). The next example is well-known (see \cite[Ex. 3.2, Ex. 11.1]{GreenStev} where this example was treated in the ungraded situation). The author has learned the idea behind the proof from a talk by G. Stevenson. 
\begin{exmp}[BGG correspondence for Artinian-Coherent pairs]\label{coherentBGG} Let $A$ be an Artinian, Koszul, Gorenstein $k$-algebra. Let $B = {\rm Ext}_A^*(k, k)$ be the Koszul dual, and assume that $B^{op}$ is a coherent graded algebra.\bigskip 
	
Consider the full subcategory $S \subseteq {\rm D}^b(\grmod A) $ defined by $S = \{ k(-i)[i] \ | \ i \in \Z \}.$ One can think of the small subcategory $S$ as abstractly equivalent to a category with object set given by $\Z$, with morphisms given by 
	\[
		{\rm Hom}(i, j) = {\rm Ext}^{j-i}_{gr A}(k, k(i-j))
	\]
	and composition given by that of the Ext algebra. The subcategory $S$ is a tilting subcategory, in the sense that: 
	\begin{enumerate}[a)]
	\item $S$ classically generate $ {\rm D}^b(\grmod A)$, i.e. $\thick(S) = {\rm D}^b(\grmod A) $,
	\item $S$ has no non-trivial extensions, meaning that ${\rm Ext}_{gr A}^n(s, s') = 0$ for $n \neq 0$ for any objects $s, s' \in S$. 
\end{enumerate}
The first fact holds since $A$ is Artinian and the second since $A$ is Koszul. Moreover, $S^{op}$ is a tilting subcategory for the opposite category $ {\rm D}^b(\grmod A)^{op}$. By Keller's Tilting Theorem \ref{derivedmorita}, we have an equivalence of triangulated categories 
\[
	{\rm RHom}(-, S): {\rm D}^b(\grmod A)^{op} \xrightarrow{\cong} {\rm D}^{\perf}(\Modsf S^{op})
\]
onto the subcategory of perfect DG-modules ${\rm D}^{\perf}(\Modsf S^{op}) \subseteq {\rm D}(\Modsf S^{op})$. This is simply the multi-object version of Keller's Theorem which we have repeatedly used throughout this thesis, but let us unpack the notation further.\bigskip

Let $\mathcal{S} \subseteq {\rm D}^b_{\mathsf{dg}}(\grmod A)$ be the small DG category obtained from $S$ with same objects as $S$ but with morphisms computed in a fixed DG enhancement $\D := {\rm D}^b_{\mathsf{dg}}(\grmod A)$ of ${\rm D}^b(\grmod A)$. Since $S$ is tilting, there is a quasi-equivalence of DG categories $\mathcal{S} \simeq {\rm H}^0(\mathcal{S}) = S$ inducing an equivalence of (perfect) derived categories
\[
	{\rm D}^{\perf}(\Modsf \mathcal{S}) \cong {\rm D}^{\perf}(\Modsf S)
\]
and similarly for opposite categories $\mathcal{S}^{op}$, $S^{op}$. Next, note that any object $M \in \mathcal{D}$ induces a left $\mathcal{S}$-module by restricting ${\rm Hom}_\D(M, -)$ to $\mathcal{S} \subseteq \mathcal{D}$. The equivalence ${\rm RHom}(-, S)$ is then obtained by composing
\[
	{\rm Hom}_\D(-, \mathcal{S}): {\rm D}^b(\grmod A) \xrightarrow{\cong} {\rm D}^{\perf}(\Modsf \mathcal{S}^{op}) 
\]
with ${\rm D}^{\perf}(\Modsf \mathcal{S}^{op}) \cong {\rm D}^{\perf}(\Modsf S^{op})$. Now, a DG module over $S^{op}$ is nothing but a complex of graded modules over $B^{op} = {\rm Ext}^*_A(k, k)^{op}$, so that ${\rm D}^{\perf}(\Modsf S^{op}) = {\rm D}^{\perf}(\Grmod B^{op})$. Since $A$ is Artinian, we have that $\gldim B^{op} < \infty$, and since $B^{op}$ is coherent we deduce that ${\rm D}^{\perf}(\Grmod B^{op}) = {\rm D}^b(\grmod B^{op}) $. The above equivalence can then be rewritten as 
\[
	{\rm RHom}(-, S): {\rm D}^b(\grmod A)^{op} \xrightarrow{\cong} {\rm D}^b(\grmod B^{op}).
\] 
Given $M \in {\rm D}^b(\grmod A)$, the complex ${\rm RHom}(M, S)$ is a complex of graded $B^{op}$-modules, with cohomology given by 
\begin{eqnarray}\label{cohomologycomputation}
	{\rm H}^j({\rm RHom}(M, S))_i \cong {\rm Ext}_{gr A}^{j}(M, k(-i)[i]) = {\rm Ext}_{gr A}^{i+j}(M, k(-i))
\end{eqnarray} 
and the standard action of $B^{op} = {\rm Ext}^{*}_A(k, k)^{op}$ by post-composition. We deduce that ${\rm RHom}(A, S) \cong B \in {\rm D}^b(\grmod B^{op})$, and so that ${\rm RHom}(-, S)$ sends perfect complexes onto complexes with finite length cohomology. This equivalence then descends to 
\[ 
	\underline{ {\rm MCM}}^\Z (A)^{op} \xrightarrow{\cong} {\rm D}^b(\qgr B^{op}).
\]
Finally, from the isomorphism \ref{cohomologycomputation}, we see that the category $\grmod B^{op} \subseteq {\rm D}^b(\grmod B^{op}) $ corresponds to complexes with linear minimal projective resolution under ${\rm RHom}(-, S)$, and so the pair $t^{\mathsf{lin}}$ on $ \underline{{\rm MCM}}^{\mathbb{Z}}(A)$ arises as the pullback of the standard t-structure on $ {\rm D}^b(\qgr B^{op}) $. \qed 
\end{exmp}\bigskip

One advantage of the argument by Morita theory is that it makes clear the formal nature of the argument, so long as $A$ is Artinian and $B$ is coherent. However, equivalences of this type hold more generally than this type of argument would initially lead one to believe. 

\begin{exmp}[BGG correspondence after Buchweitz {\cite[Appendix]{BEH}}]\label{BuchweitzBGG} 
	Let $R = k[ \underline{x}]/( \underline{q})$ be a complete intersection of quadrics, as studied in Chapter 2. Let $C = R^! = {\rm Ext}_R^*(k, k)$ be its Koszul dual, a homogeneous Clifford algebra over a polynomial subalgebra $k[ \underline{\eta}]$ with basis $\{ \underline{\eta}\}$ dual to $\{ \underline{q} \}$. Then $(R, C)$ is a pair of Koszul dual Noetherian Koszul Gorenstein algebras. Letting $A$ be either $R$ or $C$, and modifying slightly the approach in \cite[Appendix]{BEH}, we have seen that there is an equivalence
	\[
		{\rm D}^b(\grmod A)^{op} \leftrightarrows {\rm D}^b(\grmod A^!) 
	\]
	exchanging perfect complexes and complexes with Artinian cohomology, and so descending to equivalences 
	\[
		\underline{{\rm MCM}}^{\mathbb{Z}}(A)^{op} \leftrightarrows {\rm D}^b(\qgr A^!)
	\]
	and such that $t^{\mathsf{lin}}$ arises on $ \underline{{\rm MCM}}^{\mathbb{Z}}(A)$ via the pullback of the standard t-structure under this equivalence. \qed
	\end{exmp}
	\begin{rem}
		Explicit use is made in \cite[Appendix]{BEH} of the Noetherianity of ${\rm Ext}_R^*(k, k)$, and complete intersections of quadrics are the only commutative Koszul algebras with Noetherian Ext algebra. The proof does not immediately extend to the case of ${\rm Ext}_R^*(k, k)$ coherent, although an extension beyond the Noetherian case may be possible by involving new ideas, see the discussion in Sect. 6.3. 
\end{rem}\bigskip

Now, instead of imposing finiteness conditions on our algebras $A$ or their Koszul dual ${\rm Ext}^*_A(\Bbbk, \Bbbk)$, we will instead make use of criteria more internal to the stable category $ \underline{{\rm MCM}}^{\mathbb{Z}}(A)$. We begin with reviewing a parallel story. \bigskip 

\subsection*{Absolutely Koszul Algebras} 
As part of their study of graded modules over the exterior algebra $A = \bigwedge_k^*(y_0, \dots, y_n)$ through the BGG correspondence, Eisenbud, Fl{\o}ystad and Schreyer proved the following result. 
\begin{thm}[Eisenbud-Fl{\o}ystad-Schreyer \cite{EFS}]\label{EFS} Let $A = \bigwedge_k^*(y_0, \cdots, y_n)$ and let $M \in {\rm D}^b(\grmod A) $ with minimal free resolution $F_* = (F_*, \partial) \xrightarrow{\sim} M$. Then the linear part of the differential on $F$ dominates; that is, expressing $\partial$ as matrices with entries in $A$ and removing the entries of degree $\geq 2$ yields a complex ${\rm lin}^A(F_*)$ with at most bounded cohomology. 
	\end{thm}

	This notion was further analysed by Herzog-Iyengar \cite{HI} and R\"omer \cite{RO}, who introduced the notions below for (graded-)commutative algebras and whose results form the basis of this chapter. Note that closely related notions were studied by Mart\'inez-Villa and Zacharia in \cite{MVZ} (see also \cite{Madsen}). \bigskip
	
	Let $S$ be a coherent graded $\Bbbk$-algebra and let $M \in {\rm D}^b(\grmod S) $ with minimal projective resolution $F_* = (F_*, \partial) \xrightarrow{\cong} M$: 
\[
	\dots \to F_{n+1} \xrightarrow{\partial} F_n \xrightarrow{\partial} F_{n-1} \to \cdots 
\]
Writing $\mathsf{r}_S = S_{\geq 1}$, the minimal complex $F_*$ satisfies $\partial(F_{n+1}) \subseteq  F_{n} \mathsf{r}_S$, and we define a filtration by subcomplexes $\F = \{\mathcal{F}^kF_* \}_{k \geq 0}$ by setting $(\mathcal{F}^kF_*)_n = F_n\mathsf{r}_S^{k-n}$, where by convention $\mathsf{r}_S^k = S$ for $k \leq 0$. 

\begin{defn}[Herzog-Iyengar \cite{HI}] The associated graded complex ${\rm lin}^S(F_*) := gr_\F(F_*)$ is called the \emph{linear part} of $F$.
\end{defn} 

The definition in terms of the filtration $\F$ makes sense in fair generality (e.g. local Noetherian rings), but this complex takes a simpler form in our setting. Each projective module has a canonical form $F_{n+1} \cong V_{n+1} \otimes_{\Bbbk} S$ for a unique graded $\Bbbk$-module $V_{n+1}$, and so the differential $\partial: F_{n+1} \to F_n$ is uniquely determined by its restriction
\[
	\partial: V_{n+1} \to V_n \otimes_{\Bbbk} S .
\]
Expand $\partial = \partial_1 + \partial_2 + \dots$, where $\partial_i(V_{n+1}) \subseteq V_n \otimes_\Bbbk S_i$. Then $\partial^2 = 0$ implies $\partial_1^2 = 0$ for degree reasons, and one has
\[
	{\rm lin}^S(F_*) = (F_*, \partial_1)
\]
similar to the statement of Thm. \ref{EFS}. 
\begin{defn}[Herzog-Iyengar \cite{HI}] Let $S$ be a coherent graded $\Bbbk$-algebra and $M \in {\rm D}^b(\grmod S)$. Keeping the above notation, the linearity defect of $M$ is defined as $${\rm ld}_S(M) = \sup \{ n \ | \ {\rm H}_n({\rm lin}^S(F_*)) \neq 0 \}.$$ 
\end{defn}
\begin{defn}[Herzog-Iyengar \cite{HI}] Let $S$ be a coherent graded $\Bbbk$-algebra and $M \in \grmod S$. We say that $M$ is \emph{Koszul} if ${\rm ld}_S(M) = 0$. 
\end{defn}
\begin{exmp} Any module $M$ for which $F_*$ has linear differential $\partial = \partial_1$ is Koszul. \qed
\end{exmp}

\begin{exmp}[Herzog-Iyengar \cite{HI}] Let $S = k[x,y]$ and $M$ have minimal resolution $F_*$ of length two
	\[
		\xymatrix{
			0 \ar[r] & S(-1) \ar[r]^-{\begin{bsmallmatrix} x^2 \\ y \end{bsmallmatrix}} & S \oplus S(1) \ar[r] & 0.
	}
	\]
	Then ${\rm lin}^S(F_*)$ is given by
	\[
		\xymatrix{
			0 \ar[r] & S(-1) \ar[r]^-{\begin{bsmallmatrix} 0 \\ y \end{bsmallmatrix}} & S \oplus S(1) \ar[r] & 0.
	}
	\]
	and so ${\rm H}_n({\rm lin}^S(F_*)) = 0$ for $n > 0$ and $M$ is Koszul. \qed
\end{exmp}
Finiteness of linearity defect is a useful property with strong consequences. The following proposition was proved in the commutative case by Herzog-Iyengar. Recall that the Hilbert function $H_S(t) = \sum_{n \geq 0} \dim_k\! S_n\  t^n$ of a standard graded, commutative $k$-algebra $S$ is rational, of the form
\[
	H_S(t) = \frac{h_S(t)}{(1-t)^{\dim S}}.
\]
\begin{prop}[Herzog-Iyengar, {\cite[Prop. 1.8]{HI}}]\label{uniformrationality} Let $S$ be a commutative, graded connected $k$-algebra, and let $M$ be a graded $S$-module with ${\rm ld}_S(M) < \infty$. Then the Poincar\'e series $P_S^M(t) := \sum\limits_{n \geq 0} \dim_k {\rm Ext}^n_{S}(M, k)t^n \in \Z[[t]]$ is rational of the form 
	\[
		P_S^M(t) = \frac{Q^M_S(t)}{h_S(-t)}
	\]
	for some polynomial $Q^M_S(t) \in \Z[t]$.
\end{prop}
In other words, modules with finite linearity defect have rational Poincar\'e series with uniform denominator depending only on the structure of $S$. We make the following definition. 
\begin{defn}\label{abskoszul} Let $S$ be a coherent graded $\Bbbk$-algebra. We say that $S$ is absolutely Koszul if
	\begin{enumerate}[i)]
	\item $S$ is Koszul;
	\item 	${\rm ld}_S(M) < \infty$ for all finitely presented modules $M$. Equivalently, every such $M$ has a Koszul syzygy.
\end{enumerate}
\end{defn}
Herzog and Iyengar introduced this notion for commutative graded $k$-algebras $S$ (and more generally for Noetherian local rings); they actually showed that $i)$ follows from $ii)$ by using the Avramov-Eisenbud-Peeva \cite{AE, AP} characterisation of commutative Koszul algebras as graded algebras $S$ such that the Castlenuovo-Mumford regularity
\[
	\mathsf{reg}_S(M) = \sup \{ r \ | \ {\rm Tor}^S_{n}(M, \Bbbk)_{n+r} \neq 0,\ n \in \N \}
\]
is finite for each $M \in \grmod S$. It is not known to the author whether the finiteness of $\mathsf{reg}_S(M)$ for all $M \in \grmod S$ implies (or even follows from) Koszulity of $S$ in our generality, and we will simply impose condition $i)$ for now. Let us at least record a lemma. 
\begin{lem}\label{finitereg} Let $S$ be a coherent graded $\Bbbk$-algebra, and let $M, N \in \grmod S$.  
	\begin{enumerate}[a)]
	\item Assume that $M$ is Koszul, with generator degrees in the interval $[n, n']$. Then ${\rm Tor}_i^S(M, \Bbbk)_{i+j} = 0$ for $j \notin [n, n']$. 
	\item Assume that $N$ has finite linearity defect. Then $N$ has finite regularity. 
\end{enumerate}
\end{lem}
\begin{proof} To see $a)$, let $F_*$ be the minimal resolution of $M$. Since $M$ is Koszul, note that it has the same graded Betti numbers as ${\rm H}_0( {\rm lin}^S(F_*))$, which breaks down as a direct sum of $j$-linear modules for $n \leq j \leq n'$. Part $b)$ then follows since $N$ has a Koszul syzygy, which has finite regularity by $a)$. 
\end{proof}\bigskip

The prototypical examples of absolutely Koszul algebras are as follows. 
\begin{exmp}The theorem of Eisenbud-Fl{\o}ystad-Schreyer (Thm. \ref{EFS}) states that the exterior algebra $\bigwedge^*_k(y_0, \cdots, y_n)$ is absolutely Koszul. \qed
\end{exmp}
\begin{exmp}[Herzog-Iyengar, {\cite[Cor. 5.10]{HI}}\footnote{The cited result concerns the local case but applies equally to the graded setting.}] Complete intersections of quadrics are absolutely Koszul.\qed
\end{exmp}
\noindent In general, we have a proper containment 
\[
	\{ \textup{Absolutely Koszul algebras} \} \subsetneq \{ \textup{Koszul algebras} \}
\]
even when restricted to commutative $k$-algebras. The matter of comparing the relative sizes of each class has a long history in commutative algebra, and is closely connected to rationality problems for Poincar\'e series. Given a local Noetherian commutative ring $(R, \m, k)$, the Serre-Kaplansky Conjecture asked for a proof of the rationality of the Poincar\'e series
\begin{align*}
	P_R^k(t) = \sum\limits_{n \geq 0} \dim_k {\rm Ext}^n_R(k, k) t^n.
\end{align*}
Anick gave the first counterexample in his 1982 thesis \cite{AN}, and his constructions have since found many applications and refinements. When $S$ is a commutative Koszul $k$-algebra, the Poincar\'e series $P_S^k(t) = H_{S^!}(t)$ is always rational by the standard identity
\[
	H_{S^!}(t) = \frac{1}{H_S(-t)} = \frac{(1+t)}{h_S(-t)}^{\! \dim S}
\] but Jacobsson \cite{JA} gave an example of a module $M$ over such a ring for which $P_S^M(t)$ is transcendental. The construction methods of Anick and Jacobsson were later taken up by Roos \cite{RO}, who introduced the following definition in general:
\begin{defn} A Koszul $\Bbbk$-algebra $S$ is good in the sense of Roos if all finitely presented modules $M$ have rational Poincar\'e series
	\[
		P^M_S(t) = \frac{Q^M_S(t)}{d_S(t)}
	\]
	with $Q^M_S(t), d_S(t) \in \Z[t]$ and denominator $d_S(t)$ independent of $M$.
\end{defn} 
Prop. \ref{uniformrationality} states that absolutely Koszul commutative $k$-algebras are good, with uniform denominator $d_S(t) = h_S(-t)$. Roos gave examples of 'bad' (i.e. not good) Koszul algebras. For instance: 

\begin{thm}[Roos, {\cite[Thm 2.4b(A)]{RO}}]\label{badkoszulalgebra} Let $S' = k[x_1,x_2]/(x_1, x_2)^2$. Then $S := S' \otimes_k S'$ is a bad Koszul algebra in the sense of Roos. More specifically, there is a sequence of finitely generated graded $S$-modules $\{M_\alpha\}_{\alpha \in \N}$ with $P^{M_\alpha}_{S}(t)$ rational, expressed in reduced form as
	\[
		P^{M_\alpha}_{S}(t) = \frac{Q_{\alpha}(t)}{d_\alpha(t)} 
	\]
	with $\limsup (\deg\ \! d_\alpha(t)) = \infty$ as $\alpha$ runs over $\N$. \qed
\end{thm}
More precisely in the above theorem, given any candidate denominator $d_S(t) \in \Z[t]$, we see that $d_\alpha(t)$ cannot divide $d_S(t)$ for infinitely many $\alpha \in \N$ as the degree of $d_\alpha(t)$ can be arbitrarily high. In particular, taking $d_S(t) = h_S(-t)$ we deduce that ${\rm ld}_S(M_\alpha) = \infty$ for infinitely many $\alpha \in \N$ by Prop. \ref{uniformrationality}.

\begin{rem}
	The `bad' Koszul algebra $S' \otimes_k S'$ is particularly interesting, and will play a role in our construction of a counterexample to the conjectures of Minamoto and Bondal in Chapter 6.
\end{rem} \bigskip

\subsection*{Existence of t-structures} 
We are now almost in position to prove the first part of Theorem $A$. Recall that Gorenstein algebras are implicitly taken to be two-sided coherent, which ensures the existence of complete resolutions for finitely presented MCM modules, but the remaining notions refer to right modules. We also denoted $t^{\mathsf{lin}} = ( \underline{{\rm MCM}}^{\leq 0}(A),\ \underline{{\rm MCM}}^{\geq 0}(A))$. \bigskip

Let us discuss the structure of this subsection. Our main aim is to prove the following result:
\begin{prop*}[Proposition \ref{thmApart1}] Let $A$ be a Koszul Gorenstein algebra. The following are equivalent:
	\begin{enumerate}[i)]
	\item $A$ is absolutely Koszul.
	\item $t^{\mathsf{lin}}$ forms a bounded t-structure for $ \underline{{\rm MCM}}^{\mathbb{Z}}(A)$. 
\end{enumerate}
\end{prop*} 
The main tool will be an elegant characterisation of Koszul modules over a coherent Koszul algebra $S$ due to R\"omer. Recall that $M \in \grmod S$ is $n$-linear if ${\rm Tor}^S_i(M, \Bbbk)_{i+j} = 0$ for $j \neq n$ and all $i \geq 0$, and $M$ is linear if it is $0$-linear. Let $M_{\langle n \rangle} \subseteq M$ be the submodule generated by $M_n$ for fixed $n \in \Z$. We say that $M$ is \emph{componentwise linear} if $M_{\langle n \rangle}$ is $n$-linear for each $n \in \Z$. We also denote the initial degree of $M$ by
\[
	\mathsf{indeg}(M) = \inf\{ n \ | \ M_n \neq 0 \}.
\]
The following is due to R\"omer \cite[Lemma 3.2.2, Thm. 3.2.8]{RO}, \cite[Sect. 5]{IR} in commutative algebra. 
\begin{thm*}[Thm. \ref{kosmod}] Let $S$ be a coherent Koszul $\Bbbk$-algebra and $M \in \grmod S$ a module with finite regularity, with $\mathsf{indeg}(M) = n$. The following are equivalent:
	\begin{enumerate}[i)]
		\item $M$ is componentwise linear; 
		\item $M_{\langle n \rangle}$ is $n$-linear and $M/M_{\langle n \rangle}$ is componentwise linear;
		\item $M$ is Koszul.
	\end{enumerate}
\end{thm*}
While R\"omer worked over a commutative or graded-commutative, standard graded Noetherian algebras $S$, the proof works in our generality and for completeness we will reproduce the arguments, with at most superficial changes to cover our situation. We will follow the exposition in \cite[Sect. 5]{IR}. We first need a few lemmas. 

\begin{lem}[{\cite[Lemma 5.4]{IR}}]\label{Lemma5.4} Let $S$ be a coherent Koszul $\Bbbk$-algebra and let $M \in \grmod S$. If $M$ is $n$-linear, then $M \mathsf{r}_S$ is $(n+1)$-linear.
\end{lem}
\begin{proof}
	The module $M/M\mathsf{r}_S \cong \bigoplus_{i \in I} \Bbbk(-n)$ has an $n$-linear resolution since $S$ is Koszul. Consider the short exact sequence 
	\[
		0 \to M\mathsf{r}_S \to M \to M/M\mathsf{r}_S \to 0.
	\] The long exact sequence of ${\rm Tor}^S(-, \Bbbk)$ preserves internal degrees, and one deduces that ${\rm Tor}^S_i(M\mathsf{r}_S, \Bbbk)_{i+j} = 0$ for $j > n+1$. Since $M\mathsf{r}_S$ is generated in degree $n+1$, ${\rm Tor}^S_i(M\m_S, \Bbbk)_{i+j}= 0$ for $j < n + 1$ by properties of minimal resolutions, and it follows that $M\mathsf{r}_S$ is $(n+1)$-linear. 
\end{proof}
Note that as corollary we see that $n$-linear modules are componentwise linear, since if $M$ is $n$-linear and $i \geq n$, we have $M_{\langle i \rangle} = M\mathsf{r}_S^{i-n}$ is $i$-linear. \bigskip

Next, we have:
\begin{lem}[{\cite[Lemma 5.5]{IR}}]\label{Lemma5.5} Let $S$ be a coherent Koszul $\Bbbk$-algebra and let $M \in \grmod S$, with $\mathsf{indeg}(M) = n$. The following are equivalent:
	\begin{enumerate}[i)]
	\item $M$ is componentwise linear;
	\item $M_{\langle n \rangle}$ is $n$-linear and $M/M_{\langle n \rangle}$ is componentwise linear. 
\end{enumerate}
\end{lem}
\begin{proof} For $i \geq n$, we have $M_{\langle n \rangle \langle i \rangle} = M_{\langle n \rangle}\mathsf{r}_S^{i-n}$, and one sees that the sequence
	\[
		0 \to M_{\langle n \rangle \langle i \rangle} \to M_{\langle i \rangle} \to (M/M_{\langle n \rangle})_{\langle i \rangle} \to 0
	\]
	is exact. Under both hypotheses, $M_{\langle n \rangle}$ is $n$-linear and so by Lemma \ref{Lemma5.4} the module $M_{\langle n \rangle \langle i \rangle}$ is $i$-linear. Inspecting the long exact sequence of Tor modules, one sees that $M_{\langle i \rangle}$ is $i$-linear if and only if $\left(M/M_{\langle n \rangle} \right)_{\langle i \rangle}$ is $i$-linear, and the statement follows. 
\end{proof}\bigskip

We next introduce an important subcomplex of the minimal resolution, which refines the information contained in the linear part. Let $S$ be a coherent Koszul $\Bbbk$-algebra as before and $M \in \grmod S$, with minimal resolution $F_* = (F_*, \partial) \xrightarrow{\sim} M$. Let $n = \mathsf{indeg}(M)$. The minimal resolution $F_*$ has terms of the form
\[
F_i = V_i \otimes_\Bbbk S 
\]
with $V_i = \bigoplus_{j \geq n} V_{i, i+j}$ a bigraded finite $\Bbbk$-module. We define the (zeroth) \emph{linear strand} $F_*^{(n)} \subseteq F_*$ to be the subcomplex with terms of the form 
\[
	F^{(n)}_i = V_{i,i+n} \otimes_\Bbbk S. 
\]
That this is a subcomplex follows since we have $\partial_1(F_i^{(n)}) \subseteq F_{i-1}^{(n)}$ and $\partial_{\geq 2}(F_*^{(n)}) = 0$ for degree reasons. In particular $F_*^{(n)}$ is a linear subcomplex of $F_*$. \bigskip

The definition of the linear strand $F^{(n)}_*$ makes sense for any minimal, lower bounded complex of finite graded projective modules. It applies in particular to the complex $F_*/F^{(n)}_*$, whose zeroth linear strand we denote $F_*^{(n+1)}$, and call it the first linear strand of $F_*$. Continuing this way, we produce linear complexes $F_*^{(n)}, F_*^{(n+1)}, F_*^{(n+2)}, \dots $, called the zeroth, first, second, $\dots$, linear strands of $F_*$. Since the underlying bigraded $S$-module of the linear strand is a direct summand of the underlying module of $F_*$, comparing differentials one obtains a direct sum decomposition of the lineart part of $F_*$ 
\[
	{\rm lin}^S(F_*) = (F_*, \partial_1) = \bigoplus_{i \geq n} F^{(n+i)}_*
\]
in terms of the linear strands. Note that the largest $r$ indexing a non-zero summand $F_*^{(n+r)} \neq 0$ computes the regularity of $M$, via
\[
	\mathsf{reg}_S(M) - \mathsf{indeg}(M) = r
\]
and in particular that the regularity of $M$ is finite if and only if the decomposition is finite. \bigskip

The inclusions of subcomplexes $F^{(n)}_* \subseteq F_*$, $F^{(n+1)}_* \subseteq F_*/F_*^{(n)}$, $\dots$, keep track of the way one reconstructs $F_*$ from the summands of ${\rm lin}^S(F_*)$. The way this information is reflected at the level of modules is encapsulated by the following lemma. 

\begin{lem}\label{linearstrandlemma} Let $S$ be a coherent Koszul $\Bbbk$-algebra and $M \in \grmod S$, with minimal resolution $F_* \xrightarrow{\sim} M$. Let $n = \mathsf{indeg}(M)$, and assume that $M_{\langle n \rangle}$ is $n$-linear. Then the linear strand $F_*^{(n)} \subseteq F_*$ is a minimal resolution of $M_{\langle n \rangle}$, and we have a short exact sequence of resolutions
\[
\xymatrix{
	0 \ar[r] & F_*^{(n)} \ar[r] \ar[d]^-{\simeq} & F_* \ar[r] \ar[d]^-{\simeq} & F_*/F_*^{(n)} \ar[r] \ar[d]^-{\simeq} & 0 \\
	0 \ar[r] & \widetilde{M}_{\langle n \rangle} \ar[r] & \widetilde{M} \ar[r] & \widetilde{M}/\widetilde{M}_{\langle n \rangle} \ar[r] & 0.
}
\]
\end{lem}
\begin{proof}
Let $E_* \xrightarrow{\sim} M_{\langle n \rangle}$ be the minimal resolution of the $n$-linear module $M_{\langle n \rangle}$, and let $\alpha: E_* \to F_*$ be a comparison morphism lifting the inclusion map $\iota: M_{\langle n \rangle} \hookrightarrow M$. We will show that $\alpha$ is injective and that ${\rm im}(\alpha) = F_*^{(n)}$ is the zeroth linear strand of $F_*$. \bigskip
	
	Let $G_* = {\rm Cone}(\alpha)$, with distinguished triangle 
	\begin{eqnarray}\label{EFGtrian}
		E_* \to F_* \to G_* \to E_*[1].
	\end{eqnarray}
	The long exact sequence of homology shows that ${\rm H}_i(G_*) = 0$ for $i < 0$ and $i > 1$, and moreover that ${\rm H}_0(G_*) \cong M/M_{\langle n \rangle}$, fitting into an exact sequence 
		\[
		\xymatrix{ 
			0 \ar[r] &{\rm H}_1(G_*) \ar@{=}[d] \ar[r] & {\rm H}_0(E_*) \ar[r] \ar@{=}[d] & {\rm H}_0(F_*)  \ar[r] \ar@{=}[d] & {\rm H}_0(G_*) \ar[r] \ar@{=}[d] &  0 \\
		0 \ar[r] &{\rm H}_1(G_*) \ar[r] & M_{\langle n \rangle} \ar[r]^-{\iota} & M \ar[r] & M/M_{\langle n \rangle}\ar[r] &  0 } 
	\]
	Since $\iota$ is injective, we must have ${\rm H}_1(G_*) = 0$. It follows that $G_*$ is a (not necessarily minimal) resolution of $M/M_{\langle n \rangle}$. \bigskip
	
	We must then have ${\rm Tor}_i^S(G_*, \Bbbk)_{i+j} \cong {\rm Tor}_i^S(M/M_{\langle n \rangle}, \Bbbk)_{i+j} = 0$ for $j < n+1$ and all $i \geq 0$. Taking the long exact sequence of Tor associated to the distinguished triangle \ref{EFGtrian}, one sees that the boundary map ${\rm Tor}_i^S(G_*, \Bbbk)_{i+j} \to {\rm Tor}_{i-1}^S(E_*, \Bbbk)_{i+j}$ vanishes for all $j \in \Z$, since either the domain or codomain is zero by combining the above vanishing with the $n$-linearity of $E_*$. The long exact sequence breaks down into short exact sequences of totalised Tor modules 
		\[
		\xymatrix{ 
			0 \ar[r] & {\rm Tor}^S_*(E_*, \Bbbk) \ar[r] \ar@{=}[d] & {\rm Tor}^S_*(F_*, \Bbbk)  \ar[r] \ar@{=}[d] & {\rm Tor}^S_*(G_*, \Bbbk) \ar[r] \ar@{=}[d] &  0 \\ 
		0 \ar[r] & E_*\otimes_S \Bbbk \ar[r]^-{\alpha\otimes_S \Bbbk} & F_*\otimes_S \Bbbk  \ar[r] & {\rm Tor}^S_*(G_*, \Bbbk) \ar[r] &  0.}
	\] 
	The second line holds since $F_*$ and $E_*$ are minimal. By the Nakayama lemma, $\alpha: E_* \to F_*$ must then be injective. Since $E_*$ is $n$-linear, we see that ${\rm im}(\alpha) \subseteq F_*^{(n)}$ lands into the linear strand, and from the degree vanishing ${\rm Tor}_i^S(G_*, \Bbbk)_{i+j} = 0$ for $j < n+1$, we see that ${\rm im}(\alpha) = F_*^{(n)}$. We have shown that $\alpha: E_* \xrightarrow{\cong} F_*^{(n)}$ is an isomorphism. In particular $F_*^{(n)} \subseteq F_*$ resolves the submodule $M_{\langle n \rangle} \subseteq M$, as we claimed. \bigskip 

	It remains to consider the short exact sequence of complexes
	\[
		0 \to F_*^{(n)} \to F_* \to F_*/F_*^{(n)} \to 0
	\]
	Repeating the above argument with $F_*/F_*^{(n)}$ instead of $G_*$, the long exact sequence of homology shows that ${\rm H}_0(F_*/F_*^{(n)}) \cong M/M_{\langle n \rangle}$ and ${\rm H}_i(F_*/F_*^{(n)}) = 0$ for $i \neq 0$, and so the above is a short exact sequence of minimal resolutions. 
\end{proof}

In order to prove R\"omer's Theorem, we will also need the next simple lemma on the behavior of linear strands. 
\begin{lem}\label{linearstrandhomology} Let $S$ be a coherent graded $\Bbbk$-algebra and $X_*$ be a complex of finite graded projectives $S$-modules, such that $X_i = 0$ for $i < 0$ and with $\mathsf{indeg}(X_0) = n$. Assume that: 
	\begin{enumerate}[i)]
	\item $X_*$ has finitely many linear strands, that is the decomposition below is finite:
		\[
			{\rm lin}^S(X_*) = \bigoplus_{i \geq n} X_*^{(n+i)}.
		\]
	\item ${\rm H}_k({\rm lin}^S(X_*)) = 0$ for $k > 0$.
\end{enumerate}
Then ${\rm H}_k(X_*) = 0$ for $k > 0$ as well. 
\end{lem}
\begin{proof}
	Under the hypothesis, $X_*$ is a finite iterated extension of complexes $X^{(n+i)}_*$ with ${\rm H}_{k}(X^{(n+i)}) = 0$ for $k > 0$, and so $X_*$ inherits the same property. 
\end{proof}

We now turn to the proof of R\"omer's Theorem. 
\begin{thm}[R\"omer, {\cite[Lemma 3.2.2, Thm. 3.2.8]{RO}}, {\cite[Sect. 5]{IR}}]\label{kosmod} Let $S$ be a coherent Koszul $\Bbbk$-algebra and $M \in \grmod S$ a module with finite regularity, with $\mathsf{indeg}(M) = n$. The following are equivalent:
	\begin{enumerate}[i)] 
		\item $M$ is componentwise linear;
		\item $M_{\langle n \rangle}$ is $n$-linear and $M/M_{\langle n \rangle}$ is componentwise linear;
		\item $M$ is Koszul.
	\end{enumerate}
\end{thm}
\begin{proof}
The equivalence of $i)$ and $ii)$ is given by Lemma \ref{Lemma5.5}, and we will use it implicitly in the proof. We will establish the equivalence of $i)$ and $iii)$ by induction on $r(M) := \mathsf{reg}_S(M) - \mathsf{indeg}(M) \geq 0$, which is the number of linear strands in the minimal resolution $F_*$ of $M$. When $r(M) = 0$, $M$ is $n$-linear and so both Koszul and componentwise linear, and there is nothing to prove. We now assume that $r(M) > 0$. \bigskip

We will make repeated use of the following observation: consider the short exact sequence
\[
	0 \to F_*^{(n)} \to F_* \to F_*/F_*^{(n)} \to 0
\]
associated to the inclusion of the zeroth linear strand. The underlying morphism of bigraded $S$-modules is split exact, and this extends to a splitting of linear parts
\begin{eqnarray}\label{splitlinearpart}
	{\rm lin}^S(F_*) = {\rm lin}^S(F_*^{(n)}) \oplus {\rm lin}^S(F_*/F_*^{(n)}).
\end{eqnarray} 

\noindent {\bf The implication $iii) \implies i)$}. Assume that $M$ is Koszul, meaning that ${\rm H}_i({\rm lin}^S(F_*)) = 0$ for $i > 0$. By the lineart part direct sum decomposition (\ref{splitlinearpart}), the same statement holds for the zeroth linear strand $F_*^{(n)}$. We claim that the natural map ${\rm incl}_*: {\rm H}_0(F_*^{(n)}) \to {\rm H}_0(F_*) \cong M$ coming from ${\rm incl}: F_*^{(n)} \hookrightarrow F_*$ is injective with image $M_{\langle n \rangle}$; it will follow that $M_{\langle n \rangle}$ is $n$-linear in particular. \bigskip 

The above map is part of the long exact sequence of homology
\[
	\cdots \to {\rm H}_1(F_*/F^{(n)}_*) \to {\rm H}_0(F^{(n)}) \xrightarrow{{\rm incl}_*} {\rm H}_0(F_*) \to {\rm H}_0(F_*/F_*^{(n)}) \to 0. 
\]
It is clear by construction of $F^{(n)}$ that the image of ${\rm incl}_*$ is sent onto $M_{\langle n \rangle}$ under the isomorphism ${\rm H}_0(F_*) \cong M$, and we claim that ${\rm H}_1(F_*/F^{(n)}_*) = 0$. Indeed $F_*/F^{(n)}_*$ has finitely many linear strands and ${\rm H}_i({\rm lin}^S(F_*/F^{(n)}_*)) = 0$ for $i > 0$ by (\ref{splitlinearpart}) and the assumption that $M$ is Koszul, hence ${\rm H}_k(F_*/F^{(n)}_*) = 0$ for $i > 0$ by Lemma \ref{linearstrandhomology}. We conclude that ${\rm incl}_*: {\rm H}_0(F^{(n)}) \cong M_{\langle n \rangle}$ is an isomorphism. \bigskip 

Since $M_{\langle n \rangle}$ is $n$-linear, we can use Lemma \ref{linearstrandlemma} to conclude that $F_*/F_*^{(n)}$ resolves $M/M_{\langle n \rangle}$, and from the decomposition of linear parts (\ref{splitlinearpart}) we see that $M/M_{\langle n \rangle}$ is Koszul and that $r(M/M_{\langle n \rangle}) < r(M)$. By induction, we obtain that $M/M_{\langle n \rangle}$ is componentwise linear, and since $M_{\langle n \rangle}$ is $n$-linear the equivalence $i) \iff ii)$ shows that $M$ itself is componentwise linear. This proves $iii) \implies i)$. \bigskip 

\noindent {\bf The implication $i) \implies iii)$}. Assuming $i)$, then $M_{\langle n \rangle}$ is $n$-linear and $M/M_{\langle n \rangle}$ is componentwise linear. We again use Lemma \ref{linearstrandlemma} to conclude that $F_*^{(n)}$ resolves $M_{\langle n \rangle}$ and $F_*/F_*^{(n)}$ resolves $M/M_{\langle n \rangle}$, and from the linear part decomposition (\ref{splitlinearpart}) we conclude that $M$ is Koszul if and only $M/M_{\langle n \rangle}$ is Koszul, and that $r(M/M_{\langle n \rangle}) < r(M)$. Since $M/M_{\langle n \rangle}$ is componentwise linear, it is Koszul by induction and we are done. 
\end{proof}\bigskip 

We can finally prove the main result of this subsection, Prop. \ref{thmApart1} below. To alleviate notation, we will use the standard $\Omega^m M := M[-m]$ for $m \in \Z$ to denote (co)syzygies in $ \underline{{\rm MCM}}^{\mathbb{Z}}(A)$, and we will write $M_{\langle \leq n \rangle} \subseteq M$ for the submodule generated by $M_{\leq n}$. \bigskip

As we will be working in the stable category $ \underline{{\rm MCM}}^{\mathbb{Z}}(A)$, we point out that the construction $M_{\langle n \rangle}$ is not an invariant of the stable isomorphism class of $M \in \underline{{\rm MCM}}^{\mathbb{Z}}(A)$, as for instance $(M \oplus P)_{\langle n \rangle} = M_{\langle n \rangle} \oplus P_{\langle n \rangle}$ need not be stably isomorphic to $M_{\langle n \rangle}$. However, the construction $M_{\langle \leq n \rangle}$ is stably invariant, as indecomposable projective modules are cyclic and thus generated in single degree, so that $M_{\langle \leq n \rangle} \oplus P_{\langle \leq n \rangle}$ is stably isomorphic to $M_{\langle \leq n \rangle}$ for any finite graded projective $P$. \bigskip

Similarly, the statement that $M$ is $n$-linear is not stably invariant, however we will use it to mean that the unique module representative $[M]$ without projective summand is $n$-linear. Moreover the statement that $M$ is Koszul (or componentwise linear) is perfectly invariant. In light of these observations, for $M \in \underline{{\rm MCM}}^{\mathbb{Z}}(A)$ we will define $M_{\langle n \rangle} := [M]_{\langle n \rangle}$, and note that $M_{\langle n \rangle} = M_{\langle \leq n \rangle}$ whenever $\mathsf{indeg}([M]) = n$. \bigskip

Recall that the Betti table $\beta(M)$ has entry in the $i$-th column and $j$-th row given by $\beta_{i, i+j}$ as below: 
\[
		\begin{tabular}{c|c c c c c c c c}
		& $\cdots$ & 0 & 1 & $\cdots$ & i-1 & i & i+1 & $\cdots$ \\
		\hline
		\vdots & & \vdots & \vdots & & \vdots & \vdots & \vdots\\ 
		-2 & $\cdots$ &$\beta_{0, -2}$&$\beta_{1,-1}$ & $\cdots$ & $\beta_{i-1, i-3}$ & $\beta_{i, i-2}$ & $\beta_{i+1, i-1}$ & $\cdots$ \\ 
		-1 &$\cdots$ & $\beta_{0, -1}$&$\beta_{1,0}$ & $\cdots$ & $\beta_{i-1, i-2}$ & $\beta_{i, i-1}$ & $\beta_{i+1, i}$ & $\cdots$ \\ 
		0 &$\cdots$ & $\beta_{0, 0}$&$\beta_{1,1}$ & $\cdots$ & $\beta_{i-1, i-1}$ & $\beta_{i, i}$ & $\beta_{i+1, i+1}$ & $\cdots$ \\ 
		1 &$\cdots$ & $\beta_{0, 1}$&$\beta_{1,2}$ & $\cdots$ & $\beta_{i-1, i}$ & $\beta_{i, i+1}$ & $\beta_{i+1, i+2}$ & $\cdots$ \\ 
		\vdots & & \vdots & \vdots & & \vdots & \vdots & \vdots\\ 
		\end{tabular}
	\]
	and that subcategories $ \underline{{\rm MCM}}^{\leq n}(A)$ and $ \underline{{\rm MCM}}^{\geq n}(A)$ consists of modules whose Betti table is eventually supported at or above the $n$-th row (resp. at or below the $n$-th row) for $i \gg 0$. Recall also that in the presence of a t-structure, so that cohomology objects $\mathcal{H}^n(M)$ are defined, the \emph{amplitude} of $M$ is defined by 
	\[
		\sup \{ |i - i'| \ | \ \mathcal{H}^i(M) \neq 0 \textup{ and } \mathcal{H}^{i'\!}(M) \neq 0 \}.
	\]\bigskip

\begin{prop}\label{thmApart1} Let $A$ be a Koszul Gorenstein $\Bbbk$-algebra. The following are equivalent: 
	\begin{enumerate}[i)]
	\item $A$ is absolutely Koszul.
	\item $t^{\mathsf{lin}}$ forms a bounded t-structure for $ \underline{{\rm MCM}}^{\mathbb{Z}}(A)$. 
\end{enumerate}
\end{prop} 

\begin{proof} We will want to make use of R\"omer's Theorem \ref{kosmod}, but we first need to show that every $M \in \underline{{\rm MCM}}^{\mathbb{Z}}(A)$ has finite regularity under each of the above hypotheses. For $i)$ this is the statement of Lemma \ref{finitereg}. For $ii)$, note that the cohomology objects $\mathcal{H}^i(M)[-i]$ are eventually $i$-linear and so have finite regularity, and if $n = \inf \{ i \ | \ \mathcal{H}^i(M) \neq 0 \}$ then we have a distinguished triangle
\[
	\tau^{\leq n}M \to M \to \tau^{\geq n+1} M \to \tau^{\leq n}M[1].
\]
The t-cohomology long exact sequence shows that $\mathcal{H}^i(\tau^{\leq n}M) = \mathcal{H}^iM$ for $i \leq n$ and $\mathcal{H}^i(\tau^{\geq n+1}M) = \mathcal{H}^iM$ for $i \geq n+1$. Since $\tau^{\leq n}M \cong \mathcal{H}^n(M)[-n]$ is $n$-linear and $\tau^{\geq n+1}M$ has finite regularity by induction on amplitude, we conclude that $M$ has finite regularity by half-exactness of $ \underline{{\rm Ext}}^i_{gr A}(-, \Bbbk(-i-j))$. We can therefore assume that regularity is finite throughout and apply R\"omer's Theorem. \bigskip 

\noindent {\bf The implication $i) \implies ii)$}. Assume that $A$ is absolutely Koszul. We will show that $t^{\mathsf{lin}}$ defines a bounded t-structure. Recall the axioms of a t-structure:
\begin{enumerate}[T1.]
	\item $\underline{{\rm MCM}}^{\leq 0}(A) \subseteq \underline{{\rm MCM}}^{\leq 1}(A)$ and $ \underline{{\rm MCM}}^{\geq 1}(A) \subseteq \underline{{\rm MCM}}^{\geq 0}(A)$;
	\item ${\rm \underline{Hom}}_{gr A}( \underline{{\rm MCM}}^{\leq 0}(A),\ \! \underline{{\rm MCM}}^{\geq 1}(A)) = 0$;
	\item For all $M$ in $ \underline{{\rm MCM}}^{\mathbb{Z}}(A)$ there is a distinguished triangle $M^{\leq 0} \to M \to M^{\geq 1} \to M^{\leq 0}[1]$ with $M^{\leq 0}$ in $ \underline{{\rm MCM}}^{\leq 0}(A)$ and $M^{\geq 1}$ in $ \underline{{\rm MCM}}^{\geq 1}(A)$.
\end{enumerate}
That T1 holds is immediate and T2 is easily verified by considering the degrees of generators of sufficiently far syzygies. We are left with considering T3. Let $M$ in $ \underline{{\rm MCM}}^{\mathbb{Z}}(A)$ and take $m \gg 0$ so that $\widetilde{M} = \Omega^m(M)(m)$ is Koszul; such an $m$ exists by absolute Koszulity of $A$. Consider the submodule $\widetilde{M}_{\langle \leq 0 \rangle} \subseteq \widetilde{M}$ and, abusing notation by taking $\widetilde{M} \in \grmod A$ to be any representative, consider the short exact sequence in $\grmod A$ 
\[
	0 \to \widetilde{M}_{\langle \leq 0 \rangle} \to \widetilde{M} \to \widetilde{M}/\widetilde{M}_{\langle \leq 0 \rangle} \to 0.
\]
We will construct the distinguished triangle of T3
\[
	M^{\leq 0} \to M \to M^{\geq 1} \to M^{\leq 0}[1]
\]
in the following way: 
\begin{enumerate}[a)]
\item First, we will see that the graded module $\widetilde{M}_{\langle \leq 0 \rangle}$ is also Koszul, so that ${\rm Tor}^A_{i}(\widetilde{M}_{\langle \leq 0 \rangle}, \Bbbk)_{i+j} = 0$ for $j > 0$ and all $i \geq 0$ by Cor. \ref{finitereg}; this point will be established below. 
\item Secondly, one notes that ${\rm Tor}^A_{i}(\widetilde{M}/\widetilde{M}_{\langle \leq 0 \rangle}, \Bbbk)_{i+j} = 0$ for $j < 1$ and all $i \geq 0$ automatically as the module is generated in degrees $\geq 1$. 
\item Third, it is a priori unclear that the above two modules are MCM, but we can replace them by their MCM approximation. As this is an exact functor ${\rm D}^b(\grmod A) \twoheadrightarrow \underline{{\rm MCM}}^{\mathbb{Z}}(A)$, we obtain a distinguished triangle 
\[
	\left(\widetilde{M}_{\langle \leq 0 \rangle}\right)^{st} \to \widetilde{M} \to \left(\widetilde{M}/\widetilde{M}_{\langle \leq 0 \rangle}\right)^{st} \to \widetilde{M}_{\langle \leq 0 \rangle}^{st}[1]
\]
in $ \underline{{\rm MCM}}^{\mathbb{Z}}(A)$. Since MCM approximation leaves the tail of the minimal resolution unchanged, we deduce from the points a) and b) that $\left(\widetilde{M}_{\langle \leq 0 \rangle}\right)^{st} \in \underline{{\rm MCM}}^{\leq 0}(A)$ and $\left(\widetilde{M}/\widetilde{M}_{\langle \leq 0 \rangle}\right)^{st} \in \underline{{\rm MCM}}^{\geq 1}(A)$. 
\item Finally, since $\Omega^{\pm 1}(-)(\pm 1)$ translates Betti tables horizontally, we obtain the required distinguished triangle 
\[
	M^{\leq 0} \to M \to M^{\geq 1} \to M^{\leq 0}[1]
\]
by applying $\Omega^{-m}(-)(-m)$ to the previous triangle. 
	\end{enumerate}
It remains to estalish a), which we now turn to. \bigskip 

We claim that $\widetilde{M}_{\langle \leq 0 \rangle}$ is Koszul and that this follows from R\"omer's Theorem. We show this by induction on the distance of the initial degree $\mathsf{indeg}(\widetilde{M})$ to $0$, and assume that $\mathsf{indeg}(\widetilde{M}) \leq 0$, as otherwise the claim is vacuous.\bigskip 

When $\mathsf{indeg}(\widetilde{M}) = 0$, then $\widetilde{M}_{\langle \leq 0 \rangle} = \widetilde{M}_{\langle 0 \rangle}$ is linear since $\widetilde{M}$ is Koszul. Next assume that $\mathsf{indeg}(\widetilde{M}) < 0$. Define the module $N := \widetilde{M}/\widetilde{M}_{\langle \mathsf{indeg}(\widetilde{M}) \rangle}$, which is also Koszul, but with $\mathsf{indeg}(\widetilde{M}) < \mathsf{indeg}(N) \leq 0$. The module $N_{\langle \leq 0 \rangle}$ is then Koszul by induction. Now, writing $\mathsf{indeg} := \mathsf{indeg}(\widetilde{M}_{\langle \leq 0 \rangle}) = \mathsf{indeg}(\widetilde{M})$ for short, note that: 
\begin{enumerate}[1)] 
\item the module $\widetilde{M}_{\langle \leq 0 \rangle \langle \mathsf{indeg} \rangle} = \widetilde{M}_{\langle \mathsf{indeg} \rangle}$ is $\mathsf{indeg}$-linear since $\widetilde{M}$ is Koszul;
\item the module $N_{\langle \leq 0 \rangle} \cong \widetilde{M}_{\langle \leq 0 \rangle}/\widetilde{M}_{\langle \mathsf{indeg} \rangle} = \widetilde{M}_{\langle \leq 0 \rangle}/\widetilde{M}_{\langle \leq 0 \rangle \langle \mathsf{indeg} \rangle}$ is also Koszul as per the above paragraph.
	\end{enumerate}
Since $\widetilde{M}_{\langle \leq 0 \rangle \langle \mathsf{indeg} \rangle}$ is $\mathsf{indeg}$-linear and $\widetilde{M}_{\langle \leq 0 \rangle}/\widetilde{M}_{\langle \leq 0 \rangle \langle \mathsf{indeg} \rangle}$ is Koszul, it follows from R\"omer's Theorem \ref{kosmod} that $\widetilde{M}_{\langle \leq 0 \rangle}$ is also Koszul, as we wanted to show. This completes the induction, and we have established point a). \bigskip

Thus we have shown that $t^{\mathsf{lin}}$ forms a t-structure on $ \underline{{\rm MCM}}^{\mathbb{Z}}(A)$. We only need to show that it is bounded. The above construction gives the following formula for the truncation functors $\tau^{\leq n}$, $\tau^{\geq n}$ applied to $M \in \underline{{\rm MCM}}^{\mathbb{Z}}(A)$: 
\begin{align*}
	\tau^{\leq n} M &= \Omega^{-m}\left( \left(\Omega^m(M)(m)_{\langle \leq n \rangle}\right)^{st} \right)(-m) \\
	\tau^{\geq n} M &= \Omega^{-m}\left( \left(\Omega^m(M)(m)/\Omega^m(M)(m)_{\langle \leq n-1 \rangle} \right)^{st} \right)(-m) 
\end{align*}
for any $m \gg 0$ such that $\Omega^m(M)(m)$ is Koszul. It is immediate that $\tau^{\leq n} M = 0$ for all $n \ll 0$, and from the above formula we also obtain that $\tau^{\geq n}M = 0$ for all $n > \mathsf{reg}_A(M)$. Hence the t-structure is bounded, and we have shown the implication $i) \implies ii)$. \bigskip
	
\noindent {\bf The implication $ii) \implies i)$}. Assume that $t^{\mathsf{lin}} = ( \underline{{\rm MCM}}^{\leq 0}(A),\ \underline{{\rm MCM}}^{\geq 0}(A))$ forms a bounded t-structure. We will show ${\rm ld}_A(M) < \infty$ for each $M$ in $ \underline{{\rm MCM}}^{\mathbb{Z}}(A)$. Since the t-structure is bounded, we work by induction on the amplitude of $M$. When the amplitude is zero, using the truncation sequences one sees that $M \cong \mathcal{H}^{n}M[-n]$ for some $n \in \Z$, and so $M$ has an eventually $n$-linear minimal resolution and so ${\rm ld}_A(M) <\infty$. When $M$ has positive amplitude, there is a smallest interval $[n, n']$ for which $\mathcal{H}^iM = 0$ for $i \notin [n, n']$. We have a distinguished triangle
\[
	\tau^{\leq n}M \to M \to \tau^{\geq n+1}M \to \tau^{\leq n}M[1]
\]
Then $\tau^{\leq n}M \cong (\mathcal{H}^n M)[-n]$ has an eventually $n$-linear minimal free resolution and $\tau^{\geq n+1}M$ has a Koszul syzygy module by induction. We may apply the functor $\Omega^{m}(-)(m)$ for $m \gg 0$ to obtain 
\[
	\Omega^m(\tau^{\leq n}M)(m) \to \Omega^mM(m) \to \Omega^m(\tau^{\geq n+1}M)(m) \to \Omega^m(\tau^{\leq n}M)(m)[1]
\]
with the first module $n$-linear and the third module Koszul, or equivalently componentwise linear. Now, since $\Omega^m(-)(m)$ preserves the categories $ \underline{{\rm MCM}}^{\leq i}(A)$ and $ \underline{{\rm MCM}}^{\geq i}(A)$ for each $i$, it is easy to see that $\Omega^m(-)(m)$ commutes with the truncations $\tau^{\leq i}$, $\tau^{\geq i}$. Writing $\widetilde{M} = \Omega^m(M)(m)$ for short, we may rewrite the above triangle as
\[
	\tau^{\leq n} \widetilde{M} \to \widetilde{M} \to \tau^{\geq n+1} \widetilde{M} \to \tau^{\leq n} \widetilde{M}[1]
\]
Our aim is to deduce that $\widetilde{M}$ is Koszul from condition $(ii)$ in R\"omer's Theorem \ref{kosmod}. That is, we would like to compare the above to the short exact sequence 
\[
	0 \to \widetilde{M}_{\langle \mathsf{indeg}(\widetilde{M}) \rangle} \to \widetilde{M} \to \widetilde{M}/\widetilde{M}_{\langle \mathsf{indeg}(\widetilde{M}) \rangle} \to 0	
\]
and deduce that the first term is $\mathsf{indeg}(\widetilde{M})$-linear with third term componentwise linear. \bigskip

Possibly increasing $m$ and removing projective summands from $\widetilde{M}$, we will show three points: 
\begin{enumerate}[a)]
\item We have $\mathsf{indeg}(\widetilde{M}) = n = \inf \{  i \ | \ \mathcal{H}^i(M) \neq 0 \}$. 
\item The modules $\widetilde{M}_{\langle n \rangle}$ and $\widetilde{M}/\widetilde{M}_{\langle n \rangle}$ are MCM. 
\item There is a comparison map of distinguished triangles in $ \underline{{\rm MCM}}^{\mathbb{Z}}(A)$
\[
	\xymatrix{ \tau^{\leq n}\widetilde{M} \ar[r]^-{\iota} \ar[d]^{\widetilde{\iota}}_{\cong} & \widetilde{M} \ar@{=}[d] \ar[r] & \tau^{\geq n+1}\widetilde{M} \ar[d]^{\varphi}_{\cong} \ar[r] & \tau^{\leq n}\widetilde{M}[1]\ar[d] \\
	\widetilde{M}_{\langle n \rangle} \ar[r] & \widetilde{M}  \ar[r] & \widetilde{M}/\widetilde{M}_{\langle n \rangle} \ar[r] &\widetilde{M}_{\langle n \rangle}[1]  }
\]
which is an isomorphism. 
	\end{enumerate}
	Once this is established, we will deduce that $\widetilde{M}_{\langle n \rangle}$ is $n$-linear and $\widetilde{M}/\widetilde{M}_{\langle n \rangle}$ is Koszul, and so that $\widetilde{M}$ itself is Koszul by R\"omer's Theorem \ref{kosmod}. Since $M$ has the Koszul syzygy $\widetilde{M} = \Omega^m(M)(m)$, we obtain that ${\rm ld}_A(M) \leq m < \infty$. \bigskip 

The proof of a) is straightforward but b) and c) will be more involved. To see a), note that the table $\beta(\widetilde{M})$ is supported between rows $n, n'$ as $i \gg 0$ and $[n,n']$ is the smallest interval with that property. Take $m$ large enough so that this holds at $i \geq 0$ and removing potential projective summands from $\widetilde{M}$, we have $\widetilde{M}_{< n} = 0$ and $\widetilde{M}_{\langle n \rangle} \neq 0$. \bigskip 

We will prove b) and c) at the same time by playing them off each others, by gradually constructing approximations to the comparison map of c), possibly increasing $m$ along the way. The primary hurdle is that $\widetilde{M}_{\langle n \rangle}$ and $\widetilde{M}/\widetilde{M}_{\langle n \rangle}$ are not a priori MCM, and so we will need to take some of our arguments outside of the category $ \underline{{\rm MCM}}^{\mathbb{Z}}(A)$ and into the larger stable category $ \underline{\mathsf{mod}}^{\Z}\ \! A$, and eventually replace them by large syzygies. \bigskip 

We first note the standard isomorphism of distinguished triangles 
\[
	\xymatrix{
		\tau^{\leq n}\widetilde{M} \ar@{=}[d] \ar[r]^-{\iota} & \widetilde{M} \ar@{=}[d] \ar[r] & \tau^{\geq n+1}\widetilde{M} \ar[d]^-{\cong} \ar[r] & \tau^{\leq n}\widetilde{M}[1] \ar@{=}[d]\\ 
	\tau^{\leq n}\widetilde{M} \ar[r]^-{\iota} & \widetilde{M} \ar[r] & {\rm Cone}(\iota) \ar[r] & \tau^{\leq n}\widetilde{M}[1].}
\]
We can model the cone triangle as coming from a short exact sequence in ${\rm MCM}^{\Z}(A)$, of the form 
\[
	0 \to [\tau^{\leq n}\widetilde{M}] \xrightarrow{\iota} \widetilde{M} \oplus P \to {\rm Cone}(\iota) \to 0
\]
where $P$ is some projective module, and recall that we write $[N]$ for the unique module representative of $N \in \underline{{\rm MCM}}^{\mathbb{Z}}(A)$ without projective summands. Since $[\tau^{\leq n}\widetilde{M}]$ is $n$-linear and so generated in degree $n$, we have a comparison map of short exact sequences 
\begin{equation}\label{comparisonshortexactseq}
\xymatrix{ 
	0 \ar[r] & [\tau^{\leq n}\widetilde{M}] \ar[d]^-{\widetilde{\iota}} \ar[r]^-{\iota} & \widetilde{M}\oplus P \ar@{->>}[d]^-{\pi} \ar[r] & {\rm Cone}(\iota) \ar[d] \ar[r] & 0\\ 
	0 \ar[r] & \widetilde{M}_{\langle n \rangle} \ar[r] & \widetilde{M} \ar[r] & \widetilde{M}/\widetilde{M}_{\langle n \rangle} \ar[r] & 0\\ 
}
\end{equation}
where $\pi \circ \iota$ factors through $\widetilde{\iota}$ for degree reasons and the third horizontal map follows from surjectivity of $\pi$. \bigskip

Next, for any $Y \in \grmod A$, consider the $k$-linear functor ${\rm \underline{Hom}}_{gr A}(-, Y): (\grmod A)^{op} \to \Modsf k$. We claim that ${\rm \underline{Hom}}_{gr A}(-, Y)$ is half-exact on $\grmod A$: by results of Auslander-Bridger \cite[Thm. 1.40]{ABBook} (see also bottom of p.47 in \cite{ABBook}), for any $X, Y \in \grmod A$ we have isomorphisms 
\[
	{\rm \underline{Hom}}_{gr A}(X, Y) \cong {\rm Nat}\left({\rm Ext}^1_{gr A}(Y, -),\ \! {\rm Ext}^1_{gr A}(X, -) \right),
\] 
natural in $X, Y$, onto the space of natural transformations between Ext functors. Using that ${\rm Nat}({\rm Ext}^1_{gr A}(Y, -), -)$ is right-exact on the abelian category of $k$-linear functors ${\rm Fun}\left((\grmod A)^{op},\ \! \Modsf k \right)$ (with termwise abelian structure), combined with half-exactness of $X \mapsto {\rm Ext}^1_{gr A}(X, -)$, we obtain that ${\rm \underline{Hom}}_{gr A}(-, Y)$ is half-exact. \bigskip

Taking $Y = \Bbbk(-n)$, from (\ref{comparisonshortexactseq}) we obtain a commutative diagram of exact sequences 
\[
	\xymatrix@C=10pt{ {\rm \underline{Hom}}_{gr A}({\rm Cone}(\iota), \Bbbk(-n)) \ar[r] &  
		{\rm \underline{Hom}}_{gr A}(\widetilde{M}\oplus P, \Bbbk(-n)) \ar[r] &  
		{\rm \underline{Hom}}_{gr A}(\tau^{\leq n}\widetilde{M}, \Bbbk(-n)) \ar[r] &  
		{\rm \underline{Ext}}_{gr A}^1({\rm Cone}(\iota), \Bbbk(-n)) \\
		{\rm \underline{Hom}}_{gr A}(\widetilde{M}/\widetilde{M}_{\langle n \rangle}, \Bbbk(-n)) \ar[r] \ar[u] &  
		{\rm \underline{Hom}}_{gr A}(\widetilde{M}, \Bbbk(-n)) \ar[r]^-{\mathsf{res}} \ar[u]^-{\cong}_-{\pi^*} &  
		{\rm \underline{Hom}}_{gr A}(\widetilde{M}_{\langle n \rangle}, \Bbbk(-n)) \ar[u]_-{\widetilde{\iota}^*} 
	}
\]
The top row is exact as it is part of the long exact sequence of Tate cohomology, while the bottom row is exact by Auslander-Bridger. We claim that the three corner are zero: 
\[
	{\rm \underline{Hom}}_{gr A}({\rm Cone}(\iota), \Bbbk(-n)) = {\rm \underline{Ext}}_{gr A}^1({\rm Cone}(\iota), \Bbbk(-n)) = {\rm \underline{Hom}}_{gr A}({\rm Cone}(\iota), \Bbbk(-n)) = 0.
\] 
That the first two vanish follows from noting that:
\begin{enumerate}[1)]
\item ${\rm Cone}(\iota) \cong \tau^{\geq n+1}\widetilde{M}$ in $ \underline{{\rm MCM}}^{\mathbb{Z}}(A)$;
\item The syzygy level $m$ of $\widetilde{M} = \Omega^m(M)(m)$ was taken sufficiently high so that $\tau^{\geq n+1}\widetilde{M}$ is Koszul, hence has Betti table supported supported at or below the $(n+1)$-th row for all $i \geq 0$ by Lemma \ref{finitereg}.
	\end{enumerate}
The vanishing on the bottom row is clear. We conclude by a simple diagram chase that the maps $\mathsf{res}$ and $\widetilde{\iota}^* = {\rm \underline{Hom}}_{gr A}(\widetilde{\iota},\ \! \Bbbk(-n))$ are isomorphisms. \bigskip 

We claim that $\widetilde{\iota}: \tau^{\leq n}\widetilde{M} \to \widetilde{M}_{\langle n \rangle}$ must then be a stable isomorphism. Recall that by Lemma \ref{unstablelemma}, for any $N \in \underline{\mathsf{mod}}^{\Z} \ \! A$ we have 
\[
	\underline{{\rm Hom}}_A(N, \Bbbk) = \bigoplus_{j \in \Z} \underline{{\rm Hom}}_{gr A}(N, \Bbbk(j))
	= \bigoplus_{j \in \Z} {\rm Hom}_{gr A}([N], \Bbbk(j))
	= {\rm Hom}_{A}([N], \Bbbk). 
\]
Since $\widetilde{\iota}: \tau^{\leq n}\widetilde{M} \to \widetilde{M}_{\langle n \rangle}$ is a homomorphism between $n$-generated modules, it induces isomorphisms 
\[
	\xymatrix@R=10pt@C=40pt{
		\underline{{\rm Hom}}_{gr A}(\widetilde{M}_{\langle n \rangle}, \Bbbk(-n)) \ar@{=}[d] \ar[r]_-{\cong}^-{\widetilde{\iota}^*} & \underline{{\rm Hom}}_{gr A}(\tau^{\leq n}\widetilde{M}, \Bbbk(-n)) \ar@{=}[d]  \\
	\underline{{\rm Hom}}_A(\widetilde{M}_{\langle n \rangle}, \Bbbk) \ar[r]_-{\cong}^-{ {\rm \underline{Hom}}_A(\widetilde{\iota}, \Bbbk)} & \underline{{\rm Hom}}_A(\tau^{\leq n}\widetilde{M}, \Bbbk) \\
{\rm Hom}_A\left([\widetilde{M}_{\langle n \rangle}], \Bbbk \right) \ar[r]_-{\cong}^-{ {\rm Hom}_{A}([\widetilde{\iota}], \Bbbk)} \ar@{=}[u] & {\rm Hom}_A\left([\tau^{\leq n}\widetilde{M}], \Bbbk \right) \ar@{=}[u] \\
		{\rm Hom}_\Bbbk([\widetilde{M}_{\langle n \rangle}]\otimes_A \Bbbk, \Bbbk) \ar[r]_-{\cong}^-{ {\rm Hom}_{\Bbbk}([\widetilde{\iota}], \Bbbk)} \ar@{=}[u] & {\rm Hom}_\Bbbk\left([\tau^{\leq n}\widetilde{M}]\otimes_A \Bbbk, \Bbbk \right) \ar@{=}[u] 
}
\]
By the Nakayama Lemma, the map 
\[
[\widetilde{\iota}]:  [\tau^{\leq n}\widetilde{M}] \to [\widetilde{M}_{\langle n \rangle}]  
\]induced from $\widetilde{\iota}$ is an isomorphism and thus $\widetilde{\iota}$ is an isomorphism in $ \underline{\mathsf{mod}}^{\Z}\ \! A$. As consequence, we obtain that: 
\begin{enumerate}[1)]
\item $\widetilde{M}_{\langle n \rangle}$ is an MCM module, as it is stably isomorphic to one;
\item $\widetilde{M}_{\langle n \rangle}$ is $n$-linear, as it is $n$-generated and stably isomorphic to an $n$-linear module. 
\end{enumerate}
In particular we have established half of claim b). We then turn to the construction of the comparison map of claim c). \bigskip 

Next, the comparison map of short exact sequence (\ref{comparisonshortexactseq}) gives rise to a comparison map of distinguished triangles in $ {\rm D}^b(\grmod A)$, and taking MCM approximation gives rise to a comparison map
\[
	\xymatrix{ \tau^{\leq n}\widetilde{M} \ar[r]^-{\iota} \ar[d]^{\widetilde{\iota}}_{\cong} & \widetilde{M} \ar@{=}[d] \ar[r] & \tau^{\geq n+1}\widetilde{M} \ar[d]^{\varphi} \ar[r] & \tau^{\leq n}\widetilde{M}[1]\ar[d]_-{\cong} \\
	\widetilde{M}_{\langle n \rangle} \ar[r] & \widetilde{M}  \ar[r] & \left(\widetilde{M}/\widetilde{M}_{\langle n \rangle}\right)^{st} \ar[r] &\widetilde{M}_{\langle n \rangle}[1] }
\]
in $ \underline{{\rm MCM}}^{\mathbb{Z}}(A)$, where we write $\varphi$ for the remaining vertical map. We have shown that $\widetilde{\iota}$ is an isomorphism, and so $\varphi$ also must be one. We deduce that $\left(\widetilde{M}/\widetilde{M}_{\langle n \rangle}\right)^{st} \cong \tau^{\geq n+1}\widetilde{M}$ is Koszul; in particular, any sufficiently far syzygy of $\widetilde{M}/\widetilde{M}_{\langle n \rangle}$ is both Koszul and MCM. This is close to our wanted claim, that $\widetilde{M}/\widetilde{M}_{\langle n \rangle}$ be Koszul and MCM, and so we will aim to replace $\widetilde{M}/\widetilde{M}_{\langle n \rangle}$ by a large syzygy of the same form. \bigskip

We will fix this by increasing the syzygy level $m$ of $\widetilde{M} = \Omega^m(M)(m)$ once more. Let $F_* \xrightarrow{\sim} \widetilde{M}$ be the minimal resolution of $\widetilde{M}$. Recall that by Lemma \ref{linearstrandlemma}, the minimal graded projective resolutions of $\widetilde{M}_{\langle n \rangle} \subseteq \widetilde{M}$ is given by the $n$-th linear strand $F_*^{(n)} \subseteq F_*$, and so we have a short exact sequence of minimal resolutions 
\[
\xymatrix{
	0 \ar[r] & F_*^{(n)} \ar[r] \ar[d] & F_* \ar[r] \ar[d] & F_*/F_*^{(n)} \ar[r] \ar[d] & 0 \\
	0 \ar[r] & \widetilde{M}_{\langle n \rangle} \ar[r] & \widetilde{M} \ar[r] & \widetilde{M}/\widetilde{M}_{\langle n \rangle} \ar[r] & 0.
}
\]
Truncating the resolutions at homological degree $m' \gg 0$ and twisting by $(m')$, we obtain a short exact sequence of resolutions of the $m'$-th syzygies
\[
\xymatrix{
	0 \ar[r] & F_{\geq m'}^{(n)}(m') \ar[r] \ar[d] & F_{\geq m'}(m') \ar[r] \ar[d] & F_{\geq m'}(m')/F_{\geq m'}^{(n)}(m') \ar[r] \ar[d] & 0 \\
0 \ar[r] & \Omega^{m'}\widetilde{M}_{\langle n \rangle}(m') \ar[r] & \Omega^{m'}\widetilde{M}(m') \ar[r] & \Omega^{m'}\left(\widetilde{M})/\widetilde{M}_{\langle n \rangle}\right)(m') \ar[r] & 0.
}
\]
We make a few observations: 
\begin{enumerate}[1)]
\item Since $\widetilde{M}_{\langle n \rangle}$ is $n$-linear, so is $\Omega^{m'}(\widetilde{M}_{\langle n \rangle})(m')$. 
\item Since $F_{\geq m'}^{(n)}(m')$ is the $n$-th linear strand of $F_{\geq m'}(m')$ by construction, and we deduce that 
	\[
		\Omega^{m'}(\widetilde{M}_{\langle n \rangle})(m') \cong \left(\Omega^{m'}(\widetilde{M})(m')\right)_{\langle n \rangle}.
	\]
\item From the short exact sequence above we obtain that 
	\[
		\Omega^{m'}\left(\widetilde{M}/\widetilde{M}_{\langle n \rangle}\right)(m') \cong \left(\Omega^{m'}(\widetilde{M})(m')\right)/\left(\Omega^{m'}(\widetilde{M})(m')\right)_{\langle n \rangle}.
	\] 
	Note that for $m' \gg 0$ sufficiently large, the left-hand side is MCM and Koszul. 
	\end{enumerate}

Finally, replacing replacing $\Omega^m(M)(m)$ by $\Omega^{m+m'}(M)(m+m')$ (so that we henceforth write $\widetilde{M} = \Omega^{m+m'}(M)(m+m')$), the points 1)-2)-3) above show that $\widetilde{M}_{\langle n \rangle}$ is an $n$-linear MCM module and that $\widetilde{M}/\widetilde{M}_{\langle n \rangle}$ is a Koszul MCM module, establishing our previous claim b). This is enough to conclude that $\widetilde{M}$ is Koszul by R\"omer's Theorem \ref{kosmod}. Note that we also obtain claim c), since the comparison map of short exact sequence (\ref{comparisonshortexactseq}) (at our increased syzygy level) actually gives rise an isomorphism of distinguished triangles 
\[
	\xymatrix{ \tau^{\leq n}\widetilde{M} \ar[r]^-{\iota} \ar[d]^{\widetilde{\iota}}_{\cong} & \widetilde{M} \ar@{=}[d] \ar[r] & \tau^{\geq n+1}\widetilde{M} \ar[d]^-{\varphi}_-{\cong} \ar[r] & \tau^{\leq n}\widetilde{M}[1]\ar[d]_-{\cong} \\ 
	\widetilde{M}_{\langle n \rangle} \ar[r] & \widetilde{M}  \ar[r] & \widetilde{M}/\widetilde{M}_{\langle n \rangle} \ar[r] &\widetilde{M}_{\langle n \rangle}[1]. }
\]

As we have shown that $\widetilde{M} = \Omega^{m+m'}(M)(m+m')$ is Koszul, $M$ has a Koszul syzygy and so ${\rm ld}_A(M) \leq m+m' < \infty$. This establishes the implication $ii) \implies i)$. 
\end{proof} \bigskip

\begin{rem}
When $A$ is commutative, since finiteness of regularity for $M \in \underline{{\rm MCM}}^{\mathbb{Z}}(A)$ follows from Koszulity of $A$, condition $ii)$ can be replaced by the weaker conditition that $t^{\mathsf{lin}}$ forms a t-structure, which will then necessarily be bounded.
\end{rem}

This proposition provides a categorical interpretation of absolute Koszulity. In particular, this shows that $t^{\mathsf{lin}}$ fails to form a t-structure for sufficiently pathological Gorenstein Koszul algebras, notably if the algebra is bad in the sense of Roos. We will construct Koszul Gorenstein algebras which fail to be absolutely Koszul in Chapter 6. 

\subsection*{Properties of the t-structure} Let $A$ be absolutely Koszul Gorenstein. We have shown that $t^{\mathsf{lin}}$ forms a bounded t-structure on $ \underline{{\rm MCM}}^{\mathbb{Z}}(A)$ in the proof of Prop. \ref{thmApart1}. In this subsection, we simply record some of its main properties, which were either implicit or sometimes explicit in the proof of Prop. \ref{thmApart1}. \bigskip 

For any $M \in \underline{{\rm MCM}}^{\mathbb{Z}}(A)$, the truncation functors $\tau^{\leq n}$, $\tau^{\geq n}$ are defined by passing to a Koszul syzygy $\widetilde{M} = \Omega^m(M)(m)$ for $m \geq {\rm ld}_A(M)$, and applying $\Omega^{-m}(-)(-m)$ to the distinguished triangle
\begin{eqnarray}\label{pretruncationtriangle}
\left(\widetilde{M}_{\langle \leq n \rangle}\right)^{st} \to \widetilde{M} \to \left(\widetilde{M}/\widetilde{M}_{\langle \leq n \rangle}\right)^{st} \to \left(\widetilde{M}_{\langle \leq n \rangle}\right)^{st}[1] 
\end{eqnarray}
to obtain
\begin{eqnarray}
	\tau^{\leq n}M \to M \to \tau^{\geq n+1}M \to \tau^{\leq n}M[1]
\end{eqnarray}
with $\tau^{\leq n}M \in \underline{{\rm MCM}}^{\leq n}(A)$ and $\tau^{\geq n+1}M \in \underline{{\rm MCM}}^{\geq n+1}(A)$. This is independent of the choice of Koszul syzygy, since $\tau^{\leq n}$ and $\tau^{\geq n}$ compute the right and left adjoints to the corresponding inclusions of full subcategories (see \cite[Prop. 8.1.4]{HTT}). \bigskip

Note that by taking $m \geq {\rm ld}_A(M) + \dim A$ instead we may guarantee that both $\widetilde{M}_{\langle n \rangle}$ and $\widetilde{M}/\widetilde{M}_{\langle n \rangle}$ are MCM (as shown in the proof of the implication $ii) \implies i)$ in Prop. \ref{thmApart1}) and thus omit MCM approximation in the distinguished triangle \ref{pretruncationtriangle}. \bigskip 

Let $\T$ be a triangulated category with a $t$-structure $t = (\T^{\leq 0}, \T^{\geq 0})$ and a triangulated functor $F: \T \to \T$. We say that $F$ is $t$-exact (or exact when $t$ is understood) if $F(\T^{\leq 0}) \subseteq \T^{\leq 0}$ and $F(\T^{\geq 0}) \subseteq \T^{\geq 0}$. It follows that $F$ restricts to an exact functor on the abelian heart $\T^{\leq 0} \cap \T^{\geq 0}$, and commutes with taking cohomology (see \cite[Sect. 8.1]{HTT})
\[
	F \mathcal{H}^n(X) \cong \mathcal{H}^n F(X).
\]

Next, for any $M \in \underline{{\rm MCM}}^{\mathbb{Z}}(A)$ let $F_* \xrightarrow{\sim} M$ denote its minimal resolution. We have seen that the linear part of $F_*$ decomposes into linear strands  
\[
	{\rm lin}^A(F_*) = \bigoplus_{i \in \Z} F_*^{(i)}
\]
where $F_*^{(i)}$ is an $i$-linear complex, meaning with $n$-th term generated in degree $n+i$. Note that since ${\rm ld}_A(M) < \infty$ each linear strand has bounded cohomology. The next proposition follows readily from what we have shown.
\begin{prop}\label{propertiesoftlin} Let $A$ be absolutely Koszul Gorenstein. The following properties hold for $t = t^{\mathsf{lin}}$:
	\begin{enumerate}[i)]
	\item The autoequivalence $\Omega^1(-)(1)$ is $t$-exact;
	\item Let $M$ be a Koszul MCM module over $A$. Then ${\rm H}_0( {\rm lin}^A(F_*)) \cong \mathcal{H}^*(M) = \bigoplus_{n \in \Z} \mathcal{H}^n(M)[-n]$;
	\item Let $M$ be an MCM module over $A$. Then we have $\mathcal{H}^i(M)[-i] \cong (F_*^{(i)})^{st}$, the MCM approximation of the corresponding linear strand. In particular $\mathcal{H}^*(M) \cong {\rm lin}^A(F_*)^{st}$. 
\end{enumerate}
\end{prop}\bigskip

We say that $M \in \underline{{\rm MCM}}^{\mathbb{Z}}(A)$ is $t$-formal if $M \cong \mathcal{H}^*(M)$. Hence we see from Prop. \ref{propertiesoftlin} iii) that $t$-formality is equivalent to $M$ admitting a minimal resolution $F_* = (F_*, \partial)$ with eventually purely linear differential, that is $\partial = \partial_1$ on $F_{\geq m}$ for $m \gg 0$. 

\subsection*{The realisation functor}
Now assume that $A$ is an absolutely Koszul Gorenstein $\Bbbk$-algebra, so that $\mathcal{H}^{\mathsf{lin}}(A) \subseteq \underline{{\rm MCM}}^{\mathbb{Z}}(A)$ is the heart of a bounded t-structure. The category $\mathcal{H}^{\mathsf{lin}}(A)$ is then an abelian category and it is natural to expect the existence of a realisation functor between triangulated categories 
\[
	\mathsf{real}: {\rm D}^b(\mathcal{H}^{\mathsf{lin}}(A)) \to \underline{{\rm MCM}}^{\mathbb{Z}}(A)
\]
restricting to the inclusion on $\mathcal{H}^{\mathsf{lin}}(A)$. Given a triangulated category $\D$ with a t-structure $t = (\D^{\leq 0}, \D^{\geq 0})$, the construction of such a realisation functor was given by Beilinson in terms of a choice of \emph{filtered triangulated category} $(\D F, \theta)$ over $\D$, generalising the notion of the filtered derived category of an abelian category. We review this construction and show that it applies to the case at hand. Once shown to exist, standard criteria will then apply to show that $\mathsf{real}$ is fully faithful and essentially surjective, and therefore an equivalence of triangulated categories. We follow the exposition of \cite[Sect. 3]{PJ}, based on results of Beilinson \cite[Appendix]{BeilinsonPS} and Beilinson-Bernstein-Deligne \cite{BBD}. We write $\D$ for a general triangulated category throughout.

\begin{defn}[Beilinson] A filtered triangulated category $\D F$ (f-category for short) is a triangulated category equipped with two full triangulated subcategories $\D F(\leq 0)$ and $\D F(\geq 0)$, an exact automorphism $s: \D F \to \D F$ (called $f$-shift) and a natural transformation $\alpha: {\rm id} \to s$ such that $\D F(\leq n) := s^n \D F(\leq 0)$ and $\D F(\geq n) := s^n\D F(\geq 0)$ satisfy the following axioms:
	\begin{enumerate}[i)]
	\item $\D F (\geq 1) \subseteq \D F (\geq 0)$, $\D F(\leq 1) \supseteq \D F(\leq 0)$ and $\bigcup_{n \in \Z} \D F(\geq n) = \bigcup_{n \in \Z} \D F(\leq n) = \D F$.
	\item For $X \in \D F$, we have $\alpha_X = s(\alpha_{s^{-1}X})$.
	\item For $X \in \D F(\geq 1)$ and $Y \in \D F(\leq 0)$, we have ${\rm Hom}(X, Y) = 0$ and bijections ${\rm Hom}(Y, X) \cong {\rm Hom}(Y, s^{-1}X) \cong {\rm Hom}(sY, X)$ induced by $\alpha$.
	\item For any $X \in \D F$, there is a distinguished triangle $Y \to X \to Z \to Y[1]$ with $Y \in \D F(\geq 1)$ and $Z \in \D F(\leq 0)$.
\end{enumerate}
\end{defn}

\begin{defn} Let $\D$ be a triangulated category. An $f$-category $(\D F, \theta)$ over $\D$ is the data of an f-category $\D F$ along with an equivalence of triangulated categories $\theta: \D \xrightarrow{\cong} \D F(\leq 0) \cap \D F(\geq 0)$.
\end{defn}

\begin{exmp}[Filtered derived category {\cite[Ex. A.2]{Bei}}] Let $\A$ be an abelian category. Define $CF(\A)$ to be the category of complexes $X = (X, F)$ in $\A$ equipped with a finite decreasing filtration by subcomplexes 
	\[
		F:\quad  X = F^p X \supseteq F^{p+1} X \supseteq \dots \supseteq F^{q-1} X \supseteq F^{q}X = 0, \quad p,q \in \Z \textup{ with } p \leq q.
	\]
	The morphisms in $CF(\A)$ are the chain-maps respecting the filtration. Define ${\rm gr}^i_F(X) = F^iX/F^{i+1}X$ and ${\rm gr}_F(X) = \bigoplus_{i \in \Z} {\rm gr}^i_F(X)$. A morphism $f: X \to Y$ is a filtered quasi-isomorphism if the morphism on associated graded complexes ${\rm gr}(f): {\rm gr}(X) \xrightarrow{\sim} {\rm gr}(Y)$ is a quasi-isomorphism. The filtered derived category $\D F(\A)$ is the localisation of $CF(\A)$ at the filtered quasi-isomorphism. \bigskip 

	The category $\D F(\A)$ is an f-category over $\D(\A)$: we have $s(X,F) = (X, sF)$ where $(sF)^iX = F^{i-1}X$, and $\alpha_X: X \to sX$ is induced by $F^iX \subseteq F^{i-1}X$. One then takes 
	\[
		\D F(\A)(\leq n) = \{ (X, F) \ | \ {\rm gr}^i_F(X) = 0 \textup{ for } i > n \}
	\]
	\[
		\D F(\A)(\geq n) = \{ (X, F) \ | \ {\rm gr}^i_F(X) = 0 \textup{ for } i < n \}
	\]
	and the equivalence $\theta: \D(\A) \xrightarrow{\cong} \D F(\A)(\leq 0) \cap \D F(\A) (\geq 0)$ sends a complex $X \in \D(\A)$ to $(X, Tr)$ equipped with the trivial filtration
	\[
		Tr:\quad  X = F^0 X \supseteq F^{1}X = 0.
	\]
The verifications that $i) - iii)$ hold is straightforward. For $iv)$, we have truncation functors
\begin{align*} 
	&\sigma_{\leq n}: \D F(\A) \to \D F(\A)(\leq n)\\
	&\sigma_{\geq n}: \D F(\A) \to \D F(\A)(\geq n)
\end{align*}	
defined by
\begin{align*}
	&\sigma_{\leq n}X := X/F^{n+1}X\\
	&\sigma_{\geq n}X := F^nX 
\end{align*}
with the induced filtrations. We then have a triangle
\[
	\sigma_{\geq 1}X \to X \to \sigma_{\leq 0}X \to \sigma_{\geq 1}X[1]
\]
verifying $iv)$.\qed
 \end{exmp}
Note that since filtered quasi-isomorphisms are quasi-isomorphisms (\cite[Thm. 5.5.11]{Weibel}), we obtain a forgetful functor $\omega: \D F(\A) \to \D(\A)$. 
	 
\begin{exmp}[Filtered homotopy category] One can analogously define the filtered triangulated category $\mathcal{K}F (\A)$ over $\mathcal{K} \A$ consisting of complexes equipped with finite decreasing filtrations, and filtered chain-maps up to filtration-preserving homotopies. The remaining data is as above, and $\mathcal{K}F(\A)$ is a filtered triangulated category over $\mathcal{K}(\A)$.\qed
 \end{exmp}
	 
 For the question at hand, we will make use of the following. Take $\A = \grmod R$ for a two-sided coherent graded Gorenstein $\Bbbk$-algebra $R$, and let ${\rm proj}^\Z\ \! R$ be the full subcategory of finitely generated projective graded modules. 
 \begin{exmp}[Filtered singularity category]\label{filtsingcat} Define $\mathcal{K}F_{ac}({\rm proj}^\Z\ \! R) \subseteq \mathcal{K}F(\grmod R)$ as the full subcategory of acyclic complexes $C = (C, F)$ of finite projectives equipped with finite decreasing filtrations such that the associated graded $\mathsf{gr}_F(C)$ complex is also acyclic. \bigskip 
	 
	 Then $\mathcal{K}F_{ac}({\rm proj}^\Z \ \!R)$ is an $f$-category over $\mathcal{K}_{ac}({\rm proj}^\Z \ \!R)$, with same data as above. The verifications are the same as above, with only difference the question of whether $\sigma_{\leq n}C = C/F^{n+1}C$ and $\sigma_{\geq n}C = F^nC$ are acyclic. But this holds since both complexes are equipped with finite filtrations whose associated graded is acyclic. \qed
 \end{exmp}\bigskip
 Note that by transport of structure under the equivalence $ \underline{{\rm MCM}}^{\mathbb{Z}}(R) \cong \mathcal{K}_{ac}({\rm proj}^\Z \ \! R)$, we obtain an $f$-category over $ \underline{{\rm MCM}}^{\mathbb{Z}}(R)$. \bigskip

 The usefulness of $f$-categories is summarised by the following proposition. Recall from \cite[Sect. 8.1]{HTT} that an exact functor $F: \D \to \D'$ between triangulated categories equipped with $t$-structures $(\D^{\leq 0}, \D^{\geq 0})$ and $(\D'^{\leq 0}, \D'^{\geq 0})$ is $t$-exact if $F(\D^{\leq n}) \subseteq \D'^{\leq n}$ and $F(\D^{\geq n}) \subseteq \D'^{\geq n}$ for each $n \in \Z$.
 \begin{prop}[{\cite[Prop. 3.3]{PJ}}] Let $(\D F, \theta)$ be an $f$-category over $\D$. 
	 \begin{enumerate}[i)]
	 \item There is an exact functor $\omega: \D F \to \D$ such that $\omega \circ \alpha_X: \omega X \xrightarrow{\cong} \omega sX$ is an isomorphism (Forgetful functor).
	 \item Given a $t$-structure $t = (\D^{\leq 0}, \D^{\geq 0})$, there is a unique $t$-structure $\widetilde{t} = (\D F^{\leq 0}, \D F^{\geq 0})$ for which $\theta$ is a $t$-exact functor and $s\D F^{\leq 0} \subseteq \D F^{\leq -1}$. 
		 
 \end{enumerate}
 Moreover, if $\CC = \D^{\leq 0} \cap \D^{\geq 0}$ is the heart of this t-structure, then there is a canonical equivalence of the new heart $\D F^{\leq 0} \cap \D F^{\geq 0} \cong C^b(\CC)$ with the category of bounded complexes over $\CC$.
 \end{prop}
 See \cite[Rem. 3.4]{PJ} for an explicit construction of this t-structure and the corresponding equivalence. \bigskip

 We can then construct the realisation functor of a t-structure $t = (\D^{\leq 0}, \D^{\geq 0})$ on $\D$. Let $(\D F, \theta)$ be an $f$-category over $\D$, and let $G: C^b(\CC) \xrightarrow{\cong} \D F^{\leq 0} \cap \D F^{\geq 0} \subseteq \D F$ be the embedding induced from the equivalence of the previous proposition. Let $Q: C^b(\CC) \to {\rm D}^b(\CC)$ be the localisation functor. 

 \begin{thm}[{\cite[Thm 3.11]{PJ}}]\label{realthm} Given the above data, there is a unique functor $\mathsf{real}_{\D F}$ fitting into a commutative diagram
	 \[
		 \xymatrix{ 
			 C^b(\CC) \ar[d]_G \ar[r]^Q &{\rm D}^b(\CC)\ar[ddl]^{\mathsf{real_{\D F}}} \\
			 \D F\ar[d]_{\omega} \\
		 \D }
	 \]
	 Moreover, $\mathsf{real}_{\D F}: {\rm D}^b(\CC) \to \D$ is an exact functor satisfying the following properties:
	 \begin{enumerate}[i)]
	 \item $\mathsf{real}_{\D F}$ restricts to the identity on $\CC$ and is $t$-exact with respect to the standard $t$-structure on ${\rm D}^b(\CC)$ and $t = (\D^{\leq 0}, \D^{\geq 0})$ on $\D$. 
	 \item $\mathsf{real}_{\D F}$ induces isomorphisms ${\rm Hom}_{ {\rm D}^b(\CC)}(X, Y[n]) \cong {\rm Hom}_{\D}(X, Y[n])$ for any $X, Y \in \CC$ and $n \leq 1$. 
	\item The following are equivalent:
		\begin{enumerate}[a)]
		\item $\mathsf{real}_{\D F}$ is fully faithful;
		\item $\mathsf{real}_{\D F}$ induces isomorphisms ${\rm Hom}_{ {\rm D}^b(\CC)}(X, Y[n]) \cong {\rm Hom}_{\D}(X, Y[n])$ for any $X, Y \in \CC$ and all $n \geq 2$.
		\item (Effaceability criterion) For any $X, Y \in \CC \subseteq \D$, $n \geq 2$ and morphism $f: X \to Y[n]$, there is an object $W \in \CC$ and epimorphism $g: W \twoheadrightarrow X$ such that $fg: W \to Y[n]$ is zero.
		\item (Coeffaceability criterion) For any $X, Y \in \CC \subseteq \D$, $n \geq 2$ and morphism $f: X \to Y[n]$, there is an object $Z \in \CC$ and monomorphism $e: Y \hookrightarrow Z$ such that $e[n]f: X \to W[n]$ is zero.
	\end{enumerate}
\item The essential image of $\mathsf{real}_{\D F}$ is contained in the bounded part
	\[
		\D^b = \bigcup_{a \leq b} \D^{[a,b]} := \bigcup_{a \leq b} \D^{\leq b} \cap \D^{\geq a}
	\] 
	with agreement whenever $\mathsf{real}_{\D F}$ is fully faithful.
 \end{enumerate}
 \end{thm}

 We now apply these results to the absolutely Koszul Gorenstein $\Bbbk$-algebra $A$. By Ex. \ref{filtsingcat}, the t-structure $t^{\mathsf{lin}}$ admits a realisation functor 
 \[
	 \mathsf{real}: {\rm D}^b(\mathcal{H}^{\mathsf{lin}}(A)) \to \underline{{\rm MCM}}^{\mathbb{Z}}(A)
 \]
 where we hide the dependence on the f-category over $ \underline{{\rm MCM}}^{\mathbb{Z}}(A)$ in the notation. We will prove: 

 \begin{prop}[Prop. \ref{realequiv}] The functor $\mathsf{real}$ is an equivalence of triangulated categories. 
 \end{prop}
Since the $t$-structure $t^{\mathsf{lin}}$ is bounded, we see by Thm \ref{realthm} iv) that it suffices to prove that $\mathsf{real}$ is fully faithful. In order to use one of the (co)effaceability criteria, we begin by studying the monomorphisms and epimorphisms in $\mathcal{H}^{\mathsf{lin}}(A)$. We first recall the abelian category structure on $\mathcal{H}^{\mathsf{lin}}(A)$ (see \cite[Thm. 8.1.9]{HTT} for details).\bigskip

Let $f:X \to Y$ be a morphism in $\mathcal{H}^{\mathsf{lin}}(A)$ and set $Z = {\rm Cone}(f)$. It follows from the long exact sequence of $t$-cohomology that $Z \in \underline{{\rm MCM}}^{[-1, 0]}(A) = \underline{{\rm MCM}}^{\leq 0}(A) \cap \underline{{\rm MCM}}^{\geq -1}(A)$. The abelian category structure on $\mathcal{H}^{\mathsf{lin}}(A)$ is given by
\begin{align*}
	{\rm coker}(f) &= \mathcal{H}^0(Z) \cong \tau^{\geq 0} Z\\
	{\rm ker}(f) &= \mathcal{H}^{-1}(Z) \cong \tau^{\leq 0}(Z[-1]).
\end{align*}
In any abelian category, we of course have that $f$ is a monomorphism if and only if ${\rm ker}(f) = 0$, and $f$ is an epimorphism if and only if ${\rm coker}(f) = 0$. We now record a general characterisation of monomorphisms and epimorphisms in $\mathcal{H}^{\mathsf{lin}}(A)$. Some of these characterisations will be mainly of use in the next section. 

 \begin{lem}\label{monoepicharacterisation} Let $f: X \to Y$ be a morphism in $\mathcal{H}^{\mathsf{lin}}(A)$. 
	 \begin{enumerate}[i)]
	 \item The following are equivalent:
	\begin{enumerate}[a)]
	\item $f$ is a monomorphism, that is $fg = 0$ implies $g = 0$ in $\mathcal{H}^{\mathsf{lin}}(A)$;
	\item $Z = {\rm Cone}(f)$ is eventually linear;
	\item The induced map $f_m^*: \underline{{\rm Ext}}^m_{gr A}(Y, \Bbbk^{st}(-m)) \to \underline{{\rm Ext}}^m_{gr A}(X, \Bbbk^{st}(-m))$ is surjective for all $m \gg 0$;
	\item The induced map $f_m: [\Omega^m X(m)]_0 \to [\Omega^m Y(m)]_0$ is injective for all $m \gg 0$.
\end{enumerate}\bigskip

	 \item The following are equivalent:
	\begin{enumerate}[a)]
	\item $f$ is an epimorphism, that is $gf = 0$ implies $g = 0$ in $\mathcal{H}^{\mathsf{lin}}(A)$;
	\item $Z = {\rm Cone}(f)$ is eventually $(-1)$-linear;
	\item The induced map $ f_m^*: \underline{{\rm Ext}}^m_{gr A}(Y, \Bbbk^{st}(-m)) \to \underline{{\rm Ext}}^m_{gr A}(X, \Bbbk^{st}(-m))$ is injective for all $m \gg 0$;
	\item The induced map $f_m: [\Omega^m X(m)]_0 \to [\Omega^m Y(m)]_0$ is surjective for all $m \gg 0$.
\end{enumerate}
 \end{enumerate}
\end{lem} 
\begin{proof}
The equivalence of $a)$ and $b)$ simply rephrases vanishing of the kernel or cokernel of $f$. To see that $b)$ is equivalent to $c)$, consider the distinguished triangle
	\[
		X \to Y \to Z \to X[1].
	\]
	For any $m \in \N$, taking Tate cohomology with coefficients in $\Bbbk(-m)$ gives rise to a long exact sequence 
	\[
		\cdots \to \underline{{\rm Ext}}^{m+n}_{gr A}(X, \Bbbk(-m)) \to \underline{{\rm Ext}}^{m+n+1}_{gr A}(Z, \Bbbk(-m)) \to \underline{{\rm Ext}}^{m+n+1}_{gr A}(Y, \Bbbk(-m))  \xrightarrow{f^*} \underline{{\rm Ext}}_{gr A}^{m+n+1}(X, \Bbbk(-m)) \to \cdots 
	\]
	Since $X$ and $Y$ are eventually linear, there is an $n_0 \in \N$ such that for all $m \geq n_0$, we have
\[
	\underline{{\rm Ext}}_{gr A}^{m+i}(X \oplus Y, \Bbbk(-m)) = 0 \textup{ for all } i \neq 0.
\]
Assuming $m \geq n_0$, then for $n \notin \{ -1, 0\}$, we have $ \underline{{\rm Ext}}_{gr A}^{m+n+1}(Z, \Bbbk(-m)) = 0$ by vanishing of the outer terms. For $n = -1, 0$ the long exact sequence breaks down into a four-term exact sequence
	\[
		0 \to \underline{{\rm Ext}}^{m}_{gr A}(Z, \Bbbk(-m)) \to \underline{{\rm Ext}}^{m}_{gr A}(Y, \Bbbk(-m))  \xrightarrow{f_m^*} \underline{{\rm Ext}}_{gr A}^{m}(X, \Bbbk(-m)) \to \underline{{\rm Ext}}^{m+1}_{gr A}(Z, \Bbbk(-m))  \to 0.
	\]
	We obtain that $Z$ is eventually linear if and only if $f_m^*$ is surjective for all $m \gg 0$, and dually $Z$ is eventually $(-1)$-linear if and only if $f_m^*$ is injective for all $m \gg 0$. \bigskip 

Next, we show that $c)$ is equivalent to $d)$. By Lemma \ref{unstablelemma}, we have commutative diagrams 
\[
\xymatrix{ 
	\underline{{\rm Ext}}_{gr A}^m(Y, \Bbbk(-m)) \ar@{=}[d] \ar[r]^-{f_m^*} & \underline{{\rm Ext}}_{gr A}^m(X, \Bbbk(-m)) \ar@{=}[d]\\
	{\rm \underline{Hom}}_{gr A}(\Omega^m(Y)(m), \Bbbk) \ar@{=}[d] \ar[r]^-{f_m^*} & {\rm \underline{Hom}}_{gr A}(\Omega^m(X)(m), \Bbbk) \ar@{=}[d] \\
	{\rm Hom}_{gr A}([\Omega^m(Y)(m)], \Bbbk) \ar@{=}[d] \ar[r]^-{f_m^*} & 	{\rm Hom}_{gr A}([\Omega^m(X)(m)], \Bbbk) \ar@{=}[d]\\
	{\rm Hom}_{\Bbbk}([\Omega^m(Y)(m)]_0, \Bbbk) \ar[r]^-{f_m^*} & 	{\rm Hom}_{\Bbbk}([\Omega^m(X)(m)]_0, \Bbbk)
}
\]
and as ${\rm Hom}_{\Bbbk}(-, \Bbbk)$ is a duality on finite length $\Bbbk$-modules we are done.
\end{proof}

We now introduce a special monomorphism attached to any $M \in \mathcal{H}^{\mathsf{lin}}(A)$. Let $C(M)$ be the minimal complete resolution of $M$, which in degree $n$ is given by
\[
	C_n(M) = V_n \otimes_\Bbbk A
\]
where $V_n$ is some finite length graded $\Bbbk$-module. Since $M$ is eventually linear, $V_n$ is concentrated in degree $n$ for all $n \gg 0$. Since the complex is minimal, the quotient map $C_n(M) = V_n \otimes_\Bbbk A \twoheadrightarrow V_n \otimes_\Bbbk \Bbbk \cong V_n$ factors through the (reduced) $n$-th syzygy
	\[
		\xymatrix@R=5pt{ C_n(M) \ar[dr] \ar[rr] && V_n \\
			& [\Omega^n(M)] \ar[ur]
}
	\] 
	and, applying degree shift by $n$, restricts to an isomorphism on degree zero components 
	\[
		[\Omega^n(M)(n)]_0 \xrightarrow{\cong} V_n(n)_0 .
	\]
	Note that by definition $V_n(n)$ is concentrated in degree zero for any $n \gg 0$ such that $\Omega^n(M)(n)$ is linear, and in this case we simply write $V_n(n) = V_n(n)_0$. Taking MCM approximation, we obtain a map in $ \underline{{\rm MCM}}^{\mathbb{Z}}(A)$ 
	\[
		\xymatrix@R=15pt{ \Omega^n(M) \ar[dr] \ar@{..>}[r] & V_n^{st}\ar[d]^-{\pi} \\ 
			& V_n
}
	\] 
	such that composing with $\pi$ induces ${\rm \underline{Hom}}_{gr A}(\Omega^n(M), V_n^{st}) \cong {\rm \underline{Hom}}_{gr A}(\Omega^n(M), V_n)$. We will write
	\[
		\iota_n: M \to \Omega^{-n}(V_n^{st}) = V_n^{st}[n]
	\]
	for the corresponding morphism. Note that when $\Bbbk = k$ is a field, for $n \gg 0$ we have $V_n^{st}[n] = (k^{st}(-n)[n])^{\beta_{n,n}(M)}$. In general we have $V_n \in \mathsf{add}(\Bbbk(-n))$ and $V_n^{st} \in \mathsf{add}(\Bbbk^{st}(-n))$, where $\mathsf{add}(X)$ is the full subcategory whose objects are summands of finite direct sums of $X$. Note also that $V_n^{st}$ depends on $M$, although this is not reflected in the notation. 
	\begin{prop}\label{canonicalmono} The morphism $\iota_n: M \to V_n^{st}[n]$ is a monomorphism in $\mathcal{H}^{\mathsf{lin}}(A)$ for each $n \gg 0$. 
\end{prop}
\begin{proof} The proof will proceed directly from the definition of monomorphism. Fix $n \gg 0$ such that $\Omega^n(M)(n)$ is linear. Since $\Omega^n(-)(n)$ is an autoequivalence of $\mathcal{H}^{\mathsf{lin}}(A)$, the map $\iota_n: M \to V_n^{st}[n]$ is a monomorphism if and only if $\Omega^n(\iota_n)(n): \Omega^n(M)(n) \to V_n^{st}(n)$ is a monomorphism. For simplicity, we may replace $M$ by $\Omega^n(M)(n)$ to reduce to the case of $M$ linear, so that we may set $n = 0$ and simply write $\iota: M \to V_0^{st}$. \bigskip

	Replacing $M$ by $[M]$ if necessary, we may assume that $M$ has no projective summands and simply write $M = [M]$. By abuse of notation, we write $\iota: M \to V_0^{st}$ for any representative morphism in $\grmod A$. Then by construction, since $M$ is linear, the composite map
	\[
		\pi \iota: M \to V_0^{st} \to V_0 
	\]
	becomes an isomorphism upon passing to degree zero
	\[
		\pi \iota_0: M_0 \xrightarrow{\cong} V_0.
	\]
	Now let $g: N \to M$ be a morphism in $\mathcal{H}^{\mathsf{lin}}(A)$ such that $\iota g: N \to V_0^{st}$ is zero, interpreted as vanishing in the stable category $\mathcal{H}^{\mathsf{lin}}(A) \subseteq \underline{\mathsf{grmod}} \ \! A$. Since we have $V_0^{st} \in \mathsf{add}( \Bbbk^{st})$ and $V_0 \in \mathsf{add}( \Bbbk)$, by Lemma \ref{unstablelemma} we have
	\[
		{\rm \underline{Hom}}_{gr A}(N, V_0^{st}) \cong {\rm \underline{Hom}}_{gr A}(N, V_0) \cong {\rm Hom}_{gr A}([N], V_0)
	\]
	with first isomorphism given by composing with $\pi$. Again, possibly after removing projective summands from $N$ we may assume that $N = [N]$. Since $\iota g = 0$ in $ \underline{\mathsf{grmod}}\ \! A$, it follows that $\pi \iota g = 0$ in $\grmod A$. Now, since $N$ is eventually linear, it must be generated in degree $\leq 0$ and we let $N_{\langle < 0 \rangle} \subseteq N$ be the submodule generated by $N_{<0}$. Since $M$ is generated in degree zero, we have factorisations 
	\[
		\xymatrix{N \ar[dr] \ar[rr]^-{g} && M \ar[dr]_-{\pi \iota} \ar[r]^-{\iota} & V_0^{st}\ar[d]^-{\pi}\\
			& N/N_{\langle < 0 \rangle} \ar[ur]_-{\overline{g}}&& V_0.
}
	\]
	Then $\pi\iota g = 0$ implies $\pi \iota \overline{g} = 0$. Since $N/N_{\langle < 0 \rangle}$ is generated in degree zero, $\overline{g}$ is determined by its degree zero component. But $\pi \iota_0$ is an isomorphism, hence $(\pi \iota \overline{g})_0 = (\pi \iota)_0 \overline{g}_0 = 0$ implies $\overline{g}_0 = 0$, and so $g = 0$. This shows that $\iota$ is a monomorphism. 
\end{proof}

We can now prove that the realisation functor $\mathsf{real}: {\rm D}^b(\mathcal{H}^{\mathsf{lin}}(A)) \to \underline{{\rm MCM}}^{\mathbb{Z}}(A)$ gives an equivalence of triangulated categories.
\begin{prop}\label{realequiv} The functor $\mathsf{real}$ is an equivalence of triangulated categories.
 \end{prop}
\begin{proof} We will apply the coeffaceability criterion in Thm. \ref{realthm} iii). Let $X, Y \in \mathcal{H}^{\mathsf{lin}}(A)$ and $n \geq 2$, and let $f: X \to Y[n]$ be a morphism in $ \underline{{\rm MCM}}^{\mathbb{Z}}(A)$.
Since $X$ is eventually linear, there is an $n_0 \in \N$ such that for all $m \geq n_0$, we have
\[
	\underline{{\rm Ext}}_{gr A}^{m+n}(X, \Bbbk^{st}(-m)) = \underline{{\rm Ext}}_{gr A}^{m+n}(X, \Bbbk(-m)) = 0 \textup{ for all } n \neq 0.
\]
Consider the morphism $\iota_m: Y \to V_m^{st}[m]$ of Prop. \ref{canonicalmono}. Taking $m \gg n_0$ large enough so that $\iota_m$ is a monomorphism, we may take $e = \iota_m$. Then the composite morphism $e[n]f: X \to Y[n] \to V_m^{st}[m+n]$ is zero by our assumption that $m \geq n_0$, since $V_m^{st} \in \mathsf{add}(\Bbbk^{st}(-m))$. \bigskip

It follows that $\mathsf{real}$ is fully faithful, and by Thm. \ref{realthm} iv) it must be essentially surjective since $t^{\mathsf{lin}}$ is bounded. 
\end{proof}

Putting everything together, we obtain the proof of Theorem \ref{thmA}. 
\begin{thm}[Theorem \ref{thmA}] Let $A$ be a Koszul Gorenstein $\Bbbk$-algebra. The following are equivalent:
	\begin{enumerate}[i)]
	\item $A$ is absolutely Koszul;
	\item $t^{\mathsf{lin}}$ forms a bounded t-structure.
\end{enumerate}
When either of these equivalent conditions hold, the natural realisation functor
\[
	\mathsf{real}: {\rm D}^b(\mathcal{H}^{\mathsf{lin}}(A)) \xrightarrow{\cong} \underline{{\rm MCM}}^{\mathbb{Z}}(A) 
\]
is an equivalence of triangulated categories.
\end{thm}
\begin{proof}
	Combine Prop. \ref{thmApart1} and Prop. \ref{realequiv}.
\end{proof}

%% file: chapters/tstructure/thmB.tex
\section{The Artin-Zhang-Polishchuk noncommutative section ring and Theorem B}

Given an absolutely Koszul Gorenstein algebra $A$, the abelian category $\mathcal{H}^{\mathsf{lin}}(A)$ is an interesting $k$-linear, Hom-finite abelian category, and it is of interest to find multiple descriptions of this category. \bigskip

We will show that the opposite abelian category $\mathcal{H}^{\mathsf{lin}}(A)^{op}$ contains an ample sequence in a suitable sense; from this we will deduce that the (opposite) Koszul dual $(A^!)^{op} = {\rm Ext}^*_A(\Bbbk, \Bbbk)^{op}$ is coherent as a graded algebra, and obtain a description of $\mathcal{H}^{\mathsf{lin}}(A)^{op}$ in terms of $(A^!)^{op}$. 

\subsection*{Noncommutative projective schemes after Artin-Zhang}
In the foundational paper \cite{AZ}, Artin and Zhang introducing the noncommutative projective scheme associated to a Noetherian $\N$-graded noncommutative $k$-algebra $B$, and established a useful recognition theorem for its category of coherent sheaves amongst abelian categories. 

Let us begin with the classical setting. Let $X$ be a $k$-scheme, and let $\OO(1)$ be an ample line bundle on $X$, where we write $\OO(m) := \OO(1)^{\otimes m}$ for its tensor powers and $\F(m) = \F \otimes \OO(m)$. By classical results of Serre, the sequence $\{ \OO(m) \}_{m \in \Z}$ detects various properties of the category $\coh X$ of coherent sheaves: 

\begin{enumerate}[a)]
\item (Global generation) For every $\F \in \coh X$ and $m \gg 0$, the coherent sheaf $\F(m)$ is globally generated; that is, the natural morphism
	\[
		\Gamma(X, \F(m)) \otimes_k \OO \xrightarrow{ev} \F(m)
	\]
	is an epimorphism. 
\item (Detecting epimorphisms). For every epimorphism $f: \F \to \G$, the induced map on sections 
	\[
		\Gamma(f_m): \Gamma(X, \F(m) ) \to \Gamma(X, \G(m)) 
	\]
	is surjective for all $m \gg 0$.
	\end{enumerate}\bigskip 

Property $a)$ is taken as the definition of ampleness in Hartshorne, and $b)$ follows from Serre's ampleness criterion by setting $\mathcal{K} = {\rm ker}(f)$, as we have vanishing of sheaf cohomology ${\rm H}^i(X, \mathcal{K}(m)) = 0$ for $i > 0$ for all $m \gg 0$. Moreover, the section ring $S = \bigoplus_{n \geq 0} \Gamma(X, \OO(n))$ is a Noetherian $\N$-graded $k$-algebra and the truncated sections
\[
	\Gamma_{\geq m}(X, \F) := \bigoplus_{n \geq m} \Gamma(X, \F(n)) 	
\]
form finitely generated modules over $S$. When $\OO(1)$ is very ample so that $S$ is generated by $S_1$ over $S_0$, the above properties are instrumental in the proof of Serre's Theorem, which recovers the category of coherent sheaves via the Serre quotient $\coh X = \qgr S$. \bigskip

Now consider a triple $(\CC, \OO, s)$, with $\CC$ a Noetherian, $k$-linear, Hom-finite abelian category, where Noetherian means that every object $X \in \CC$ is Noetherian; $\OO$ is a distinguished object of $\CC$; and $s: \CC \xrightarrow{\cong} \CC$ is an autoequivalence of $\CC$, which we think of formally as $- \otimes \OO(1)$. We may formally write $\OO(m) = s^m(\OO)$, and we can mimic taking global sections of any $M \in \CC$ by defining $\Gamma(M(m)) := {\rm Hom}(\OO(-m), M)$. As before, we obtain a section ring 

\begin{align*}
	S = \bigoplus_{n \geq 0} \Gamma(\OO(n)) &= \bigoplus_{n \geq 0}{\rm Hom}(\OO(-n), \OO) \\
	&\cong \bigoplus_{n \geq 0} {\rm Hom}(\OO, \OO(n))
\end{align*}
with graded algebra structure given by composing $g \in {\rm Hom}(\OO, \OO(m))$ and $f \in {\rm Hom}(\OO, \OO(n))$ as
\[
	f * g := s^m(f) \circ g.
\]
Similarly we obtain homogeneous and truncated section functors $\Gamma_*, \Gamma_{\geq m}: \CC \to \Grmod S$ by setting
\[
	\Gamma_*(M) = \bigoplus_{n \in \Z} \Gamma(M(n)) = \bigoplus_{n \in \Z}{\rm Hom}(\OO(-n), M)
\]
\[
	\Gamma_{\geq m}(M) = \bigoplus_{n \geq m} \Gamma(M(n)) = \bigoplus_{n \geq m}{\rm Hom}(\OO(-n), M)
\]
with $S$ operating on the right by composition. 

\begin{defn}[{\cite[4.2.1]{AZ}}]\label{ampleAZ} The autoequivalence $s$ is ample if the following conditions are satisfied in $\CC$:
	\begin{enumerate}[a)]
	\item (Global generation) For every $M \in \CC$ and $m \in \Z$, there exists indices $i_1, \dots, i_s \geq m$ for which there exists an epimorphism
		\[
			\bigoplus_{j = 1}^s \OO(-i_j) \twoheadrightarrow M.
		\]
	\item (Detecting epimorphisms) For every epimorphism $f: M \to N$ in $\CC$, the induced map 
		\[
			\Gamma(M(m)) \to \Gamma(N(m)) 
		\]
		is surjective for all $m \gg 0$.
\end{enumerate}
\end{defn}

Artin and Zhang proved the following extension of Serre's Theorem. We refer to \cite{AZ} for the definition of the $\chi_1$ condition. 
\begin{thm}[{\cite[Thm 4.5]{AZ}}] Let $(\CC, \OO, s)$ be a triple as above, and assume that $s$ is ample. Then the section ring $S = \bigoplus_{n \geq 0} \Gamma(\OO(n))$ is a Noetherian $\N$-graded algebra (satisfying the $\chi_1$ condition), and there is an equivalence of abelian categories
	\[
		\CC \xrightarrow{\cong} \qgr S
	\]
	given by sending $M \mapsto \Gamma_{\geq m}(M)$ for any $m \in \Z$. \bigskip

	Conversely, if $B$ is a Noetherian $\N$-graded algebra (satisfying the $\chi_1$ condition), then $(\qgr B, \pi B, (1))$ forms a triple as above, meaning that $\qgr B$ is a Noetherian abelian category with distinguished object $\pi B$ the image of $B$, and autoequivalence $s = (1)$ descencending from the degree shift functor. Moreover, $(1)$ is ample, and the natural morphism of graded algebras $B \to \Gamma_{\geq 0}(\pi B)$ is an isomorphism in large enough degrees.
\end{thm}
\begin{rem}
	The Artin-Zhang Theorem extends Serre's Theorem not only to the noncommutative setting but also improves on the commutative case. Set $(\CC, \OO, s) = (\coh X, \OO_X, - \otimes \mathcal{L})$ for $X$ a $k$-scheme and $\mathcal{L}$ an ample, but not very ample, line bundle. Then $s = - \otimes \mathcal{L}$ is ample in the sense of Def. \ref{ampleAZ}, and its section ring $S = \bigoplus_{n \geq 0} \Gamma(X, \mathcal{L}^{\otimes n})$ is finitely generated but not necessarily in degree $0$ and $1$. Serre's Theorem requires $S$ to be generated by $S_1$ over $S_0$, but by the above theorem, the functor $\Gamma_{\geq 0}(X, -)$ still induces an equivalence $\coh X \cong \qgr S$. In other words, Serre's Theorem still holds for finitely generated graded $k$-algebras arising as the section ring of an ample line bundle on some projective variety. This was verified directly in the case of line bundle of degree $1$ and $2$ on an elliptic curve in \cite[Chp. 7]{Pavlov}. \bigskip

	Note that the general statement $\coh \proj S \cong \qgr S$ fails to hold for arbitrary finitely generated commutative $k$-algebras $S$, see \cite[Prop. 2.17]{Or1} for correct behavior there. 
\end{rem}
\subsection*{Noncommutative projective schemes after Polishchuk}
Many natural examples of noncommutative projective schemes occur in mathematics beyond those covered by the Artin-Zhang schemes. The crux of the matter is that many abelian categories $\CC$ which ought to be realised as categories of coherent sheaves fail to be Noetherian (see \cite{Pol} for natural examples). This was remedied in \cite{Pol} by systematically working with coherent algebras instead, and we present the relevant results here. \bigskip 

We have seen previously that attached to an $\N$-graded locally finite coherent $k$-algebra $B = B_0 \oplus B_1 \oplus \dots$, we obtain an abelian category $\qgr B = \grmod B / \mathsf{grmod}_0 \ \! B$. We now look to establish a recognition theorem for such categories, following Polishchuk \cite{Pol}. \bigskip

Let $\CC$ be a Hom-finite $k$-linear abelian category, not necessarily Noetherian. As previously, we consider a triple $(\CC, \OO, s)$ where $\OO$ is a distinguished object and $s$ an autoequivalence of $\CC$, and define $\OO(n) := s^n(\OO)$. \bigskip

Polishchuk works in the setting of connected graded $k$-algebras, and so one introduces the connected section ring $S = k \oplus S_{\geq 1} = k \oplus \Gamma_{\geq 1}(\OO)$, where for any $M \in \CC$ the sections $\Gamma_{\geq m}(M)$ are defined as previously; the latter inherits as before the structure of a graded $S$-module. The notion of ampleness of $s$ then takes a different form: 

\begin{defn}[{\cite[Sect. 2]{Pol}}] The autoequivalence $s$ is ample if the following properties hold:
	\begin{enumerate}[i)]
	\item (Detecting epimorphisms) For every epimorphism $f: M \to N$ in $\CC$, the induced map 
		\[
			\Gamma(M(m)) \to \Gamma(N(m)) 
		\]
		is surjective for all $m \gg 0$.
	\item (Finite generation) For every $M \in \CC$ and $m \in \Z$, the module $\Gamma_{\geq m}(M)$ is finitely generated over $S$.
	\item (Global generation) For every $M \in \CC$ and $m \in \Z$, there exists indices $i_1, \dots, i_s \geq m$ for which there exists an epimorphism
		\[
			\bigoplus_{j = 1}^s \OO(-i_j) \twoheadrightarrow M.
		\]
\end{enumerate}
\end{defn}
\begin{rem} The terminology is a minor departure from Polishchuk's \cite{Pol}, and it is good to highlight the discrepancy. Polishchuk works more generally with a sequence of objects $\E = (E_i)_{i \in \Z}$, to which he associates a $\Z$-algebra $A = A(\E)$, defined as a bigraded $k$-algebra $A = \bigoplus_{i \leq j} A_{ij}$ with $A_{ii} = k$ and $A_{ij} = {\rm Hom}_\CC(E_i, E_j)$ for $i < j$, and with multiplication having sole component $A_{jk} \otimes A_{ij} \to A_{ik} $ given by composition in $\CC$, where $A_{ii} = k \hookrightarrow {\rm End}_\CC(E_i)$ is identified with scalar multiples of the identity. Global sections are defined by $\Gamma_*(M) = \bigoplus_{i \in \Z} {\rm Hom}_\CC(E_i, M)$, which form graded modules over $A$. Setting $E_i = \OO(-i)$, one obtains the $\Z$-algebra $A = A(\E)$ from $S$ via $A = \bigoplus_{i\leq j} A_{ij} = \bigoplus_{i \leq j} S_{j-i}$. The notions readily translate from there, see \cite[Rem. 2.1]{Pol}.
\end{rem}\bigskip

Polishchuk then extends the Artin-Zhang Theorem to non-Noetherian abelian categories. 
\begin{prop}[{\cite[Prop. 2.3]{Pol}}] Let $(\CC, \OO, s)$ be a triple as above with $s$ ample. Then:
	\begin{enumerate}[i)]
	\item For every $M \in \CC$ and $m \in \Z$, the graded module $\Gamma_{\geq m}(M)$ is a coherent $S$-module. 
	\item The section ring $S$ is a coherent graded algebra.
\end{enumerate}
\end{prop}
\begin{thm}[{\cite[Thm. 2.4]{Pol}}] Let $(\CC, \OO, s)$ be a triple as above with $s$ ample. If $S$ denotes the section ring, then there is an equivalence of abelian categories 
\[
	\CC \xrightarrow{\cong} \qgr S
\]
given by sending $M \mapsto \Gamma_{\geq m}(M)$ for any $m \in \Z$. 
\end{thm}
When $s$ is ample, we will more generally call $\{ \OO(n) \}_{n \in \Z}$ an ample sequence. \bigskip

\subsection*{Ample sequences in the category $\mathcal{H}^{\mathsf{lin}}(A)^{op}$}

Fix an absolutely Koszul Gorenstein $\Bbbk$-algebra $A$, and consider the triple $(\CC, \OO, s) := (\mathcal{H}^{\mathsf{lin}}(A)^{op}, \Bbbk^{st}, \Omega^{1}(-)(1))$. We aim to show that $s = \Omega^{1}(-)(1) = (-)(1)[-1]$ is ample, and then compare the section ring with ${\rm Ext}^*_A(\Bbbk, \Bbbk)^{op}$. We first need a minor assumption and observation, only relevant to the case where $\Bbbk \neq k$. \bigskip

The definition of the section ring $S = k \oplus S_{\geq 1} = k \oplus \Gamma_{\geq 1}(\OO)$ can be changed to allow any intermediate subalgebra $k \subseteq \widetilde{S}_0 \subseteq S_0$ in degree zero, since $S_0 = {\rm End}_\CC(\OO)$ is finite dimensional over $k$. This was noted by Minamoto, whose work on coherence of higher preprojective algebras as discussed in Chapter 4 led him to consider section rings of the form $S = \Pi(\Lambda)$
\[
	S_0 \oplus S_1 \oplus \dots = \Lambda \oplus \Pi_1(\Lambda) \oplus \dots
\]
for $\Lambda$ a Fano finite dimensional algebra (see \cite[Sect. 3]{Mina}). He notes that Polishchuk's results go through with the same proof in this more general case, and we will use this implicitly. \bigskip

We now impose the following assumption: 
\begin{assump}\label{homologicallyhomogeneous} We assume that $A$ is homologically homogeneous, of infinite global dimension. That is, for every simple summand $S$ of $\Bbbk$, we assume that $\pdim S = \infty$. 
\end{assump} 

The stabilisation functor ${\rm st}: {\rm D}^b(\grmod A) \to \underline{{\rm MCM}}^{\mathbb{Z}}(A)$ induces an algebra morphism $\varphi: \Bbbk = {\rm End}_{gr A}(\Bbbk) \to {\rm \underline{End}}_{gr A}(\Bbbk^{st})$. By our assumption (\ref{homologicallyhomogeneous}), no summand of $\Bbbk$ is perfect over $A$ and so the summands of $\Bbbk$ survive stabilisation. Since $\Bbbk$ is semisimple, the kernel of $\varphi$ is generated by idempotents, but since non-zero idempotents are sent to non-zero idempotents the map $\varphi$ must be injective, giving a sequence of intermediate algebras
\[
	k \subseteq \Bbbk \subseteq {\rm \underline{End}}_{gr A}(\Bbbk^{st})
\]
Working in the opposite abelian category $\mathcal{H}^{\mathsf{lin}}(A)^{op}$, we shall set $\widetilde{S}_0 = \Bbbk^{op}$, we will use the section ring
	\begin{align*}
	S &= \Bbbk^{op} \oplus \bigoplus_{n \geq 1} {\rm Hom}_{\mathcal{H}^{\mathsf{lin}}(A)^{op}}(\Bbbk^{st}, \Omega^n(\Bbbk^{st})(n))\\
	&= \Bbbk^{op} \oplus \bigoplus_{n \geq 1} {\rm Hom}_{\mathcal{H}^{\mathsf{lin}}(A)}(\Bbbk^{st}, \Omega^{-n}(\Bbbk^{st})(-n))^{op}\\
	&= \Bbbk^{op} \oplus \bigoplus_{n \geq 1} \underline{{\rm Ext}}_{gr A}^{n}(\Bbbk^{st}, \Bbbk^{st}(-n))^{op}
\end{align*}
and denote it by $\underline{{\rm Ext}}^{\Delta_+}_{gr A}(\Bbbk^{st}, \Bbbk^{st})^{op}$. Of course, this agrees with Polishchuk's definition when $\Bbbk = k$. The truncated global sections of $M \in \mathcal{H}^{\mathsf{lin}}(A)^{op}$ can be described
\[
	\Gamma_{\geq m}(M) = \bigoplus_{n \geq m} \underline{{\rm Ext}}^{n}_{gr A}(M, \Bbbk^{st}(-n)) 
\]
with its natural right $S$-module structure by post-composition in $ \underline{{\rm MCM}}^{\mathbb{Z}}(A)$. To help analyse its structure, we have the following well-known lemma: 
\begin{lem}[{\cite[Thm. 6.3 (4)]{SMI}}]\label{linearfinitegeneration} Let $N \in \grmod A$ be a linear module. Then ${\rm Ext}^*_A(N, \Bbbk)$ is a linear module over ${\rm Ext}^*_A(\Bbbk, \Bbbk)^{op}$. Therefore it is generated by ${\rm Ext}^0_A(N, \Bbbk)$, and in particular ${\rm Ext}^*_A(N, \Bbbk)$ is finitely generated. 
\end{lem}
\begin{rem}
	The paper \cite{SMI} works over a field $k$ instead of the more general semisimple base $\Bbbk$, but the proof of \cite[Thm. 6.3 (4)]{SMI} goes through without change. 
\end{rem}

We finally get to the main point. 
\begin{prop}\label{ThmBPart1} Let $A$ be absolutely Koszul Gorenstein, satisfying hypothesis (\ref{homologicallyhomogeneous}). Consider the triple $(\mathcal{H}^{\mathsf{lin}}(A)^{op}, \Bbbk^{st}, \Omega^{1}(-)(1))$. Then $\Omega^{1}(-)(1)$ is an ample autoequivalence. 
\end{prop}
\begin{proof}
We need to verify three conditions: 
	\begin{enumerate}[i)]
	\item (Detecting epimorphisms) For every epimorphism $f: M \twoheadrightarrow N$ in $\mathcal{H}^{\mathsf{lin}}(A)^{op}$, the induced map on sections 
		\[
			\Gamma(f_m): \underline{{\rm Ext}}^{m}_{gr A}(M, \Bbbk^{st}(-m)) \to \underline{{\rm Ext}}^m_{gr A}(N, \Bbbk^{st}(-m)) 
		\]
		is surjective for all $m \gg 0$.
	\item (Finite generation) For every $M \in \mathcal{H}^{\mathsf{lin}}(A)^{op}$ and $m \in \Z$, the truncated global sections $\Gamma_{\geq m}(M) = \bigoplus_{n \geq m} \underline{{\rm Ext}}^{n}_{gr A}(M, \Bbbk^{st}(-n))$ is finitely generated over $S$.
	\item (Global generation) For every $M \in \mathcal{H}^{\mathsf{lin}}(A)^{op}$ and $m \in \Z$, there exists indices $i_1, \dots, i_s \geq m$ for which one has an epimorphism
\[
	\bigoplus_{j=1}^s \Omega^{-i_j}(\Bbbk^{st})(-i_j) \twoheadrightarrow M.
\]
\end{enumerate}

To see $i)$, note that epimorphisms $f: M \twoheadrightarrow N$ in $\mathcal{H}^{\mathsf{lin}}(A)^{op}$ are simply monomorphisms $f: N \hookrightarrow M$ in $\mathcal{H}^{\mathsf{lin}}(A)$, in which case the result follows from the characterisation of monomorphisms (Lemma \ref{monoepicharacterisation}) . For $iii)$, note that we have monomorphisms
\[
	\iota_m: M \hookrightarrow V_m[m] = \Omega^{-m} V_m
\]
for all $m \gg 0$. Since $\Omega^{-m}V_m \in \mathsf{add}(\Omega^{-m}(\Bbbk^{st})(-m))$, we have a split monomorphism $\Omega^{-m}V_m \hookrightarrow \Omega^{-m}(\Bbbk^{st})(-m)^{\oplus N}$ for some $N \gg 0$. Passing to the opposite category, we obtain an epimorphism by composition
\[
	\Omega^{-m}(\Bbbk^{st})(-m)^{\oplus N} \twoheadrightarrow \Omega^{-m}V_m \twoheadrightarrow M
\]
and $iii)$ follows. \bigskip

It remains to prove finite generation of $\Gamma_{\geq m}(M) = \bigoplus_{n \geq m} \underline{{\rm Ext}}^n_{gr A}(M, \Bbbk^{st}(-n))$ over $S = \underline{{\rm Ext}}^{\Delta_+}_A(\Bbbk^{st}, \Bbbk^{st})^{op}$. Since $\Gamma_{\geq i}(M)/\Gamma_{\geq i+1}(M)$ has finite length over $\Bbbk$, it suffices to prove finite generation of $\Gamma_{\geq m}(M)$ for $M$ fixed and $m \gg 0$. There exists an $n_0 \in \N$ such that for any $m \geq n_0$, the syzygy $\Omega^m(M)(m)$ is linear. Since we have $\Gamma_{\geq m}(M) = \Gamma_{\geq 0}(\Omega^m(M)(m))$, it is sufficient to prove that $\Gamma_{\geq 0}(M)$ is finitely generated for every linear MCM module $M$. \bigskip

By Lemma \ref{linearfinitegeneration}, for $M$ linear ${\rm Ext}^*_A(M, \Bbbk)$ is finitely generated over ${\rm Ext}^*_A(\Bbbk, \Bbbk)$. After possibly removing projective summands from $M$ (which does not affect the result), the stabilisation functor induces a bijection
\[
	{\rm Ext}^*_A(M, \Bbbk) \cong \underline{{\rm Ext}}^{\geq 0}_A(M, \Bbbk^{st}) 
\]
and the induced right ${\rm Ext}^*_A(\Bbbk, \Bbbk)^{op}$-module structure on the latter factors through the map ${\rm Ext}^*_A(\Bbbk, \Bbbk)^{op} \to \underline{{\rm Ext}}^{\Delta_+}_A(\Bbbk^{st}, \Bbbk^{st})^{op} = S$. It follows that 
\[
	\Gamma_{\geq 0}(M) = \bigoplus_{n \geq 0} \underline{{\rm Ext}}^n_{gr A}(M, \Bbbk^{st}(-n)) = \underline{{\rm Ext}}^{\geq 0}_A(M, \Bbbk^{st})
\]
is finitely generated over $S$, as we wanted. 
\end{proof}

The proposition shows that $S = \underline{{\rm Ext}}^{\Delta_+}_A(\Bbbk^{st}, \Bbbk^{st})^{op}$ is coherent as a graded $\Bbbk$-algebra, and we want to compare it with the classical Ext algebra. Recall that for any graded algebra $B$ and $n \in \N$, we denote by $B^{(n)} = \bigoplus_{i \in \Z} B_{ni}$ be the $n$-th Veronese subalgebra of $B$, and likewise for any graded $B$ module $N$ we write $N^{(n)} = \bigoplus_{i \in \Z} N_{ni}$.

Fixing such an $n \in \N$ and replacing $\Omega^1(-)(1)$ by $\Omega^n(-)(n)$, note that the corresponding section ring is replaced by $S^{(n)} = (\underline{{\rm Ext}}^{\Delta_+}_A(\Bbbk^{st}, \Bbbk^{st})^{op})^{(n)}$. We will also need: 

\begin{lem}\label{amplepower} For any $n \in \N$, the autoequivalence $\Omega^n(-)(n)$ of $\mathcal{H}^{\mathsf{lin}}(A)^{op}$ is also ample. 
\end{lem}
\begin{proof}
The verifications of $i)$ and $iii)$ are immediate from the previous proposition, with the only subtle point the finite generation of $\Gamma_{\geq m}(M)^{(n)} = \bigoplus_{ni \geq m} \underline{{\rm Ext}}^{ni}_{gr A}(M, \Bbbk^{st}(-ni))$ over $S^{(n)}$. We may similarly reduce to the case of $M$ linear. \bigskip 

By Lemma \ref{linearfinitegeneration}, the module ${\rm Ext}^*_A(M, \Bbbk)$ is generated by ${\rm Ext}^0_A(M, \Bbbk)$ over ${\rm Ext}^*_A(\Bbbk, \Bbbk)^{op}$. Since ${\rm Ext}^*_A(\Bbbk, \Bbbk)^{op}$ is generated by ${\rm Ext}^1_A(\Bbbk, \Bbbk)$ as a $\Bbbk$-algebra, this survives passage to Veroneses and ${\rm Ext}^*_A(M, \Bbbk)^{(n)}$ is generated by ${\rm Ext}^0_A(M, \Bbbk)$ over $({\rm Ext}_A^*(\Bbbk, \Bbbk)^{op})^{(n)}$. The rest of the argument is the same.
\end{proof}
\begin{rem} Given a triple $(\CC, \OO, s)$ with $s$ ample and $n \in \N$, it isn't a priori clear if $s^n$ is also ample. It is easy to see that this holds for arbitrary $n$ if the section ring $S$ of $s$ is generated by $S_1$ over $S_0$. The above proof is a workaround for our situation. 
\end{rem}

By Lemma \ref{amplepower}, it follows that for any $n \in \N$, the Veronese $S^{(n)} = (\underline{{\rm Ext}}^{\Delta_+}_A(\Bbbk^{st}, \Bbbk^{st})^{op})^{(n)}$ is also coherent. This is useful in light of the following:
\begin{prop}[Polishchuk]\label{Veronese} Let $B$ be a graded $\Bbbk$-algebra, generated by $B_1$ over $\Bbbk$ and with finitely many relations. Let $n \in \N$. Then $B$ is coherent if and only if $B^{(n)}$ is coherent. When both conditions hold, we have an equivalence of abelian categories
	\[
		\qgr B \cong \qgr B^{(n)}
	\]
	sending $M \in \qgr B$ to $M^{(n)} = \bigoplus_{i \in \Z} M_{ni}$. 
\end{prop}
\begin{proof}
	This is \cite[Prop. 2.6]{Pol}. Again, Polishchuk works over $k$ but the argument goes through unchanged over $\Bbbk$. 
\end{proof}

Combining everything, we obtain a proof of Theorem \ref{thmB}.
\begin{thm}[Theorem \ref{thmB}] Let $A$ be absolutely Koszul Gorenstein, satisfying (\ref{homologicallyhomogeneous}). Then $(A^!)^{op} = {\rm Ext}_A^*(\Bbbk, \Bbbk)^{op}$ is coherent, and we have a contravariant equivalence of abelian categories
	\[
		\mathcal{H}^{\mathsf{lin}}(A)^{op} \cong \qgr (A^!)^{op} 
	\]
	sending $M$ to ${\rm Ext}^{*}_A(M, \Bbbk)$. \bigskip

	The converse holds in the Artinian case: if $A$ is an Artinian Koszul Gorenstein algebra with $(A^!)^{op}$ is coherent, then $A$ is absolutely Koszul. 
\end{thm}
\begin{proof} 
	Assume that $A$ is absolutely Koszul. Since $A$ is Koszul, ${\rm Ext}_A^*(\Bbbk, \Bbbk)^{op}$ is generated in degree one over $\Bbbk$ with finitely many (quadratic) relations. By Prop. \ref{Veronese}, ${\rm Ext}^*_A(\Bbbk, \Bbbk)^{op}$ is coherent if and only if one of its $n$-th Veronese $({\rm Ext}^*_A(\Bbbk, \Bbbk)^{op})^{(n)}$ is coherent. Taking $n \gg \idim(A)$, we have an isomorphism of Veroneses subalgebras
	\[
		({\rm Ext}^*_A(\Bbbk, \Bbbk)^{op})^{(n)} \cong \Bbbk^{op} \oplus \big(\bigoplus_{i \geq 1} \underline{{\rm Ext}}^{ni}_{gr A}(\Bbbk^{st}, \Bbbk^{st}(-ni))\big)^{op} = S^{(n)} 
	\]
	since $\Bbbk^{st}$ is eventually linear and Tate cohomology eventually agrees with Ext. We have noted that this last algebra is coherent by Lemma \ref{amplepower}, and therefore so is $(A^!)^{op} = {\rm Ext}^*_A(\Bbbk, \Bbbk)^{op}$.\bigskip

	Next, consider the composition
		\[
			\mathcal{H}^{\mathsf{lin}}(A)^{op} \xrightarrow{ {\rm Ext}^*_A(-, \Bbbk)} \qgr (A^!)^{op} \xrightarrow{(-)^{(n)}} \qgr S^{(n)}.
	\] 
	Note that ${\rm Ext}^{\geq m}_A(M, \Bbbk) \hookrightarrow {\rm Ext}^*_A(M, \Bbbk)$ is an isomorphism in $\qgr (A^!)^{op}$ for any $m \geq 0$, and that ${\rm Ext}^{\geq 1}_A(-, \Bbbk)$ is well-defined on the stable category $\mathcal{H}^{\mathsf{lin}}(A) \subseteq \underline{{\rm MCM}}^{\mathbb{Z}}(A)$, hence we may compute ${\rm Ext}^*_A(-, \Bbbk)$ as ${\rm Ext}^{\geq m}_A(-, \Bbbk)$ for any choice of $m \geq 1$. For any $M \in \mathcal{H}^{\mathsf{lin}}(A)$, there is an $n_0 \geq \idim A$ such that for any $m \geq n_0$, we have 
	\[
		{\rm Ext}^{\geq m}_A(M, \Bbbk) \cong \bigoplus_{i \geq m} \underline{{\rm Ext}}^{i}_{gr A}(M, \Bbbk(-i)) \cong \bigoplus_{i \geq m} \underline{{\rm Ext}}^{i}_{gr A}(M, \Bbbk^{st}(-i)) 
	\]
	as $(A^!)^{op} = {\rm Ext}_A^*(\Bbbk, \Bbbk)^{op}$-modules. Composing with $(-)^{(n)}$ sends $M$ to
	\[
		\Gamma_{\geq m}(M)^{(n)} = \bigoplus_{ni \geq m} \underline{{\rm Ext}}^{ni}_{gr A}(M, \Bbbk^{st}(-ni))
	\]
	with the notation from Lemma \ref{amplepower}. We obtain a commutative diagram
	\[
		\xymatrix@R=30pt@C=50pt{ \mathcal{H}^{\mathsf{lin}}(A)^{op} \ar[dr]_-{\Gamma_{\geq m}(-)^{(n)}}^-{\cong} \ar[r]^-{ {\rm Ext}^*_A(-, \Bbbk)} & \qgr (A^!)^{op} \ar[d]^-{(-)^{(n)}}_-{\cong} \\
			& \qgr S^{(n)}
}
	\]
	and so the top row is an equivalence. This proves the main claim. \bigskip

	For the converse, let $A$ be Artinian Koszul Gorenstein and assume that $E = (A^!)^{op} = {\rm Ext}_A^*(\Bbbk, \Bbbk)^{op}$ is coherent. We have seen that the BGG correspondence already holds in this case (Ex. \ref{coherentBGG}), meaning that we have an equivalence of triangulated categories 
	\[
		\underline{{\rm MCM}}^{\mathbb{Z}}(A)^{op} \cong {\rm D}^b(\qgr E)
	\]
	such that pulling back the standard bounded t-structure on the right hand side gives rise to $t^{\mathsf{lin}}$, which is therefore also a bounded t-structure. Therefore $A$ is absolutely Koszul by Theorem \ref{thmA}.
\end{proof}

As corollary, we obtain Theorem \ref{thmC}. 
\begin{thm}[Theorem \ref{thmC}] Let $A$ be absolutely Koszul Gorenstein, and homologically homogeneous. Let $E = (A^!)^{op} = {\rm Ext}^*_A(\Bbbk, \Bbbk)^{op}$. Then there is an equivalence of triangulated categories 
	\[
		\underline{{\rm MCM}}^{\mathbb{Z}}(A)^{op} \cong {\rm D}^b(\qgr E)
	\]
	such that $t^{\mathsf{lin}}$ arises as the pullback of the standard t-structure on the right hand side.
\end{thm}

%% file: chapters/tstructure/singdim.tex
\section{The virtual dimension of a Koszul Gorenstein algebra}
In this subsection, unless otherwise stated, $A$ will denote an absolutely Koszul Gorenstein $\Bbbk$-algebra, of Gorenstein dimension $d$. It is often the case that the triangulated category $\underline{{\rm MCM}}^{\mathbb{Z}}(A)$ admits a Serre functor $\Se_A$ of \emph{standard form}
\[
	\Se_A(-) = - \otimes_A \omega_A [d-1]
\]
for some invertible $A$-bimodule  $\omega_A =$ $_1A_\alpha (a)$, where $\alpha$ an automorphism of $A$ of degree zero. This holds for instance when $A$ is Frobenius, or when $A$ is commutative with isolated singularities. From now on, we will assume that $A$ admits a Serre functor $\Se_A$ of standard form. \bigskip

Define $\omega_\mathcal{H} := \omega_A[-a]$ so that the category $\mathcal{H} := \mathcal{H}^{\mathsf{lin}}(A)$ of eventually linear modules is stable under $M \mapsto M \otimes_A \omega_\mathcal{H} := M(a)[-a]$. Setting $\nu := d-1 + a$, we write 
\[
	\Se_\mathcal{H}(-) = - \otimes \omega_\mathcal{H}[\nu]
\]
for the induced autoequivalence on $ {\rm D}^b(\mathcal{H}) $. The following is an immediate consequence of Theorem A. 
\begin{prop}\label{serredualityheart} The autoequivalence $\Se_\mathcal{H}$ is a Serre functor for $ {\rm D}^b(\mathcal{H})$. In particular for any $X, Y \in \mathcal{H}$ we have 
	\[
		{\rm Ext}^i_\mathcal{H}(X, Y) \cong {\rm Ext}^{\nu-i}_{\mathcal{H}}(Y, X \otimes \omega_\mathcal{H})^*
	\]
	and so $\gldim \mathcal{H} = \nu$. 
\end{prop}
\begin{rem} 
	Since the category $\mathcal{H}^{op} = \mathcal{H}^{\mathsf{lin}}(A)^{op} \cong \qgr (A^!)^{op}$ sometimes arises as the category of coherent sheaves on some projective variety, it fails in general to contain enough projectives or injectives, and so we interpret ${\rm Ext}$ through the derived category. 
\end{rem}

\begin{defn} For $A$ a Koszul Gorenstein $\Bbbk$-algebra, we call $\nu = d-1 + a$ the \emph{virtual} dimension of $A$.
\end{defn}
Note that when $A$ is absolutely Koszul Gorenstein (and $\gldim A = \infty$) we have an immediate inequality $\nu = \gldim \mathcal{H} \geq 0$. When $A$ is commutative, the invariant $\nu$ is easily computed from numerical data attached to $A$, and this inequality holds more generally.
\begin{prop}\label{hilbertdegree} Assume that $A$ is commutative Gorenstein graded connected over $k$, so that its Hilbert series is rational of the form
	\[
	H_A(t) = \frac{h_A(t)}{(1-t)^d}.
	\]
	Then $\deg h_A(t) = \nu + 1$. In particular $\nu \geq 0$ unless $A$ is a polynomial algebra. 
\end{prop}
\noindent We will use a well-known lemma, which follows from a prime avoidance argument. 
\begin{lem} Let $R = (R, \m, k)$ be a local or graded local Noetherian commutative ring. If $\m$ contains a non-zerodivisor $r$, then there is a non-zerodivisor $ r' \in \m \setminus \m^2$.
\end{lem}
\begin{proof}[Proof of Prop. \ref{hilbertdegree}]
	Since $\depth A = \dim A$ for $A$ in the proposition, the lemma guarantees that a regular sequence ${\bf x} = (x_1, \dots, x_d)$ of linear forms $x_i \in A_1$ exists, and both sides of the equality $\deg h_A(t) = \nu + 1$ are stable under the reduction $A \mapsto A/( {\bf x})$. The proof of Prop. \ref{hilbertdegree} then immediately reduces to the case $d = 0$ where it holds since both sides equal the socle degree. Note that $\nu = -1$ corresponds to $h_A(t)$ constant, in which case $A/({\bf x}) = k$ and $A$ is regular, hence a polynomial algebra. 
\end{proof}

Next, again in the commutative case, we have $\nu = \dim A - 1 + a = \dim X + a$ for $X = \proj A$, in which case the possible inequalities 
\begin{align*}
&\nu < \dim X\\
&\nu = \dim X\\
&\nu > \dim X
\end{align*}
correspond to the three cases of Orlov's semiorthogonal decomposition theorem, and dictate the direction of the embedding. \bigskip

Now, from Prop. \ref{serredualityheart} we immediately obtain:
\begin{cor} Let $A$ be an absolutely Koszul Gorenstein $\Bbbk$-algebra with a Serre functor of standard form. The following holds:
	\begin{enumerate}[1)]
	\item When $\nu = 0$, the category $ \underline{{\rm MCM}}^{\mathbb{Z}}(A)$ is semisimple. 
	\item When $\nu \leq 1$, every indecomposable $M \in \underline{{\rm MCM}}^{\mathbb{Z}}(A)$ is eventually $n$-linear for some $n \in \Z$.
\end{enumerate}
\end{cor}
\begin{proof}
	Indeed $\gldim \mathcal{H} = 0$ implies that ${\rm D}^b(\mathcal{H}) \cong \underline{{\rm MCM}}^{\mathbb{Z}}(A)$ is semisimple, while for $\gldim \mathcal{H} \leq 1$ every complex $X \in {\rm D}^b(\mathcal{H}) $ is formal by \cite[Sect. 5.2.5]{HandbookTilting} (see also \cite[Sect. 6, Thm. 3.1]{HandbookTilting}), and so each indecomposable $M \in \underline{{\rm MCM}}^{\mathbb{Z}}(A)$ has the form $M \cong \mathcal{H}^n(M)[-n]$ for some unique $n \in \Z$.
\end{proof} 

\begin{exmp} Let $A = k[x_0, \dots, x_d]/(Q)$ be a quadric hypersurface with isolated singularities at the origin. Then $H_A(t) = \frac{1+t}{(1-t)^d}$ and so $\nu = 0$, and $ \underline{{\rm MCM}}^{\mathbb{Z}}(A)$ is semisimple as first shown in \cite{BEH}. 
\end{exmp}

%% file: chapters/applications/applications.tex
\chapter{Applications of absolute Koszulity}
\input{chapters/applications/BGG}
\input{chapters/applications/coherence}

\input{chapters/applications/discussion}

%% file: chapters/applications/BGG.tex
\section{Application: BGG correspondence beyond complete intersections}

In this subsection we collect known constructions and examples from the literature to which we can apply our results. These will be mostly pulled from commutative algebra, and so from now on we let $R, S$ stand for Noetherian, standard graded commutative $k$-algebras over a field $k$. All such examples and the methods leading to them follow either explicitly or implicitly by work of Conca-Iyengar-Nguyen-R\"omer \cite{CINR}. \bigskip

We have already encountered the important example of short Gorenstein rings. The next proposition follows from the literature \cite{AL, AIS}. Moreover, a complete understanding of the Betti tables of indecomposable stable modules over short Gorenstein rings has been achieved by Avramov-Gibbons-Wiegand in \cite{AGW}.

We will give a separate proof of the following proposition, distinct from the existing literature, as to exemplify the results and methods of chapter 4 and 5. 
\begin{prop}[Short Gorenstein rings]\label{shortgorensteinring} Let $R$ be an Artinian Gorenstein algebra of embedding dimension $e \geq 2$ with $\m^3 = 0$, hence of socle degree $a = 2$. Then $R$ is absolutely Koszul with $\nu = a-1 = 1$. 
\end{prop}
\begin{proof}
	$R$ is absolutely Koszul by \cite{AIS}, where it is shown that the only indecomposable non-free graded $R$-module which are not Koszul are the cosyzygies of the residue field $k$, and so every module has an eventually Koszul syzygy. Note that for $a = 2$, $e \geq 2$ is necessary since $e = 1$ gives $R \cong k[x]/(x^3)$, which is not Koszul. \bigskip 

	We can also give a direct representation theoretic proof via tilting theory. Since $a = 2$, there is a full exceptional collection $\underline{\mathsf{mod}}^\Z\ \! R = \big\langle k, k(-1)[1] \big\rangle$, which can be verified directly as in the proof of Prop. \ref{helicity}, or alternatively follows from Orlov's Theorem. It is immediate that this collection is strong, and we have 
\[
	\underline{{\rm End}}_{gr R}(k \oplus k(-1)[1]) \cong \begin{pmatrix} {\rm \underline{End}}_{gr R}(k) & {\rm \underline{Ext}}^1_{gr R}(k, k(-1))\\ 0 & {\rm \underline{End}}_{gr R}(k(-1)[1]) \end{pmatrix} = \begin{pmatrix} k & R_1^* \\ 0 & k \end{pmatrix}
\]
which is isomorphic to the path algebra $kQ_e$ of the $e$-Kronecker quiver 
	\[
		\xymatrix@C=30pt{ 0 \ar@<3ex>[r]^-{y_1} \ar@<3ex>@{}[r]_{\vdots} \ar@<-3ex>[r]_-{y_e}  & 1}
	\]
	In particular $kQ_e$ is representation infinite for $e \geq 2$ and representation finite otherwise. Making use of the opposite tilting object $T^{op} = k \oplus k(-1)[1] \in (\underline{\mathsf{mod}}^{\Z}\ \! R)^{op}$, we have a contravariant equivalence
	\[
		G: (\underline{\mathsf{mod}}^{\Z}\ \! R)^{op} \xrightarrow{\cong} {\rm D}^b(\modsf kQ_e^{op}) 
	\] 
	sending $R$-modules to complexes of left $kQ_e$-modules, or covariant $Q_e$-representations, and we have $G(k) = P(0)$ and $G(k(-1)[1]) = P(1)$. Next, writing $\tau_R = \Se_R \circ [-1] = (2)[-2]$ for the Auslander-Reiten translate of $R$ and $\tau_{Q_e} = \Se_{kQ_e} \circ [-1] = D(kQ_e) \otimes^{\mathbb{L}}_{kQ_e}(-)[-1]$ for that of $kQ_e$, these are related by
	\[
		G \circ \tau_R \cong \tau^{-1}_{Q_e} \circ G.
	\]
	We also have $\tau^{-1}_{Q_e} = {\rm RHom}_{kQ_e}(D(kQ_e), -)[1]$. If $X$ is a $kQ_e^{op}$-module, we write $X = (X_0, X_1)$ for the corresponding quiver representation. For any $R$-module $M$, we then have
\begin{align*}
	\underline{{\rm Ext}}_{gr R}^{2n+i}(M, k(-2n)) &= {\rm \underline{Hom}}_{gr R}(M, \tau_R^{-2n}k[i])\\ 
	&\cong {\rm Hom}_{{\rm D}^b(kQ_e^{op})}(G(\tau_R^{-2n}k[i]), G(M))\\
	&\cong {\rm Hom}_{{\rm D}^b(kQ_e^{op})}(\tau_{Q_e}^{2n}G(k)[-i], G(M))\\
	&\cong {\rm Hom}_{{\rm D}^b(kQ_e^{op})}(P(0), \tau_{Q_e}^{-2n}G(M)[i])\\
	&\cong {\rm H}^{i}(\tau^{-2n}_{Q_e} G(M))_0\\
	\\
	\underline{{\rm Ext}}_{gr R}^{2n+1+i}(M, k(-2n-1)) &= {\rm \underline{Hom}}_{gr R}(M, \tau_R^{-2n}k(-1)[1][i])\\
	&\cong {\rm Hom}_{{\rm D}^b(kQ_e^{op})}(G(\tau_R^{-2n}k(-1)[1][i]), G(M))\\
	&\cong {\rm Hom}_{{\rm D}^b(kQ_e^{op})}(\tau_{Q_e}^{2n}G(k(-1)[1])[-i], G(M))\\
	&\cong {\rm Hom}_{{\rm D}^b(kQ_e^{op})}(P(1), \tau_{Q_e}^{-2n}G(M)[i])\\
	&\cong {\rm H}^{i}(\tau^{-2n}_{Q_e} G(M))_1.
\end{align*}
Setting $M = k$ so that $G(M) = P(0)$, we conclude that $R$ is Koszul if and only if the complex $\tau^{-2n}_{Q_e}P(0)$ is supported in cohomological degree zero (i.e. is quasi-isomorphic to a module) for all $n \geq 0$. This is well-known to characterise representation infinite quivers (see Appendix A.1), and so holds if and only if $e \geq 2$.

More generally this shows that $M$ has a linear resolution if and only if $G(M) \in \mod kQ_e^{op} \subseteq {\rm D}^b(kQ_e^{op}) $ is a module, and remains a module under the iteration of $\tau^{-1}_{Q_e}$. Assuming $e \geq 2$, the only indecomposable $kQ_e$-modules without this property are the preinjective modules $\{\tau^{n}I(0), \tau^n I(1)\}_{n \geq 0}$, since regular modules are closed under $\tau^{\pm 1}_{Q_e}$. Define an $R$-module $M$ to be completely linear if $\beta_{i,j}(M) = 0$ whenever $i \neq j$ for all $i, j \in \Z$. Then $G$ restricts to a bijection on the following classes, up to isomorphism: 
\[
	\{ \textup{Completely linear indecomposable stable $R$-modules}\} \leftrightarrow \{ \textup{Regular indecomposable $kQ_e$-modules} \}.
\]
Since $kQ_e$ is hereditary, complexes in $ {\rm D}^b(kQ_e^{op}) $ are formal, and translating through $G$ one sees that all indecomposable $R$-modules are $n$-linear except for the cosyzygies of $k$, which are sent onto the preinjective modules (up to suspension). It follows that $R$ is absolutely Koszul.
\end{proof}
\begin{rem} Applying Theorem \ref{thmC} and Minamoto's Theorem \ref{Minamototheorem}, we obtain equivalences of triangulated categories 
	\[
		{\rm D}^b(\qgr ((R^!)^{op})) \cong {\rm D}^b(\mathcal{H}^{\mathsf{lin}}(R))^{op} \cong (\underline{\mathsf{mod}}^\Z \ \! R)^{op} \cong {\rm D}^b(kQ_e^{op}) \cong {\rm D}^b(\qgr \Pi(Q^{op}_e)) 
	\]
	which are compatible with the t-structures, and which induce equivalence of abelian categories 
	\[
		\qgr ((R^!)^{op}) \cong \mathcal{H}^{\mathsf{lin}}(R)^{op} \cong \mathcal{H}^{\omega^{-1}}(kQ_e^{op}) \cong \qgr \Pi(Q^{op}_e). 
	\] 

	The two presentations by $\qgr ((R^!)^{op})$ and $\qgr \Pi(Q^{op}_e)$ correspond to picking different choices of ample sequences in $\mathcal{H}^{\mathsf{lin}}(R)^{op}$. 
\end{rem}\bigskip

The absolute Koszulity of short Gorenstein rings allows us to produce many more interesting examples of absolutely Koszul Gorenstein algebras, using a theorem of Conca-Iyengar-Nguyen-R\"omer: 
\begin{thm}[{\cite[Thm. 2.4]{CINR}}]\label{finiteflatdim} Let $R, S$ be standard graded commutative algebras, and let $\varphi: R \to S$ be a $k$-algebra morphism of finite flat dimension. If $S$ is absolutely Koszul, then so is $R$. 
\end{thm}

When $R$ is Gorenstein, we can apply this to the zero-dimensional reduction $R \to \overline{R} := R/({\bf x})$ by a regular sequence of linear forms ${\bf x} = (x_1, \dots, x_d)$. 
\begin{cor}\label{shortreduction} Let $R$ be Gorenstein of codimension $\geq 2$ with $\nu = d - 1 + a = 1$. Then $R$ is absolutely Koszul (and in particular Koszul).
\end{cor}
\begin{proof} 
	The codimension and $\nu$ invariant are preserved under passing to $\overline{R}$, in which case $\overline{R}$ is an Artinian Gorenstein algebra of embedding dimension $e \geq 2$ and socle degree $a = \nu + 1 = 2$. By Prop. \ref{shortgorensteinring}, $\overline{R}$ is absolutely Koszul and by Thm. \ref{finiteflatdim} so is $R$.
\end{proof}

Using this proposition, we can provide a great deal of interesting examples. 
\begin{cor}\label{examplenu1} Let $k$ be an algebraically closed field. The following algebras are absolutely Koszul Gorenstein with $\nu = 1$:
	\begin{enumerate}[a)] 
	\item The coordinate ring $R_{E, d}$ of an elliptic curve $E \subseteq \mathbb{P}^{d-1}$ of degree $d \geq 4$. 
	\item The coordinate ring $R_{X_d}$ of an anticanonically embedded smooth del Pezzo surface $X_d \subseteq \mathbb{P}^d$ of degree $d \geq 4$. 
	\item More generally, the coordinate ring $R_{X_d}$ of a smooth variety $X_d \subseteq \mathbb{P}^{d+n-2}$ of dimension $n \geq 2$ and degree $d \geq 4$. These include: 
		\begin{enumerate}[i)]
		\item The coordinate ring of a smooth complete intersection of two quadrics $X_4 = V(Q_1, Q_2) \subseteq \mathbb{P}^{n+2}$.
	\item The Pl\"ucker coordinate ring of the Grassmannian $X_5 = {\rm Gr}(2, 5) \subseteq \mathbb{P}^9$. 
	\item The coordinate ring of the Segre variety $X_6 = \mathbb{P}^2 \times \mathbb{P}^2 \subseteq \mathbb{P}^8$. 
	\item The coordinate ring of the Segre variety $X_6 = \mathbb{P}^1 \times \mathbb{P}^1 \times \mathbb{P}^1 \subseteq \mathbb{P}^7$. 
	\item The coordinate ring of the degree $7$ threefold $X_7 \subseteq \mathbb{P}^8$.
	\item The coordinate ring of the $2$nd Veronese embedding $X_8 = \mathbb{P}^3 \subseteq \mathbb{P}^{14}$.
\end{enumerate}
	\item $R/({\bf x})$ where $R$ is one of the above and ${\bf x} = (x_1, \dots, x_c)$ is a regular sequence of linear forms $x_i \in R_1$ of length $1 \leq c \leq d$.
\end{enumerate} 
\end{cor}
\noindent Let us put those examples in context before giving the proof. Given an smooth Fano variety $X$ of dimension $n \geq 2$, we define the index of $X$ by $\mathsf{ind}(X) = \sup \{\ r \ |\ -K_X \sim rH \textup{ for some } H \in {\rm Div}(X)\}$. We call $H$ a fundamental divisor if $\mathsf{ind}(X)H \sim -K_X$, and define the degree $d = \deg X = H^{n}$ via intersection theory. There is an upper bound $\mathsf{ind}(X) \leq \dim X + 1$, and so we define the coindex 
	\[
		\mathsf{coind}(X) = \dim X + 1 - \mathsf{ind}(X).
	\] 
	The cases $\mathsf{coind}(X) = 0, 1$ characterise projective space and quadric hypersurfaces, respectively, and the next case is more interesting. We say that $X = (X, H)$ is a (polarised) del Pezzo variety if $H$ is a fundamental divisor, $\mathsf{coind}(X) = 2$ and $H$ satisfies
	\[
		{\rm H}^n(X, \OO_X(mH)) = 0 \textup{ for all } 0 < n < \dim X \textup{ and } m \in \Z.
	\] 
The smooth polarised del Pezzo varieties were classified by Fujita in characteristic zero \cite[Thm. 3.3.1]{PS}, and the divisor $H$ is very ample as soon as $d \geq 3$. The corresponding embedding gives cubic hypersurfaces for $d =3$, and the examples $b), c)$ are the sole examples of degree $d \geq 4$ up to taking linear sections, at least in characteristic zero. \bigskip 

When $R$ is a graded Gorenstein algebra with $X = \proj R$ a smooth Fano variety, then polarised by $\OO_X(1) = \widetilde{R(1)}$, the coindex and the invariant $\nu$ are related. We have $\omega_X \cong \widetilde{\omega_R} \cong \OO_X(a)$ for $a < 0$, so that $-a \leq \mathsf{ind}(X)$, and so $\mathsf{coind}(X) \leq \nu + 1$, with equality when $\OO_X(1)$ comes from a fundamental divisor. Setting $\nu = 1$ then naturally leads to del Pezzo varieties.\bigskip

We now prove the corollary. 
\begin{proof}[Proof of Cor. \ref{examplenu1}] Case $a)$: It is well-known that $R_{E,d}$ is Gorenstein, and one can show this as follows. By \cite[Cor. 6.18]{EisSyzygy}, the ring $R_{E,d}$ is a normal Cohen-Macaulay domain, and the so natural morphism to the section ring $R_{E,d} \to \bigoplus_{n \geq 0}{\rm H}^0(E, \OO_E(n))$ is bijective, where $\OO_E(1)$ is the restriction to $E$ of $\OO_{\mathbb{P}^{d-1}}(1)$. That $R_{E,d}$ is Gorenstein follows from Stanley's numerical criterion \cite[Thm. 4.4]{Stanley} which says that a graded Cohen-Macaulay domain $R$ is Gorenstein if and only if its Hilbert series ${\rm H}_R(t)$ is symmetric, in that
\[ 
	{\rm H}_R(t^{-1}) = (-1)^{\dim R}t^s H_R(t)
\]
for some $s \in \Z$. One easily computes the Hilbert series of the section ring to be 
\[
	{\rm H}_{R_{E,d}}(t) = \frac{1 + (d-2)t + t^2}{(1-t)^2} 
\]
and symmetry follows. Hence $R_{E,d}$ is Gorenstein. We have $\nu = \dim E + a = 1+0 = 1$, and the requirement $\mathsf{codim}\ \! R_{E,d} \geq 2$ to use Cor. \ref{shortreduction} forces $d \geq 4$, with $d = 3$ corresponding to plane cubic curves. \bigskip 

\noindent Case $b)$: We can prove that $R_{X_d}$ is Gorenstein by reducing to the previous case. A generic linear section $X_{d,t}$ of $X_{d} \subseteq \mathbb{P}^d$ is a smooth irreducible curve, which has trivial canonical bundle by the adjunction formula (using that $\omega_{X_d}^{-1} = \OO_{X_d}(1)$) and so $X_{d,t}$ is an elliptic curve in $\mathbb{P}^{d-1}$. We then have $R_{X_d}/l = R_{X_{d,t}}$ for the corresponding regular linear form $l \in R_{X_d,1}$ and so $R_{X_d}$ must be Gorenstein by the previous case. We have $\nu = \dim X_d + a = 2 -1 = 1$, and again $\mathsf{codim}\ \! R(X_d) \geq 2$ forces $d \geq 4$, with $d = 3$ corresponding to a cubic surface.\bigskip

\noindent Case $c)$: Let $X_d \subseteq \mathbb{P}^{d+n-2}$ be a smooth variety of degree $d \geq 4$ and dimension $n \geq 2$. By Bertini's Theorem, a generic linear section $X_d' = X_d \cap \mathbb{P}^{d+n-3} \subseteq \mathbb{P}^{d+n-3}$ is smooth and irreducible so long as $\dim X_d \geq 2$, and $X_d' \subseteq \mathbb{P}^{d+n-3}$ is another smooth variety of degree $d$ and dimension $n-1$. Iterating, we construct a smooth del Pezzo surface $\widetilde{X}_d \subseteq X_d$ of the same degree, with coordinate ring $R_{\widetilde{X}_d} = R_{X_d}/({\bf x})$ a quotient by a regular sequence of linear forms ${\bf x} = (x_1, \dots, x_{n-2})$. Since $R_{\widetilde{X}_d} = R_{X_d}/({\bf x})$ is Gorenstein, so is $R_{X_d}$, and since $R_{\widetilde{X}_d}$ is absolutely Koszul, so is $R_{X_d}$ by Thm. \ref{finiteflatdim}. \bigskip 

\noindent Case $d)$: Finally, the property of $R$ being Gorenstein with $\nu = \dim R - 1 + a = 1$ and codimension $\geq 2$ is unchanged under reduction modulo ${\bf x}$, and Cor. \ref{shortreduction} shows that $R/({\bf x})$ is absolutely Koszul. 
\end{proof} 

There also are a few interesting examples with $\nu = 2$. For the next proof, recall that the multiplicity $e(R)$ of $R$ in dimension $d > 0$ is given in terms of the Hilbert polynomial of $R$ by the coefficient $e_0$ in the expansion \cite[Defn. 4.1.5]{BHBook} 
\[
	P_R(n) = \sum_{i=0}^{d-1} (-1)^{d-1-i} e_{d-1-i} \binom{n+i}{i}
\]
and that the Hilbert polynomial of $R$ can be computed from $X = \proj R$ as $P_R(n) = \chi(\OO_X(n))$. The following examples are due to Conca-Iyengar-Nguyen-R\"omer. 
\begin{prop}[Conca-Iyengar-Nguyen-R\"omer, {\cite[Thm. 5.1]{CINR}}]\label{canonicalcurves} Let $k$ be an algebraically closed field of characteristic zero. The following algebras are absolutely Koszul Gorenstein algebras with $\nu = 2$:
	\begin{enumerate}[a)] 
	\item The coordinate ring $R_C$ of the canonical embedding of a non-hyperelliptic smooth projective curve $C \subseteq \mathbb{P}^{g-1}$ of genus $g \geq 3$, which is neither trigonal nor isomorphic to a plane quintic. 
	\item The coordinate ring $R_X$ of a smooth projective variety $X \subseteq \mathbb{P}^{g+n-2}$ of dimension $n \geq 2$ whose section by a generic linear subspace $X \cap \mathbb{P}^{g-1}$ is a canonical curve of the above form. 
\end{enumerate} 
\end{prop}
\begin{proof}
	The ring $R_C$ is a Gorenstein domain with isolated singularity, with $\nu = \dim C + a = 1 + 1 = 2$ since $C$ is canonically embedded. In \cite[Thm. 5.1]{CINR}, it is shown\footnote{Our presentation is somewhat out of order, since they show Thm. 5.1 by reducing to the case of the coordinate ring of a canonical curve. However we point out that the curves arising this way are precisely those covered by Petri's Theorem.} that any graded Gorenstein domain with isolated singularity in characteristic zero with $2 e(R) = 2 \mathsf{codim}\ \! R + 2$ is absolutely Koszul, as soon as it has quadratic relations. We can compute the Hilbert polynomial of $R_C$ via Riemann-Roch as $P_{R_C}(n) = \chi(\OO_C(n)) = \chi(\omega_C^{\otimes n}) = (2g-2)n + 1 - g$, which we write 
	\[
		P_{R_C}(n) = -(2-2g)(n+1) + g - 1 
	\]
	to obtain $e(R_C) = e_0 = g-1$. Since $\mathsf{codim}\ \! R_C = g-2$, the equality holds, while the condition of having quadratic relations is equivalent to the stated conditions by Petri's Theorem. The second case follows from Bertini's Theorem as before. 
\end{proof}

Next, we collect various classes of absolutely Koszul algebras from the literature. Note that some references work in the setting of local Noetherian rings, but the arguments immediately adapt to the graded case. \bigskip 

We say that an algebra $R$ has the Backelin-Roos property if there is a Golod map $\varphi: Q \twoheadrightarrow R$ with $Q$ a complete intersection, see \cite[Sect. 3]{CINR} for the definition of Golod maps. Such algebras are widespread and were originally introduced in connections with rationality questions for Poincar\'e series. 
\begin{exmp}\label{literature} The following graded $k$-algebras are absolutely Koszul. 
	\begin{enumerate}[a)]
	\item Koszul algebras with the Backelin-Roos property \cite[Thm. 5.9]{HI}, \cite[Prop. 3.4]{CINR}; 
	\item Retracts of absolutely Koszul algebras \cite[Prop. 2.3(3)]{CINR}. 
	\item Artinian Gorenstein algebras of embedding dimension $e \geq 3$ and socle degree $a = 3$ with an exact pair of linear zero divisors; in particular a generic Artinian Gorenstein algebras with $e \geq 3$ and $a = 3$ \cite[Prop. 4.1, Thm. 4.2]{HS}. 
	\item Koszul algebras $S$ of small embedding codepth, that is $\edim S - \depth S \leq 3$, or Koszul Gorenstein algebras $R$ of codimension $\leq 4$ \cite[Thm. 6.4, Prop. 6.3]{AKM}. 
	\item Algebras of the form $S/(I+L)$ with $S = k[x_1, \dots, x_n]$, $I, L$ quadratic monomial ideals with $I$ generated by a regular sequence and $L$ having a $2$-linear resolution over $S$ \cite[Thm. 4.1]{CINR}.
	\item The $c$-th Veronese subalgebra $S^{(c)}$ of $S = k[x_1, \dots, x_n]$ in characteristic zero, for $n, c$ taking values \cite[Cor. 5.4]{CINR}:
		\begin{enumerate}[i)]
		\item $n \leq 3$ and all $c$;
		\item $n \leq 4$ and $c \leq 4$;
		\item $n \leq 6$ and $c = 2$.
	\end{enumerate}
	\item The Segre product $S_{m,n}$ of $k[x_1, \dots, x_m]$ with $k[x_1, \dots, x_n]$ in characteristic zero, with $m \leq n$ taking values \cite[Prop 5.9]{CINR}:
		\begin{enumerate}[i)]
		\item $m \leq 2$;
		\item $m = 3$ and $n \leq 5$;
		\item $m = n = 4$. 
	\end{enumerate}
\end{enumerate}
\end{exmp} 
Note that most of $c)-g)$ are special cases of $a)$, as shown in the respective references. We make a special mention of the Koszul Gorenstein algebras of codimension $\leq 4$, which should give rise to a large class of interesting examples by the Buchsbaum-Eisenbud structure theorem in codimension $\leq 3$ \cite{BE} (see also \cite{ReidGorenstein}). \bigskip 

We now take a more detailed look at some of the previous examples. 
\subsection*{The cone over an elliptic curve of degree $d \geq 4$}
Let $R_{E,d}$ be the coordinate ring of $E \subseteq \mathbb{P}^{d-1}$ from Cor. \ref{examplenu1}, and we assume that $\OO_E(1) \cong \OO_E(dx)$ for some closed point $x \in E$. Applying Orlov's Theorem\footnote{We fix a choice of cut-off $i = 0$ throughout this subsection.} and Theorem A, we obtain equivalences of triangulated categories.
\[
	{\rm D}^b(\coh E) \cong \underline{{\rm MCM}}^{\mathbb{Z}}(R_{E, d}) \cong {\rm D}^b(\mathcal{H}^{\mathsf{lin}}(R_{E,d})).
\]
In his thesis \cite{Pavlov}, Pavlov characterised the images of indecomposable linear MCM modules in $ {\rm D}^b(\coh E) $ under Orlov's equivalence \cite[Chp. 6, Sect. 6.2]{Pavlov}. His method leads more generally to a description of the induced t-structure with hereditary heart $\mathcal{H}^{\mathsf{lin}}(R_{E, d}) \subseteq  {\rm D}^b(\coh E) $. We recall a few results and facts from \cite{Pavlov}. \bigskip

Recall that any spherical object $\F \in {\rm D}^b(\coh E) $ gives rise to an autoequivalence $T_\F \in {\rm Aut}( {\rm D}^b(\coh E) )$, the associated spherical twist, also called Thomas-Seidel twist (see \cite{Pavlov}, \cite{Huy} for definitions). Given our closed point $x \in E$ above, the skyscraper sheaf $k(x)$ is a spherical object, and $\OO_E$ is spherical since $E$ is Calabi-Yau. This defines two autoequivalences of $ {\rm D}^b(\coh E) $ which we denote 
\begin{align*}
	\mathbb{A} &:= T_{\OO_E} \\
	\mathbb{B} &:= T_{k(x)}
\end{align*}
We have $\mathbb{B} = T_{k(x)} \cong - \otimes_{\OO_E} \OO_E(x)$ (\cite[Lemma 2.4.3]{Pavlov}), but $\mathbb{A} = T_{\OO_E}$ has no such simple description. Moreover, $\mathbb{A}$ and $\mathbb{B}$ satisfy the braid relations $\mathbb{A}\mathbb{B}\mathbb{A} \simeq \mathbb{B}\mathbb{A}\mathbb{B}$, see \cite[Prop. 2.13]{ST}.\bigskip 

Next, consider the degree shift autoequivalence $M \mapsto M(1)$ of $ \underline{{\rm MCM}}^{\mathbb{Z}}(R_{E,d})$, and write $\sigma$ for the corresponding autoequivalence of $ {\rm D}^b(\coh E) $ under Orlov's equivalence $ {\rm D}^b(\coh E) \cong \underline{{\rm MCM}}^{\mathbb{Z}}(R_{E,d})$. 
\begin{prop}[Pavlov, {\cite[Thm. 4.1.2, Lemma 2.4.4]{Pavlov}}] In the above situation, we have the following identifications: 
\begin{enumerate}[1)]
\item The suspended sheaf $\OO_E[1] \in {\rm D}^b(\coh E) $ corresponds to $k^{st} \in \underline{{\rm MCM}}^{\mathbb{Z}}(R_{E,d})$ under Orlov's equivalence. 
\item The autoequivalence $\sigma \in {\rm D}^b(\coh E) $ is given by $\sigma = \mathbb{B}^d \circ \mathbb{A}$.
	\end{enumerate}\qed
\end{prop}
\begin{rem}
	We can rewrite $\mathbb{B}^d \cong - \otimes_{\OO_E} \OO_E(x)^{\otimes d} = - \otimes_{\OO_E} \OO_E(dx) = - \otimes_{\OO_E} \OO_E(1)$ for $\OO_E(1)$ coming from the embedding $E \subseteq \mathbb{P}^{d-1}$, and thus rewrite $\sigma = -\otimes_{\OO_E} \OO_E(1) \circ \mathbb{A}$. 
\end{rem}

For any $M \in \underline{{\rm MCM}}^{\mathbb{Z}}(R_{E, d})$, write $\F_M$ for the corresponding complex of sheaves in $ {\rm D}^b(\coh E) $. We have $\omega_E \cong \OO_E$, and so Serre duality gives
\[
	{\rm Hom}_{ {\rm D}^b(E) }(\F, \OO_E[1]) \cong {\rm Hom}_{ {\rm D}^b(E) }(\OO_E, \F)^*
\]
for any $\F \in {\rm D}^b(\coh E)$. In particular, the Betti numbers of $M$ can be computed as (\cite[Thm. 4.1.2]{Pavlov}) 
\begin{align*}
	\beta_{i,j}(M) = \dim_k \underline{{\rm Ext}}^i_{gr R_{E,d}}(M, k^{st}(-j))&= \dim_k \underline{ {\rm Hom}}_{gr R_{E,d}}(M, k^{st}(-j)[i]) \\
	&= \dim_k \underline{ {\rm Hom}}_{gr R_{E,d}}(M(j)[-i], k^{st}) \\
	&= \dim_k {\rm Hom}_{ {\rm D}^b(E) }(\sigma^j\F_M[-i], \OO_E[1]) \\
	&= \dim_k {\rm Hom}_{ {\rm D}^b(E) }(\OO_E, \sigma^j\F_M[-i] )\\
	&= \dim_k {\rm H}^{0}(E, \sigma^j\F_M[-i]).
\end{align*}
Define a new autoequivalence $\upgamma = \sigma [-1] = \mathbb{B}^d \circ \mathbb{A} \circ [-1]$. This last quantity can be rewritten 
\begin{align*}
	\beta_{i,j}(M) = \dim_k {\rm H}^{0}(E, \sigma^j\F_M[-i]) &= \dim_k {\rm H}^0(E, \upgamma^j \F_M [j-i]) \\
	&= \dim_k {\rm H}^{j-i}(E, \upgamma^j \F_M).
\end{align*}
Finally, recall that the bounded t-structure $t^{\mathsf{lin}}$ is defined as 
\begin{align*}
	\underline{{\rm MCM}}^{\leq 0}(R_{E,d}) &= \{ M \ | \ \beta_{i,j}(M) = 0 \textup{ for } j-i > 0 \textup{ whenever } i \gg 0 \}\\
	&= \{ M \ | \ \beta_{i,j}(M) = 0 \textup{ for } j-i > 0 \textup{ whenever } j \gg 0 \}\\
	\underline{{\rm MCM}}^{\geq 0}(R_{E,d}) &= \{ M \ | \ \beta_{i,j}(M) = 0 \textup{ for } j-i < 0 \textup{ whenever } i \gg 0 \}\\
	&= \{ M \ | \ \beta_{i,j}(M) = 0 \textup{ for } j-i < 0 \textup{ whenever } j \gg 0 \}.
\end{align*}
We see that this gives rise to a bounded t-structure $t^{\upgamma} = ({\rm D}^{\leq 0, \upgamma}, {\rm D}^{\geq 0, \upgamma})$ on $ {\rm D}^b(\coh E) $ defined by 
\begin{align*}
	{\rm D}^{\leq 0, \upgamma} &= \{ \F \ | \ {\rm H}^{n}(E, \upgamma^j\F) = 0 \textup{ for } n > 0 \textup{ whenever } j \gg 0 \}\\
	{\rm D}^{\geq 0, \upgamma} &= \{ \F \ | \ {\rm H}^{n}(E, \upgamma^j\F) = 0 \textup{ for } n < 0 \textup{ whenever } j \gg 0 \}.
\end{align*}
The hereditary heart $\mathcal{H}^{\upgamma}(E) := {\rm D}^{\leq 0, \upgamma} \cap {\rm D}^{\geq 0, \upgamma} \subseteq {\rm D}^b(\coh E) $ then consists of complexes eventually without non-zero sheaf cohomology 
\[
	\mathcal{H}^{\upgamma}(E) = \{ \F \in {\rm D}^b(\coh E) \ | \ {\rm H}^i(E, \upgamma^j \F) = 0 \textup{ for } i \neq 0 \textup{ for all } j \gg 0\}. 
\]
This category can likely be described further using Pavlov's work, see in particular \cite[Sect. 6.2]{Pavlov} for a description of the indecomposables $\F_M$ corresponding to indecomposable linear MCM modules $M$ under $\mathcal{H}^{\upgamma}(E) \cong \mathcal{H}^{\mathsf{lin}}(R_{E, d})$. \bigskip 

Lastly, our work shows that $ {\rm D}^b(\mathcal{H}^{\upgamma}(E)) \cong {\rm D}^b(\coh E) $, and so for each indecomposable sheaf $\F \in \coh E$, there is a unique $n = n_\F \in \Z$ for which $\F[n] \in \mathcal{H}^{\upgamma}(E)$. It would be interesting to have a description of $n_\F$ in terms of $\F$. 

\subsection*{Cones over smooth del Pezzo varieties of degree $d \geq 4$}
The hereditary heart $\mathcal{H}^{\mathsf{lin}}(R_{X_d})$ for $X_d$ a smooth del Pezzo variety of degree $d \geq 4$ are examples of hereditary Ext-finite abelian categories with Serre duality, and such categories tend to be rather special. Along this line, let us recall a famous structure theorem of Happel. Recall that an abelian category $\CC$ is connected if it is not of the form $\CC = \CC_1 \times \CC_2$ for two non-zero orthogonal full subcategories $\CC_{1}, \CC_2$. 

\begin{thm}[Happel \cite{HappelTilting}] Let $k$ be an algebraically closed field and $\A$ an hereditary, Ext-finite, connected $k$-linear abelian category. Assume that $ {\rm D}^b(\A) $ admits a tilting object. Then, up to derived equivalence $\A$ is of the form 
	\[
		{\rm D}^b(\A) \cong \begin{cases} {\rm D}^b(\modsf kQ) & \textup{ with Q a finite acyclic quiver;} \\ {\rm D}^b(\coh \X) & \textup{ with }\X = \mathbb{P}^1(p_1, \dots, p_t) \textup{ a weighted projective line}. \end{cases} 
	\] 
\end{thm}

Note that $\A$ itself need not be equivalent to $\modsf kQ$ or $\coh \X$, as it can fail to be Noetherian. As immediate special case, we obtain:
\begin{cor} Let $k$ be algebraically closed and $R$ an absolutely Koszul, commutative, Gorenstein $k$-algebra with isolated singularities and $\nu = 1$. Assume that $\mathcal{H}^{\mathsf{lin}}(R)$ is connected and that $ \underline{{\rm MCM}}^{\mathbb{Z}}(R) \cong {\rm D}^b(\mathcal{H}^{\mathsf{lin}}(R)) $ admits a tilting object. Then, up to derived equivalence $\mathcal{H}^{\mathsf{lin}}(R)$ has one of two forms: 
	\[
		{\rm D}^b(\mathcal{H}^{\mathsf{lin}}(R)) \cong \begin{cases} {\rm D}^b(\modsf kQ) & \textup{ with Q a finite acyclic quiver;} \\ {\rm D}^b(\coh \X) & \textup{ with }\X = \mathbb{P}^1(p_1, \dots, p_t) \textup{ a weighted projective line}. \end{cases} 
	\]
\end{cor}
Note that a priori the connectedness assumption can fail to hold\footnote{In the case $\nu = 0$, consider $R = k[x,y]/(xy)$. It is easy to see that $\mathcal{H}^{\mathsf{lin}}(R) \cong \modsf k \times \modsf k$, with two simple objects corresponding to $R/(x)$ and $R/(y)$.}, however Happel's Theorem will describe the connected summands. \bigskip

Happel's Theorem applies to the coordinate rings $R_{X_d}$ of smooth del Pezzo surfaces $X_d \subseteq \mathbb{P}^d$ of degree $d \geq 4$ of Cor. \ref{examplenu1}. To see this, we will show that $ \underline{{\rm MCM}}^{\mathbb{Z}}(R_{X_d})$ contains a tilting object (essentially the same argument was already used in Chapter 2 for the cubic surface in $\mathbb{P}^3$). Recall that $a = a(R_{X_d}) = -1$ since $X_d$ is anticanonically embedded, and so Orlov's Theorem yields 
\[
	{\rm D}^b(\coh X_d ) = \big\langle \OO_{X_d}, \Phi_0( \underline{{\rm MCM}}^{\mathbb{Z}}(R_{X_d})) \big\rangle 
\]
Recall that $X_d$ is abstractly isomorphic either to the blowup $\pi: Bl_m\mathbb{P}^2 \twoheadrightarrow \mathbb{P}^2$ of $m = 9-d$ points on $\mathbb{P}^2$ in general position, or additionally to the variety $X_8 = \mathbb{P}^1 \times \mathbb{P}^1 \hookrightarrow \mathbb{P}^8$ in degree $8$. \bigskip 

In the first case $X_d \cong Bl_m\mathbb{P}^2 \twoheadrightarrow \mathbb{P}^2$, let $E_j = \pi^{-1}(p_j)$ denote the exceptional divisors. Then by Orlov's blow-up formula \cite[Sect. 11.2]{Huy} (see also references in \cite{AKO}), the derived pullback provides a fully faithful embedding $\mathbb{L}\pi^*: {\rm D}^b(\coh \mathbb{P}^2) \hookrightarrow {\rm D}^b(\coh Bl_m \mathbb{P}^2) $ with semiorthogonal decomposition 
\begin{align*}
	{\rm D}^b(\coh X_d ) 	&= \big\langle \mathbb{L}\pi^* {\rm D}^b(\coh \mathbb{P}^2), \OO_{E_1}, \cdots, \OO_{E_{m}} \big\rangle \\
				&= \big\langle \pi^*\OO_{\mathbb{P}^2}, \pi^*(\OO_{\mathbb{P}^2}(1)), \pi^*(\OO_{\mathbb{P}^2}(2)), \OO_{E_1}, \cdots, \OO_{E_{m}} \big\rangle \\
				&= \big\langle \OO_{X_d}, \pi^*(\OO_{\mathbb{P}^2}(1)), \pi^*(\OO_{\mathbb{P}^2}(2)), \OO_{E_1}, \cdots, \OO_{E_{m}} \big\rangle. 
\end{align*} 
By \cite[Thm. 2.5]{AKO}, the latter forms a full strong exceptional collection of sheaves in $ {\rm D}^b(\coh X_d) $. One can perform calculations directly: by adjunction and using that $\pi$ has rational fibres, we have
\begin{align*}
	{\rm Hom}_{ {\rm D}^b(X_d) }(\pi^*(\OO(n)), \OO_{E_j}) & \cong {\rm Hom}_{ {\rm D}^b(\mathbb{P}^2)}(\OO(n), {\rm R}\pi_* (\OO_{E_j}) ) \\
	& \cong {\rm Hom}_{ {\rm D}^b(\mathbb{P}^2)}(\OO(n), \OO_{p_j} )
\end{align*} 
and one sees that its endomorphism algebra $\Lambda \cong k\widetilde{Q}/I$ is given by the quiver $\widetilde{Q} = \widetilde{Q}^{(1)}_m$ below 
		\[
	\xymatrix@R=5pt@C=20pt{ 
		&&&\bullet\\
		&&&\bullet\\
		\bullet \ar@<1ex>[r]^{ \{x_i \}} \ar[r] \ar@<-1ex>[r] & \bullet \ar@<1ex>[r]^{ \{x_i \}}  \ar[r] \ar@<-1ex>[r] &\bullet \ar[uur]^{q_1} \ar[ur]_{\ q_2}  \ar[ddr]_{q_m}& \vdots\\
		&&&\vdots\\
		&&&\bullet
	}
\]
with $\{x_i \}$ a basis of sections for ${\rm H}^0(\mathbb{P}^2, \OO_{\mathbb{P}^2}(1))$ and with relations $x_ix_{i'} = x_{i'}x_i$ and $q_j l_j(x_0,x_1,x_2) = 0 = q_j l'_j(x_0, x_1, x_2)$, where the $j$-th point $p_j = V(l_j, l'_j) \in \mathbb{P}^2$ is cutout by said linear forms $l_j, l'_j$. By Orlov's Theorem $ \underline{{\rm MCM}}^{\mathbb{Z}}(R_{X_d})$ inherits a full strong exceptional collection with endomorphism algebra $kQ/I$ with $Q = Q_m^{(1)}$ 
		\[
	\xymatrix@R=5pt@C=20pt{ 
		&&\bullet\\
		&&\bullet\\
		\bullet \ar@<1ex>[r]^{ \{x_i \}}  \ar[r] \ar@<-1ex>[r] &\bullet \ar[uur]^{q_1} \ar[ur]_{\ q_2}  \ar[ddr]_{q_m}& \vdots\\
		&&\vdots\\
		&&\bullet
	}
\]
with the same relations.\footnote{There is a strong similarity with the Squid algebra of Chapter 3.} \bigskip

Next let $X_8 = \mathbb{P}^1 \times \mathbb{P}^1 \subseteq \mathbb{P}^8$. Let $\E = \OO_{\mathbb{P}^1} \oplus \OO_{\mathbb{P}^1}(1)$, and define $\OO(i, j) := \OO_{\mathbb{P}^1}(i) \boxtimes \OO_{\mathbb{P}^1}(j)$. It is well-known that $\E \boxtimes \E$ is a tilting bundle on $\mathbb{P}^1 \times \mathbb{P}^1$, whose summands form a full strong exceptional collection of line bundles
\[
	{\rm D}^b(\coh X_8) = \big\langle \OO, \OO(0, 1), \OO(1, 0), \OO(1, 1)  \big\rangle.
\]
with corresponding quiver path algebra $k\widetilde{Q}/I$ given by $\widetilde{Q} = \widetilde{Q}^{(2)}$
\[
	\xymatrix{\bullet \ar@<-.5ex>[d]_u \ar@<.5ex>[d]^v  \ar@<-.5ex>[r]_x \ar@<.5ex>[r]^y & \bullet \ar@<-.5ex>[d]_u \ar@<.5ex>[d]^v  \\ \bullet \ar@<-.5ex>[r]_x \ar@<.5ex>[r]^y  & \bullet}
\]
with commuting square relations $I = (ux - xu, vx-xv, uy-yu, vy-yv)$. Similarly, Orlov's Theorem gives 
\[
	{\rm D}^b(\coh X_d ) = \big\langle \OO_{X_d}, \Phi_0( \underline{{\rm MCM}}^{\mathbb{Z}}(R_{X_8})) \big\rangle
\]
and again $ \underline{{\rm MCM}}^{\mathbb{Z}}(R_{X_8})$ inherits a full strong exceptional sequence, with endomorphism algebra $kQ$ where $Q = Q^{(2)}$ is 
\[
	\xymatrix{& \bullet \ar@<-.5ex>[d]_x \ar@<.5ex>[d]^y  \\ \bullet \ar@<-.5ex>[r]_u \ar@<.5ex>[r]^v  & \bullet}
\]
In particular the path algebra $kQ^{(2)}$ has no relations. \bigskip 

Because the quiver path algebras $kQ^{(1)}/I$ and $kQ^{(2)}$ are connected, $\mathcal{H}^{\mathsf{lin}}(R_{X_d})$ is a connected category, and we have just shown that $ {\rm D}^b(\mathcal{H}^{\mathsf{lin}}(R)) $ contains a tilting complex. Happel's Theorem then implies an equivalence 
	\[
		{\rm D}^b(\mathcal{H}^{\mathsf{lin}}(R_{X_d})) \cong \begin{cases} {\rm D}^b(\modsf kQ) & \textup{ with Q a finite acyclic quiver;} \\ {\rm D}^b(\coh \X) & \textup{ with }\X = \mathbb{P}^1(p_1, \dots, p_t) \textup{ a weighted projective line}. \end{cases} 
	\]
	It would be very interesting to know which of these occur. The following table contains the extent of the author's knowledge (recall that $m = 9-d$): 
\[
	\begin{tabular}{|c|c|c|c|c|c|c|c|}
	\hline
$d$	& $d = 4$ & $d=5$ & $d=6$ & $d=7$ & $d=8$ & $d=8$ & $d=9$ \\
	$X_d$ & $Bl_5\mathbb{P}^2$ & $Bl_4\mathbb{P}^2$  & $Bl_3\mathbb{P}^2$  & $Bl_2\mathbb{P}^2$  & $Bl_1\mathbb{P}^2$ & $\mathbb{P}^1 \times \mathbb{P}^1$ & $\mathbb{P}^2$ \\
	\hline
	${\rm D}^b(\mathcal{H}^{\mathsf{lin}}(R_{X_d}))$&   	${\rm D}^b(\coh \mathbb{P}^1(2,2,2,2,2))$ & ?? & ??& ?? & ?? & ${\rm D}^b(\modsf kQ^{(2)})$ & ${\rm D}^b(\modsf kQ_0^{(1)})$ \\
	\hline
\end{tabular}
\] 
The case $d = 4$ corresponds to a complete intersection of two quadrics $Bl_5\mathbb{P}^2 \cong X_4 = V(Q_1, Q_2) \subseteq \mathbb{P}^4$, and we have already given an exposition of this case in chapter 2, using methods of Buchweitz and Kuznetsov. In particular $\coh \mathbb{P}^1(2,2,2,2,2) \cong \coh \OO$ for the hereditary order $\OO$ on $\mathbb{P}^1$ associated to the pencil of quadrics. The second case $d = 8$ follows since the quiver $Q^{(2)}$ above had no relation. \bigskip 

The quiver $kQ_0^{(1)}$ for $d = 9$ is the $3$-Kronecker quiver, and this case was studied by Iyama-Yoshino \cite{IY} (see also \cite{KMVdB}) who classified the rigid indecomposables MCM modules over the third Veronese\footnote{Technically they worked with the invariant ring of the natural $\mu_3$ action by cube roots of unity, which agrees with the third Veronese algebra in characteristic $\neq 3$.} of $S = k[x_0, x_1, x_2]$ by use of Kac's Theorem. Since the result of Kac's Theorem has been extended to hereditary categories of the form $\coh \X$ by Crawley-Boevey \cite{CB}, one could contemplate looking for similar classifications for all rings $R_{X_d}$, $d \geq 4$, at least once the above table has been filled. This will be investigated in later work. \bigskip

%% file: chapters/applications/coherence.tex
\section{Application: The Coherence Conjectures of Minamoto and Bondal}
We finally return after a long digression to the original motivation for the study of the t-structure $t^{\mathsf{lin}}$ on $ \underline{{\rm MCM}}^{\mathbb{Z}}(A)$. In this section we take $A$ to be a Koszul Frobenius $\Bbbk$-algebra of top degree $a \geq 1$, and for simplicity work over $\Bbbk = k$, where $k$ is our fixed ground field. \bigskip

We have seen in chapter 4 that one can attach a finite dimensional $k$-algebra $\Lambda$ to $A$ by the construction
\[
	\Lambda = \underline{{\rm End}}_{gr A}(\bigoplus_{i=0}^{a-1} \Omega^i k(i)) 
\]
such that both $\Lambda$ and $\Lambda^{op}$ are $d$-representation infinite algebras, of global dimension $d = a-1$. As these are endomorphism algebras of tilting objects, these come equipped with equivalence of triangulated categories
\[
	\underline{\mathsf{mod}}^\Z \ \! A \xrightarrow{\cong} {\rm D}^b(\Lambda) 
\]
\[
	(\underline{\mathsf{mod}}^\Z \ \! A)^{op} \xrightarrow{\cong} {\rm D}^b(\Lambda^{op}).
\]
and we inherit a pair of $(d+1)$-preprojective algebras $\Pi(\Lambda)$ and $\Pi(\Lambda^{op})$, both of global dimension $d+1=a$. Moreover, since $A$ is Frobenius, its Koszul dual $A^! = {\rm Ext}_A^*(k, k)$ is an Artin-Schelter regular algebra \cite[Thm. 5.10]{SMI} of global dimension $a$. Now, the Koszul and Frobenius properties are left-right symmetric and it isn't hard to see that $(A^{op})^! \cong (A^!)^{op}$, and so $(A^!)^{op}$ is also Artin-Schelter regular of global dimension $a$. 
The next theorem summarises what we have shown in the previous two chapters. Recall that all notions refer to right modules. We set $\Pi = \Pi(\Lambda^{op})$ and $E = (A^!)^{op}$. 
\begin{thm}\label{summary6.2} Let $A$ be a Koszul Frobenius algebra, and let $\Pi$ and $E$ be the graded algebras defined above. The following are equivalent: 
	\begin{enumerate}[i)]
	\item $A$ is absolutely Koszul.
	\item $\Pi$ is coherent.
	\item $E$ is coherent. 
\end{enumerate}
Moreover, when these conditions hold, we have equivalences of triangulated categories 
\[
{\rm D}^b(\qgr E) \cong (\underline{\mathsf{mod}}^\Z \ \! A)^{op} \cong {\rm D}^b(\Lambda^{op}) \cong {\rm D}^b(\qgr \Pi )
\]
preserving the relevant t-structures, so that this descends to an equivalence of abelian categories
\[
\qgr E \cong \mathcal{H}^{\mathsf{lin}}(A)^{op} \cong \mathcal{H}(\Lambda^{op}) \cong \qgr \Pi.
\]
\end{thm}
\begin{proof} 
The equivalence $i) \iff ii)$ follows from Prop. \ref{summary4.2} and Thm. A. That of $i) \iff iii)$ follows from Thm. \ref{thmB}, and the last statement is a combination of Thm. \ref{thmC}, Prop. \ref{pullbacktstructure} and Minamoto's Theorem \ref{Minamototheorem}.
\end{proof}

We now construct counterexamples to coherence. 
\begin{thm}\label{abscounterexamples} For each $n \geq 4$, there is a (commutative) Koszul Frobenius $k$-algebra $R_n$ of socle degree $n$ which is not absolutely Koszul. 
\end{thm} 
\noindent The associated higher preprojective algebra $\Pi_{n} = \Pi(\Lambda_{R_n}^{op})$ and AS-regular algebra $E_n = (R^!_{n})^{op}$ are then counterexamples to Conjectures \ref{Minamotoconjecture} and \ref{Bondalconjecture} in all global dimension $n \geq 4$. \bigskip

We first need a lemma, which was briefly stated in the last section. 
\begin{lem}[{\cite[Prop. 2.3(3)]{CINR}}]\label{abskoszulretract} Let $R, S$ be standard graded commutative $k$-algebras. Let $\varphi: R \to S$ be a retract of $k$-algebras with section $\sigma: S \to R$, $\varphi \sigma = id_S$. If $R$ is absolutely Koszul, then so is $S$. 
\end{lem} 

\begin{con} We begin by taking a bad Koszul algebra in the sense of Roos. Let $S = k[x,y]/(x,y)^2$, and consider $S_m = S^{\otimes m}$. Then $S_2$ is the bad Koszul algebra of Roos from Thm. \ref{badkoszulalgebra}; in particular it cannot be absolutely Koszul by Prop. \ref{uniformrationality}. Moreover, $S_m$ retracts onto $S_2$ for all $m \geq 2$ and so these are never absolutely Koszul. However the $S_m$ are not Frobenius. \bigskip 

	Let $R$ = $S \ltimes (S(2))^*$ be the trivial algebra extension by the symmetric $S$-bimodule $(S(2))^* = S^*(-2)$, with multiplication
	\[
		(r, \varphi) \cdot (r', \varphi') = (rr', r\varphi' + r'\varphi).
	\]
	Then $R$ is commutative graded Frobenius, with Hilbert series $H_R(t) = 1 + 4t + t^2$. In particular $R$ satisfies $R_{\geq 3} = 0$ and so $R$ is Koszul by Prop. \ref{shortgorensteinring} (in fact, absolutely Koszul). In fact, setting $V = {\rm span}_k \{x,y\}$, the algebra $R$ is simply the algebra 
	\[
		R = k \oplus (V \oplus V^*) \oplus k 
	\]
	with multiplication induced from the self duality pairing on $V \oplus V^*$. The Koszul Frobenius algebra $R = S \ltimes (S(2))^*$ retracts onto $S$ by construction, with $\varphi: R \to S$ the projection and section $\sigma: S \to R$ the natural inclusion. Define
	\[
		R_n = \begin{cases} R^{\otimes m} & n = 2m\\ R^{\otimes m} \otimes k[\epsilon] & n = 2m+1  \end{cases}
	\]
	where $k[\epsilon] = k[x]/(x^2)$. Then $R_n$ is a commutative Koszul Frobenius algebra of socle degree $n$. For all $m \geq 2$, we have $k$-algebra retractions induced from $\varphi, \sigma$
	\[
		S_2 \leftrightarrows S_m \leftrightarrows R^{\otimes m} \leftrightarrows R^{\otimes m} \otimes k[\epsilon]
	\]
	and so $R_n$ retracts onto $S_2$ for all $n \geq 4$, and cannot be absolutely Koszul then by the previous lemma. This establishes Thm. \ref{abscounterexamples}. \qed
\end{con}
\begin{rem}
	Since $R_n$ is a commutative (graded) local $k$-algebra, we have that $R_n^! = {\rm Ext}^*_{R_n}(k, k) \cong \U \pi^*(R_n)$ is the universal enveloping algebra of a graded Lie algebra $\pi^*(R_n)$ supported in degrees $\geq 1$ called the \emph{Homotopy Lie algebra} of $R_n$ \cite[Chp 10]{AvramovBook}. In particular, ${\rm Ext}_{R_n}^*(k, k)$ is a Hopf algebra, and so there is anti-isomorphism $\sigma: {\rm Ext}_{R_n}^*(k, k) \xrightarrow{\cong} {\rm Ext}^*_{R_n}(k, k)^{op}$ given by the antipode map. This shows more clearly that coherence fails on both sides simultaneously. 
\end{rem}\newpage

%% file: chapters/applications/discussion.tex
\section{Discussion and conjectures}

\subsection*{The structure of complete resolutions} 
We record here a few thoughts and conjectures on the structure of complete resolutions over absolutely Koszul Gorenstein algebras, inspired by the equivalence of Theorem \ref{thmA}. \bigskip 

Let $A = \Bbbk \oplus A_1 \oplus \dots $ be a two-sided coherent, Koszul Gorenstein $\Bbbk$-algebra, and recall that $\dim A$ refers to the Gorenstein dimension of $A$. We define the \emph{width} of $M \in \grmod A$ and the \emph{global width} of $A$ by 
\begin{align*}
	&\mathsf{width}(M) = \sup\{\ |j_1 - j_2| \ | \ (M \otimes_A \Bbbk)_{j_k}  \neq 0,\ k = 1,2 \}.	\\ 
	&\mathsf{gl.width}(A) = \sup \{\ \mathsf{width}(M) \ | \ M \in {\rm MCM}^{\mathbb{Z}}(A) \textup{ indecomposable} \}. 
\end{align*} 
In plain terms, the width of $M$ is the largest difference between degrees of generators. When $A$ is absolutely Koszul, referring to $t^{\mathsf{lin}}$-cohomology recall that we define the \emph{amplitude} of $M \in \underline{{\rm MCM}}^{\mathbb{Z}}(A)$ (and analogously the \emph{global amplitude} of $A$) by
\begin{align*}
&\mathsf{amp}(M) = \sup\{\ |j_1 - j_2| \ | \ \mathcal{H}^{j_k}(M) \neq 0,\ k = 1,2 \}.	\\ 
&\mathsf{gl.amp}(A) = \sup \{\ \mathsf{amp}(M) \ | \ M \in \underline{{\rm MCM}}^{\mathbb{Z}}(A) \textup{ indecomposable} \}. 
\end{align*}
These notions are closely related to the regularity, since for any Koszul module $K$ with initial degree $\mathsf{indeg}(K) = 0$, we have
\[
	\mathsf{amp}(K) = \mathsf{width}(K) = \mathsf{reg}_A(K).
\] 
We have an immediate inequality 
\begin{align*}
	&\mathsf{gl.amp}(A) \leq \mathsf{gl.width}(A). 
\end{align*} 
Now, let us assume additionally that $A$ is two-sided absolutely Koszul. We say that $N \in \grmod A$ is coKoszul if $N^*$ is Koszul. In this situation, we have a notion of `middle part' of the complete resolution of MCM modules. 
\begin{defn} Let $A$ be two-sided coherent, two-sided absolutely Koszul Gorenstein. Let $M \in \underline{{\rm MCM}}^{\mathbb{Z}}(A)$, and let $C$ be a complete resolution of $M$. The \emph{middle part} of $C$ consists of the terms $\{ C_n \}$ for which
	\[
		{\rm coker}(C_{n-1} \to C_n) = \Omega^n M 
	\]
	is neither Koszul nor coKoszul. 
\end{defn}
Since $\Omega^n M$ is Koszul for any $n \gg 0$ and coKoszul for any $n \ll 0$ in the situation above, the middle part of $C$ consists of a finite number of terms in the complete resolution. The global width of $A$ attempts to measure the complexity of the middle parts of all indecomposable complete resolutions, while the global amplitude measures the complexity of the tails. \bigskip 

The complexity of the Betti tables of complete resolutions appears to be sensitive to the invariant $\nu = \dim A - 1 + a$. More precisely, there is some evidence for the following dichotomies. 
\begin{conj}[Dichotomy for width]\label{glwidthdichotomy} Let $A$ be a two-sided coherent, Koszul, Gorenstein $\Bbbk$-algebra, and let $\nu = \dim A -1 + a$ be the virtual dimension of $A$. Then: 
	\begin{enumerate}[1)]
	\item If $\nu \leq 1$, then $\mathsf{gl.width}(A) \leq N < \infty$ for some $N \in \N$. That is, there is a uniform bound $N$ such that every indecomposable $M \in \underline{{\rm MCM}}^{\mathbb{Z}}(A)$ in generated in at most $N$ degrees. 
	\item If $\nu \geq 2$, then $\mathsf{gl.width}(A) = \infty$. That is, for every $N \in \N$ there is an indecomposable $M \in \underline{{\rm MCM}}^{\mathbb{Z}}(A)$ such that $\mathsf{width}(M) > N$. 
\end{enumerate}
\end{conj}

\begin{conj}[Dichotomy for amplitude]\label{glampdichotomy} Let $A$ be a two-sided coherent, two-sided absolutely Koszul, Gorenstein $\Bbbk$-algebra, and let $\nu = \dim A -1 + a$ be the virtual dimension of $A$. Then: 
	\begin{enumerate}[1)]
	\item If $\nu \leq 1$, then $\mathsf{gl.amp}(A) \leq 1$. That is, every indecomposable $M \in \underline{{\rm MCM}}^{\mathbb{Z}}(A)$ is eventually $n$-linear for some $n \in \Z$. 
	\item If $\nu \geq 2$, then $\mathsf{gl.amp}(A) = \infty$. That is, for every $N \in \N$ there is an indecomposable $M \in \underline{{\rm MCM}}^{\mathbb{Z}}(A)$ such that $\mathsf{amp}(M) > N$. 
\end{enumerate}
\end{conj} 
We have shown some partial cases already, such as Conj. \ref{glampdichotomy} 1). The general case will be studied in subsequent work.

%% file: chapters/appendix/appendix.tex
\chapter{Appendix}
\input{chapters/appendix/artin}

\input{chapters/appendix/morita}
\input{chapters/appendix/orlov}

%% file: chapters/appendix/artin.tex
\section{Representation theory of quivers and finite dimensional algebras}\label{artinappendix}
This appendix is a quick summary of the relevant facts on the representation theory of algebras used implicitly in this thesis. The reader is referred to \cite{ARS, HappelBook, SimSkow} for details. 

\subsection*{Finite dimensional algebras and quivers}
Throughout this section $k$ will stand for an algebraically closed field, and $\Lambda$ will refer to a finite dimensional $k$-algebra. Since $\Lambda$ is finite dimensional and thus Artinian, $\Lambda /\rad\Lambda$ is semisimple and so isomorphic to a product of matrix algebras $M_{e_1}(k) \times \dots \times M_{e_r}(k)$ by the Artin-Wedderburn theorem. 
	\begin{defn}
		A finite-dimensional algebra $\Lambda$ is basic if $\Lambda / \rad \Lambda \cong k \times \dots \times k$. Equivalently, the decomposition of $\Lambda_\Lambda = \bigoplus_{i = 1}^r P(i)$ into indecomposable projectives $P(i)$ is multiplicity free, i.e. the $P(i)$ are pairwise non-isomorphic. 
	\end{defn}
	Every finite dimensional algebra is Morita equivalent to a basic algebra, and two Morita equivalent basic algebras are isomorphic. \bigskip

	A quiver $Q = (Q_0, Q_1)$ is a directed graph, with vertex set $Q_0$ and arrow set $Q_1$. Given an arrow $a \in Q_1$, we let $s(a), t(a)$ denote its source and target vertices. 
	{\tiny{
	\begin{align*} 
	\xymatrix{ &\\
	\bullet \ar@/^/[r] & \bullet \ar@/^/[l] }
		&& \xymatrix{ && \bullet \\
			\bullet \ar@<1ex>[r] \ar@<-1ex>[r] & \bullet \ar[r] \ar[dr] \ar[ur] & \bullet \ar[d] \\
	&& \bullet} \\ 
\end{align*}
\begin{align*}
	\xymatrix@R=15pt@C=15pt{ 
		& \bullet \ar[r] & \bullet \ar[dr] \\ 
		\bullet \ar[ur]	& && \bullet \ar[d] \\ 
		\bullet \ar@{.}[u] &&& \bullet \ar[dl] \\
		& \bullet \ar[ul] & \ar[l] \bullet
	}
	\end{align*} 
}}
Let $kQ_d$ be the $k$-span of $Q_d$ for $d = 0, 1$, and assume for simplicity that both sets are finite. Then $kQ_0 = \prod_{i \in Q_0} k$ is a semisimple $k$-algebra with standard basis of idempotents $\{e_i \}$ and $kQ_1$ is a $kQ_0$-bimodule for which $e_j kQ_1 e_i$ is the $k$-span of arrows $\{ a: i \to j \}$. We define the path algebra $$kQ = T_{kQ_0}(kQ_1) = \bigoplus_{d \geq 0} kQ_d $$ as the tensor algebra of $kQ_1$ over $kQ_0$, with $kQ_d = (kQ_1)^{\otimes d}$ the $k$-space of paths of length $d$ for all $d \in \N$. Equivalently, $kQ$ is the $k$-span of the graded set $\{Q_d\}$ of paths of length $d \in \N$, with multiplication given by ``function composition'' for $a, b \in Q_1$ whenever possible:
	\[
		m_2(a,b) = \begin{cases} ab & s(a) = t(b) \\
		0 & \textup{ else. }\end{cases} 
	\]
	Note that $kQ$ is an augmented $kQ_0$-algebra, and is the free $kQ_0$-algebra on the bimodule $kQ_1$. The Jacobson radical is given by the arrow ideal $\rad kQ = (kQ_1) = kQ_{\geq 1}$. \bigskip

	Given $\Lambda$ basic, write $1 = e_1 + \dots + e_r$ for a full set of primitive idempotents $\{e_i\}$, so that $P(i) = e_i \Lambda$. We construct the \emph{ordinary quiver} $Q = Q_\Lambda$ as follows. Take $Q_0 = \{e_i \}$. Next, let $Q_1$ be a set of arrows $\{a: i \to j\}$ in bijection with a $k$-basis $e_j (\rad \Lambda/ \rad^2\Lambda) e_i$.

	\begin{prop}[Gabriel's theorem]Let $\Lambda$ be basic with ordinary quiver $Q = Q_\Lambda$. Then there is a surjective homomorphism $kQ \twoheadrightarrow \Lambda$ whose kernel $I$ consists of decomposable elements, meaning $I \subseteq kQ_{\geq 2} = \rad^2kQ$, so that $kQ/I \cong \Lambda$. The quiver $Q$ is uniquely determined by this property so long as $I \subseteq \rad^2 kQ$.
	\end{prop}
	A generating set for $I = (\rho_1, \dots, \rho_c)$ can always be taken as linear combinations of paths with same source and target. When $I$ is in this form, we call the pair $(Q, I)$ a bound quiver. Note that $kQ/I$ is finite dimensional if and only if $kQ_{\geq d} \subseteq I$ for $d \gg 0$. 

	A representation of a bound quiver $(Q, I)$ consists of a set of vector spaces $\{V_i\}_{i \in Q_0}$ and operators $\varphi_a: V_i \to V_j$ for every $a: i \to j$, such that $\{ \varphi_a \}_{a \in Q_1}$ satisfy the relations $I$. Morphisms of representations are simply natural transformations. We denote the category of representations of $(Q, I)$ by $\mathsf{Rep}(Q, I)$ (resp. $\mathsf{rep}(Q, I)$ for the subcategory of representations with finite dimensional $V_i$).

	Given a quiver representation $\left( \{V_i\}, \{ \varphi_a \} \right)$, the path algebra $kQ/I$ operates on the left of $V = \bigoplus_{i \in Q_0} V_i$ via the $\varphi_a$, giving a functor to left $kQ/I$-modules. 

	\begin{prop} This functor gives rise to an equivalence of categories
		\[
		\mathsf{Rep}(Q, I) \cong \Modsf (kQ/I^{op}).
		\]
		The inverse sends $V \mapsto \{V_i\}$ with $V_i = e_iV$. When $kQ/I$ is finite dimensional, this restricts to
		\[
			\mathsf{rep}(Q, I) \cong \modsf (kQ/I^{op}).
		\]
	\end{prop}

	\subsection*{Trichotomy} Fix a finite dimensional algebra $\Lambda$ and restrict attention to finitely generated modules. For an algebra $A$, we say that a $k$-linear functor
	\[
F: \modsf A \to \modsf \Lambda
	\]
	is a \emph{representation embedding} if it sends indecomposables to indecomposables and reflects isomorphisms, that is $F(M) \cong F(N)$ implies $M \cong N$. Given $\Lambda$, the classification problem for indecomposable $\Lambda$-modules can be a various complexity. We say that the representation theory of $\Lambda$ is of
	\begin{itemize}
	\item \emph{finite} type, if there are finitely many indecomposable $\Lambda$-modules;
	\item \emph{tame} type, if there are infinitely many indecomposables and, for any dimension $d$, there are finitely many $k[T]-\Lambda$-bimodules $\F_1, \dots, \F_{\mu(d)}$ such that $\F_i$ are free of finite rank as $k[T]$-modules and all but finitely many indecomposable $\Lambda$-modules of dimension $d$ are of the form $M_{i, \lambda} = \F_i/(T-\lambda)\F_i$, for some $i = 1, \dots, \mu(d)$ and $\lambda \in k$;
	\item \emph{wild} type, if for any finitely generated $k$-algebra $A$ there is a representation embedding $F: \modsf A \to \modsf \Lambda$.
\end{itemize}

By Drozd's Theorem, the representation type of a finite dimensional algebra $\Lambda$ must be either finite, tame or wild. 

\subsection*{Grothendieck groups and Euler forms}  
Let $\Lambda = kQ/I$ be a finite dimensional algebra, presented as in Gabriel's theorem. Given a finite dimensional quiver representation $V = \{V_i\}$, we define the dimension vector $\underline{\dim}\ \! V := (\dim V_i)_{i \in Q_0}$. Letting $K_0(\Lambda) = K_0(\modsf \Lambda)$, we have an isomorphism of abelian groups
\[
	\underline{\dim}: K_0(\Lambda) \xrightarrow{\cong} \Z^{|Q_0|}.
\]
Letting $S(i)$ be the simple top of the indecomposable projective $P(i)$, $S(i)$ can be represented by the unique quiver representation with $\dim V_j = \delta_{ij}$, and so $[S(i)]$ forms a finite $\Z$-basis of $K_0(\Lambda)$. The dimension vector $ \underline{\dim}\ \! V$ gives the Jordan-H\"older multiplicity of each $S(i)$ in $V$. When $\gldim \Lambda < \infty$, the Euler form 
\[
	\langle X, Y \rangle = \sum_{i \in \Z} \dim {\rm Ext}^i_{\Lambda}(X, Y)
\]
is defined and descends to a pairing $\langle -, - \rangle: K_0(\Lambda) \otimes K_0(\Lambda) \to \Z$. One can show that this pairing is perfect.  

\subsection*{Hereditary algebras} An algebra $\Lambda$ is hereditary if $\gldim \Lambda = 1$. The finite dimensional hereditary algebras have historically been the best studied finite dimensional algebras, and their module category exhibits a rich structure. 

\begin{prop} A finite dimensional basic hereditary algebra $\Lambda$ is isomorphic to $kQ$ for some acyclic quiver $Q$, and any such path algebra is hereditary.
\end{prop}

\begin{prop}
	Let $Q$ be an acyclic quiver. The representation type of $kQ$ is determined as follows:
	\begin{enumerate}[i)]
	\item $kQ$ is of finite type if and only if the underlying graph of $Q$ is a simply-laced (ADE) Dynkin diagram.
	\item $kQ$ is of tame type if and only if the underlying graph of $Q$ is of affine ADE type. 
	\item $kQ$ is wild otherwise. 
\end{enumerate}
\end{prop}

Define a bilinear form $\langle -, - \rangle_Q$ on $\Z^{|Q_0|}$ by 
\[
	\langle d, d' \rangle_Q = \sum_{i \in Q_0} d_id'_i - \sum_{a: i \to j} d_id'_j.
\]

\begin{prop}
	The isomorphism $ \underline{\dim}: K_0(kQ) \xrightarrow{\cong} \Z^{|Q_0|}$ sends the Euler form to the above. That is
	$$  \langle X, Y \rangle = \langle \underline{\dim}X, \underline{\dim}Y \rangle_Q.$$
\end{prop}

The underlying graph of $Q$ determines a symmetrizable Cartan matrix, to which we attach a Kac-Moody algebra $\mathfrak{g} = \mathfrak{g}_Q$. Recall that $\mathfrak{g}$ is graded by a finite free abelian group $\Gamma$ called the root lattice, with $\Z$-basis of simple roots $\{\epsilon_i \}$. We let $\Delta = \{ \alpha \in \Gamma \ | \ \mathfrak{g}_\alpha \neq 0 \}$ be the sets of roots of $\mathfrak{g}$. Now denote by $(-, -)_Q$ the symmetrization of the above pairing, meaning
\[
	(d, d')_Q = \langle d, d' \rangle_Q + \langle d', d \rangle_Q.
\]
We can identify the lattice $\Z^{|Q_0|} = \Gamma$ with the root lattice of $\mathfrak{g}$, sending the dimension vector of the simple representation $S(i)$ to the simple root $\epsilon_i$. This identifies $(-, -)_Q$ with the Weyl-invariant symmetric bilinear form on $\Gamma$. 

\begin{prop}[Kac's Theorem] The above identification induces a bijection between the dimension vectors of indecomposable representations of $Q$ and the positive roots $\Delta^{+}$ of $\mathfrak{g}$. Furthermore:
	\begin{enumerate}[i)]
	\item There is a unique indecomposable $X$ with $ \underline{\dim}X \mapsto \alpha$ for each positive real root $\alpha$.
	\item There are infinitely many indecomposables $Y$ with $ \underline{\dim} Y \mapsto \beta$ for each positive imaginary root $\beta$.
\end{enumerate}
\end{prop}

\subsection*{Auslander-Reiten theory} Let $\Lambda$ be a finite dimensional algebra, and let 
\[
\xi: 0 \to X \to Y \to Z \to 0
\]
be a short exact sequence in $\modsf \Lambda$.
\begin{defn}
	A short exact sequence $\xi$ is almost-split if $0 \neq \xi \in {\rm Ext}^1_{\Lambda}(Z, X)$, $X, Z$ are indecomposable and for every non-isomorphism $t: W \to Z$ with $W$ indecomposable factors as 
	\[
		\xymatrix@R1pc{ &&& W \ar[d]^t \ar@{..>}[dl] \\
			0 \ar[r]& X \ar[r]& Y \ar[r]& Z \ar[r]& 0.
}
	\]
\end{defn}

Almost-split sequences were introduced by Auslander-Reiten and their existence is fundamental to the representation theory of finite dimensional algebras. Every indecomposable object sits inside an almost-split sequence. Let $X$ be indecomposable, and pick a minimal projective presentation
\[
	P_1 \to P_0 \to X \to 0.
\]
Writing $(-)^{\vee} = {\rm Hom}_{\Lambda}(-, \Lambda)$, we define the transpose ${\rm Tr}(X) \in \modsf \Lambda^{op}$ by the dual presentation
\[
	P_0^{\vee} \to P_1^{\vee} \to {\rm Tr}(X) \to 0.
\]
Define the functors 
\[
\tau = {\rm DTr}(-): \modsf \Lambda \to \modsf \Lambda
\] 
\[
	\tau^{-1} = {\rm TrD}(-): \modsf \Lambda \to \modsf \Lambda.
\] 
Define the projectively stable (resp. injectively stable) categories $\underline{\modsf}\Lambda$ (resp. $ \overline{\modsf}\Lambda $) with homomorphism space
\[
	\underline{ {\rm Hom}}(X, Y) = {\rm Hom}_{\Lambda}(X, Y)/\mathcal{P}(X, Y)
\]
and
\[
	\overline{ {\rm Hom}}(X, Y) = {\rm Hom}_{\Lambda}(X, Y)/\mathcal{I}(X, Y)
\]
where $\mathcal{P}(X, Y)$ (resp. $\mathcal{I}(X, Y)$) is the ideal of morphisms factoring through a projective object (resp. injective object).

\begin{prop} The functors ${\rm Tr}, \tau, \tau^{-1}$ descend to equivalences of stable categories
	\begin{align*}
	{\rm Tr}: &\underline{\modsf}\Lambda \leftrightarrows \underline{\modsf}\Lambda^{op}: {\rm Tr} \\
	\tau: &\overline{\modsf}\Lambda \leftrightarrows \underline{\modsf}\Lambda^{op}: \tau^{-1}
\end{align*}
 \end{prop}

 \begin{thm}[Auslander-Reiten] Let $\Lambda$ be a finite dimensional $k$-algebra. Then every non-injective indecomposable $X$ sits in a unique almost-split sequence
	 \[
		0 \to X \to Y \to \tau^{-1} X \to 0
	 \]
	and every non-projective indecomposable $Z$ sits in a unique almost-split sequence
	\[
		0 \to \tau Z \to Y \to Z \to 0.
	\]
 \end{thm}

 \begin{defn} Let $P(i) = e_i\Lambda$ and $I(i) = {\rm D}(\Lambda e_i) = e_i {\rm D}(\Lambda)$ be the indecomposable projective and injective $\Lambda$ modules, as $\{e_i\}$ runs through a full set of primitive idempotents. 
	 \begin{enumerate}[i)]
	 \item Indecomposables of the form $\tau^{-n}P(i)$ for $n \geq 0$ are called preprojective.
	 \item Indecomposables of the form $\tau^{n}I(i)$ for $n \geq 0$ are called preinjective. 
 \end{enumerate}
 \end{defn}

 We can encode the structure of almost-split sequences in a combinatorial structure called the Auslander-Reiten (AR) quiver $\Gamma = \Gamma \modsf \Lambda$. The quiver $\Gamma = (\Gamma_0, \Gamma_1)$ is defined by
 \begin{enumerate}[a)]
 \item $\Gamma_0 = \{ X \textup{ indecomposable } \}/\cong$.
 \item $\Gamma_1$ has a set of arrows $\{ a: [X] \to [Y] \}$ in bijection with a $k$-basis of the space of irreducible maps ${\rm Irr}(X, Y)$.
	 \end{enumerate}
	 See \cite{HappelBook} for the definition of irreducible maps. 

	 \begin{prop}The AR quiver is related to almost-split sequences as follows. Let
		 \[
			 \xi: 0 \to X \xrightarrow{f} Y \xrightarrow{g} Z \to 0
		 \]
		 be an almost-split sequence, and decompose $Y = \bigoplus_{i = 1}^r Y_i^{\oplus e_i}$ into indecomposables. 
		 \begin{enumerate}[i)]
		 \item The classes of the components $f_i: X \to Y_i$ give a full set of arrows coming out of $X$ in $\Gamma$, and we have $\dim {\rm Irr}(X, Y_i) = e_i$.
		 \item The classes of the components $g_i: Y_i \to Z$ give a full set of arrows coming into $Z$ in $\Gamma$, and we have $\dim {\rm Irr}(Y_i, Z) = e_i$.
	 \end{enumerate}
	 \end{prop}
	 We can picture the above in $\Gamma$ as
\[
	\xymatrix@R.3pc{	& [Y_1] \ar[ddr]^{e_1}& \\ 
		& [Y_2] \ar[dr]_{e_2} & \\
				[X] \ar[uur]^{e_1} \ar[ur]_{e_2} \ar[dr]^{e_{r-1}} \ar[ddr]_{e_r} & \vdots & [Z]\\
				& [Y_{r-1}] \ar[ur]^{e_{r-1}} & \\  
			& [Y_r] \ar[uur]_{e_r} & }
\]
writing $e_i$ for the multiplicity of arrows. \bigskip

Now let $\Lambda$ be hereditary, say basic so that $\Lambda = kQ$ for some $Q$ acyclic. Denote by $\mathcal{P}$ and $\mathcal{I}$ the subcategories of preprojective and preinjective modules. We say that an indecomposable is regular if it is neither in $\mathcal{P}$ nor $\mathcal{I}$, and denote by $\mathcal{R}$ the subcategory of regular indecomposables. The categories $\mathcal{P}, \mathcal{I}$ each form a component in the AR quiver $\Gamma$, called the preprojective and preinjective components. The category $\mathcal{R}$ breaks down into connected components, each closed under $\tau^{\pm}$, called regular components. 

\begin{prop} The following are equivalent:
	\begin{enumerate}[i)]
	\item $kQ$ is of finite representation type (i.e. $Q$ is Dynkin).
	\item $\mathcal{R} = \emptyset$.
	\item $\mathcal{P} = \mathcal{I}$.
\end{enumerate}
\end{prop}

When $kQ$ is representation infinite, the quiver $\Gamma$ is directed in that the nonzero morphisms in $\modsf kQ$ go from left to right in $(\mathcal{P}, \mathcal{R}, \mathcal{I})$.
\begin{prop} When $kQ$ is representation infinite, we have:
	\begin{enumerate}[i)]
	\item ${\rm Hom}_{\Lambda}(\mathcal{R}, \mathcal{P}) = 0$.
	\item ${\rm Hom}_{\Lambda}(\mathcal{I}, \mathcal{R}) = 0$.
	\item ${\rm Hom}_{\Lambda}(\mathcal{I}, \mathcal{P}) = 0$.
\end{enumerate}
\end{prop} 

\begin{prop} Let $kQ$ be representation infinite. The following are equivalent:
	\begin{enumerate}[i)]
	\item $kQ$ is of tame representation type (i.e $Q$ is affine Dynkin).
	\item The functor $\tau$ has finite order on $\mathcal{R}$.
\end{enumerate}
\end{prop}
The regular components in the tame case are tubes $\{ \mathcal{T}_{\lambda} \}_{\lambda \in \mathbb{P}^1}$ parameterized by $\mathbb{P}^1$. The category of modules whose indecomposable summands belong to $\mathcal{T}_\lambda$ is a serial category, with finitely many simple objects forming the base of the tube. We call the number of simples the \emph{rank} of the tube. All but finitely many tubes have rank one.
\begin{prop}
	We have ${\rm rk}\ \! K_0(\Lambda) = 2 + \sum_{\lambda} ({\rm rk}(\mathcal{T}_\lambda) - 1)$. 
\end{prop}

\subsection*{Derived categories of finite-dimensional algebras}

We let $ {\rm D}^b(\Lambda) := {\rm D}^b(\modsf \Lambda) \cong {\rm D}_{\modsf \Lambda}^b(\Modsf \Lambda)$, and ${\rm D}^{\perf}(\Lambda) \subseteq {\rm D}^b(\Lambda)$ for the subcategory of perfect complexes, meaning those quasi-isomorphic to bounded complexes of projectives. \bigskip

Given a finite dimensional algebra $\Lambda$, we are primarily interested in the bounded derived category $ {\rm D}^b(\Lambda)$ and the perfect derived category ${\rm D}^{\perf}(\Lambda)$. We recall some of their properties. 
	\begin{prop}[Happel] The following hold. 
		\begin{enumerate}[i)]
		\item $ {\rm D}^b(\Lambda) $ is an Hom-finite Krull-Schmidt category. 
		\item $ {\rm D}^b(\Lambda) $ admits a Serre functor if and only if $\gldim \Lambda < \infty$. In this case it is given by$$ \Se_\Lambda(-) = - \otimes^{\mathbb{L}}_{\Lambda} D\Lambda$$ with inverse given by the right adjoint $$\Se^{-1}_{\Lambda}(-) = - \otimes_{\Lambda}^{\mathbb{L}} {\rm RHom}_{\Lambda}(D\Lambda, \Lambda) \simeq {\rm RHom}(D\Lambda, -).$$
		\item ${\rm D}^{\perf}(\Lambda)$ admits a Serre functor if and only if $\Lambda$ is Gorenstein, whence it is also given by the above.  
	\end{enumerate}
	\end{prop}
	By \cite{ReitVdBalmostsplit}, existence of a Serre functor $\mathbb{S}$ for the triangulated category $\mathcal{T} = {\rm D}^b(\Lambda)$ is equivalent to the existence of almost-split triangles, with Auslander-Reiten translate $\tau = \mathbb{S} \circ [-1]$. Since the module category $\modsf \Lambda$ admits almost-split short exact ssequence and a translate $\tau = {\rm DTr}$, one may ask how whether these are compatible. Let us temporarily write $\tau_{\Delta} = \mathbb{S} \circ [-1]$ for the translate in $ {\rm D}^b(\Lambda)$ and $\tau$ for the translate in $\modsf \Lambda$. 

	\begin{prop} Let $\Lambda$ be hereditary. Consider a short exact sequence in $\modsf \Lambda$ 
		\[
\xi: 0 \to X \to Y \to Z \to 0
		\]
		with associated triangle in $ {\rm D}^b(\Lambda) $
		\[
			\xi^{\Delta}:	X \to Y \to Z \to X[1].
		\]
		Then $\xi$ is almost-split in $\modsf \Lambda$ if and only if $\xi^{\Delta}$ is almost-split in $ {\rm D}^b(\Lambda) $. Furthermore, we have:
		\begin{enumerate}[i)]
		\item $\tau_{\Delta}Z = \tau Z = {\rm Tor}^{\Lambda}_1(Z, D\Lambda)$ for each indecomposable non-projective module $Z$.
		\item $\tau_{\Delta}^{-1}X = \tau^{-1}X = {\rm Ext}^1_{\Lambda}(D\Lambda, X)$ for each indecomposable non-injective module $X$.
		\item $\tau_{\Delta} P(i) = I(i)[-1]$.
	\end{enumerate}
	\end{prop}

	\begin{rem}
		The reader is warned that this result is special to the hereditary case; the translates $\tau$ and $\tau_\Delta$ usually differ in nature for $\gldim \Lambda \geq 2$. We will primarily work with Serre functors on triangulated categories, and so $\tau$ will always refer to $\tau_{\Delta}$ if any ambiguity arises. 
\end{rem}

Happel has worked out the structure of the Auslander-Reiten quiver of ${\rm D}^b(kQ)$. First, a standard fact.

	\begin{prop}\label{hereditaryformal}
		Let $H$ be an abelian category with ${\rm Ext}_{\mathcal{H}}^2(-, -) = 0$. Then all objects in $ {\rm D}^b(H) $ are formal. That is, for each $X$ we have an isomorphism $X \cong \bigoplus_{n \in \Z} {\rm H}^n(X)[-n]$ in ${\rm D}^b(H)$, and so each object is the shifted sum of its cohomology objects. 
	\end{prop}

	\begin{thm}[Happel \cite{HappelBook}] The AR components of ${\rm D}^b(kQ)$ look as follows.
		\begin{enumerate}[i)]
		\item The preprojective component $\mathcal{P}$ and shifted preinjective components $\mathcal{I}[-1]$ form one connected component in ${\rm D}^b(kQ)$, called the \emph{transjective} component $\mathcal{PI}$. We have 
	\[
		\mathcal{PI} = \{ \tau^{n} P(i) \ | \ P(i) \textup{ indecomposable projective}\ n \in \Z \}.
	\]
\item The connected components are either shifts of regular components in $\modsf kQ$ or shifts of $\mathcal{PI}$.
	\end{enumerate}
	\end{thm}

%% file: chapters/appendix/morita.tex
\section{Derived Morita theory and Tilting theory}\label{moritaappendix}

We briefly review and give references for the general Morita theorem of Keller in the context of algebraic triangulated categories, based on the articles \cite{KellerDG, KellerTilting}. The reader is referred to these articles for more details. \bigskip 

Let $k$ be a commutative ring. A $k$-linear triangulated category $\T$ is called \emph{algebraic} if it arises as the stable category $ \underline{\E}$ of a Frobenius category $\E$, see \cite[Sect. 3.6]{KellerDG} for definitions. We note that all triangulated categories appearing in this thesis are algebraic. \bigskip 

Closely related is the notion of dg category. 
\begin{defn}
	We say that an additive $k$-linear category $\CC$ is a differential graded (dg) category if $\CC$ has the following properties: 
	\begin{enumerate}[1)]
	\item The Hom objects ${\rm Hom}_\CC(X, Y)$ are complexes for any $X, Y$.
	\item The identity $1_X$ is closed for each $X$, i.e. $d(1_X) = 0$.
	\item Composition of morphisms
		\[
			{\rm Hom}_\CC(Y, Z) \otimes_k {\rm Hom}_\CC(X, Y) \to {\rm Hom}_\CC(X, Z)
		\]
		is a chain-map of complexes. 
\end{enumerate}
\end{defn}
To any dg category $\CC$ we can attach a $k$-linear category ${\rm H}^0 \CC$ with the same objects as $\CC$ and morphism space ${\rm Hom}_{ {\rm H}^0\CC}(X, Y) = {\rm H}^0 {\rm Hom}_\CC(X, Y)$. A dg functor between dg categories $F: \CC \to \mathcal{D}$ is an additive functor such that the induced linear map $F_{X, Y}: {\rm Hom}_{\CC}(X, Y) \to {\rm Hom}_{\mathcal{D}}(F(X), F(Y))$ is a morphism of complexes; such functors descend to functors ${\rm H}^0(F): {\rm H}^0(\CC) \to {\rm H}^0(\mathcal{D})$. A quasi-equivalence $F: \CC \xrightarrow{\sim} \mathcal{D}$ is a dg functor $F$ such that ${\rm H}^0(F)$ is an equivalence; equivalently, the morphisms $F_{X,Y}$ on Hom complexes are quasi-isomorphisms. \bigskip

Given any algebraic triangulated category $\T$, one can find a dg category $\CC$ along with an equivalence $\T \cong {\rm H}^0(\CC)$. We will set the notation 
\[
	{\rm RHom}_\T(X, Y) := {\rm Hom}_\CC(X, Y).
\]
Note that ${\rm RHom}_\T(X, X)$ is naturally a dg algebra over $k$ for any object $X \in \T$. \bigskip

Define an object $X \in \T$ to be compact if ${\rm Hom}_\T(X, -)$ commutes with arbitrary direct sums in $\T$. Let $\T^c$ be the subcategory of compact objects. Given a set of objects $S \subseteq \T$, we define $\thick(S) \subseteq \T$ (respectively $\mathsf{loc}(S) \subseteq \T$) to be the smallest triangulated subcategory of $\T$ closed under finite sums and summands (respectively arbitrary sums; direct summands will follow, see \cite{NeemanBook}). Next, we say that a set of objects $S \subseteq \T$ is: 
\begin{enumerate}[i)]
\item A set of \emph{classical generators}, if $\thick(S) = \T$.
\item A set of \emph{compact generators}, if $S \subseteq \T^c$ consists of compact objects and $\mathsf{loc}(S) = \T$. 
	\end{enumerate}
	We can now state Keller's derived Morita theorem. We refer to \cite{KellerDG} for the module and derived category of a small DG category. 
	\begin{thm}[Keller {\cite[Thm. 3.8]{KellerDG}}]\label{derivedmorita} Let $\T$ be an algebraic triangulated category. Assume that $\T$ is idempotent closed. Let $\CC$ be the associated dg category such that ${\rm H}^0(\CC) \cong \T$. Let $S \subseteq \T$ be a small full subcategory and $\mathcal{S}$ the DG category obtained from $S$, with same objects as $S$ but with morphisms given by ${\rm Hom}_\CC(s, s')$ for any $s, s' \in S$. The following holds:
		\begin{enumerate}[i)] 
		\item If $S \subseteq \T$ is a set of classical generators, then we have an exact equivalence of categories
			\[
				{\rm RHom}(S, -): \T \xrightarrow{\cong} {\rm D}^{\perf}(\Modsf \mathcal{S}).
			\]
		\item If $\T$ is closed under arbitrary direct sums and $S \subseteq \T$ is a set of compact generators, then we have an exact equivalence of categories
			\[
				{\rm RHom}(S, -): \T \xrightarrow{\cong} {\rm D}(\Modsf \mathcal{S}).
			\]\qed
	\end{enumerate}
	\end{thm}
\noindent This theorem further specialises. We say that $S \subseteq \T$ is tilting if 
	\begin{enumerate}[i)]
	\item $S$ is a set of classical generators;
	\item ${\rm Hom}_\T(s, s'[n]) = 0$ for $n \neq 0$ and all $s, s' \in S$.
\end{enumerate}
We note the special case when $S$ consists of a single object $S = \{ T \}$. The object $T$ is then called a tilting object. \bigskip

When $S$ is a tilting subcategory with $\mathcal{S}$ as in the above theorem, letting $\tau^{\leq 0}$ be the left truncation functor for standard t-structure on $ {\rm D}^b(k) $, we have quasi-isomorphisms of Hom complexes 
\[
	{\rm Hom}_\CC(s, s') \xleftarrow{\sim} \tau^{\leq 0}{\rm Hom}_\CC(s, s') \xrightarrow{\sim} {\rm Hom}_\T(s, s')
\]
for any $s, s' \in S$. Writing $\tau^{\leq 0}\mathcal{S}$ for the DG category with morphism complexes as above, we obtain quasi-equivalences 
\[
	\mathcal{S} \xleftarrow{\sim} \tau^{\leq 0}\mathcal{S} \xrightarrow{\sim} {\rm H}^0\mathcal{S} = S.
\]
Lastly, we note that quasi-equivalent dg categories $\mathcal{S} \simeq {\rm H}^0(\mathcal{S}) = S$ have equivalent derived categories ${\rm D}(\Modsf \mathcal{S}) \cong {\rm D}(\Modsf S)$ (and perfect derived categories $({\rm D}^{\perf}(\Modsf \mathcal{S}) \cong {\rm D}(\Modsf S)$). Keller's Tilting theorem then specialises to

	\begin{cor}[Keller] Assume that $S \subseteq \T$ is a tilting subcategory. Then there is an equivalence of triangulated categories
		\[
			\T \xrightarrow{\cong} {\rm D}^{\perf}(\Modsf S). 
		\]
		As a special case, when $S = \{T \}$ consists of a single tilting object $T$, we have an equivalence of triangulated categories
		\[
			\T \xrightarrow{\cong} {\rm D}^{\perf}(\Modsf {\rm End}_\T(T)). 
		\]
	\end{cor}
	This last statement is \cite[Thm. 8.7]{KellerTilting}. 

%% file: chapters/appendix/orlov.tex
\section{Semiorthogonal decompositions and Orlov's Theorem}\label{orlovappendix}

In this section we review Orlov's semiorthogonal decomposition theorem and standard background notions leading up to it. Everything in here is in \cite{Or1}, except for the description of the crucial left adjoint to MCM approximation, which is due to Buchweitz.\bigskip

Let $\T$ be a $k$-linear triangulated category over some field $k$, and $\A \subseteq \T$ a full triangulated subcategory. 
\begin{defn} We say that $\A$ is right admissible (resp. left admissible) if the embedding $\iota: \A \hookrightarrow \T$ has a right adjoint $Q: \T \to \A$ (resp. left adjoint). We say that $\A$ is admissible when it is both left and right admissible. 
\end{defn}

Define the right orthogonal category $\A^{\perp} = \{ X \in \T \ | \ {\rm Hom}_\T(A, X) = 0 \textup{ for all } A \in \A \}$, and similarly define the left orthogonal category $\ \!^{\perp}\A$. By making use of the counit map of the adjunction and the long exact sequence of Hom spaces, one sees that $\A \subseteq \T$ is right admissible if and only if for all $X \in \T$, there is a distinguished triangle
\[
	X_\A \to X \to X_{\A^{\perp}} \to X_\A[1]
\]
with $X_\A \in \A$ and $X_{\A^{\perp}} \in \A^{\perp}$. Analogously, $\A \subseteq \T$ is left admissible if and only if for all $X \in \T$, there is a distinguished triangle
\[
X_{\ \!^{\perp}\A} \to X \to X_\A \to X_{\ \!^{\perp}\A}[1]
\]
with $X_{\ \!^{\perp}\A} \in \ \!^{\perp}\A$ and $X_\A$ in $\A$.

\begin{defn} Let $\A, \B \subseteq \T$ be full triangulated subcategories. We say that $\T$ has a (weak) semiorthogonal decomposition $\T = \langle \A, \B \rangle$ if ${\rm Hom}_\T(\B, \A) = 0$ and for all $X \in \T$, there is a distinguished triangle
	\[
		X_\B \to X \to X_\A \to X_\B[1].
	\]
	with $X_\A \in \A$ and $X_\B$ in $\B$. Equivalently, $\A$ is left admissible (and then $\B = \ \!^{\perp}\A$). Again equivalently, $\B$ is right admissible (and then $\A = \B^{\perp}$). We say that $\T = \langle \A, \B \rangle$ is a semiorthogonal decomposition if $\A$ (equivalently, $\B$) are admissible. 
\end{defn}

Given a weak semiorthogonal decomposition, we can refine it by decomposing $\A$ or $\B$ further.

\begin{defn} A sequence of full triangulated subcategories $(\A_1, \dots, \A_n)$ in $\T$ is a weak semiorthogonal decomposition $\T = \langle \A_1, \dots, \A_n \rangle$ if there is a sequence of left admissible subcategories $\T_1 = \A_1 \subset \T_2 \subset \dots \subset \T_n = \T$ such that $\A_{k}$ is the left orthogonal of $A_{k-1}$ in $T_k$, for all $k = 2, \dots, n$. If all $A_k$ are admissible, then $\T = \langle \A_1, \dots, \A_n \rangle$ is a semiorthogonal decomposition.
\end{defn}

Semiorthogonal decompositions occur often in algebraic geometry and representation theory. At one extreme, the derived category $ {\rm D}^b(X) $ of a smooth projective Calabi-Yau variety $X$ does not admit any semiorthogonal decomposition. At the other extreme, some varieties (e.g. $X = \mathbb{P}^n$) admit semiorthogonal decompositions $ {\rm D}^b(X) = \langle \A_1, \dots, \A_n \rangle$ with $\A_i \cong {\rm D}^b(pt)$, so that the pieces $\A_i$ are as simple as possible. 

\begin{defn} An object $E \in \T$ is exceptional if 
	\[
		{\rm Hom}_\T(E, E[n]) = \begin{cases} k & n = 0 \\ 0 & n \neq 0.\end{cases}
	\]
	A sequence $\sigma = (E_1, \dots, E_n)$ is an exceptional sequence if the $E_i$ are exceptional objects, and ${\rm Hom}_\T(E_i, E_j[n]) = 0$ for all $n \in \Z$ whenever $j > i$. The sequence $\sigma$ is \emph{strong} if furthermore ${\rm Hom}_\T(E_i, E_j[n]) = 0$ for all $n \neq 0$ and any $i,j$, and is \emph{full} if $\thick(E_1, \dots, E_n) = \T$.
\end{defn}
Given $\sigma = (E_1, \dots, E_n)$ a full strong exceptional sequence, the sum $T = \bigoplus_{i=1}^n E_i$ is a tilting object for $\T$ in the sense of Appendix \ref{moritaappendix}. 

\begin{exmp} Let $A = kQ/I$ be a finite dimensional algebra with $I \subseteq \mathsf{rad}^2 kQ$ with $Q$ a finite acyclic quiver. Order $Q_0 = \{1, \dots, n\}$ in which the arrow directions $a:i \to j$ are increasing, and let $e_i$ be the idempotent path at the $i$-th vertex. Then the indecomposable right projective modules $P(i) = e_iA$ form a full strong exceptional collection
	\[
		{\rm D}^b(A) = \langle P(1), \dots, P(n) \rangle. 
	\]
\end{exmp}

\subsection*{Mutations of exceptional collections}
Assume that $\T$ is Ext-finite. Given objects $A, B, C$ in $\T$, denote by ${\rm Hom}^{\bullet}(A, B) = \bigoplus_{n \in \Z} {\rm Hom}_\T(A, B[n])[-n]$ the object in ${\rm D}^b(k)$, and define
\[
	{\rm Hom}^{\bullet}(A, B)\otimes_k C = \bigoplus_{n \in \Z}{\rm Hom}_\T(A, B[n])\otimes_k C[-n].
\]
Interpret the dual ${\rm Hom}^{\bullet}(A, B)^*$ in ${\rm D}^b(k)$ accordingly. Following Gorodentsev, we have canonical distinguished triangles in $\T$
\[
	L_A(B)[-1] \to	{\rm Hom}^{\bullet}(A, B) \otimes_k A \xrightarrow{ev} B \to L_A(B)
\]
\[
	R_A(B) \to A \xrightarrow{coev} {\rm Hom}^{\bullet}(A, B)^* \otimes_k B \to R_A(B)[1]
\]
which uniquely define $L_A(B)$, $R_A(B)$. Now, an exceptional collection $(E,F)$ of length two is called an exceptional pair. We have the standard result.
\begin{prop}[Gorodentsev {\cite{Rud}}]\label{braidaction} Let $(E, F)$ be an exceptional pair in $\T$. The operations $L, R$ descend to an action on the set of exceptional pairs 
	\begin{align*}
		R\!:\ & (E, F) \mapsto (F, R_E(F))\\
		L\!:\ & (E, F) \mapsto (L_E(F), E)
	\end{align*}
	called right and left mutations. Moreover, $R, L$ are inverses in that we have isomorphisms of pairs $L \circ R (E, F) \cong (E, F)$ and $R \circ L (E, F) \cong (E, F)$.
\end{prop}

We can extend these to mutations on isomorphism classes of exceptional collections of any lengths.

\subsection*{Orlov's semiorthogonal decomposition theorem}
Let $A = \bigoplus_{i \geq 0} A_i$ be a two-sided Noetherian, graded connected $k$-algebra throughout. We say that $A$ is Artin-Schelter Gorenstein if $A$ is Gorenstein of dimension $d$ (meaning $\idim(_A A) = \idim(A_A) = d < \infty$) and $A$ satisfies the additional \emph{Gorenstein condition} 
\[
	{\rm Ext}_A^i(k, A) = \begin{cases} 0 & i \neq d \\ k(-a) & i =d. \end{cases}
\]
interpreted as isomorphism of graded modules. This latter condition follows for free in the commutative case \cite{BHBook}, but has to be imposed otherwise. The integer $a$ is the $a$-invariant (and $-a$ is the Gorenstein parameter\footnote{The reader is warned that $`a'$ is often used to denote the Gorenstein parameter in the literature, notably in \cite{Or1}.}). Since $A$ is Noetherian, we can define the abelian category 
\[
	\qgr A := \grmod A / \mathsf{grmod}_0 \ \! A
\]
as the Serre quotient category of finitely generated graded $A$-modules by the Serre subcategory of finite length modules. The significance of this construction is due to a classical theorem of Serre.
\begin{thm}[Serre] Let $R = \bigoplus_{n \geq 0} R_n$ be a standard graded, Noetherian commutative $k$-algebra. Let $X = \proj R$ be its projective scheme. Then the sheafification functor $\widetilde{(-)}: \grmod R \to \coh X$ descends to an exact equivalence of abelian categories 
	\[
		\widetilde{(-)}: \qgr R \xrightarrow{\cong} \coh X.
	\]
	Its inverse sends $\F$ to the class of the finitely generated graded $R$-module
	\[
		\Gamma_{\geq i}(X, \F) = \bigoplus_{n \geq i} \Gamma(X, \F(n))
	\]
	for any $i \in \Z$. 
\end{thm}
Note that the image of the homogeneous section functor $\Gamma_*(X, \F) = \bigoplus_{n \in \Z}\Gamma(X, \F(n))$ may fail to be finitely generated as an $R$-module, e.g. for $\F$ a skyscraper sheaf.\bigskip

Now in general, for any choice of cut-off $i \in \Z$, denote by $\mathsf{grmod}_{\geq i}\ \! A$ the full abelian subcategory of graded $A$-modules with $M_{j} = 0$ for $j < i$. The quotient functor $\pi: \grmod A \twoheadrightarrow \qgr A$ restricts to an essentially surjective exact functor $\pi_i: \mathsf{grmod}_{\geq i}\ \! A \twoheadrightarrow \qgr A$. \bigskip

Assume from this point on that $A$ is Artin-Schelter Gorenstein. Then Orlov observes \cite[Sect. 2]{Or1} that $\pi_i$ admits a right adjoint $\omega_i: \qgr A \to \mathsf{grmod}_{\geq i}\ \! A$ given by
\[
	\omega_i(M) = \bigoplus_{n \geq i}{\rm Hom}_{\qgr A}(\pi A, M(n)).
\]
Moreover, we have $\pi_i \omega_i \cong {\rm id}$ and so $\omega_i$ is fully faithful. This extends to an adjoint pair on derived categories
\[
	\pi_i: 	{\rm D}^b(\mathsf{grmod}_{\geq i}\ \! A) \leftrightarrows {\rm D}^b(\qgr A): {\rm R}\omega_i
\]
with ${\rm R}\omega_i$ fully faithful and $\pi_i$ essentially surjective. \bigskip

Since $A$ is Gorenstein, we have another Verdier quotient ${\rm st}: {\rm D}^b(\grmod A) \twoheadrightarrow {\rm D}_{sg}^{\Z}(A)$ given by stabilisation. Using semiorthogonal decomposition techniques, Orlov observes that ${\rm st}$ restricted to $ {\rm D}^b(\mathsf{grmod}_{\geq i}\ \! A)$ stays essentially surjective and gains a left adjoint which we shall denote by $X \mapsto X_{[\geq i]}$
\[
	(-)_{[\geq i]}: {\rm D}^{\Z}_{sg}(A) \rightleftarrows {\rm D}^b(\mathsf{grmod}_{\geq i}\ \! A): {\rm st}
\]
Moreover, we have ${\rm st}\circ  (-)_{[\geq i]} \cong {\rm id}$, and so $(-)_{[\geq i]}$ is fully faithful. We may then compose both adjoints as shown below:
\[
	\xymatrix{ & {\rm D}^b(\mathsf{grmod}_{\geq i}\ \!A) \ar@<.5ex>[dl]^-{st} \ar@<-.5ex>[dr]_-{\pi_i} & \\ {\rm D}_{sg}^{\Z}(A) \ar@<.5ex>[ur]^-{(-)_{[\geq i]}} \ar@<-.5ex>[rr]_{\Phi_i} && {\rm D}^b(\qgr A) \ar@<-.5ex>[ll]_{\Psi_i} \ar@<-.5ex>[ul]_-{ {\rm R}\omega_i} 
	}
\] 
This yields an adjoint pair $(\Phi_i, \Psi_i)$ with $\Phi_i = \pi \circ (-)_{[\geq i]}$ and $\Psi_i = {\rm st} \circ {\rm R}\omega_i$. Note that while the categories $ {\rm D}^b(\qgr A)$ and ${\rm D}^{\Z}_{sg}(A)$ do not depend on the resulting cutoff $i$, both functors $(\Phi_i, \Psi_i)$ do and will generally differ as $i$ varies. \bigskip 

The following is Orlov's semiorthogonal decomposition theorem.
\begin{thm}[{\cite[Thm. 2.5]{Or1}}]\label{OrlovSOD} Let $A$ be an Artin-Schelter Gorenstein $k$-algebra with $a$-invariant $a \in \Z$. The above functors and triangulated categories are related as follows:
	\begin{enumerate}[i)]
	\item (Fano case) if $a < 0$, there is a semiorthogonal decomposition
		\[
			{\rm R}\omega_i\left({\rm D}^b(\qgr A)\right) = \big\langle A(-i + a +1), A(-i+a+2), \dots, A(-i), {\rm D}^{\Z}_{sg}(A)_{[\geq i]} \big\rangle.
		\]
		Applying $\pi_i$, this descends to a semiorthogonal decomposition 
		\[
		{\rm D}^b(\qgr A) = \big\langle \pi A(-i + a + 1), \pi A(-i+a+2), \dots, \pi A(-i), \Phi_i{\rm D}^{\Z}_{sg}(A)) \big\rangle.
		\]
	\item (Calabi-Yau case) if $a = 0$, the essential images of both embeddings in $ {\rm D}^b(\mathsf{grmod}_{\geq i}\ \! A)$ are equal 
		\[
			{\rm D}^{\Z}_{sg}(A)_{[\geq i]} = {\rm R}\omega_i {\rm D}^b(\qgr A)  
		\]
		hence $(\Phi_i, \Psi_i)$ give inverse equivalences
		\[
			\Phi_i: {\rm D}^{\Z}_{sg}(A) \cong {\rm D}^b(\qgr A) :\Psi_i.
		\]
	\item (General type case) if $a > 0$, there is a semiorthogonal decomposition
		\[
			{\rm D}^{\Z}_{sg}(A)_{[\geq i]} = \big\langle k(-i), k(-i-1), \dots, k(-i-a+1), {\rm R}\omega_{i+a} {\rm D}^b(\qgr A) \big\rangle.
		\]
		Applying stabilisation ${\rm st}$, this descends to a semiorthogonal decomposition
		\[
			{\rm D}^{\Z}_{sg}(A) = \big\langle k^{st}(-i), k^{st}(-i-1), \dots, k^{st}(-i-a+1), \Psi_{i+a} {\rm D}^b(\qgr A) \big\rangle.
		\]
\end{enumerate}
\end{thm}
\begin{rem}
	The above theorem is usually just stated in terms of $\Psi_i$ and $\Phi_i$ and the resulting semiorthogonal decompositions in $ {\rm D}^b(\qgr A)$ and ${\rm D}_{sg}^{\Z}(A)$. However the above statement is what Orlov shows as part of the proof, and this stronger version has many uses.
\end{rem}\bigskip

The existence of the left adjoint $(-)_{[\geq i]}: {\rm D}_{sg}^{\Z}(A) \to {\rm D}^b(\mathsf{grmod}_{\geq i}\ \! A)$ follows by abstract nonsense arguments involving admissible subcategories, and it is important in applications to have a concrete description. The following is due to Buchweitz.\bigskip

First apply the equivalence ${\rm D}^{\Z}_{sg}(A) \cong \underline{{\rm MCM}}^{\mathbb{Z}}(A)$. We will construct a functor 
\[
	(-)_{[\geq i]}: \underline{{\rm MCM}}^{\mathbb{Z}}(A) \to {\rm D}^b(\mathsf{grmod}_{\geq i}\ \! A)
\]
which we also denote $(-)_{[\geq i]}$ by abuse of notation. We will then give a direct proof that $(-)_{[\geq i]}$ is a left adjoint to MCM approximation, giving a model for the sought-after left adjoint implicit in Orlov's result. \bigskip

Let $M$ be a graded MCM $A$-module with complete resolution $C$. Let $C_{[<i]} \subseteq C$ be the graded submodule whose terms are graded free $A$-modules generated in degree $<i$. The short exact sequence 
\[
\xi:	0 \to C_{[<i]} \to C \to C/C_{[<i]} \to 0
\]
is split as a sequence of $A$-modules and we have ${\rm Hom}_{gr A}(C_{[<i]}, C/C_{[<i]}) = 0$. Since the differential $d: C_n \to C_{n-1}$ is $A$-linear and homogeneous of degree zero, we see that $C_{[<i]}$ is a subcomplex, and we interpret $\xi$ as a short exact sequence of complexes. Applying the same idea shows that $C \mapsto C_{[<i]}$ and $C/C_{[<i]}$ is natural in $C$ and preserves homotopy equivalences.

The quotient complex $C_{[\geq i]} = C/C_{[<i]}$ has terms in $\mathsf{grmod}_{\geq i}\ \! A$. Since $A$ is graded connected, up to homotopy one can replace $C$, and thus $C_{[\geq i]}$, by a minimal complex, which shows that $C_{[\geq i]}$ has bounded cohomology as the degrees of generators of $C_n$ strictly increase as $n \to \infty$ (and respectively decrease as $n \to -\infty$). Hence we obtain a functor $ \underline{{\rm MCM}}^{\mathbb{Z}}(A) \to {\rm D}^b(\mathsf{grmod}_{\geq i}\ \! A)$ given by $M \mapsto C_{[\geq i]}$. We now consider the pair of functors 
\[
	(-)_{[\geq i]}: \underline{{\rm MCM}}^{\mathbb{Z}}(A) \leftrightarrows {\rm D}^b(\mathsf{grmod}_{\geq i}\ \! A): {\rm st}.
\]
\begin{prop}[Buchweitz]\label{BuchweitzAdjoint} The pair $\left((-)_{[\geq i]},\ \! {\rm st} \right)$ is an adjoint pair. 
\end{prop} 
\begin{proof} Let $\mathcal{K} = \mathcal{K}({\rm proj}^{\Z}\ \! A)$ be the homotopy category of complexes of finite graded projectives. Then $\xi$ gives a triangle in $\mathcal{K}$ 
\[
	C_{[<i]} \to C \to C_{[\geq i]} \to C_{[<i]}[1]
\]
natural in $C$. Let $F \in {\rm D}^b(\mathsf{grmod}_{\geq i}\ \! A)$ with projective resolution $P_* \xrightarrow{\sim} F$, taken so that all terms $P_n$ have generators in degree $\geq i$. Then ${\rm Hom}_{\mathcal{K}}(C_{[<i]}, P_*) = 0$ and so 
\[
	{\rm Hom}_{\mathcal{K}}(C_{[\geq i]}, P_*) \cong {\rm Hom}_{\mathcal{K}}(C, P_*).
\]
Next, writing $P_*^t = P_{\geq n}$ as the tail truncation and $P_*^h = P_{\leq n-1}$ for the head for $n \gg 0$, we have a distinguished triangle
\[
	P_*^{{\rm h}} \to P_* \to P_*^{{\rm t}} \to P_*^{{\rm h}}[1]
\]
with the head $P_*^{{\rm h}}$ a perfect complex and the tail $P_*^{{\rm t}}$ a shifted resolution of an MCM module. Since we have ${\rm Hom}_{\mathcal{K}}(C, P_*^{{\rm h}}) = 0$, the long exact sequence of Hom spaces gives
\[
	{\rm Hom}_{\mathcal{K}}(C, P_*) \cong {\rm Hom}_{\mathcal{K}}(C, P_*^{{\rm t}}).
\]
As $P_*^{{\rm t}}$ is a shifted resolution of an MCM module, it extends to a complete resolution $D$ of the MCM approximation $F^{st}$. Since each morphism $C \to P_*^{{\rm t}}$ has a unique lift to $D$ up to homotopy, the map $D \to P_*^{{\rm t}}$ induces an isomorphism
\[
	{\rm Hom}_{\mathcal{K}}(C, D) \cong {\rm Hom}_{\mathcal{K}}(C, P_*^t).
\]
Combining the above, we obtain natural isomorphisms
\begin{align*}
	{\rm Hom}_{{\rm D}^b(\mathsf{grmod}_{\geq i}A)}(C_{[\geq i]}, F) &\cong {\rm Hom}_{\mathcal{K}}(C_{[\geq i]}, P_*)\\
	&\cong {\rm Hom}_{\mathcal{K}}(C, P_*)\\
	&\cong {\rm Hom}_{\mathcal{K}}(C, D)\\
	&\cong \underline{{\rm Hom}}_{gr A}(M, F^{st}). 
\end{align*}
\end{proof}

By uniqueness of adjoints, it follows that Orlov's left adjoint is given on MCM modules by $M_{[\geq i]} = C_{[\geq i]}$. This has some immediately interesting consequences:

\begin{cor}\label{a>0} Let $A$ be Artin-Schelter Gorenstein with $a > 0$. Let $C$ be a complete resolution of the MCM module $k^{st}$. Then $C_{[\geq 0]}$ is a projective resolution of $k$. In other words, $k^{st}$ can be ``unstabilised''.
\end{cor}
We note that this typically fails for $a \leq 0$, as can be seen already from the Tate resolution of the residue field over a hypersurface ring.